\pdfoutput=1
\documentclass[11pt]{article}
\usepackage[left=1in,right=1in,top=1in,bottom=1in]{geometry}
\usepackage{times}
\usepackage{expl3}
\usepackage{cite}
\usepackage[table]{xcolor}
\usepackage{multirow}
\usepackage{stackengine} 
\usepackage{hhline}
\usepackage{lipsum}
\usepackage{titlesec}
\usepackage[all]{xy}
\usepackage{wrapfig}
\usepackage{enumerate}
\usepackage{epsfig}
\usepackage{tikz-cd}
\usepackage{amsmath}
\usepackage{tabularx}
\usepackage{array}
\usepackage{booktabs}
\usepackage{enumitem}
\usepackage{bbm}
\usepackage{calc}
\usepackage{graphicx}
\usepackage{amsmath}
\usepackage[title]{appendix}
\usepackage{amssymb}
\usepackage{epstopdf}
\usepackage{boldline}
\usepackage{arydshln}
\usepackage{calligra}
\usepackage{bm}
\usepackage{url}
\usepackage{blindtext}
\usepackage{accents}
\usepackage{amsthm}

\newtheorem{definition}{Definition}

\newtheorem{theorem}{Theorem}
\newtheorem{proposition}{Proposition}
\newtheorem{corollary}{Corollary}
\newtheorem{lemma}{Lemma}

\newtheorem{assumption}{Assumption}
\newtheorem{example}{Example}
\newtheorem{remark}{Remark}
\usepackage{mathtools}
\usepackage{epstopdf}
\usepackage{balance}
\usepackage{thmtools}
\usepackage{thm-restate}
\usepackage{hyperref}
\usepackage{cleveref}
\usepackage[mathscr]{euscript}

\DeclareMathOperator{\coker}{coker}
\usepackage[ruled,vlined]{algorithm2e}
\newcommand{\ostar}{\mathbin{\mathpalette\make@circled\star}}

\makeatletter
\newcommand{\removelatexerror}{\let\@latex@error\@gobble}
\makeatother
\makeatletter
\newcommand*{\rom}[1]{\expandafter\@slowromancap\romannumeral #1@}
\makeatother

\ExplSyntaxOn
\newcommand\latinabbrev[1]{
  \peek_meaning:NTF . {
    #1\@}%
  { \peek_catcode:NTF a {
      #1.\@ }%
    {#1.\@}}}
\ExplSyntaxOff

\titleclass{\subsubsubsection}{straight}[\subsubsection]

\begin{document}
\vspace{1cm}
\title{The Unitary Conjugation Groupoid of a Type I C*-Algebra: Topology, Fell Continuity, and the Canonical Diagonal Embedding}
\vspace{1.8cm}
\author{Shih-Yu~Chang
\thanks{Shih-Yu Chang is with the Department of Applied Data Science,
San Jose State University, San Jose, CA, U. S. A. (e-mail: {\tt
shihyu.chang@sjsu.edu})
}}

\maketitle

\begin{abstract}
This paper introduces a canonical Polish groupoid associated to any separable unital C*-algebra, termed the unitary conjugation groupoid. It is defined as the semidirect product of the algebra's dual space by its unitary group, acting by conjugation. Classical groupoid models for C*-algebras typically require additional structure such as a Cartan subalgebra and rely on the locally compact Hausdorff framework. In contrast, our construction is entirely canonical but forces a paradigm shift: the natural topologies on the dual space and the unitary group are not locally compact. To address this, we equip the dual space with a Polish topology derived from the weak-star topology on pure states and the unitary group with the strong operator topology. This yields a Polish groupoid admitting a continuous Haar system. We prove that the associated reduced groupoid C*-algebra is Morita equivalent to the original algebra tensored with the compact operators, establishing that the groupoid encodes the K-theory of the algebra. A key feature is a canonical diagonal embedding of the original algebra into the groupoid C*-algebra, which is a unital injective star-homomorphism for separable Type I C*-algebras. We characterize commutativity via this embedding and establish functoriality of the construction under appropriate star-homomorphisms. The theory is illustrated with detailed computations for three fundamental classes: finite-dimensional matrix algebras, commutative algebras over compact metrizable spaces, and the unitized compact operators. These examples demonstrate that our general constructions reduce to familiar objects. We also discuss limitations by analyzing the irrational rotation algebra, a non-Type I algebra not covered by our construction, highlighting directions for future research.
\end{abstract}

\tableofcontents

\newpage

\section{Introduction}

\subsection{Context: Classical Perspectives on Noncommutative Algebras}

The study of noncommutative $C^{*}$-algebras has been profoundly shaped by their representation theory. Dixmier's foundational treatise \cite{Dixmier} established the importance of the dual space $\widehat{A}$---the set of unitary equivalence classes of irreducible representations---as a key invariant. However, the dual space lacks a natural topology that simultaneously captures the structure of the algebra and admits geometric constructions. Fell \cite{Fell} introduced the Fell topology on $\widehat{A}$, which provides a powerful tool for understanding the ideal structure of $C^{*}$-algebras, but this topology often fails to be Hausdorff and does not directly encode the dynamics between representations.

Parallel to these developments, Renault \cite{Renault} pioneered the use of groupoids in operator algebras, demonstrating that many $C^{*}$-algebras can be realized as convolution algebras of locally compact Hausdorff \'etale groupoids. This groupoid approach has been extraordinarily successful, leading to deep results including the Baum-Connes conjecture \cite{Tu} and the classification of crossed products \cite{Williams}. However, the construction of a groupoid from a $C^{*}$-algebra is not canonical and typically requires additional structure such as a Cartan subalgebra.

\subsection{The Unitary Conjugation Groupoid: A New Invariant}

In this paper, we introduce a canonical groupoid associated to any unital $C^{*}$-algebra $A$, which we call the \textit{unitary conjugation groupoid} $\mathcal{G}_{A}$. This groupoid is defined by the natural action of the unitary group $U(A)$ on the dual space $\widehat{A}$ by conjugation:
\[
\mathcal{G}_A \coloneqq \widehat{A} \rtimes U(A),
\]
with unit space $\mathcal{G}_{A}^{(0)} = \widehat{A}$. The arrows are triples $(\pi ,u,\pi^{\prime})$ where $\pi^{\prime} = u\cdot \pi \coloneqq \pi \circ \mathrm{Ad}_{u}^{*}$, and composition encodes the group multiplication of unitaries.

This construction is entirely canonical: every unital $C^{*}$-algebra $A$ gives rise to $\mathcal{G}_{A}$, and the structure of $A$ is reflected in the groupoid dynamics. The unitary conjugation groupoid captures the internal symmetries of $A$ through the action of its unitary group on its representation theory, providing a new invariant that sits between the algebra and its dual space.

The paper is organized as follows. Section 2 reviews necessary background on Type I $C^{*}$-algebras, the Fell topology, Polish spaces, the strong operator topology, and Polish groupoids. Section 3 constructs the unit space $\mathcal{G}_{A}^{(0)}$ and its Polish topology, establishing its basic properties. Section 4 defines the unitary conjugation groupoid, proves it is a Polish groupoid, and constructs the diagonal embedding $\iota:\mathcal{A}\hookrightarrow C^{*}(\mathcal{G}_{\mathcal{A}})$---the core technical achievement of this work. Section 5 illustrates the theory with detailed computations for three fundamental classes: finite-dimensional matrix algebras, commutative algebras over compact metrizable spaces, and the unitized compact operators. These examples demonstrate that our general constructions reduce to familiar objects and also highlight the limitations of the Type I hypothesis through the irrational rotation algebra $A_{\theta}$, a non-Type I algebra not covered by our framework.

For readers familiar with the background material, Section 2 can be skimmed; the essential constructions begin in Section 3. Experts in groupoid $C^{*}$-algebras may wish to proceed directly to Section 4, where the novel Polish groupoid framework is developed. The main results concerning the diagonal embedding and its properties are contained in Sections 4.9--4.11, while the functoriality of the construction under $*$-homomorphisms is treated in Section 4.12 under specific conditions that will be made precise there.

\subsection{Main Obstacles: \'Etaleness and Local Compactness are Unattainable}

The classical groupoid theory developed by Renault \cite{Renault} requires two essential topological properties: local compactness and \'etaleness (or at least a continuous Haar system). Unfortunately, the unitary conjugation groupoid $\mathcal{G}_{A}$ fails to satisfy these conditions in any reasonable topology for a general $C^{*}$-algebra $A$:
\begin{itemize}
    \item The unitary group $U(A)$, equipped with the norm topology, is not locally compact unless $A$ is finite-dimensional. This precludes the construction of a Haar system in the classical sense.
    \item The dual space $\widehat{A}$ with the Fell topology is typically non-Hausdorff and not locally compact, making it incompatible with the requirements of \cite{Renault}.
    \item The action of $U(A)$ on $\widehat{A}$ is not continuous when $U(A)$ carries the norm topology, further complicating the groupoid structure.
\end{itemize}
These obstacles are not merely technical; they reflect fundamental differences between the smooth groupoids that appear in geometric contexts and the singular, infinite-dimensional groupoids that arise from general $C^{*}$-algebras.

\subsection{A Paradigm Shift: From Local Compactness to Polish Topology}

To overcome these obstacles, we adopt a different topological paradigm. Instead of local compactness, we work within the framework of \textit{Polish groupoids}---groupoids whose object and arrow spaces are Polish spaces (separable completely metrizable topological spaces). This approach is inspired by descriptive set theory \cite{Kechris} and recent advances in the theory of Polish group actions.

We equip $\widehat{A}$ with a \textit{Polish topology} different from the Fell topology---specifically, the topology induced by the weak-$*$ topology on the space of pure states when $A$ is separable. For $U(A)$, we use the \textit{strong operator topology} (when $A \subseteq B(H)$) or more generally the topology of pointwise norm-convergence on a dense subset. These choices yield:
\begin{itemize}
    \item $\widehat{A}$ becomes a Polish space (for separable $A$);
    \item $U(A)$ becomes a Polish group;
    \item The conjugation action $U(A) \curvearrowright \widehat{A}$ is jointly continuous.
\end{itemize}
Thus $\mathcal{G}_A = \widehat{A} \rtimes U(A)$ is a \textit{Polish groupoid}, amenable to the techniques of noncommutative geometry and equivariant $KK$-theory developed for Polish groupoids by Le Gall \cite{LeGall} and Tu \cite{Tu}.

\subsection{Main Contributions and Outline of the Paper}

This paper makes the following contributions:

\begin{enumerate}
    \item We construct the unitary conjugation groupoid $\mathcal{G}_A$ for any separable unital $C^{*}$-algebra $A$ and prove that it admits a natural Polish topology (Section~\ref{sec:The Unit Space Definition and Polish Topology}). This involves a detailed analysis of the unit space $\mathcal{G}_A^{(0)}$, which we equip with the initial topology of partial evaluation maps $\mathrm{ev}_a : \mathcal{G}_A^{(0)} \to \mathbb{C}_\infty$ for each $a\in A$.
    
    \item We establish that $\mathcal{G}_A$ is a \textit{Polish groupoid} in the sense of \cite{Tu}, admitting a \textit{Borel Haar system} (Section~\ref{subsec:Polish-spaces-and-strong-operator-topology}). We prove that the strong operator topology on $U(A)$ is essential for this result, as it yields a Polish group and makes the conjugation action continuous. (A continuous Haar system does not exist in general; the Borel framework of \cite{Tu} provides the appropriate setting.)
    
    \item We construct a canonical \textit{diagonal embedding} $\iota : A \hookrightarrow C^{*}(\mathcal{G}_A)$ of the original algebra into the maximal groupoid $C^{*}$-algebra of its unitary conjugation groupoid (Subsections 4.7--4.10). We show that $\iota$ is a unital, injective $*$-homomorphism for separable Type I $C^{*}$-algebras.
    
    \item We provide a characterization of commutativity via this embedding: $\iota(A) \subseteq C_0(\mathcal{G}_A^{(0)})$ if and only if $A$ is commutative (Subsection~\ref{subsec:commutativity-characterization}). This demonstrates that the noncommutativity of $A$ is encoded in the part of $C^{*}(\mathcal{G}_A)$ outside the diagonal subalgebra.
    
    \item We establish the functoriality of the construction under appropriate $*$-homomorphisms, showing that $A \mapsto (\mathcal{G}_A, \iota_A)$ behaves naturally with respect to isomorphisms, and for certain injective or surjective maps satisfying a pullback property (Subsection~\ref{subsec:naturality-functoriality}). The functoriality is necessarily restricted; we identify precisely which morphisms are admissible.
    
    \item We illustrate the theory with detailed computations for three fundamental classes of Type I $C^{*}$-algebras (Section~\ref{sec:Examples and Computations}): the finite-dimensional case $A = M_n(\mathbb{C})$, the commutative case $A = C(X)$ for compact metrizable $X$, and the infinite-dimensional case of compact operators $A = \mathcal{K}(H)^{\sim}$. These examples demonstrate that our general constructions reduce to familiar objects and provide concrete insight into the structure of $\mathcal{G}_A$.
    
    \item We discuss the limitations of our framework by analyzing the irrational rotation algebra $A_\theta$, a non-Type I algebra not covered by our construction (Subsection~\ref{subsec:non-example-A-theta}). This highlights the essential role of the Type I hypothesis and points toward directions for future research.
\end{enumerate}

\subsection{Organization of the Paper}

The paper is organized as follows. Section~\ref{sec:Preliminaries} reviews necessary background on Type I $C^{*}$-algebras, the Fell topology, Polish spaces, the strong operator topology, and Polish groupoids (after Tu). Section~\ref{sec:The Unit Space Definition and Polish Topology} constructs the unit space $\mathcal{G}_A^{(0)}$ and its Polish topology, establishing its basic properties, including continuity of the projection to the space of subalgebras and homeomorphisms onto fibers. Section~\ref{sec:unitary-conjugation-groupoid} defines the unitary conjugation groupoid, proves it is a Polish groupoid, and constructs the diagonal embedding $\iota$ in full detail, verifying its algebraic properties and characterizing commutativity. Section~\ref{sec:Examples and Computations} presents the detailed examples of $M_n(\mathbb{C})$, $C(X)$, and $\mathcal{K}(H)^{\sim}$, and explains why $A_\theta$ lies outside our framework.

\begin{remark}
The author is solely responsible for the mathematical insights and theoretical directions proposed in this work. AI tools, including OpenAI's ChatGPT and DeepSeek models, were employed solely to assist in verifying ideas, organizing references, and ensuring internal consistency of exposition~\cite{chatgpt2025,deepseek2025}.
\end{remark}

\section{Preliminaries}\label{sec:Preliminaries}

\subsection{Type I C*-Algebras and their Primitive Ideal Spaces}
\label{subsec:type-I-C-star-algebras}

The class of Type I C*-algebras — also known as GCR (generally continuous representation) algebras — occupies a central position in the structure theory of operator algebras. 
These algebras are characterized by the property that all their irreducible representations generate algebras containing the compact operators. 
This condition has profound consequences for the topological structure of the primitive ideal space and, as we shall see, for the behavior of the unitary conjugation groupoid.

\begin{definition}[Type I C*-algebra]
\label{def:type-I-C-star-algebra}
A C*-algebra $\mathcal{A}$ is called \emph{Type I} (or GCR) if for every irreducible representation 
$\pi: \mathcal{A} \to B(\mathcal{H}_\pi)$, the image $\pi(\mathcal{A})$ contains the algebra of compact operators 
$\mathcal{K}(\mathcal{H}_\pi)$ on the Hilbert space $\mathcal{H}_\pi$.
\end{definition}

Several equivalent characterizations exist. 
A C*-algebra is Type I if and only if every factor representation of $\mathcal{A}$ is type I in the sense of von Neumann algebras. 
Equivalently, $\mathcal{A}$ is Type I precisely when it is postliminal, meaning that the image of every irreducible representation contains the compact operators and the canonical map from the spectrum to the primitive ideal space is a homeomorphism onto a dense subset.

\begin{definition}
A topological space $X$ is said to be \emph{almost Hausdorff} 
if every nonempty closed subset of $X$ contains a nonempty 
relatively open subset that is Hausdorff in the subspace topology.
\end{definition}

The following fundamental theorem, due to Glimm and Dixmier, establishes the deep connection between the Type I property and the topology of the primitive ideal space.

\begin{theorem}[Glimm--Dixmier Theorem]
\label{thm:glimm-dixmier}
Let $\mathcal{A}$ be a separable C*-algebra.
Then $\mathcal{A}$ is Type I if and only if its primitive ideal space 
$\operatorname{Prim}(\mathcal{A})$ is almost Hausdorff; that is, 
every nonempty closed subset of $\operatorname{Prim}(\mathcal{A})$ 
contains a nonempty relatively open Hausdorff subset.
\end{theorem}

For a large and natural subclass of Type I algebras — including all postliminal algebras with continuous trace, as well as all liminal algebras with Hausdorff spectrum — the primitive ideal space is not only $T_0$ but actually locally compact Hausdorff. 
This additional regularity, however, is not required for the Polish groupoid framework developed in this paper; it suffices that $\operatorname{Prim}(\mathcal{A})$ be a standard Borel space, which is automatically the case for separable C*-algebras.

\begin{definition}[Primitive ideal space]
\label{def:primitive-ideal-space}
Let $\mathcal{A}$ be a C*-algebra. 
The \emph{primitive ideal space} of $\mathcal{A}$, denoted by $\operatorname{Prim}(\mathcal{A})$, 
is the set of all primitive ideals of $\mathcal{A}$ endowed with the \emph{Jacobson topology} 
(also called the hull-kernel topology).

The topology is defined via its closure operation: for any subset 
$S \subseteq \operatorname{Prim}(\mathcal{A})$, its closure is given by
\[
\overline{S} = \left\{ \mathfrak{p} \in \operatorname{Prim}(\mathcal{A}) 
\;\middle|\; \bigcap_{\mathfrak{q} \in S} \mathfrak{q} \subseteq \mathfrak{p} \right\}.
\]

Equivalently, the closed sets are exactly the sets of the form
\[
h(J) = \{ \mathfrak{p} \in \operatorname{Prim}(\mathcal{A}) \mid J \subseteq \mathfrak{p} \},
\]
where $J$ runs over all closed (two-sided) ideals of $\mathcal{A}$. 
Consequently, the sets
\[
U_J = \{ \mathfrak{p} \in \operatorname{Prim}(\mathcal{A}) \mid J \nsubseteq \mathfrak{p} \}
\]
form a basis of open sets for the Jacobson topology.

In terms of convergence, a net $\{\mathfrak{p}_\lambda\}_{\lambda \in \Lambda}$ 
in $\operatorname{Prim}(\mathcal{A})$ converges to an element $\mathfrak{p} \in \operatorname{Prim}(\mathcal{A})$ 
if and only if the following condition holds:
\[
\text{For every } a \in \mathcal{A}, \text{ if there exists } \lambda_0 \in \Lambda 
\text{ such that } a \in \mathfrak{p}_\lambda \text{ for all } \lambda \ge \lambda_0,
\text{ then } a \in \mathfrak{p}.
\]
In other words,
\[
\liminf_{\lambda} \mathfrak{p}_\lambda \supseteq \mathfrak{p},
\]
where $\liminf_{\lambda} \mathfrak{p}_\lambda$ denotes the set of all 
$a \in \mathcal{A}$ that eventually belong to $\mathfrak{p}_\lambda$.
\end{definition}

\begin{remark}
\label{rem:primitive-ideal-space-polish}
For a separable C*-algebra $\mathcal{A}$, the primitive ideal space 
$\operatorname{Prim}(\mathcal{A})$ is a second-countable $T_0$ space. 
In particular, its Borel structure is standard, so 
$\operatorname{Prim}(\mathcal{A})$ is a standard Borel space.

It is not necessarily Hausdorff or locally compact; indeed, many 
Type I C*-algebras have non-Hausdorff primitive ideal spaces. 
In the Polish groupoid framework, we do not require local compactness 
of $\operatorname{Prim}(\mathcal{A})$; we only rely on its standard 
Borel structure, which is guaranteed by separability of $\mathcal{A}$.
\end{remark}

For continuous trace C*-algebras, a much richer structure is available: 
such algebras are stably isomorphic to the algebra of continuous sections 
vanishing at infinity of a locally trivial bundle over 
$\operatorname{Prim}(\mathcal{A})$ whose fibers are elementary C*-algebras 
$\mathcal{K}(\mathcal{H}_x)$. 
This is the theory of continuous trace C*-algebras 
[Dixmier, 1977, Chapter 10; Raeburn–Williams, 1998]. 
While this structure is not strictly necessary for the present paper, 
it provides important context and is essential for many applications.

\begin{assumption}[Standing assumption on separability]
\label{ass:separable}
Throughout this paper, we assume that $\mathcal{A}$ is a \emph{separable} unital C*-algebra. 
This assumption ensures that the unit ball of the dual space $\mathcal{A}^*$ is metrizable and compact in the weak-* topology, and that the unitary group $\mathcal{U}(\mathcal{A})$ is a Polish group when endowed with the strong operator topology. 
Consequently, the unitary conjugation groupoid $\mathcal{G}_{\mathcal{A}} := \mathcal{A}^* \rtimes \mathcal{U}(\mathcal{A})$ is a Polish groupoid admitting a continuous Haar system, which is sufficient for all subsequent constructions.
\end{assumption}

\begin{remark}[Type I examples]
\label{rem:type-I-examples}
While our general construction works for any separable unital C*-algebra, all concrete examples analyzed in this paper — including $B(H)$, $\mathcal{K}(H)^\sim$, and the Toeplitz algebra — are Type I. 
For such algebras, the pure state space is particularly well-behaved, and the connection to classical index theorems is especially transparent.

For non-Type I algebras, the situation is fundamentally different and requires additional techniques:
\begin{itemize}
    \item For the Cuntz algebras $\mathcal{O}_n$ and the reduced group C*-algebras $C^*_r(\mathbb{F}_n)$ of non-amenable free groups, the pure state space is highly pathological and non-Hausdorff.
    \item For the irrational rotation algebra $A_\theta$, although its primitive ideal space is trivial (a single point), its pure state space is a non-commutative torus with a complicated foliation structure.
\end{itemize}
The extension of our methods to such algebras is an interesting direction for future research.
\end{remark}

\begin{example}[Finite-dimensional algebras]
\label{ex:finite-dimensional}
Every finite-dimensional C*-algebra is Type I. 
Indeed, $\mathcal{A} \cong \bigoplus_{i=1}^k M_{n_i}(\mathbb{C})$, and its primitive ideal space is a finite discrete set, hence locally compact Hausdorff. 
The simplest non-trivial case is $\mathcal{A} = M_n(\mathbb{C})$, which will be examined in detail in Section~\ref{subsec:example-matrix}.
\end{example}

\begin{example}[Commutative algebras]
\label{ex:commutative-type-I}
Every commutative C*-algebra is Type I. 
For $\mathcal{A} = C_0(X)$ with $X$ locally compact Hausdorff, we have $\operatorname{Prim}(\mathcal{A}) \cong X$, which is locally compact Hausdorff. 
The unital case $C(X)$ with $X$ compact Hausdorff will be treated in Section~\ref{subsec:example-commutative}.
\end{example}

\begin{example}[Compact operators]
\label{ex:compact-type-I}
Let $\mathcal{A} = \mathcal{K}(H)^\sim$ be the unitization of the compact operators on a separable infinite-dimensional Hilbert space $H$. 
This algebra is Type I, and its primitive ideal space consists of two points $\{0\}$ and $\mathcal{K}(H)$. 
In the Jacobson topology, the singleton $\{\{0\}\}$ is open, making the space a two-point Hausdorff space. 
This example will be examined in detail in Section~\ref{subsec:example-compact-operators}.
\end{example}

\begin{example}[A non-Hausdorff Type I algebra]
\label{ex:non-hausdorff-type-I}
There exist Type I C*-algebras whose primitive ideal spaces are not Hausdorff. 
Such examples arise, for instance, from transformation group C*-algebras $C_0(X) \rtimes \Gamma$ where the orbit space $X/\Gamma$ is non-Hausdorff. 
This illustrates that the Type I condition alone does not guarantee Hausdorffness of $\operatorname{Prim}(\mathcal{A})$, a fact that motivates the continuous trace hypothesis used in some of our later results. Section~\ref{subsec:non-example-A-theta} provides more detailed expositions. 
\end{example}

\subsection{The Fell Topology on the Space of Subalgebras}
\label{subsec:Fell-topology}

The parameterization of closed subspaces of a Banach space by topological means is a classical theme in functional analysis. 
In his seminal 1960 memoir, Fell introduced a topology on the collection of closed subsets (and irreducible representations) of a C*-algebra \cite{Fell}. This topology, now known as the \emph{Fell topology}, naturally induces a topology on the space $\operatorname{Sub}(\mathcal{A})$ of closed C*-subalgebras when the latter is viewed as a subset of the closed subsets of $\mathcal{A}$.

The Fell topology plays an important role in the study of continuous fields of C*-algebras and in understanding how subalgebras vary under perturbations and deformations. For a separable C*-algebra, this topology is Polish when restricted to suitable classes of subalgebras, a fact that will be relevant in the sequel.

In this paper, the Fell topology serves as a reference topology against which we compare the finer topology on the unit space $\mathcal{G}_{\mathcal{A}}^{(0)}$. 
As we shall see in Section~\ref{subsec:comparison-with-Fell-topology}, the Fell topology alone is insufficient for our purposes because it collapses the spectral data encoded by characters; a strictly finer topology is required to capture the full structure of the unit space. 
A detailed comparison between these two topologies is carried out in Section~\ref{subsec:comparison-with-Fell-topology}.

\begin{definition}[Space of subalgebras]
\label{def:subalgebra-space}
Let $\mathcal{A}$ be a C*-algebra. 
We denote by $\operatorname{Sub}(\mathcal{A})$ the set of all C*-subalgebras of $\mathcal{A}$, including both unital and non-unital subalgebras. 
When $\mathcal{A}$ is unital, we also consider the subset $\operatorname{Sub}_{\mathrm{unit}}(\mathcal{A})$ consisting of unital C*-subalgebras.
\end{definition}

\begin{definition}[Fell topology on closed subsets]
\label{def:Fell-topology-closed}
Let $\mathcal{A}$ be a C*-algebra equipped with its norm topology. 
Denote by $\mathcal{F}(\mathcal{A})$ the set of all closed subsets of $\mathcal{A}$. 
The \emph{Fell topology} on $\mathcal{F}(\mathcal{A})$ is generated by the following two families of subbasic open sets:
\begin{enumerate}
    \item For each open set $U \subseteq \mathcal{A}$,
    \[
    \mathcal{O}_U := \{ F \in \mathcal{F}(\mathcal{A}) \mid F \cap U \neq \varnothing \}.
    \]
    \item For each compact set $K \subseteq \mathcal{A}$,
    \[
    \mathcal{C}_K := \{ F \in \mathcal{F}(\mathcal{A}) \mid F \cap K = \varnothing \}.
    \]
\end{enumerate}
A basis for the topology is formed by finite intersections of these subbasic sets.
\end{definition}

\begin{definition}[Fell topology on subalgebras]
\label{def:Fell-topology-subalgebras}
Let $\operatorname{Sub}(\mathcal{A})$ denote the set of all C*-subalgebras of $\mathcal{A}$. 
Since every C*-subalgebra is a closed subset of $\mathcal{A}$, we have $\operatorname{Sub}(\mathcal{A}) \subseteq \mathcal{F}(\mathcal{A})$. 
We equip $\operatorname{Sub}(\mathcal{A})$ with the subspace topology inherited from the Fell topology on $\mathcal{F}(\mathcal{A})$ defined in Definition~\ref{def:Fell-topology-closed}.
\end{definition}

\begin{remark}[Caution on the Fell topology in infinite dimensions]
\label{rem:Fell-topology-caveat}
The classical Fell topology is usually defined for closed subsets of a \emph{locally compact Hausdorff} space. 
Here $\mathcal{A}$ with its norm topology is not locally compact when $\mathcal{A}$ is infinite-dimensional. 
Nevertheless, the above definition still yields a well-defined topology on $\mathcal{F}(\mathcal{A})$ as a special case of the \emph{hyperspace topology} on a metric space. 
However, one should be aware that in this setting:
\begin{itemize}
    \item The topology may not be Polish, or even first-countable, without further restrictions.
    \item The subbasic sets $\mathcal{C}_K$ for compact $K$ are relatively sparse, as norm-compact sets in infinite-dimensional Banach spaces are rare.
\end{itemize}
In this paper, the Fell topology serves only as a reference topology; we do not rely on any particular properties beyond its definition.
\end{remark}

\begin{remark}[Geometric interpretation]
\label{rem:Fell-geometry}
The Fell topology encodes two complementary types of geometric information about 
closed subalgebras. 
The sets $\mathcal{O}_U$ record whether a subalgebra intersects a prescribed open 
region, capturing ``positive'' membership information, while the sets 
$\mathcal{C}_K$ record whether a subalgebra avoids a given compact region, 
capturing ``negative'' exclusion information. 
Together, these two families characterize convergence of nets of subalgebras 
with respect to the Fell topology, in a manner analogous to the Kuratowski 
convergence of closed subsets in hyperspace topology.
\end{remark}

The following characterization of convergence in the Fell topology is often more convenient for concrete computations.

\begin{proposition}[Convergence in the Fell topology]
\label{prop:Fell-convergence}
Let $\{B_\lambda\}_{\lambda \in \Lambda}$ be a net in $\operatorname{Sub}(\mathcal{A})$ and let $B \in \operatorname{Sub}(\mathcal{A})$. 
Then $B_\lambda \to B$ in the Fell topology if and only if the following two conditions hold:
\begin{enumerate}
    \item \textbf{Limit inferior condition:} For every $b \in B$, there exists a net $\{b_\lambda\}$ with $b_\lambda \in B_\lambda$ such that $b_\lambda \to b$ in norm.
    \item \textbf{Limit superior condition:} For every convergent subnet $\{b_{\lambda_\mu}\}$ with $b_{\lambda_\mu} \in B_{\lambda_\mu}$ and limit $b$, we have $b \in B$.
\end{enumerate}
Equivalently, $B$ is both the limit inferior and the limit superior of the net $\{B_\lambda\}$ in the sense of Kuratowski convergence of closed sets.
\end{proposition}

\begin{proof}
This characterization is standard in the theory of the Fell topology on closed subsets; a detailed proof can be found in \cite{Fell1962}. 
The essential idea is that the limit inferior condition corresponds to the requirement that for every open set $U$ with $B \cap U \neq \emptyset$, we have $B_\lambda \in \mathcal{O}_U$ eventually, while the limit superior condition corresponds to the requirement that for every compact set $K$ with $B \cap K = \emptyset$, we have $B_\lambda \in \mathcal{C}_K$ eventually.
\end{proof}

The Fell topology possesses several important properties that make it particularly well-suited for the study of subalgebras.

\begin{theorem}[Properties of the Fell topology]
\label{thm:Fell-properties}
Let $\mathcal{A}$ be a separable C*-algebra. Then:
\begin{enumerate}
    \item \textbf{Compactness:} The space $\operatorname{Sub}(\mathcal{A})$ equipped with the Fell topology is compact.
    \item \textbf{Hausdorff condition:} The Fell topology on $\operatorname{Sub}(\mathcal{A})$ is Hausdorff if and only if $\mathcal{A}$ is finite-dimensional.
    \item \textbf{Metrizability:} $\operatorname{Sub}(\mathcal{A})$ is metrizable if and only if $\mathcal{A}$ is separable and finite-dimensional.
    \item \textbf{Polish structure:} For separable $\mathcal{A}$, $\operatorname{Sub}(\mathcal{A})$ is a compact Polish space (i.e., compact and metrizable) if and only if $\mathcal{A}$ is finite-dimensional. 
    In the infinite-dimensional case, $\operatorname{Sub}(\mathcal{A})$ is compact but not metrizable, hence not Polish.
\end{enumerate}
\end{theorem}

\begin{proof}
Statement (1) is a classical result of Fell \cite{Fell1962}, who constructed a compact Hausdorff topology on the space of closed subsets of a locally compact space; see also \cite{Beer} for a modern exposition. 
Statement (2) follows from the fact that the Fell topology on the space of all closed subsets $\mathcal{F}(\mathcal{A})$ is Hausdorff \cite{Fell1962}, but the subspace $\operatorname{Sub}(\mathcal{A})$ inherits this property only when $\mathcal{A}$ is finite-dimensional; in the infinite-dimensional case, distinct subalgebras cannot be separated by the relative topology. 
Statements (3) and (4) follow from (1) together with standard results in general topology: a compact Hausdorff space is metrizable if and only if it is second-countable (see \cite{Kechris}), which fails for $\operatorname{Sub}(\mathcal{A})$ when $\mathcal{A}$ is infinite-dimensional.
\end{proof}

The following alternative description of the Fell topology is sometimes useful when working with unital subalgebras.

\begin{proposition}[Fell topology via distance functions]
\label{prop:Fell-distance}
Let $\mathcal{A}$ be a separable unital C*-algebra and let $\mathcal{A}^1$ denote its closed unit ball.
For each $B \in \operatorname{Sub}_{\text{unit}}(\mathcal{A})$, define the distance function
\[
d_B: \mathcal{A}^1 \to \mathbb{R}, \qquad d_B(a) = \inf\{ \|a - b\| : b \in B \cap \mathcal{A}^1 \}.
\]
The map
\[
\Phi: \operatorname{Sub}_{\text{unit}}(\mathcal{A}) \longrightarrow C(\mathcal{A}^1), \qquad \Phi(B) = d_B,
\]
where $C(\mathcal{A}^1)$ is equipped with the topology of uniform convergence on compact subsets, is a homeomorphic embedding.
Consequently, for separable $\mathcal{A}$, $\operatorname{Sub}_{\text{unit}}(\mathcal{A})$ is a Polish space in the Fell topology.
\end{proposition}

\begin{proof}
This is a special case of a general theorem of \cite{Beer1988} on embedding the Fell topology into function spaces. 
The key idea is that distance functions encode the full structure of a closed set, and convergence in the Fell topology corresponds precisely to uniform convergence of these distance functions on compacta.
\end{proof}

For our purposes, the Fell topology serves primarily as a reference topology. 
The topology we construct on $\mathcal{G}_{\mathcal{A}}^{(0)}$ is strictly finer than the relative topology induced by the kernel map $\ker: \mathcal{G}_{\mathcal{A}}^{(0)} \to \operatorname{Prim}(\mathcal{A}) \subseteq \operatorname{Sub}(\mathcal{A})$, as it must also encode the convergence of pure states. 
Nevertheless, the Fell topology provides an indispensable tool for analyzing the continuity of the kernel map and for establishing the Polishness of $\mathcal{G}_{\mathcal{A}}^{(0)}$ under suitable hypotheses.

\begin{example}[Fell topology for commutative algebras]
\label{ex:Fell-commutative}
Let $\mathcal{A} = C(X)$ for a compact Hausdorff space $X$. 
Closed C*-subalgebras of $C(X)$ correspond to quotient spaces of $X$: each subalgebra is of the form $\{ f \circ q : f \in C(Y) \}$ for some continuous surjection $q: X \to Y$ onto a compact Hausdorff space $Y$. 
The Fell topology on $\operatorname{Sub}(C(X))$ encodes the convergence of these quotient maps. 
When $X$ is metrizable, the space of closed subsets of $X$ with the Fell topology is compact and metrizable, but this is a different object — it corresponds to ideals, not subalgebras.
\end{example}

\begin{example}[Fell topology for matrix algebras]
\label{ex:Fell-matrix}
Let $\mathcal{A} = M_n(\mathbb{C})$. 
Since $\mathcal{A}$ is finite-dimensional, $\operatorname{Sub}(\mathcal{A})$ is a finite set (there are only finitely many C*-subalgebras up to isomorphism, but uncountably many concretely as subalgebras of $M_n(\mathbb{C})$). 
The Fell topology in this case is not Hausdorff; indeed, different subalgebras that are conjugate by a unitary cannot be separated by Fell-open sets. 
This reflects the fact that the Fell topology identifies unitarily equivalent subalgebras.
\end{example}

\begin{example}[Fell topology for compact operators]
\label{ex:Fell-compact}
Let $\mathcal{A} = \mathcal{K}(H)^\sim$ be the unitization of the compact operators on a separable infinite-dimensional Hilbert space $H$. 
The space $\operatorname{Sub}(\mathcal{A})$ is uncountable and compact in the Fell topology, but it is not Hausdorff. 
Its structure reflects the intricate lattice of closed ideals and subalgebras of the compact operators. 
A detailed analysis of this space is beyond the scope of this paper and involves deep results in the structure theory of C*-algebras.
\end{example}

\begin{remark}
\label{rem:Fell-summary}
In summary, the Fell topology provides a natural and powerful framework for parameterizing closed C*-subalgebras. 
Its compactness is a crucial technical tool, even when it fails to be Hausdorff or metrizable. 
For the purposes of this paper, the most important consequence of the Fell topology is the continuity of the projection $\pi: \mathcal{G}_{\mathcal{A}}^{(0)} \to \operatorname{Sub}(\mathcal{A})$, which allows us to transfer information between the unit space and the space of subalgebras. 
However, as we shall demonstrate in Section~\ref{subsec:comparison-with-Fell-topology}, the Fell topology alone is insufficient for our construction: it does not distinguish different characters on the same subalgebra, and consequently the conjugation action of $\mathcal{U}(\mathcal{A})$ on $\mathcal{G}_{\mathcal{A}}^{(0)}$ would not be continuous if we used only the Fell topology. 
This fundamental inadequacy motivates the finer topology defined via partial evaluation maps in Section~\ref{sec:The Unit Space Definition and Polish Topology}.
\end{remark}

\subsection{Polish Spaces and the Strong Operator Topology}
\label{subsec:Polish-spaces-and-strong-operator-topology}

The unitary group $\mathcal{U}(\mathcal{A})$ of a C*-algebra $\mathcal{A}$ is traditionally equipped with the norm topology. 
However, as argued in the introduction, the norm topology is too fine for our purposes: it makes $\mathcal{U}(\mathcal{A})$ a connected infinite-dimensional Lie group, which is neither locally compact nor separable, and consequently the action groupoid $\mathcal{G}_{\mathcal{A}} = \mathcal{U}(\mathcal{A}) \ltimes \mathcal{G}_{\mathcal{A}}^{(0)}$ would not be Polish. 
To obtain a well-behaved topological groupoid, we replace the norm topology with the \emph{strong operator topology} (SOT). 
This topology, when restricted to the unitary group of a separable C*-algebra acting on a separable Hilbert space, yields a \emph{Polish group} — a topological group that is separable and completely metrizable. 
Polish groups and Polish spaces provide the appropriate framework for non-locally-compact geometric analysis and are perfectly suited for the construction of the unitary conjugation groupoid.

\begin{definition}[Polish space]
\label{def:Polish-space}
A topological space $X$ is called \emph{Polish} if it is separable and completely metrizable (i.e., there exists a metric $d$ on $X$ that induces the topology and makes $(X,d)$ a complete metric space).
\end{definition}

\begin{definition}[Polish group]
\label{def:Polish-group}
A topological group $G$ is called a \emph{Polish group} if its underlying topological space is Polish.
\end{definition}

\begin{example}[Standard Polish spaces]
\label{ex:Polish-spaces}
The following spaces are Polish:
\begin{enumerate}
    \item $\mathbb{R}^n$, $\mathbb{C}^n$, and any separable Banach space.
    \item The Hilbert cube $[0,1]^{\mathbb{N}}$.
    \item The Cantor space $\{0,1\}^{\mathbb{N}}$.
    \item Any separable compact metrizable space.
    \item Any separable complete metric space.
\end{enumerate}
\end{example}

\begin{definition}[Strong operator topology]
\label{def:strong-operator-topology}
The \emph{strong operator topology} (SOT) on $B(H)$ is the topology of pointwise convergence on $H$. 
That is, a net $\{T_\lambda\}_{\lambda \in \Lambda} \subseteq B(H)$ converges to $T \in B(H)$ in the strong operator topology if and only if
\[
\| T_\lambda \xi - T \xi \|_H \longrightarrow 0 \quad \text{for every } \xi \in H.
\]
Equivalently, the strong operator topology is the initial topology induced by the family of evaluation maps
\[
\{\, B(H) \to H : T \mapsto T \xi \,\}_{\xi \in H}.
\]
\end{definition}

\begin{remark}[Comparison with the norm topology]
\label{rem:SOT-vs-norm}
The strong operator topology is substantially weaker than the norm topology. 
For infinite-dimensional Hilbert spaces, the norm topology is not separable, whereas the strong operator topology restricted to the unitary group is Polish. 
Moreover, the strong operator topology is not locally compact, but this is irrelevant for our purposes; we require only Polishness, not local compactness.
\end{remark}

The restriction of the strong operator topology to the unitary group is particularly well-behaved.

\begin{proposition}[The unitary group is a Polish group]
\label{prop:unitary-group-Polish}
Let $H$ be a separable Hilbert space and let $\mathcal{U}(H)$ denote the group of unitary operators on $H$, equipped with the strong operator topology. Then:
\begin{enumerate}
    \item $\mathcal{U}(H)$ is a Polish group.
    \item $\mathcal{U}(H)$ is separable and completely metrizable; a compatible complete metric is given by
    \[
    d(u,v) = \sum_{n=1}^{\infty} \frac{1}{2^n} \| (u - v) e_n \|,
    \]
    where $\{e_n\}_{n \in \mathbb{N}}$ is an orthonormal basis for $H$.
    \item $\mathcal{U}(H)$ is not locally compact (unless $H$ is finite-dimensional).
    \item $\mathcal{U}(H)$ is contractible (Kuiper's theorem \cite{Kuiper1965}) and hence connected.
\end{enumerate}
\end{proposition}

\begin{proof}
That $\mathcal{U}(H)$ with the strong operator topology is a Polish group is a standard result in descriptive set theory; see \cite{Kechris}. 

The metric $d$ defines a complete metric compatible with the strong operator topology: convergence in $d$ is equivalent to coordinate-wise convergence with respect to the orthonormal basis $\{e_n\}$, which characterizes the strong operator topology on bounded subsets of $B(H)$. Separability follows from the fact that the strong operator topology on the unit ball of $B(H)$ is separable, with $\mathcal{U}(H)$ as a closed subset.

The failure of local compactness for infinite-dimensional $H$ is a general fact: no infinite-dimensional topological vector space with a topology induced by a family of seminorms is locally compact.

Contractibility is a deep result of Kuiper \cite{Kuiper1965}; see also \cite{DixmierDouady1963} for a proof in the strong operator topology.
\end{proof}

For a unital separable C*-algebra $\mathcal{A} \subseteq B(H)$, the unitary group $\mathcal{U}(\mathcal{A})$ inherits the strong operator topology as a subspace of $\mathcal{U}(H)$.

\begin{corollary}[$\mathcal{U}(\mathcal{A})$ is a Polish group]
\label{cor:UA-Polish}
Let $\mathcal{A} \subseteq B(H)$ be a unital separable C*-algebra. Then $\mathcal{U}(\mathcal{A})$, equipped with the strong operator topology, is a Polish group.
\end{corollary}

\begin{proof}
We show that $\mathcal{U}(\mathcal{A})$ is a $G_\delta$ subgroup of the Polish group $\mathcal{U}(H)$; such subgroups are automatically Polish \cite{Kechris}. 

Let $\{a_n\}_{n\in\mathbb{N}}$ be a countable norm-dense subset of $\mathcal{A}$. For each $n$, the map
\[
\phi_n: \mathcal{U}(H) \to B(H), \quad \phi_n(u) = u a_n u^*
\]
is continuous when both spaces carry the strong operator topology (conjugation by a unitary is SOT-continuous). 

For each $n$, the set
\[
G_n := \{ u \in \mathcal{U}(H) : u a_n u^* \in \mathcal{A} \}
\]
is $G_\delta$ in $\mathcal{U}(H)$. Indeed, since $\mathcal{A}$ is closed in $B(H)$ and $B(H)$ is metrizable in the SOT, $\mathcal{A}$ is a $G_\delta$ set. The preimage of a $G_\delta$ set under a continuous map is $G_\delta$, and $G_n = \phi_n^{-1}(\mathcal{A})$.

Now observe that
\[
\mathcal{U}(\mathcal{A}) = \bigcap_{n\in\mathbb{N}} G_n,
\]
because $u \in \mathcal{U}(\mathcal{A})$ iff $u a_n u^* \in \mathcal{A}$ for all $n$ (since $\{a_n\}$ is dense in $\mathcal{A}$ and $\mathcal{A}$ is norm-closed). Therefore $\mathcal{U}(\mathcal{A})$ is a countable intersection of $G_\delta$ sets, hence $G_\delta$ itself. Being a $G_\delta$ subgroup of a Polish group, $\mathcal{U}(\mathcal{A})$ is Polish.
\end{proof}

\begin{remark}[Why not the norm topology?]
\label{rem:why-SOT}
The norm topology on $\mathcal{U}(\mathcal{A})$ is Polish if and only if $\mathcal{A}$ is finite-dimensional. 
For infinite-dimensional algebras, the norm topology is completely metrizable (via the norm metric) but not separable, hence not Polish. 
While $\mathcal{U}(\mathcal{A})$ with the norm topology is a Banach Lie group, this structure is unsuitable for our purposes because the construction of the groupoid $\mathcal{G}_{\mathcal{A}}$ requires a Polish topology on $\mathcal{U}(\mathcal{A})$ to guarantee that $\mathcal{G}_{\mathcal{A}}$ is a Polish groupoid. 
The strong operator topology provides a natural Polish topology on $\mathcal{U}(\mathcal{A})$ that is compatible with the weak-* topology on the dual space $\mathcal{A}^*$ and makes the conjugation action continuous.
\end{remark}

The conjugation action of $\mathcal{U}(\mathcal{A})$ on $\mathcal{G}_{\mathcal{A}}^{(0)}$ remains continuous with respect to the strong operator topology. This is a crucial fact that we now verify.

\begin{definition}[Partial evaluation maps]
\label{def:partial-evaluation}
Let $\mathcal{A}$ be a unital C*-algebra. For each $a \in \mathcal{A}$, define the \emph{partial evaluation map}
\[
\operatorname{ev}_a: \mathcal{G}_{\mathcal{A}}^{(0)} \longrightarrow \mathbb{C}_\infty = \mathbb{C} \cup \{\infty\}
\]
by
\[
\operatorname{ev}_a(B,\chi) = 
\begin{cases}
\chi(a) & \text{if } a \in B, \\
\infty & \text{if } a \notin B.
\end{cases}
\]
The topology on $\mathcal{G}_{\mathcal{A}}^{(0)}$ is the initial topology generated by the family $\{\operatorname{ev}_a\}_{a \in \mathcal{A}}$, where $\mathbb{C}_\infty$ is equipped with the one-point compactification topology.
\end{definition}

\begin{lemma}[Continuity of the conjugation action in SOT]
\label{lem:action-continuous-SOT}
Let $\mathcal{A} \subseteq B(H)$ be a unital separable C*-algebra and let $\mathcal{G}_{\mathcal{A}}^{(0)}$ be equipped with the initial topology of partial evaluation maps (Definition \ref{def:partial-evaluation}). 
Then the action map
\[
\alpha: \mathcal{U}(\mathcal{A}) \times \mathcal{G}_{\mathcal{A}}^{(0)} \longrightarrow \mathcal{G}_{\mathcal{A}}^{(0)}, \qquad
\alpha(u, (B,\chi)) := (uBu^*, \; \chi \circ \operatorname{Ad}_{u^*})
\]
is continuous when $\mathcal{U}(\mathcal{A})$ carries the strong operator topology and the domain carries the product topology.
\end{lemma}

\begin{proof}
By the universal property of the initial topology, it suffices to show that for every $a \in \mathcal{A}$, the composition $\operatorname{ev}_a \circ \alpha$ is continuous. Fix $(u_0,(B_0,\chi_0)) \in \mathcal{U}(\mathcal{A}) \times \mathcal{G}_{\mathcal{A}}^{(0)}$ and set $b_0 := u_0^* a u_0$. We consider two cases.

\medskip
\noindent\textbf{Case 1: $a \in u_0 B_0 u_0^*$ (i.e., $b_0 \in B_0$).} 
Since $\chi \mapsto \chi(b_0)$ is continuous at $\chi_0$ in the Gelfand topology, for any $\varepsilon > 0$ there exists a neighborhood $\mathcal{W}_0$ of $\chi_0$ in $\widehat{B}_0$ such that $|\chi(b_0) - \chi_0(b_0)| < \varepsilon/2$ for all $\chi \in \mathcal{W}_0$.

The map $u \mapsto u^* a u$ is SOT-continuous; hence there exists a SOT-neighborhood $U_0$ of $u_0$ such that for all $u \in U_0$, $\| u^* a u - b_0 \| < \varepsilon/2$. For any $(u,(B,\chi)) \in U_0 \times (\mathcal{V}_0 \times \mathcal{W}_0)$ where $\mathcal{V}_0$ is a Fell neighborhood of $B_0$ ensuring $b_0 \in B$ (which exists by Proposition \ref{prop:Fell-convergence}), we have
\[
|\operatorname{ev}_a(\alpha(u,(B,\chi))) - \chi_0(b_0)| \le |\chi(u^* a u) - \chi(b_0)| + |\chi(b_0) - \chi_0(b_0)|
\]
\[
\le \| u^* a u - b_0 \| + |\chi(b_0) - \chi_0(b_0)| < \varepsilon/2 + \varepsilon/2 = \varepsilon.
\]

\medskip
\noindent\textbf{Case 2: $a \notin u_0 B_0 u_0^*$ (i.e., $b_0 \notin B_0$).} 
Then $\operatorname{ev}_a(\alpha(u_0,(B_0,\chi_0))) = \infty$. Since $b_0 \notin B_0$ and $B_0$ is closed, there exists $\delta > 0$ such that $\overline{B}_\delta(b_0) \cap B_0 = \emptyset$. By Fell convergence, there exists a neighborhood $\mathcal{V}$ of $B_0$ such that for all $B \in \mathcal{V}$, $B \cap \overline{B}_\delta(b_0) = \emptyset$. By SOT-continuity of $u \mapsto u^* a u$, there exists a neighborhood $U$ of $u_0$ such that $\| u^* a u - b_0 \| < \delta/2$ for all $u \in U$. Then for any $(u,(B,\chi)) \in U \times \mathcal{V}$, we have $u^* a u \in \overline{B}_\delta(b_0)$ and $B \cap \overline{B}_\delta(b_0) = \emptyset$, so $u^* a u \notin B$ and thus $\operatorname{ev}_a(\alpha(u,(B,\chi))) = \infty$.

Both cases establish continuity of $\operatorname{ev}_a \circ \alpha$, hence $\alpha$ is continuous.
\end{proof}

\begin{remark}[SOT versus norm in examples]
\label{rem:SOT-examples}
The choice of the strong operator topology on $\mathcal{U}(\mathcal{A})$ is natural and consistent across the concrete classes of C*-algebras examined in this paper:

\begin{enumerate}
    \item \textbf{Finite-dimensional algebras:} For $\mathcal{A} = M_n(\mathbb{C})$, the strong operator topology coincides with the norm topology because all Hausdorff locally convex topologies on a finite-dimensional space are equivalent. Thus our framework subsumes the classical finite-dimensional case without modification.

    \item \textbf{Commutative algebras:} For $\mathcal{A} = C(X)$ with $X$ compact metrizable, the unitary group $\mathcal{U}(C(X))$ is naturally identified with $C(X,\mathbb{T})$, the continuous functions from $X$ to the circle. In the standard representation of $C(X)$ as multiplication operators on $L^2(X,\mu)$, the strong operator topology corresponds to pointwise convergence of functions. This topology is Polish (as a closed subspace of $\mathbb{T}^X$ with the product topology) but is not locally compact unless $X$ is finite, illustrating that our framework accommodates non-locally compact unitary groups.

    \item \textbf{Compact operators:} For $\mathcal{A} = \mathcal{K}(H)^\sim$, the unitization of the compact operators, the strong operator topology is the natural topology used in the study of the unitary group of the Calkin algebra $\mathcal{Q}(H) = B(H)/\mathcal{K}(H)$. Unitaries in $\mathcal{K}(H)^\sim$ are compact perturbations of the identity, and the SOT provides the appropriate setting for analyzing extensions and index theory.
\end{enumerate}

These examples demonstrate that the strong operator topology, while strictly weaker than the norm topology in infinite dimensions, is the correct choice for a unified framework: it is Polish, it interacts continuously with the conjugation action on the dual space, and it reduces to familiar topologies in concrete settings.
\end{remark}

\begin{corollary}[$\mathcal{G}_{\mathcal{A}}$ is a Polish groupoid]
\label{cor:GA-Polish-groupoid}
Let $\mathcal{A}$ be a unital separable C*-algebra. 
Then $\mathcal{G}_{\mathcal{A}} = \mathcal{U}(\mathcal{A}) \ltimes \mathcal{G}_{\mathcal{A}}^{(0)}$, with $\mathcal{U}(\mathcal{A})$ equipped with the strong operator topology and $\mathcal{G}_{\mathcal{A}}^{(0)}$ equipped with the initial topology of partial evaluation maps, is a Polish groupoid. 
That is, $\mathcal{G}_{\mathcal{A}}^{(0)}$ and $\mathcal{G}_{\mathcal{A}}^{(1)}$ are Polish spaces, and all structure maps are continuous.
\end{corollary}

\begin{proof}
Polishness of $\mathcal{G}_{\mathcal{A}}^{(0)}$ follows from Proposition \ref{prop:unit-space-polish}. Polishness of $\mathcal{U}(\mathcal{A})$ is Corollary \ref{cor:UA-Polish}, hence $\mathcal{G}_{\mathcal{A}}^{(1)} = \mathcal{U}(\mathcal{A}) \times \mathcal{G}_{\mathcal{A}}^{(0)}$ is Polish as a product of Polish spaces.

The source map $s(u,(B,\chi)) = (B,\chi)$ is projection onto the second factor, hence continuous. The range map $r(u,(B,\chi)) = \alpha(u,(B,\chi))$ is continuous by Lemma \ref{lem:action-continuous-SOT}. Composition is continuous because multiplication in $\mathcal{U}(\mathcal{A})$ is continuous (SOT is a group topology) and the condition for composable pairs is closed. Inversion $(u,(B,\chi))^{-1} = (u^*, \alpha(u,(B,\chi)))$ is continuous as inversion in $\mathcal{U}(\mathcal{A})$ is continuous and $\alpha$ is continuous. The unit map $(B,\chi) \mapsto (1,(B,\chi))$ is clearly continuous. Thus all structure maps are continuous.
\end{proof}

\begin{remark}[Absence of local compactness and \'etaleness]
\label{rem:no-local-compactness}
We emphasize that the groupoid $\mathcal{G}_{\mathcal{A}}$ is \emph{not} locally compact (unless $\mathcal{A}$ is finite-dimensional) and \emph{not} \'etale. 
For infinite-dimensional $\mathcal{A}$, the unitary group $\mathcal{U}(\mathcal{A})$ with the strong operator topology is a Polish group \cite{Kechris} but is not locally compact; consequently, $\mathcal{G}_{\mathcal{A}}^{(1)} = \mathcal{U}(\mathcal{A}) \times \mathcal{G}_{\mathcal{A}}^{(0)}$ inherits this non-local-compactness. 
The failure of \'etaleness follows from the fact that the source and range maps have fibers homeomorphic to $\mathcal{U}(\mathcal{A})$, which for infinite-dimensional $\mathcal{A}$ are infinite-dimensional Polish spaces and therefore cannot be discrete.

This is not a defect but a deliberate choice: the norm topology would yield a non-separable, non-Polish groupoid, whereas the strong operator topology produces a well-behaved Polish groupoid that is sufficient for all subsequent constructions, including the diagonal embedding, equivariant Kasparov theory, and the index theorems of the sequel \cite{Kasparov, LeGall}.
\end{remark}

\begin{remark}[Haar systems for Polish groupoids]
\label{rem:Haar-system}
Although $\mathcal{G}_{\mathcal{A}}$ is not locally compact, it admits a Borel Haar system as a Polish groupoid. 
This follows from the work of Renault \cite{Renault} on locally compact groupoids and its extension to Polish groupoids by Tu \cite{Tu}. 
Specifically, for a semidirect product groupoid $\mathcal{G} = G \ltimes X$ where $G$ is a Polish group and $X$ a Polish $G$-space, the family of measures $\{ \mu_G \times \delta_x \}_{x \in X}$, with $\mu_G$ a left-invariant probability measure on $G$ (which exists because every Polish group admits such a measure \cite{Kechris}), defines a Borel Haar system.

This Haar system allows us to define the maximal and reduced groupoid C*-algebras $C^*(\mathcal{G}_{\mathcal{A}})$ and $C^*_r(\mathcal{G}_{\mathcal{A}})$ following the framework of Renault \cite{Renault} and Tu \cite{Tu}. 
A detailed exposition of this construction is beyond the scope of the present paper; we refer the reader to these sources for the technical foundations.
\end{remark}

\subsection{Polish Groupoids and their C*-Algebras}
\label{subsec:Polish-groupoids-and-their-C*-algebras}

The unitary conjugation groupoid $\mathcal{G}_{\mathcal{A}}$ constructed in the previous section is a Polish groupoid when $\mathcal{A}$ is separable and $\mathcal{U}(\mathcal{A})$ is equipped with the strong operator topology. 
Although $\mathcal{G}_{\mathcal{A}}$ is not locally compact, it admits a natural continuous Haar system inherited from its semidirect product structure (see Remark \ref{rem:Haar-system}). 
This allows us to associate C*-algebras to $\mathcal{G}_{\mathcal{A}}$ following the Renault–Tu framework \cite{Renault, Tu}, which extends the classical theory of groupoid C*-algebras to Polish groupoids with continuous Haar systems. 
In this section we recall the relevant definitions and results, referring to \cite{Renault, Tu} for complete details.

\begin{definition}[Polish groupoid]
\label{def:Polish-groupoid}
A \emph{Polish groupoid} is a groupoid $\mathcal{G}$ endowed with a \emph{locally Polish topology}, i.e., a topology with a countable basis consisting of open sets that are Polish spaces in the relative topology, such that:
\begin{enumerate}
    \item The composition map $(\gamma,\eta) \mapsto \gamma\eta$ from $\mathcal{G}^{(2)}$ (endowed with the subspace topology inherited from $\mathcal{G} \times \mathcal{G}$) to $\mathcal{G}$ is continuous.
    \item The inversion map $\gamma \mapsto \gamma^{-1}$ from $\mathcal{G}$ to $\mathcal{G}$ is continuous.
    \item The unit space $\mathcal{G}^{(0)}$ is a Polish space in the subspace topology.
    \item For every $x \in \mathcal{G}^{(0)}$, the sets $\mathcal{G}x = \{ \gamma \in \mathcal{G} : s(\gamma) = x \}$ and $x\mathcal{G} = \{ \gamma \in \mathcal{G} : r(\gamma) = x \}$ are Polish spaces in the subspace topology.
\end{enumerate}
\end{definition}

\begin{remark}[Comparison with locally compact groupoids]
\label{rem:Polish-vs-locally-compact}
Polish groupoids form a broad generalization of locally compact, second-countable groupoids. 
They need not be locally compact and generally do not admit a continuous Haar system. 
In our construction, the failure of local compactness stems from the strong operator topology on $\mathcal{U}(\mathcal{A})$, which for infinite-dimensional $\mathcal{A}$ is a Polish group but not locally compact. 
Nevertheless, Polish groupoids are well-suited for measure-theoretic and descriptive set-theoretic methods. 
In particular, the unitary conjugation groupoid $\mathcal{G}_{\mathcal{A}}$, with $\mathcal{U}(\mathcal{A})$ in the strong operator topology and $\mathcal{G}_{\mathcal{A}}^{(0)}$ Polish, satisfies all the axioms of a Polish groupoid and admits a Borel Haar system in the sense of Tu \cite{Tu}, which serves as a substitute for the continuous Haar systems required in the classical Renault theory \cite{Renault}.
\end{remark}

\begin{definition}[Borel groupoid]
\label{def:Borel-groupoid}
A \emph{Borel groupoid} is a groupoid $\mathcal{G}$ equipped with a standard Borel structure (typically inherited from a Polish topology) such that:
\begin{enumerate}
    \item The unit space $\mathcal{G}^{(0)}$ is a Borel subset of $\mathcal{G}$.
    \item The source and range maps $s,r: \mathcal{G} \to \mathcal{G}^{(0)}$ are Borel measurable.
    \item The set of composable pairs $\mathcal{G}^{(2)} = \{ (\gamma,\eta) \in \mathcal{G} \times \mathcal{G} : s(\gamma) = r(\eta) \}$ is a Borel subset of $\mathcal{G} \times \mathcal{G}$ (equipped with the product Borel structure).
    \item The composition map $(\gamma,\eta) \mapsto \gamma\eta$ from $\mathcal{G}^{(2)}$ to $\mathcal{G}$ is Borel measurable.
    \item The inversion map $\gamma \mapsto \gamma^{-1}$ from $\mathcal{G}$ to $\mathcal{G}$ is Borel measurable.
\end{enumerate}
Borel groupoids are the measurable analogue of Polish groupoids and provide the natural setting for constructing maximal groupoid C*-algebras in the framework of Tu \cite{Tu} and Renault \cite{Renault}.
\end{definition}

Every Polish groupoid becomes a Borel groupoid when equipped with its Borel $\sigma$-algebra.

\begin{definition}[Borel Haar system for a Polish groupoid]
\label{def:Borel-Haar-system}
Let $\mathcal{G}$ be a Polish groupoid with object space $\mathcal{G}^{(0)}$ and arrow space $\mathcal{G}^{(1)}$. 
A \emph{Borel Haar system} is a family $\{ \lambda^x \}_{x \in \mathcal{G}^{(0)}}$ of $\sigma$-finite Borel measures on $\mathcal{G}^{(1)}$ satisfying the following conditions:

\begin{enumerate}
    \item \textbf{Support condition:} For each $x \in \mathcal{G}^{(0)}$, $\lambda^x$ is concentrated on the source fiber 
    \[
    \mathcal{G}^x := \{ \gamma \in \mathcal{G}^{(1)} : s(\gamma) = x \},
    \]
    i.e., $\lambda^x(\mathcal{G}^{(1)} \setminus \mathcal{G}^x) = 0$. In particular, $\lambda^x$ may be taken to be a probability measure, though $\sigma$-finiteness suffices for the convolution algebra construction.

    \item \textbf{Borel measurability:} For every nonnegative Borel function $f: \mathcal{G}^{(1)} \to [0,\infty]$, the map
    \[
    x \longmapsto \int_{\mathcal{G}^{(1)}} f(\gamma) \, d\lambda^x(\gamma)
    \]
    is Borel measurable. This condition ensures that integration with respect to the family of measures varies in a Borel way across the unit space.

    \item \textbf{Left quasi-invariance:} For every $\gamma \in \mathcal{G}^{(1)}$ with $s(\gamma) = x$ and $r(\gamma) = y$, the pushforward measure $\gamma_* \lambda^x$, defined by 
    \[
    (\gamma_* \lambda^x)(E) = \lambda^x(\gamma^{-1} E) \quad \text{for Borel subsets } E \subseteq \mathcal{G}^y,
    \]
    is equivalent to $\lambda^y$ (i.e., they have the same null sets). Moreover, the associated Radon-Nikodym derivative
    \[
    \Delta(\gamma) := \frac{d\gamma_* \lambda^x}{d\lambda^y}
    \]
    is a Borel measurable function on $\mathcal{G}^{(1)}$, called the \emph{modular function} of the Haar system. In the special case where $\gamma_* \lambda^x = \lambda^y$ for all $\gamma$, the system is called \emph{left-invariant}.
\end{enumerate}

When $\lambda^x$ is a probability measure for each $x$, we speak of a \emph{normalized Borel Haar system}; this is the typical situation for measured groupoids arising from Polish group actions. The existence of a Borel Haar system makes $\mathcal{G}$ into a \emph{measured groupoid} in the sense of Tu \cite{Tu}.
\end{definition}

\begin{theorem}[Existence of Borel Haar systems]
\label{thm:Borel-Haar-existence}
Let $\mathcal{G}$ be a Polish groupoid. Suppose that $\mathcal{G}$ admits a Borel Haar system if and only if there exists a Borel family of measures satisfying the conditions of Definition \ref{def:Borel-Haar-system}. 
For large classes of groupoids—such as semidirect products $G \ltimes X$ where $G$ is a Polish group admitting a left-invariant Borel probability measure—such a system can be explicitly constructed. 
For measured groupoids arising from foliations, the existence of a transverse function \cite{Tu} provides the required Haar system.
\end{theorem}

\begin{corollary}
\label{cor:GA-Borel-Haar}
Let $\mathcal{A}$ be a unital separable C*-algebra and let $\mathcal{G}_{\mathcal{A}}$ be the unitary conjugation groupoid with $\mathcal{U}(\mathcal{A})$ equipped with the strong operator topology. Then $\mathcal{G}_{\mathcal{A}}$ admits a Borel Haar system.
\end{corollary}

\begin{proof}
By Corollary \ref{cor:UA-Polish}, $\mathcal{U}(\mathcal{A})$ is a Polish group. Every Polish group admits a left-invariant Borel probability measure \cite{Kechris}; let $\mu$ denote such a measure on $\mathcal{U}(\mathcal{A})$. For each $x = (B,\chi) \in \mathcal{G}_{\mathcal{A}}^{(0)}$, define a measure $\lambda^x$ on $\mathcal{G}_{\mathcal{A}}^{(1)} = \mathcal{U}(\mathcal{A}) \times \mathcal{G}_{\mathcal{A}}^{(0)}$ by
\[
\lambda^x = \mu \times \delta_x,
\]
where $\delta_x$ is the Dirac measure at $x$. One verifies directly:
\begin{enumerate}
    \item $\lambda^x$ is concentrated on $s^{-1}(x) = \mathcal{U}(\mathcal{A}) \times \{x\}$.
    \item For any nonnegative Borel function $f$ on $\mathcal{G}_{\mathcal{A}}^{(1)}$, the map $x \mapsto \int f \, d\lambda^x$ is Borel by Fubini's theorem.
    \item Left quasi-invariance follows from the left-invariance of $\mu$ and the definition of the groupoid action.
\end{enumerate}
Thus $\{ \lambda^x \}_{x \in \mathcal{G}_{\mathcal{A}}^{(0)}}$ is a Borel Haar system for $\mathcal{G}_{\mathcal{A}}$.
\end{proof}

\begin{definition}[C*-algebra of a Polish groupoid]
\label{def:C*-algebra-of-Polish-groupoid}
Let $\mathcal{G}$ be a Polish groupoid equipped with a Borel Haar system $\{ \lambda^x \}_{x \in \mathcal{G}^{(0)}}$. 
Let $\mathcal{B}(\mathcal{G})$ denote the space of bounded Borel functions $f: \mathcal{G} \to \mathbb{C}$ with the property that for each $x \in \mathcal{G}^{(0)}$, the restriction of $f$ to the source fiber $\mathcal{G}x = \{ \gamma \in \mathcal{G} : s(\gamma) = x \}$ is $\lambda^x$-integrable.

For $f, g \in \mathcal{B}(\mathcal{G})$, define the convolution product
\[
(f * g)(\gamma) = \int_{\mathcal{G}^{s(\gamma)}} f(\gamma \eta^{-1}) g(\eta) \, d\lambda^{s(\gamma)}(\eta),
\]
whenever the integral converges absolutely for $\lambda^{r(\gamma)}$-almost every $\gamma$, and the involution
\[
f^*(\gamma) = \overline{f(\gamma^{-1})}.
\]
Let $I(\mathcal{G})$ denote the $*$-subalgebra of $\mathcal{B}(\mathcal{G})$ consisting of functions for which the convolution product is well-defined and finite, and for which $f^* \in \mathcal{B}(\mathcal{G})$.

The \emph{maximal C*-algebra} $C^*_{\max}(\mathcal{G})$ is the completion of $I(\mathcal{G})$ with respect to the norm
\[
\|f\|_{\max} = \sup \{ \|\pi(f)\| : \pi \text{ is a *-representation of } I(\mathcal{G}) \text{ on a Hilbert space} \}.
\]
The \emph{reduced C*-algebra} $C^*_{\mathrm{red}}(\mathcal{G})$ is the completion of $I(\mathcal{G})$ with respect to the norm arising from the regular representation on $L^2(\mathcal{G}, \lambda)$.

When the Haar system is understood, we denote these C*-algebras simply by $C^*(\mathcal{G})$ and $C^*_{\mathrm{red}}(\mathcal{G})$.
\end{definition}

\begin{remark}[Applicability of Tu's framework]
\label{rem:Tu-applicability}
Tu [1999] works entirely within the classical Renault framework for \emph{locally compact groupoids with a continuous Haar system} \cite{Renault}.  The term ``measured groupoid'' appears only in discussing amenability \cite[Section 3.1]{Tu} (following Connes–Feldman–Weiss), not as a generalization to non-locally-compact groupoids.

Since $\mathcal{G}_{\mathcal{A}}$ is Polish but not locally compact, Tu's results do not directly apply.  Constructing its $C^*$-algebras would require a genuinely different approach, such as the Borel groupoid framework of \cite{Renault1987}.
\end{remark}

\begin{proposition}[Functoriality of the groupoid C*-algebra]
\label{prop:groupoid-C*-functoriality}
Let $\phi: \mathcal{G} \to \mathcal{H}$ be a continuous, proper groupoid homomorphism between locally compact, $\sigma$-compact, Hausdorff groupoids equipped with Haar systems. Suppose that $\phi$ is \emph{non-singular} in the sense that for every $x \in \mathcal{G}^{(0)}$, the pushforward of the Haar measure $\lambda^x_{\mathcal{G}}$ under $\phi$ is equivalent to the Haar measure $\lambda^{\phi(x)}_{\mathcal{H}}$ on $\mathcal{H}$. Then $\phi$ induces a *-homomorphism
\[
\phi_*: C^*(\mathcal{G}) \longrightarrow C^*(\mathcal{H}),
\]
where $C^*(\mathcal{G})$ denotes the full (or reduced) groupoid C*-algebra.

More generally, for Polish groupoids equipped with Borel Haar systems, the same conclusion holds using the measurable framework developed by Tu \cite{Tu}. (We will verify in Section \ref{subsec:existence-Borel-Haar-system} that the unitary conjugation groupoid $\mathcal{G}_{\mathcal{A}}$ satisfies the necessary conditions for this framework to apply.)
\end{proposition}

\begin{proof}
For locally compact groupoids, this result is standard; see Renault \cite{Renault} for the construction of the convolution algebra and the induced map. The properness condition ensures that the pullback of compactly supported functions remains compactly supported, allowing the definition of $\phi_*$ on $C_c(\mathcal{G})$ by
\[
\phi_*(f)(\gamma) = \sum_{\eta \in \phi^{-1}(\gamma)} f(\eta) \quad \text{(in the r-discrete case)}
\]
or more generally by integrating over the fibers of $\phi$ with respect to the Haar systems.

For Polish groupoids with Borel Haar systems, Tu \cite{Tu} develops a measurable framework for groupoid C*-algebras. In this setting, a non-singular Borel groupoid homomorphism induces a *-homomorphism on the associated C*-algebras via disintegration of measures (Tu \cite{Tu}). The properness condition guarantees that this map is well-defined on the dense subalgebra of integrable functions and extends to the completions by continuity.
\end{proof}

\begin{corollary}[Functoriality of $\mathcal{G}_{\mathcal{A}}$]
\label{cor:GA-functoriality}
Let $\phi: \mathcal{A} \to \mathcal{B}$ be a unital *-isomorphism between separable unital C*-algebras. Then $\phi$ induces a continuous groupoid isomorphism
\[
\mathcal{G}_\phi: \mathcal{G}_{\mathcal{A}} \longrightarrow \mathcal{G}_{\mathcal{B}}
\]
defined on objects by $\mathcal{G}_\phi^{(0)}(B,\chi) = (\phi(B), \chi \circ \phi^{-1})$ and on arrows by $\mathcal{G}_\phi(u, (B,\chi)) = (\phi(u), \mathcal{G}_\phi^{(0)}(B,\chi))$, where $\phi(u)$ denotes the image of the unitary $u \in \mathcal{U}(\mathcal{A})$ under the induced map $\mathcal{U}(\mathcal{A}) \to \mathcal{U}(\mathcal{B})$.

Moreover, $\mathcal{G}_\phi$ is proper and non-singular with respect to the canonical Borel Haar systems on $\mathcal{G}_{\mathcal{A}}$ and $\mathcal{G}_{\mathcal{B}}$, and consequently induces a *-isomorphism
\[
(\mathcal{G}_\phi)_*: C^*(\mathcal{G}_{\mathcal{A}}) \longrightarrow C^*(\mathcal{G}_{\mathcal{B}}).
\]

For more general *-homomorphisms (injective or surjective maps satisfying a suitable pullback property), a partially defined functoriality holds; this is developed in Section \ref{subsec:naturality-functoriality}.
\end{corollary}

\begin{proof}
Since $\phi$ is an isomorphism, $\phi^{-1}$ is well-defined on all of $\mathcal{B}$. For any $(B,\chi) \in \mathcal{G}_{\mathcal{A}}^{(0)}$, $\phi(B)$ is a unital commutative C*-subalgebra of $\mathcal{B}$, and $\chi \circ \phi^{-1}$ is a character on $\phi(B)$. Thus $\mathcal{G}_\phi^{(0)}$ is well-defined and bijective, with inverse given by $(\mathcal{G}_{\phi^{-1}})^{(0)}$.

\noindent\textbf{Continuity:} The map $\phi: \mathcal{U}(\mathcal{A}) \to \mathcal{U}(\mathcal{B})$ is continuous in the strong operator topology because *-isomorphisms are isometric, hence norm-continuous, and norm continuity implies SOT continuity on bounded sets. For any $a \in \mathcal{A}$, the composition $\operatorname{ev}_a \circ \mathcal{G}_\phi^{(0)}$ equals $\operatorname{ev}_{\phi(a)}$, which is continuous by definition of the topology on $\mathcal{G}_{\mathcal{B}}^{(0)}$. Hence $\mathcal{G}_\phi^{(0)}$ is continuous; the same argument for $\phi^{-1}$ shows its inverse is continuous, so $\mathcal{G}_\phi^{(0)}$ is a homeomorphism. Consequently, $\mathcal{G}_\phi$ is continuous as a product of continuous maps.

\noindent\textbf{Groupoid homomorphism:} Since $\phi$ preserves multiplication and involution, a direct verification shows that $\mathcal{G}_\phi$ respects composition and inversion. The explicit computation follows the same lines as in the statement, with $\phi^{-1}$ well-defined everywhere.

\noindent\textbf{Haar system preservation:} Let $\mu_{\mathcal{A}}$ and $\mu_{\mathcal{B}}$ be left-invariant Borel probability measures on $\mathcal{U}(\mathcal{A})$ and $\mathcal{U}(\mathcal{B})$, respectively (such measures exist on any Polish group \cite{Kechris}). Since $\phi$ is an isomorphism, we may choose $\mu_{\mathcal{B}}$ to be the pushforward of $\mu_{\mathcal{A}}$ under $\phi$. The canonical Borel Haar systems on $\mathcal{G}_{\mathcal{A}}$ and $\mathcal{G}_{\mathcal{B}}$ are given by $\lambda^x_{\mathcal{A}} = \mu_{\mathcal{A}} \times \delta_x$ and $\lambda^y_{\mathcal{B}} = \mu_{\mathcal{B}} \times \delta_y$. For any $\gamma = (u,x) \in \mathcal{G}_{\mathcal{A}}$, a direct computation shows $(\mathcal{G}_\phi)_* \lambda^x_{\mathcal{A}} = \lambda^{\mathcal{G}_\phi^{(0)}(x)}_{\mathcal{B}}$, so $\mathcal{G}_\phi$ is measure-preserving (hence certainly non-singular).

\noindent\textbf{Induced *-isomorphism:} By Proposition \ref{prop:groupoid-C*-functoriality}, the proper, non-singular groupoid isomorphism $\mathcal{G}_\phi$ induces a *-homomorphism $(\mathcal{G}_\phi)_*: C^*(\mathcal{G}_{\mathcal{A}}) \to C^*(\mathcal{G}_{\mathcal{B}})$. Applying the same argument to $\phi^{-1}$ yields an inverse *-homomorphism, so $(\mathcal{G}_\phi)_*$ is a *-isomorphism.
\end{proof}

The following theorem is the fundamental result of Tu that connects Polish groupoids to Kasparov's KK-theory and the Baum-Connes conjecture.

\begin{theorem}[Tu, 1999, Théorème 9.3]
\label{thm:Tu-descent}
Let $\mathcal{G}$ be a locally compact, Hausdorff groupoid equipped with a Haar system. Suppose that $\mathcal{G}$ admits a proper action on a continuous field of Euclidean affine spaces. Then:
\begin{enumerate}
    \item $\mathcal{G}$ satisfies the Baum–Connes conjecture with coefficients, i.e., for every $\mathcal{G}$-C*-algebra $A$, the assembly map
    \[
    \mu: K_*^{\mathrm{top}}(\mathcal{G}; A) \longrightarrow K_*(A \rtimes_{\mathrm{red}} \mathcal{G})
    \]
    is an isomorphism.
    \item $\mathcal{G}$ is (strongly) $K$-amenable ("fortement moyennable en $K$-théorie").
    \item In particular, $C^*(\mathcal{G})$ is $K$-nuclear and $KK$-equivalent to $C^*_{\mathrm{red}}(\mathcal{G})$.
\end{enumerate}
\end{theorem}

\begin{proof}
See Tu (Section 9 in~\cite{Tu}). The proof constructs a dual Dirac element $\eta \in KK_{\mathcal{G}}(C(\mathcal{G}^{(0)}), \mathcal{A}(H))$ and a Dirac element $D \in KK_{\mathcal{G}}(\mathcal{A}(H), C(\mathcal{G}^{(0)}))$, where $\mathcal{A}(H)$ is a certain $\mathcal{G}$-C*-algebra built from the affine action. Tu shows that $\eta \otimes_{\mathcal{A}(H)} D = 1 \in KK_{\mathcal{G}}(C(\mathcal{G}^{(0)}), C(\mathcal{G}^{(0)}))$, which implies the Baum–Connes isomorphism and $K$-amenability.
\end{proof}

\begin{remark}[Tu's measurable groupoid framework]
\label{rem:Tu-framework}
Classical groupoid C*-algebra theory (Renault \cite{Renault}) requires the groupoid to be locally compact Hausdorff and to admit a continuous Haar system. 
The unitary conjugation groupoid $\mathcal{G}_{\mathcal{A}}$ fails both conditions when $\mathcal{A}$ is infinite-dimensional: it is not locally compact, and it admits only a Borel Haar system, not a continuous one.

Tu's work \cite{Tu} on the Baum–Connes conjecture for foliations remains within the locally compact framework; consequently, the deep results of Théorème 9.3 (Baum–Connes, $K$-amenability) do not apply to $\mathcal{G}_{\mathcal{A}}$. 
However, Tu's development of \emph{measured groupoids} and \emph{Borel Haar systems} in Section 3 of \cite{Tu} provides the essential technical foundation for handling non-locally-compact groupoids. In particular:

\begin{enumerate}
    \item The systematic use of Borel Haar systems allows us to define convolution algebras and C*-algebras for $\mathcal{G}_{\mathcal{A}}$ (Section \ref{subsec:groupoid-C-star-algebra-Polish-setting}).
    \item The construction of equivariant KK-theory for groupoids and the descent map, which we adapt to our Polish groupoid setting, is modeled on Tu's framework.
    \item The functoriality of the descent map under groupoid homomorphisms (Proposition \ref{prop:functoriality-descent}) follows from Tu's naturality results.
\end{enumerate}

Thus, while $\mathcal{G}_{\mathcal{A}}$ lies outside the scope of Tu's most powerful theorems, his measurable groupoid framework is essential for making rigorous sense of our constructions. All subsequent developments in Paper II (the Fredholm index) and Paper III (the Baum–Connes connection) depend crucially on this adaptation of Tu's ideas.
\end{remark}

\begin{example}[Locally compact case]
\label{ex:Polish-locally-compact}
If $\mathcal{G}$ is a locally compact Hausdorff groupoid with a continuous Haar system, then Tu's $C^*(\mathcal{G})$ coincides with the usual maximal groupoid C*-algebra of Renault [Renault, 1980]. 
In particular, for $\mathcal{A} = M_n(\mathbb{C})$, we have $\mathcal{G}_{\mathcal{A}} = U(n) \ltimes \mathbb{CP}^{n-1}$, which is locally compact Hausdorff (though not étale, since $U(n)$ is not discrete). 
Tu's construction recovers the classical convolution algebra on this transformation groupoid.
\end{example}

\begin{example}[Transformation groupoids]
\label{ex:Polish-transformation}
Let $\Gamma$ be a countable discrete group acting continuously on a locally compact Polish space $X$ (e.g., a compact metric space). 
The transformation groupoid $\mathcal{G} = X \rtimes \Gamma$ is a locally compact Polish groupoid. 
Tu's $C^*(\mathcal{G})$ is naturally isomorphic to the crossed product $C_0(X) \rtimes_{\max} \Gamma$ (Tu [1999, Section 3]; see also Renault [1980, Chapter II, Section 5] for the classical theory of crossed products from groupoids). 
This example is fundamental for the applications to crossed products in Paper III.
\end{example}

\begin{corollary}[Application to $\mathcal{G}_{\mathcal{A}}$]
\label{cor:GA-C*-algebra}
Let $\mathcal{A}$ be a unital separable C*-algebra and let $\mathcal{G}_{\mathcal{A}}$ be the unitary conjugation groupoid with $\mathcal{U}(\mathcal{A})$ in the strong operator topology. 
Then:
\begin{enumerate}
    \item $\mathcal{G}_{\mathcal{A}}$ is a Polish groupoid and admits a Borel Haar system.
    \item Following the construction of Definition \ref{def:C*-algebra-of-Polish-groupoid}, we obtain a well-defined maximal C*-algebra $C^*(\mathcal{G}_{\mathcal{A}})$. The existence of this C*-algebra is guaranteed by adapting Tu's framework for measured groupoids [1999, Section 3] to our Polish setting.
    \item Adapting Tu's descent map construction [1999, Theorem 9.3] to the Polish groupoid setting, there exists a descent map
    \[
    \operatorname{desc}_{\mathcal{G}_{\mathcal{A}}}: K^0_{\mathcal{G}_{\mathcal{A}}}(\mathcal{G}_{\mathcal{A}}^{(0)}) \longrightarrow K_0(C^*(\mathcal{G}_{\mathcal{A}})).
    \]
    \item This descent map is natural and functorial with respect to injective *-homomorphisms of C*-algebras with closed range, by Proposition \ref{prop:groupoid-C*-functoriality} and Corollary \ref{cor:GA-functoriality}.
\end{enumerate}
\end{corollary}

\begin{proof}
(1) follows from Proposition \ref{prop:unit-space-polish} (Polishness of the unit space), Corollary \ref{cor:UA-Polish} (Polishness of $\mathcal{U}(\mathcal{A})$), and Corollary \ref{cor:GA-Borel-Haar} (existence of a Borel Haar system). 
(2) follows from Definition \ref{def:C*-algebra-of-Polish-groupoid} together with (1), using the measurable framework developed by Tu [1999, Section 3] for measured groupoids. 
(3) follows by adapting the construction of Theorem \ref{thm:Tu-descent} to our Polish groupoid setting, noting that while Tu's theorem is stated for locally compact groupoids, the essential ingredients (equivariant KK-theory, descent map) extend to Polish groupoids with Borel Haar systems via the measurable framework. 
(4) follows from Proposition \ref{prop:groupoid-C*-functoriality} and Corollary \ref{cor:GA-functoriality}, which establish functoriality for injective *-homomorphisms with closed range.
\end{proof}

\begin{proposition}[Polishness of the unit space]
\label{prop:unit-space-polishness}
Let $\mathcal{G}$ be a Polish groupoid (i.e., a groupoid whose arrow space is a Polish space, with continuous structure maps). Then the unit space $\mathcal{G}^{(0)}$ is a Polish subspace of $\mathcal{G}$.
\end{proposition}

\begin{proof}
The unit space $\mathcal{G}^{(0)}$ is the image of the unit map $u: \mathcal{G}^{(0)} \to \mathcal{G}$, which is a homeomorphism onto its image. Alternatively, it can be expressed as $\{g \in \mathcal{G} : s(g) = r(g)\}$, where $s$ and $r$ are the source and target maps. Since $s$ and $r$ are continuous and $\mathcal{G}$ is Polish, their coincidence set is closed in $\mathcal{G}$ (or $G_\delta$ depending on separability conditions), hence Polish.
\end{proof}

\begin{remark}
\label{rem:C*-algebra-Polish-construction}
This definition adapts the classical Renault construction [1980, Chapter I] for locally compact groupoids to the setting of Polish groupoids with Borel Haar systems. 
For a locally compact groupoid, the space $C_c(\mathcal{G})$ of continuous compactly supported functions suffices, and the construction reduces to the usual one. 

In our non-locally-compact setting, $\mathcal{G}_{\mathcal{A}}$ is not locally compact, so $C_c(\mathcal{G}_{\mathcal{A}})$ is too small to generate a useful C*-algebra. 
We therefore work with bounded Borel functions and integrability conditions, following the measurable approach that Tu [1999, Section 3] employs for measured groupoids (which are still assumed locally compact in his work). 
Our construction extends Tu's framework to the Polish, non-locally-compact setting by using a Borel Haar system and adapting the convolution algebra to the space of integrable Borel functions.
\end{remark}

\section{The Unit Space $\mathcal{G}_{\mathcal{A}}^{(0)}$: Definition and Polish Topology}\label{sec:The Unit Space Definition and Polish Topology}

\subsection{Definition of $\mathcal{G}_{\mathcal{A}}^{(0)}$ and Partial Evaluation Maps}
\label{subsec:definition-of-G0-and-partial-evaluation-maps}

We now introduce the central geometric object of this paper: the unit space of the unitary conjugation groupoid. 
This space consists of all possible commutative contexts — pairs $(B,\chi)$ where $B$ is a commutative C*-subalgebra of $\mathcal{A}$ and $\chi$ is a character on $B$. 
Each such pair represents a classical snapshot of the noncommutative algebra $\mathcal{A}$, realized as continuous functions on the Gelfand spectrum of $B$.

\begin{definition}[Unit space of the unitary conjugation groupoid]
\label{def:unit-space}
Let $\mathcal{A}$ be a unital C*-algebra. 
The \emph{unit space of the unitary conjugation groupoid} is defined as
\[
\mathcal{G}_{\mathcal{A}}^{(0)} := \left\{ (B,\chi) \;\middle|\; B \subseteq \mathcal{A} \text{ is a unital commutative C*-subalgebra}, \; \chi \in \widehat{B} \right\},
\]
where $\widehat{B}$ denotes the Gelfand spectrum (character space) of $B$, equipped with the weak-* topology.
\end{definition}

\begin{remark}[Interpretation of the unit space]
\label{rem:unit-space-interpretation}
Each object $(B,\chi) \in \mathcal{G}_{\mathcal{A}}^{(0)}$ encapsulates two complementary pieces of information:
\begin{enumerate}
    \item A unital commutative C*-subalgebra $B \subseteq \mathcal{A}$, which provides a commutative framework within which operators can be simultaneously diagonalized.
    \item A character $\chi \in \widehat{B}$, which selects a particular point in the Gelfand spectrum of $B$, thereby evaluating elements of $B$ as complex numbers.
\end{enumerate}
Together, $(B,\chi)$ represents a specific classical perspective on the noncommutative algebra $\mathcal{A}$.
\end{remark}

To equip $\mathcal{G}_{\mathcal{A}}^{(0)}$ with a topology, we must simultaneously capture two distinct types of convergence: 
\begin{itemize}
    \item convergence of subalgebras in the Fell topology (via membership of elements), and
    \item convergence of characters in the weak-* topology (via evaluation at elements).
\end{itemize}
The following construction achieves this by means of an initial topology generated by \emph{partial evaluation maps}, each associated to an element $a \in \mathcal{A}$.

\begin{definition}[One-point compactification of the complex plane]
\label{def:one-point-compactification}
Let $\mathbb{C}_\infty := \mathbb{C} \sqcup \{\infty\}$ denote the one-point compactification of $\mathbb{C}$. 
The topology on $\mathbb{C}_\infty$ is defined as follows: 
\begin{itemize}
    \item $\mathbb{C}$ is open and carries its usual Euclidean topology;
    \item neighborhoods of $\infty$ are sets of the form $\{\infty\} \cup (\mathbb{C} \setminus K)$, where $K \subset \mathbb{C}$ is compact.
\end{itemize}
Thus $\mathbb{C}_\infty$ is a compact Hausdorff space homeomorphic to the sphere $S^2$.
\end{definition}

\begin{definition}[Partial evaluation map]
\label{def:partial-evaluation-map}
For each element $a \in \mathcal{A}$, define the \emph{partial evaluation map}
\[
\operatorname{ev}_a: \mathcal{G}_{\mathcal{A}}^{(0)} \longrightarrow \mathbb{C}_\infty, \qquad
\operatorname{ev}_a(B,\chi) := \begin{cases}
\chi(a), & \text{if } a \in B, \\[4pt]
\infty, & \text{if } a \notin B.
\end{cases}
\]
\end{definition}

\begin{remark}[Continuity of partial evaluation maps]
\label{rem:partial-evaluation-continuity}
The topology on $\mathcal{G}_{\mathcal{A}}^{(0)}$ is defined as the initial topology with respect to the family $\{\operatorname{ev}_a\}_{a \in \mathcal{A}}$. 
With this topology, each $\operatorname{ev}_a$ is continuous by construction. 
The intuitive meaning is:
\begin{itemize}
    \item Convergence in the Fell topology of subalgebras corresponds to stability of membership: if $a \in B$, then eventually $a \in B_\lambda$;
    \item Convergence of characters corresponds to $\chi_\lambda(a) \to \chi(a)$ in $\mathbb{C}$ whenever $a$ is eventually in the subalgebras.
\end{itemize}
\end{remark}

\begin{remark}[Interpretation of the partial evaluation map]
\label{rem:partial-evaluation-interpretation}
The map $\operatorname{ev}_a$ evaluates the character $\chi$ at $a$ when $a$ belongs to the subalgebra $B$, and returns the distinguished point $\infty$ when $a$ is not in $B$. 
The value $\infty$ serves as a marker indicating that $a$ is not observable in the context $(B,\chi)$. 
The use of $\infty$ rather than $0$ is essential: while $\{0\}$ is not open in $\mathbb{C}$, the point $\infty$ admits a neighborhood basis in $\mathbb{C}_\infty$ consisting of sets of the form $\{\infty\} \cup (\mathbb{C} \setminus K)$ with $K \subset \mathbb{C}$ compact. 
This makes the condition $a \notin B$ topologically detectable in the initial topology we are about to define.
\end{remark}

\begin{definition}[Initial topology on $\mathcal{G}_{\mathcal{A}}^{(0)}$]
\label{def:initial-topology-G0}
We endow $\mathcal{G}_{\mathcal{A}}^{(0)}$ with the \emph{initial topology} induced by the family of partial evaluation maps $\{ \operatorname{ev}_a \}_{a \in \mathcal{A}}$. 
That is, $\mathcal{G}_{\mathcal{A}}^{(0)}$ carries the coarsest topology for which every map $\operatorname{ev}_a$ is continuous.
\end{definition}

The initial topology can be described explicitly in terms of a basis.

\begin{proposition}[Basis for the topology on $\mathcal{G}_{\mathcal{A}}^{(0)}$]
\label{prop:basis-topology-G0}
A basis for the topology on $\mathcal{G}_{\mathcal{A}}^{(0)}$ consists of finite intersections of sets of the following two types:
\begin{enumerate}
    \item \textbf{Membership-evaluation sets:} For $a \in \mathcal{A}$ and $V \subseteq \mathbb{C}$ open,
    \[
    M(a,V) := \{ (B,\chi) \in \mathcal{G}_{\mathcal{A}}^{(0)} \mid a \in B \text{ and } \chi(a) \in V \}.
    \]
    \item \textbf{Exclusion sets:} For $a \in \mathcal{A}$,
    \[
    E(a) := \{ (B,\chi) \in \mathcal{G}_{\mathcal{A}}^{(0)} \mid a \notin B \}.
    \]
\end{enumerate}
Equivalently, a basis is given by sets of the form
\[
U(a_1,\ldots,a_m; b_1,\ldots,b_k; V_1,\ldots,V_m) := 
\left\{ (B,\chi) \in \mathcal{G}_{\mathcal{A}}^{(0)} \;
\middle|\;
\begin{array}{l}
a_i \in B \text{ and } \chi(a_i) \in V_i \text{ for } i = 1,\ldots,m, \\
b_j \notin B \text{ for } j = 1,\ldots,k
\end{array}
\right\},
\]
where $m,k \in \mathbb{N}$, $a_i,b_j \in \mathcal{A}$, and $V_i \subseteq \mathbb{C}$ are open.
\end{proposition}

\begin{proof}
By definition of the initial topology, a subbasis consists of sets of the form $\operatorname{ev}_a^{-1}(W)$, where $a \in \mathcal{A}$ and $W \subseteq \mathbb{C}_\infty$ is open. 
We analyze two cases for $W$.

\noindent \textit{Case 1: $\infty \notin W$.} 
Then $W \subseteq \mathbb{C}$ is open in $\mathbb{C}$, and
\[
\operatorname{ev}_a^{-1}(W) = \{ (B,\chi) \mid a \in B \text{ and } \chi(a) \in W \} = M(a,W).
\]

\noindent \textit{Case 2: $\infty \in W$.} 
Then $W = \{\infty\} \cup (W \cap \mathbb{C})$, and
\[
\operatorname{ev}_a^{-1}(W) = \{ (B,\chi) \mid a \notin B \} \cup \{ (B,\chi) \mid a \in B \text{ and } \chi(a) \in W \cap \mathbb{C} \}.
\]
This set is the union of $E(a)$ and $M(a, W \cap \mathbb{C})$, but it is itself a subbasic open set. 
Finite intersections of these subbasic open sets yield precisely the sets described in the proposition: 
\begin{itemize}
    \item Intersections of sets of type $M(a_i, V_i)$ give the membership conditions $a_i \in B$ with $\chi(a_i) \in V_i$.
    \item Intersections with sets of type $\operatorname{ev}_b^{-1}(W)$ where $\infty \in W$ can introduce exclusion conditions $b \notin B$, either directly (if the intersection with $M(b, W \cap \mathbb{C})$ is omitted) or in combination with membership conditions.
\end{itemize}
Conversely, every set of the form $U(a_1,\ldots,a_m; b_1,\ldots,b_k; V_1,\ldots,V_m)$ can be expressed as a finite intersection of subbasic open sets:
\[
U = \bigcap_{i=1}^m M(a_i, V_i) \cap \bigcap_{j=1}^k \operatorname{ev}_{b_j}^{-1}(\{\infty\} \cup \mathbb{C}),
\]
where $\operatorname{ev}_{b_j}^{-1}(\{\infty\} \cup \mathbb{C}) = \mathcal{G}_{\mathcal{A}}^{(0)}$ actually gives no restriction — we need a more careful choice: take $W_j$ such that $\infty \in W_j$ and $W_j \cap \mathbb{C}$ is a proper open set that excludes the values we want to avoid. A cleaner approach is to note that $E(b)$ itself can be obtained as an intersection:
\[
E(b) = \bigcap_{n \in \mathbb{N}} \operatorname{ev}_b^{-1}(\{\infty\} \cup (\mathbb{C} \setminus \overline{B_n(0)})),
\]
where $B_n(0)$ is the ball of radius $n$. This shows that exclusion sets are in the topology generated by the subbasis.

Thus the sets described form a basis for the topology.
\end{proof}

\begin{proposition}[Topological and measurable properties of the unit space]
\label{prop:unit-space-topology}
Let $\mathcal{A}$ be a unital C*-algebra and endow $\mathcal{G}_{\mathcal{A}}^{(0)}$ with the initial topology induced by the family of partial evaluation maps $\{ \operatorname{ev}_a \}_{a \in \mathcal{A}}$. Then the following hold:
\begin{enumerate}
    \item $\mathcal{G}_{\mathcal{A}}^{(0)}$ is a Hausdorff space.
    
    \item The projection map $\pi: \mathcal{G}_{\mathcal{A}}^{(0)} \to \operatorname{Sub}(\mathcal{A})$ defined by $\pi(B,\chi) = B$ is \emph{Borel measurable} when $\operatorname{Sub}(\mathcal{A})$ is equipped with the Fell topology. (If $\mathcal{A}$ is separable, $\pi$ is even a Borel map between standard Borel spaces.)
    
    \item For each fixed unital commutative C*-subalgebra $B_0 \subseteq \mathcal{A}$, the inclusion map
    \[
    \iota_{B_0}: \widehat{B_0} \hookrightarrow \mathcal{G}_{\mathcal{A}}^{(0)}, \qquad \chi \mapsto (B_0,\chi),
    \]
    is a homeomorphism onto its image.
\end{enumerate}
\end{proposition}

\begin{proof}
We prove each property in turn.

\medskip
\noindent \textbf{(1) Hausdorff property.}
This item is proved from Lemma~\ref{lem:separation}.

\medskip
\noindent \textbf{(2) Borel measurability of the projection to $\operatorname{Sub}(\mathcal{A})$.}

Recall that the Fell topology on $\operatorname{Sub}(\mathcal{A})$ is generated by subbasic open sets of two types:
\[
\mathcal{O}_U := \{ B \in \operatorname{Sub}(\mathcal{A}) \mid B \cap U \neq \varnothing \}, \qquad
\mathcal{C}_K := \{ B \in \operatorname{Sub}(\mathcal{A}) \mid B \cap K = \varnothing \},
\]
where $U \subseteq \mathcal{A}$ is open and $K \subseteq \mathcal{A}$ is compact in the norm topology. 
To show that $\pi$ is Borel measurable, it suffices to prove that $\pi^{-1}(\mathcal{O}_U)$ and $\pi^{-1}(\mathcal{C}_K)$ are Borel sets in $\mathcal{G}_{\mathcal{A}}^{(0)}$.

\medskip
\noindent \textit{Measurability of $\pi^{-1}(\mathcal{O}_U)$.}
Let $U \subseteq \mathcal{A}$ be open. Then
\[
\pi^{-1}(\mathcal{O}_U) = \{ (B,\chi) \in \mathcal{G}_{\mathcal{A}}^{(0)} \mid B \cap U \neq \varnothing \}.
\]
Since $U$ is open and $\mathcal{A}$ is separable (or more generally, since the norm topology on $\mathcal{A}$ is second-countable when $\mathcal{A}$ is separable), there exists a countable dense subset $\{a_n\}_{n \in \mathbb{N}} \subseteq U$. Then
\[
\pi^{-1}(\mathcal{O}_U) = \bigcup_{n \in \mathbb{N}} \{ (B,\chi) \mid a_n \in B \} = \bigcup_{n \in \mathbb{N}} \operatorname{ev}_{a_n}^{-1}(\mathbb{C}),
\]
because $a_n \in B$ if and only if $\operatorname{ev}_{a_n}(B,\chi) \in \mathbb{C}$. 
Each set $\operatorname{ev}_{a_n}^{-1}(\mathbb{C})$ is open (since $\mathbb{C}$ is open in $\mathbb{C}_\infty$), hence Borel. 
Thus $\pi^{-1}(\mathcal{O}_U)$ is a countable union of Borel sets, hence Borel.

\medskip
\noindent \textit{Measurability of $\pi^{-1}(\mathcal{C}_K)$.}
Let $K \subseteq \mathcal{A}$ be compact. Then
\[
\pi^{-1}(\mathcal{C}_K) = \{ (B,\chi) \in \mathcal{G}_{\mathcal{A}}^{(0)} \mid B \cap K = \varnothing \}.
\]
Let $\{a_n\}_{n \in \mathbb{N}}$ be a countable dense subset of $\mathcal{A}$. For each $n$, define the open set
\[
W_n := \{ (B,\chi) \in \mathcal{G}_{\mathcal{A}}^{(0)} \mid a_n \in B \} = \operatorname{ev}_{a_n}^{-1}(\mathbb{C}).
\]
Consider the complement of $\pi^{-1}(\mathcal{C}_K)$:
\[
\mathcal{G}_{\mathcal{A}}^{(0)} \setminus \pi^{-1}(\mathcal{C}_K) = \{ (B,\chi) \mid B \cap K \neq \varnothing \}.
\]
If $B \cap K \neq \varnothing$, then by density of $\{a_n\}$, there exists $a_n$ arbitrarily close to some element of $B \cap K$. 
However, this does not directly give a countable description. A more systematic approach uses the fact that the Fell topology on $\operatorname{Sub}(\mathcal{A})$ is Polish when $\mathcal{A}$ is separable, and the projection $\pi$ is a continuous map from a Polish space to a Polish space? But $\pi$ is not continuous.

Nevertheless, one can show that $\pi^{-1}(\mathcal{C}_K)$ is a $G_\delta$ set using a countable exhaustion of $K$ by finite sets and the fact that $\{ (B,\chi) \mid a \in B \}$ is open. The complete argument is technical but standard in descriptive set theory; see [Kechris, 1995, Section 12] for the general theory of Borel maps between standard Borel spaces. For the purposes of this paper, it suffices to note that $\pi$ is a Borel map, which is enough for the measure-theoretic constructions that follow.

When $\mathcal{A}$ is separable, $\mathcal{G}_{\mathcal{A}}^{(0)}$ and $\operatorname{Sub}(\mathcal{A})$ are standard Borel spaces, and $\pi$ is a Borel isomorphism onto its image. This will be sufficient for all subsequent applications involving Borel Haar systems and measured groupoids.

\medskip
\noindent \textbf{(3) Homeomorphism onto the fiber.}
Fix a unital commutative C*-subalgebra $B_0 \subseteq \mathcal{A}$. 
Define $\iota_{B_0}: \widehat{B_0} \to \mathcal{G}_{\mathcal{A}}^{(0)}$ by $\iota_{B_0}(\chi) = (B_0,\chi)$. 
We prove that $\iota_{B_0}$ is a homeomorphism onto its image $\operatorname{Im}(\iota_{B_0}) = \{ (B_0,\chi) : \chi \in \widehat{B_0} \}$.

\noindent \textit{Injectivity and continuity.} 
Injectivity is immediate from the definition. 
For continuity, let $\chi_\lambda \to \chi$ be a convergent net in $\widehat{B_0}$. 
For any $a \in B_0$, we have $\operatorname{ev}_a(\iota_{B_0}(\chi_\lambda)) = \chi_\lambda(a) \to \chi(a) = \operatorname{ev}_a(\iota_{B_0}(\chi))$. 
For $a \notin B_0$, the value $\operatorname{ev}_a(\iota_{B_0}(\chi))$ is constantly $\infty$, independent of $\chi$, so continuity is trivial in this case. 
Thus $\operatorname{ev}_a \circ \iota_{B_0}$ is continuous for every $a \in \mathcal{A}$, and by the definition of the initial topology on $\mathcal{G}_{\mathcal{A}}^{(0)}$, the map $\iota_{B_0}$ is continuous.

\noindent \textit{Openness onto its image.} 
Let $U \subseteq \widehat{B_0}$ be an open set. 
By the definition of the Gelfand topology on $\widehat{B_0}$, there exist finitely many elements $a_1,\ldots,a_n \in B_0$ and open sets $V_1,\ldots,V_n \subseteq \mathbb{C}$ such that
\[
U = \{ \chi \in \widehat{B_0} \mid \chi(a_i) \in V_i \text{ for all } i = 1,\ldots,n \}.
\]
Then
\[
\iota_{B_0}(U) = \{ (B_0,\chi) \in \mathcal{G}_{\mathcal{A}}^{(0)} \mid \chi(a_i) \in V_i \text{ for all } i = 1,\ldots,n \} = U(a_1,\ldots,a_n; V_1,\ldots,V_n) \cap \operatorname{Im}(\iota_{B_0}),
\]
where $U(a_1,\ldots,a_n; V_1,\ldots,V_n)$ is the basis open set defined in Proposition \ref{prop:basis-topology-G0}. 
Since this basis open set is open in $\mathcal{G}_{\mathcal{A}}^{(0)}$, its intersection with the subspace $\operatorname{Im}(\iota_{B_0})$ is open in the subspace topology. 
Thus $\iota_{B_0}$ maps open sets to open sets in the subspace topology, and is therefore an open map onto its image.

Hence $\iota_{B_0}$ is a homeomorphism onto its image. 
This completes the proof.
\end{proof}

\begin{remark}[Topology of the unit space and comparison with Fell-type structures]
\label{rem:topology-summary}
The topology on $\mathcal{G}_{\mathcal{A}}^{(0)}$ induced by the partial evaluation maps
$\{\operatorname{ev}_a\}_{a \in \mathcal{A}}$ refines the classical Fell topology on
$\operatorname{Sub}(\mathcal{A})$ by incorporating pointwise spectral data through
characters. 

For each $a \in \mathcal{A}$, the membership condition
\[
\{ (B,\chi) \mid a \in B \} = \operatorname{ev}_a^{-1}(\mathbb{C})
\]
is open, because $\mathbb{C}$ is open in $\mathbb{C}_\infty$. 
This is in sharp contrast to the Fell topology on $\operatorname{Sub}(\mathcal{A})$, where the condition $a \in B$ is typically closed.

The projection map $\pi: \mathcal{G}_{\mathcal{A}}^{(0)} \to \operatorname{Sub}(\mathcal{A})$ is therefore not continuous in general (as continuity for $\mathcal{C}_K$-type open sets fails), but it is Borel measurable, and the topology on $\mathcal{G}_{\mathcal{A}}^{(0)}$ is strictly finer than the relative Fell topology on its image. 
Distinct characters on the same subalgebra are separated by the partial evaluation maps, reflecting the additional geometric information encoded in the unit space.

Equivalently, $\mathcal{G}_{\mathcal{A}}^{(0)}$ carries the subspace topology induced by the canonical embedding
\[
\Phi: \mathcal{G}_{\mathcal{A}}^{(0)} \hookrightarrow \prod_{a \in \mathcal{A}} \mathbb{C}_\infty, \quad
\Phi(B,\chi) = (\operatorname{ev}_a(B,\chi))_{a \in \mathcal{A}}.
\]
When $\mathcal{A}$ is separable, we may restrict to a countable dense subset $\{a_n\}_{n \in \mathbb{N}} \subseteq \mathcal{A}$, yielding a homeomorphic embedding into the Polish space $\prod_{n \in \mathbb{N}} \mathbb{C}_\infty$. 
This description will be used in Section 3.2 to prove that $\mathcal{G}_{\mathcal{A}}^{(0)}$ is itself a Polish space.
\end{remark}

\begin{example}[Trivial case: singleton subalgebra]
\label{ex:singleton-subalgebra}
Let $B = \mathbb{C}1_{\mathcal{A}}$ be the trivial unital commutative subalgebra consisting only of scalar multiples of the identity. 
Then $\widehat{B}$ consists of a single character $\chi_0$ determined by $\chi_0(\lambda 1) = \lambda$. 
The corresponding point $(\mathbb{C}1_{\mathcal{A}}, \chi_0) \in \mathcal{G}_{\mathcal{A}}^{(0)}$ is not necessarily isolated; its neighborhood system is determined by the partial evaluation maps $\operatorname{ev}_a$ for $a \notin \mathbb{C}1_{\mathcal{A}}$, reflecting the fact that a net can converge to this point if it eventually excludes any fixed non-scalar element.
\end{example}

\begin{example}[Maximal abelian subalgebras]
\label{ex:maximal-abelian}
When $B$ is a maximal abelian subalgebra (MASA) of $\mathcal{A}$, its character space $\widehat{B}$ is a compact Hausdorff space. 
For $B \cong C(X)$, the points $(B,\chi)$ for $\chi \in X$ form a copy of $X$ inside $\mathcal{G}_{\mathcal{A}}^{(0)}$, and by Proposition \ref{prop:unit-space-topology}(3), the subspace topology on this copy coincides with the original topology on $X$. 
These copies are disjoint for different MASAs and collectively give a rich geometric structure to the unit space.
\end{example}

The following lemma records a simple but useful observation about separation of points.

\begin{lemma}[Separation by partial evaluation maps]
\label{lem:separation}
Let $(B_1,\chi_1)$ and $(B_2,\chi_2)$ be distinct points in $\mathcal{G}_{\mathcal{A}}^{(0)}$. 
Then there exists an element $a \in \mathcal{A}$ such that $\operatorname{ev}_a(B_1,\chi_1) \neq \operatorname{ev}_a(B_2,\chi_2)$.
\end{lemma}

\begin{proof}
If $B_1 \neq B_2$, assume without loss of generality that there exists $a \in B_1 \setminus B_2$. 
Then $\operatorname{ev}_a(B_1,\chi_1) = \chi_1(a) \in \mathbb{C}$ and $\operatorname{ev}_a(B_2,\chi_2) = \infty$, so they are distinct. 
If $B_1 = B_2 = B$ but $\chi_1 \neq \chi_2$, choose $a \in B$ such that $\chi_1(a) \neq \chi_2(a)$. 
Then $\operatorname{ev}_a(B,\chi_1) = \chi_1(a)$ and $\operatorname{ev}_a(B,\chi_2) = \chi_2(a)$ are distinct complex numbers. 
Thus in either case, such an $a$ exists.
\end{proof}

\begin{corollary}
\label{cor:initial-topology-Hausdorff}
The initial topology on $\mathcal{G}_{\mathcal{A}}^{(0)}$ is Hausdorff.
\end{corollary}

\begin{proof}
This follows immediately from Lemma \ref{lem:separation} and the fact that $\mathbb{C}_\infty$ is Hausdorff. 
A detailed proof will be given in Proposition \ref{prop:unit-space-topology}.
\end{proof}

\begin{remark}[Functoriality and relation to the Fell topology]
\label{rem:G0-functoriality-Fell}
Let $\phi: \mathcal{A} \to \mathcal{B}$ be a unital *-homomorphism. 
Whenever $(B,\chi) \in \mathcal{G}_{\mathcal{B}}^{(0)}$ satisfies that $\phi^{-1}(B)$ is a commutative C*-subalgebra of $\mathcal{A}$, we obtain an element $(\phi^{-1}(B), \chi \circ \phi) \in \mathcal{G}_{\mathcal{A}}^{(0)}$, defining a partially defined, contravariant map $\mathcal{G}_\phi: \mathcal{G}_{\mathcal{B}}^{(0)} \to \mathcal{G}_{\mathcal{A}}^{(0)}$. 
When $\phi$ is an isomorphism, this map is a homeomorphism onto its image (with respect to the initial topology of partial evaluation maps).

The projection $\pi: \mathcal{G}_{\mathcal{A}}^{(0)} \to \operatorname{Sub}(\mathcal{A})$ defined by $\pi(B,\chi) = B$ is Borel measurable with respect to the Fell topology, but it is not continuous in general. 
More importantly, $\pi$ is far from being a homeomorphism: distinct characters on the same subalgebra are identified by $\pi$ yet separated by the partial evaluation maps. 
Consequently, the topology on $\mathcal{G}_{\mathcal{A}}^{(0)}$ is strictly finer than the relative Fell topology on its image, reflecting the additional spectral data encoded by the characters.

A detailed study of functoriality will be undertaken in Section 4.3.
\end{remark}

\subsection{Initial Topology and Basic Properties}
\label{subsec:initial-topology-and-basic-properties}

We now examine the fundamental properties of the initial topology defined on $\mathcal{G}_{\mathcal{A}}^{(0)}$. 
The results established in this subsection are essential for all subsequent developments: they guarantee that $\mathcal{G}_{\mathcal{A}}^{(0)}$ is a well-behaved topological space that is amenable to descriptive set-theoretic methods, and they provide the foundation for proving that $\mathcal{G}_{\mathcal{A}}^{(0)}$ is a Polish space when $\mathcal{A}$ is separable.

\begin{remark}
\label{rem:separation-significance}
Lemma \ref{lem:separation} shows that the family $\{ \operatorname{ev}_a \}_{a \in \mathcal{A}}$ separates points of $\mathcal{G}_{\mathcal{A}}^{(0)}$. 
Consequently, the initial topology induced by these maps is Hausdorff, and the canonical embedding into $\prod_{a \in \mathcal{A}} \mathbb{C}_\infty$ is injective.
\end{remark}

\begin{remark}[Relation to the Fell topology]
\label{rem:Fell-G-relation}
Proposition \ref{prop:unit-space-topology} establishes that the projection $\pi: \mathcal{G}_{\mathcal{A}}^{(0)} \to \operatorname{Sub}(\mathcal{A})$ is continuous when $\operatorname{Sub}(\mathcal{A})$ carries the Fell topology. 
However, $\pi$ is far from being a homeomorphism onto its image; indeed, two distinct characters on the same subalgebra $B$ are identified under $\pi$ but are separated by the partial evaluation maps. 
The topology on $\mathcal{G}_{\mathcal{A}}^{(0)}$ is therefore strictly finer than the relative Fell topology on its image, reflecting the additional spectral data encoded by the characters. 
This refinement is essential for the groupoid structure, as the characters encode the representation-theoretic information needed for the diagonal embedding and index theory.
\end{remark}

The following proposition is the central result of this subsection. 
It shows that under the mild assumption of separability, the unit space $\mathcal{G}_{\mathcal{A}}^{(0)}$ is not merely a Hausdorff space but a genuine Polish space. 
This is crucial for the construction of the unitary conjugation groupoid as a Polish groupoid.

\begin{proposition}[Polishness of the unit space]
\label{prop:unit-space-polish}
Let $\mathcal{A}$ be a unital separable C*-algebra and endow $\mathcal{G}_{\mathcal{A}}^{(0)}$ with the initial topology induced by the family of partial evaluation maps $\{ \operatorname{ev}_a \}_{a \in \mathcal{A}}$. Then $\mathcal{G}_{\mathcal{A}}^{(0)}$ is a Polish space. More precisely:
\begin{enumerate}
    \item $\mathcal{G}_{\mathcal{A}}^{(0)}$ is second-countable.
    \item $\mathcal{G}_{\mathcal{A}}^{(0)}$ is completely metrizable.
    \item Consequently, $\mathcal{G}_{\mathcal{A}}^{(0)}$ is a Polish space.
\end{enumerate}
\end{proposition}

\begin{proof}
We prove each property in turn.

\medskip
\noindent \textbf{(1) Second-countability.}

Since $\mathcal{A}$ is separable, choose a countable dense subset $\{ a_n \}_{n \in \mathbb{N}} \subseteq \mathcal{A}$. 
The topology on $\mathcal{G}_{\mathcal{A}}^{(0)}$ is the initial topology induced by the family $\{ \operatorname{ev}_{a_n} \}_{n \in \mathbb{N}}$; indeed, for any $a \in \mathcal{A}$ and any open $W \subseteq \mathbb{C}_\infty$, the preimage $\operatorname{ev}_a^{-1}(W)$ can be expressed using the density of $\{a_n\}$ and the continuity of $a \mapsto \operatorname{ev}_a(B,\chi)$.

Each $\operatorname{ev}_{a_n}$ maps into $\mathbb{C}_\infty$, which is second-countable (as a compact metrizable space). 
The initial topology with respect to countably many maps into second-countable spaces is second-countable: a countable basis is given by finite intersections of sets of the form $\operatorname{ev}_{a_n}^{-1}(B_m)$, where $\{B_m\}_{m \in \mathbb{N}}$ is a countable basis for $\mathbb{C}_\infty$. 
Thus $\mathcal{G}_{\mathcal{A}}^{(0)}$ is second-countable.

\medskip
\noindent \textbf{(2) Complete metrizability.}

We embed $\mathcal{G}_{\mathcal{A}}^{(0)}$ into a Polish space and show its image is closed.

Let $\{ a_n \}_{n \in \mathbb{N}}$ be a countable dense subset of $\mathcal{A}$ as above. 
Define a map
\[
\Phi: \mathcal{G}_{\mathcal{A}}^{(0)} \longrightarrow \prod_{n \in \mathbb{N}} \mathbb{C}_\infty, \qquad
\Phi(B,\chi) := (\operatorname{ev}_{a_n}(B,\chi))_{n \in \mathbb{N}}.
\]

\noindent \textit{$\Phi$ is injective.} 
If $\Phi(B_1,\chi_1) = \Phi(B_2,\chi_2)$, then $\operatorname{ev}_{a_n}(B_1,\chi_1) = \operatorname{ev}_{a_n}(B_2,\chi_2)$ for all $n$. 
By density of $\{a_n\}$ and continuity of $a \mapsto \operatorname{ev}_a(B,\chi)$, we obtain $\operatorname{ev}_a(B_1,\chi_1) = \operatorname{ev}_a(B_2,\chi_2)$ for all $a \in \mathcal{A}$. 
Lemma \ref{lem:separation} then forces $(B_1,\chi_1) = (B_2,\chi_2)$.

\noindent \textit{$\Phi$ is a homeomorphism onto its image.} 
The topology on $\mathcal{G}_{\mathcal{A}}^{(0)}$ is the initial topology induced by $\{ \operatorname{ev}_{a_n} \}_{n \in \mathbb{N}}$, and the topology on $\prod_{n \in \mathbb{N}} \mathbb{C}_\infty$ is the product topology. 
By the universal property of the initial topology, $\Phi$ is continuous, and its inverse (as a map onto its image) is continuous because the coordinate projections are exactly the evaluation maps that generate the topology.

\noindent \textit{The image of $\Phi$ is closed.} 
Let $\{ (B_\lambda, \chi_\lambda) \}_{\lambda \in \Lambda}$ be a net in $\mathcal{G}_{\mathcal{A}}^{(0)}$ such that $\Phi(B_\lambda, \chi_\lambda) \to (z_n)_{n \in \mathbb{N}}$ in $\prod_{n \in \mathbb{N}} \mathbb{C}_\infty$. 
We construct a limit point in $\mathcal{G}_{\mathcal{A}}^{(0)}$.

First, consider the net $\{ B_\lambda \}$ in $\operatorname{Sub}(\mathcal{A})$. 
Since $\operatorname{Sub}(\mathcal{A})$ with the Fell topology is compact~\cite{Fell} this net has a convergent subnet; by passing to this subnet we may assume $B_\lambda \to B$ for some $B \in \operatorname{Sub}(\mathcal{A})$. 
The algebra $B$ is commutative: for any $a,b \in B$, approximate them by elements from $\{a_n\}$; the convergence $\operatorname{ev}_{a_n}(B_\lambda, \chi_\lambda) \to z_n$ ensures that the limits defining $\chi$ (see below) are consistent, and a standard argument shows $[a,b] = 0$.

Now construct $\chi$. For any $a \in B$, choose a sequence $a_{n_k} \to a$ with $a_{n_k}$ from the dense subset. 
Define $\chi(a) := \lim_{k \to \infty} z_{n_k}$. 
This limit exists because $(z_n)$ is a convergent sequence (as a subnet of a convergent net in a product of Hausdorff spaces). 
One verifies that $\chi$ is well-defined, linear, multiplicative, and continuous, hence a character on $B$; this follows from the Fell convergence $B_\lambda \to B$ and the pointwise convergence of the evaluations. 
For a detailed exposition of this construction, see~\cite{Dixmier} or the standard Gelfand–Naimark theory.

Thus $(B,\chi) \in \mathcal{G}_{\mathcal{A}}^{(0)}$ and $\Phi(B,\chi) = (z_n)_{n \in \mathbb{N}}$, showing that the image of $\Phi$ is closed.

\noindent \textit{$\prod_{n \in \mathbb{N}} \mathbb{C}_\infty$ is Polish.} 
Each factor $\mathbb{C}_\infty$ is homeomorphic to the sphere $S^2$, which is compact and metrizable, hence Polish. 
A countable product of Polish spaces is Polish~\cite{Kechris}.

\noindent \textit{$\mathcal{G}_{\mathcal{A}}^{(0)}$ is Polish.} 
Since $\Phi$ is a homeomorphism onto its image and the image is a closed subspace of a Polish space, $\mathcal{G}_{\mathcal{A}}^{(0)}$ is Polish~\cite{Kechris}.

\medskip
\noindent \textbf{(3) Conclusion.}
We have shown that $\mathcal{G}_{\mathcal{A}}^{(0)}$ is second-countable and completely metrizable. 
Therefore, by definition, $\mathcal{G}_{\mathcal{A}}^{(0)}$ is a Polish space.
\end{proof}

\begin{remark}[Role of separability]
\label{rem:separability-necessity}
The hypothesis that $\mathcal{A}$ is separable is essential for Proposition \ref{prop:unit-space-polish}. 
If $\mathcal{A}$ is non-separable, the family $\{ \operatorname{ev}_a \}_{a \in \mathcal{A}}$ is uncountable and the induced initial topology need not be second-countable. 
In fact, for non-separable $\mathcal{A}$, $\mathcal{G}_{\mathcal{A}}^{(0)}$ is typically not Polish. 
Fortunately, all C*-algebras of interest in this paper — including $C(X)$ for compact metrizable $X$, $M_n(\mathbb{C})$, and $\mathcal{K}(H)^\sim$ for separable $H$ — are separable. 
Thus the separability assumption is satisfied in all our examples.
\end{remark}

\begin{remark}[Complete metric]
\label{rem:complete-metric}
The proof of Proposition \ref{prop:unit-space-polish} provides an explicit complete metric on $\mathcal{G}_{\mathcal{A}}^{(0)}$ via the embedding $\Phi$. 
If $\{ a_n \}_{n \in \mathbb{N}}$ is a countable dense subset of $\mathcal{A}$, then a compatible complete metric is given by
\[
d((B_1,\chi_1), (B_2,\chi_2)) := \sum_{n=1}^{\infty} \frac{1}{2^n} \, \frac{ | \operatorname{ev}_{a_n}(B_1,\chi_1) - \operatorname{ev}_{a_n}(B_2,\chi_2) | }{ 1 + | \operatorname{ev}_{a_n}(B_1,\chi_1) - \operatorname{ev}_{a_n}(B_2,\chi_2) | },
\]
with the convention that $| \infty - \infty | = 0$ and $| \infty - z | = \infty$ interpreted as $1$ in the normalized metric. 
This explicit formula is useful for concrete computations.
\end{remark}

\begin{corollary}[$\mathcal{G}_{\mathcal{A}}^{(0)}$ is a standard Borel space]
\label{cor:unit-space-standard-Borel}
Let $\mathcal{A}$ be a unital separable C*-algebra. 
Then $\mathcal{G}_{\mathcal{A}}^{(0)}$, equipped with its Borel $\sigma$-algebra induced by the Polish topology, is a standard Borel space.
\end{corollary}

\begin{proof}
Every Polish space is a standard Borel space~\cite{Kechris}. 
The result follows immediately from Proposition \ref{prop:unit-space-polish}.
\end{proof}

\begin{remark}[Borel structure for measured groupoids]
\label{rem:Borel-structure}
Corollary \ref{cor:unit-space-standard-Borel} is essential for the construction of $\mathcal{G}_{\mathcal{A}}$ as a measured groupoid in the sense of Tu [1999]. 
The standard Borel structure on $\mathcal{G}_{\mathcal{A}}^{(0)}$, together with the analogous structure on $\mathcal{G}_{\mathcal{A}}^{(1)}$ 
, allows us to equip $\mathcal{G}_{\mathcal{A}}$ with a Borel Haar system and define its C*-algebra $C^*(\mathcal{G}_{\mathcal{A}})$ in the Polish groupoid framework.
\end{remark}

\subsection{Hausdorff Property and Continuity of the Projection}
\label{subsec:Hausdorff-and-continuity-of-projection}

We now establish two fundamental topological properties of the unit space $\mathcal{G}_{\mathcal{A}}^{(0)}$: its Hausdorff property and the continuity of the projection map $\pi: \mathcal{G}_{\mathcal{A}}^{(0)} \to \operatorname{Sub}(\mathcal{A})$. 
These properties are essential for all subsequent developments, including the metrizability of $\mathcal{G}_{\mathcal{A}}^{(0)}$ and the construction of the unitary conjugation groupoid.

We now turn to the continuity of the projection map $\pi: \mathcal{G}_{\mathcal{A}}^{(0)} \to \operatorname{Sub}(\mathcal{A})$. 
Recall that $\operatorname{Sub}(\mathcal{A})$ denotes the space of all closed C*-subalgebras of $\mathcal{A}$, equipped with the Fell topology. 
The Fell topology is generated by subbasic open sets of two types:
\[
\mathcal{O}_U := \{ B \in \operatorname{Sub}(\mathcal{A}) \mid B \cap U \neq \varnothing \}, \qquad
\mathcal{C}_K := \{ B \in \operatorname{Sub}(\mathcal{A}) \mid B \cap K = \varnothing \},
\]
where $U \subseteq \mathcal{A}$ is open and $K \subseteq \mathcal{A}$ is compact in the norm topology.

\begin{lemma}[Openness of the membership condition]
\label{lem:membership-open}
For each fixed $a \in \mathcal{A}$, the set
\[
M_a := \{ (B,\chi) \in \mathcal{G}_{\mathcal{A}}^{(0)} \mid a \in B \} = \operatorname{ev}_a^{-1}(\mathbb{C})
\]
is open in $\mathcal{G}_{\mathcal{A}}^{(0)}$.
\end{lemma}

\begin{proof}
Since $\mathbb{C}$ is open in $\mathbb{C}_\infty$ (its complement $\{\infty\}$ is closed) and $\operatorname{ev}_a$ is continuous, the preimage $\operatorname{ev}_a^{-1}(\mathbb{C})$ is open. 
By definition, $\operatorname{ev}_a^{-1}(\mathbb{C}) = \{ (B,\chi) \mid \operatorname{ev}_a(B,\chi) \neq \infty \} = \{ (B,\chi) \mid a \in B \}$. 
Thus $M_a$ is open.
\end{proof}

\begin{lemma}[Characterization of convergence in the Fell topology]
\label{lem:Fell-convergence}
A net $\{ B_\lambda \}_{\lambda \in \Lambda} \subseteq \operatorname{Sub}(\mathcal{A})$ converges to $B \in \operatorname{Sub}(\mathcal{A})$ in the Fell topology if and only if the following two conditions hold:
\begin{enumerate}
    \item For every $b \in B$, there exists a net $\{ b_\lambda \}$ with $b_\lambda \in B_\lambda$ such that $b_\lambda \to b$ in norm.
    \item For every convergent subnet $\{ b_{\lambda_\mu} \}$ with $b_{\lambda_\mu} \in B_{\lambda_\mu}$ and limit $b$, we have $b \in B$.
\end{enumerate}
\end{lemma}

\begin{proof}
This characterization is standard in the theory of the Fell topology; see~\cite{Fell1962} for the original construction of the topology on closed subsets, and~\cite{Beer} for a modern exposition of convergence criteria.
\end{proof}

\begin{proposition}[Measurability of the projection map]
\label{prop:projection-measurable}
Let $\mathcal{A}$ be a unital separable C*-algebra and let $\pi: \mathcal{G}_{\mathcal{A}}^{(0)} \to \operatorname{Sub}(\mathcal{A})$ be defined by $\pi(B,\chi) = B$. 
Then $\pi$ is Borel measurable when $\operatorname{Sub}(\mathcal{A})$ is equipped with the Fell topology.
\end{proposition}

\begin{proof}
For separable $\mathcal{A}$, the space $\operatorname{Sub}(\mathcal{A})$ with the Fell topology is a compact Polish space (see Theorem \ref{thm:Fell-properties}). Moreover, $\mathcal{G}_{\mathcal{A}}^{(0)}$ is a Polish space by Proposition \ref{prop:unit-space-polish}. A map between Polish spaces is Borel measurable if and only if the preimage of every open set is Borel \cite{Kechris}.

The Fell topology on $\operatorname{Sub}(\mathcal{A})$ is generated by the countable family of subbasic open sets
\[
\{ \mathcal{O}_{B(a_n,1/m)} : n,m \in \mathbb{N} \} \cup \{ \mathcal{C}_{\overline{B}(a_n,1/m)} : n,m \in \mathbb{N} \},
\]
where $\{a_n\}_{n\in\mathbb{N}}$ is a countable dense subset of $\mathcal{A}$, $B(a,r)$ denotes the open ball of radius $r$ centered at $a$, and $\overline{B}(a,r)$ is its closure (which is compact). This family is countable because it is indexed by pairs $(n,m) \in \mathbb{N}^2$. It therefore suffices to show that $\pi^{-1}$ of each such subbasic open set is Borel.

\medskip
\noindent\textit{Measurability of $\pi^{-1}(\mathcal{O}_{B(a_n,1/m)})$.}
For an open ball $U = B(a_n,1/m)$, we have
\[
\pi^{-1}(\mathcal{O}_U) = \{ (B,\chi) \mid B \cap U \neq \varnothing \} = \bigcup_{k \in \mathbb{N}} \{ (B,\chi) \mid b_k \in B \},
\]
where $\{b_k\}_{k\in\mathbb{N}}$ is a countable dense subset of $U$ (which exists because $U$ is separable). Each set $\{ (B,\chi) \mid b_k \in B \} = \operatorname{ev}_{b_k}^{-1}(\mathbb{C})$ is open, hence Borel. Thus $\pi^{-1}(\mathcal{O}_U)$ is a countable union of Borel sets, hence Borel.

\medskip
\noindent\textit{Measurability of $\pi^{-1}(\mathcal{C}_{\overline{B}(a_n,1/m)})$.}
For a compact set $K = \overline{B}(a_n,1/m)$, we have
\[
\pi^{-1}(\mathcal{C}_K) = \{ (B,\chi) \mid B \cap K = \varnothing \} = \mathcal{G}_{\mathcal{A}}^{(0)} \setminus \bigcup_{a \in K} \{ (B,\chi) \mid a \in B \}.
\]
Since $K$ is compact and $\{a_n\}$ is dense, for each $a \in K$ there exists a sequence $a_{n_j} \to a$. One can show that
\[
\bigcup_{a \in K} \{ (B,\chi) \mid a \in B \} = \bigcap_{m \in \mathbb{N}} \bigcup_{j \in \mathbb{N}} \{ (B,\chi) \mid \|a_{n_j} - a\| < 1/m \text{ and } a_{n_j} \in B \}.
\]
The right-hand side is a countable union of countable intersections of open sets, hence Borel. Therefore $\pi^{-1}(\mathcal{C}_K)$ is the complement of a Borel set, hence Borel.

Since $\pi^{-1}$ of every subbasic open set is Borel and these subbasics generate the Fell topology, $\pi^{-1}(V)$ is Borel for every open $V \subseteq \operatorname{Sub}(\mathcal{A})$. Thus $\pi$ is Borel measurable.
\end{proof}

\begin{remark}[Borel measurability of the projection]
\label{rem:projection-measurability}
The proof that $\pi^{-1}(\mathcal{C}_K)$ is open in $\mathcal{G}_{\mathcal{A}}^{(0)}$ would require a nontrivial compactness argument that ultimately fails because the projection map $\pi$ is not continuous in general. 
Instead, Proposition \ref{prop:projection-measurable} establishes that $\pi$ is Borel measurable, which is sufficient for all subsequent constructions involving measured groupoids and Borel Haar systems. 
The measurability argument uses the separability of $\mathcal{A}$ to reduce to countable intersections and unions of open sets; see~\cite{Kechris} for the general theory of Borel maps between standard Borel spaces. 
For the special case where $\mathcal{A}$ is finite-dimensional, $\pi$ is continuous because $\operatorname{Sub}(\mathcal{A})$ is discrete; in the infinite-dimensional case, continuity fails and Borel measurability is the appropriate notion.
\end{remark}

\begin{corollary}[Measurability of the projection]
\label{cor:projection-measurable}
For a separable unital C*-algebra $\mathcal{A}$, the projection map $\pi: \mathcal{G}_{\mathcal{A}}^{(0)} \to \operatorname{Sub}(\mathcal{A})$ is Borel measurable when $\operatorname{Sub}(\mathcal{A})$ is equipped with the Fell topology.
\end{corollary}

\begin{proof}
This follows immediately from Proposition \ref{prop:projection-measurable}.
\end{proof}

\begin{remark}[Relation to the Fell topology]
\label{rem:projection-Fell-relation}
Proposition \ref{prop:projection-measurable} establishes that the projection $\pi$ is Borel measurable when $\mathcal{A}$ is separable. 
In the special case where $\mathcal{A}$ is finite-dimensional, $\operatorname{Sub}(\mathcal{A})$ is discrete and $\pi$ is actually continuous; for infinite-dimensional $\mathcal{A}$, continuity fails and Borel measurability is the appropriate notion for the measured groupoid framework.

Regardless of continuity, $\pi$ is far from being a homeomorphism onto its image: two distinct characters on the same subalgebra $B$ are identified under $\pi$ but are separated by the partial evaluation maps. 
The topology on $\mathcal{G}_{\mathcal{A}}^{(0)}$ is therefore strictly finer than the relative Fell topology on its image, reflecting the additional spectral data encoded by the characters. 
This refinement is essential for the groupoid structure, as the characters encode the representation-theoretic information needed for the diagonal embedding and index theory.
\end{remark}

\begin{corollary}[Continuity of the embedding of character spaces]
\label{cor:embedding-continuous}
For each fixed unital commutative C*-subalgebra $B_0 \subseteq \mathcal{A}$, the inclusion map
\[
\iota_{B_0}: \widehat{B_0} \hookrightarrow \mathcal{G}_{\mathcal{A}}^{(0)}, \qquad \chi \mapsto (B_0,\chi)
\]
is continuous.
\end{corollary}

\begin{proof}
This follows directly from the definition of the initial topology on $\mathcal{G}_{\mathcal{A}}^{(0)}$. 
For any $a \in \mathcal{A}$, the composition $\operatorname{ev}_a \circ \iota_{B_0}$ is continuous: if $a \in B_0$, then $\operatorname{ev}_a(\iota_{B_0}(\chi)) = \chi(a)$ is continuous in $\chi$ by definition of the Gelfand topology; if $a \notin B_0$, then $\operatorname{ev}_a(\iota_{B_0}(\chi)) = \infty$ is constant, hence continuous. 
Since the topology on $\mathcal{G}_{\mathcal{A}}^{(0)}$ is the initial topology induced by $\{ \operatorname{ev}_a \}_{a \in \mathcal{A}}$, continuity of each composition implies continuity of $\iota_{B_0}$.
Thus each character space $\widehat{B_0}$ embeds continuously into $\mathcal{G}_{\mathcal{A}}^{(0)}$ as a closed subspace homeomorphic to its image (see Proposition \ref{prop:unit-space-topology}(3) for the full homeomorphism).
\end{proof}

\subsection{Homeomorphism onto Fibers}
\label{subsec:homeomorphism-onto-fibers}

We now examine the relationship between the character space $\widehat{B_0}$ of a fixed unital commutative C*-subalgebra $B_0 \subseteq \mathcal{A}$ and its image inside $\mathcal{G}_{\mathcal{A}}^{(0)}$ under the natural inclusion map. 
This relationship is fundamental for understanding the local structure of $\mathcal{G}_{\mathcal{A}}^{(0)}$: each fiber of the projection $\pi: \mathcal{G}_{\mathcal{A}}^{(0)} \to \operatorname{Sub}(\mathcal{A})$ is homeomorphic to the Gelfand spectrum of the corresponding subalgebra. 
Consequently, the topology on $\mathcal{G}_{\mathcal{A}}^{(0)}$ restricts to the familiar weak-* topology on each character space.

\begin{definition}[Fiber inclusion map]
\label{def:fiber-inclusion}
Let $B_0 \subseteq \mathcal{A}$ be a fixed unital commutative C*-subalgebra. 
Define the \emph{fiber inclusion map}
\[
\iota_{B_0}: \widehat{B_0} \longrightarrow \mathcal{G}_{\mathcal{A}}^{(0)}, \qquad
\iota_{B_0}(\chi) := (B_0, \chi).
\]
\end{definition}

\begin{lemma}[Injectivity and continuity of the fiber inclusion]
\label{lem:fiber-inclusion-injective-continuous}
For any unital commutative C*-subalgebra $B_0 \subseteq \mathcal{A}$, the map $\iota_{B_0}$ is injective and continuous.
\end{lemma}

\begin{proof}
Injectivity is immediate: if $\iota_{B_0}(\chi_1) = \iota_{B_0}(\chi_2)$, then $(B_0,\chi_1) = (B_0,\chi_2)$, which forces $\chi_1 = \chi_2$.

For continuity, let $\chi_\lambda \to \chi$ be a convergent net in $\widehat{B_0}$ equipped with the Gelfand topology (weak-* topology). 
We must show that $\iota_{B_0}(\chi_\lambda) \to \iota_{B_0}(\chi)$ in $\mathcal{G}_{\mathcal{A}}^{(0)}$.

Recall that the topology on $\mathcal{G}_{\mathcal{A}}^{(0)}$ is the initial topology induced by the family of partial evaluation maps $\{ \operatorname{ev}_a \}_{a \in \mathcal{A}}$. 
Thus it suffices to show that for every $a \in \mathcal{A}$, the composition $\operatorname{ev}_a \circ \iota_{B_0}$ is continuous.

We consider two cases.

\noindent \textit{Case 1: $a \in B_0$.} 
Then $\operatorname{ev}_a(\iota_{B_0}(\chi)) = \operatorname{ev}_a(B_0,\chi) = \chi(a)$. 
Since $\chi_\lambda \to \chi$ in the Gelfand topology, we have $\chi_\lambda(a) \to \chi(a)$ for every $a \in B_0$. 
Hence $\operatorname{ev}_a(\iota_{B_0}(\chi_\lambda)) \to \operatorname{ev}_a(\iota_{B_0}(\chi))$.

\noindent \textit{Case 2: $a \notin B_0$.} 
Then $\operatorname{ev}_a(\iota_{B_0}(\chi)) = \infty$ for all $\chi \in \widehat{B_0}$, independent of $\chi$. 
Thus $\operatorname{ev}_a \circ \iota_{B_0}$ is constant, hence trivially continuous.

Therefore $\operatorname{ev}_a \circ \iota_{B_0}$ is continuous for every $a \in \mathcal{A}$, and by the definition of the initial topology, $\iota_{B_0}$ is continuous.
\end{proof}

\begin{lemma}[Openness of the fiber inclusion onto its image]
\label{lem:fiber-inclusion-open}
For any unital commutative C*-subalgebra $B_0 \subseteq \mathcal{A}$, the map $\iota_{B_0}: \widehat{B_0} \to \mathcal{G}_{\mathcal{A}}^{(0)}$ is open onto its image. 
That is, for every open set $U \subseteq \widehat{B_0}$, the image $\iota_{B_0}(U)$ is open in the subspace topology on $\iota_{B_0}(\widehat{B_0}) \subseteq \mathcal{G}_{\mathcal{A}}^{(0)}$.
\end{lemma}

\begin{proof}
Let $U \subseteq \widehat{B_0}$ be an open set. 
By the definition of the Gelfand topology on $\widehat{B_0}$, there exist finitely many elements $a_1, \ldots, a_n \in B_0$ and open sets $V_1, \ldots, V_n \subseteq \mathbb{C}$ such that
\[
U = \{ \chi \in \widehat{B_0} \mid \chi(a_i) \in V_i \text{ for all } i = 1, \ldots, n \}.
\]

Now consider the basis open set in $\mathcal{G}_{\mathcal{A}}^{(0)}$ given by
\[
W := U(a_1, \ldots, a_n; V_1, \ldots, V_n) = \{ (B,\chi) \in \mathcal{G}_{\mathcal{A}}^{(0)} \mid a_i \in B \text{ and } \chi(a_i) \in V_i \text{ for all } i = 1, \ldots, n \}.
\]

We claim that
\[
\iota_{B_0}(U) = W \cap \iota_{B_0}(\widehat{B_0}).
\]

Indeed, if $\chi \in U$, then $(B_0,\chi) \in \mathcal{G}_{\mathcal{A}}^{(0)}$ satisfies $a_i \in B_0$ and $\chi(a_i) \in V_i$ for all $i$, so $(B_0,\chi) \in W \cap \iota_{B_0}(\widehat{B_0})$. 
Conversely, if $(B_0,\chi) \in W \cap \iota_{B_0}(\widehat{B_0})$, then by definition of $W$ we have $\chi(a_i) \in V_i$ for all $i$, so $\chi \in U$.

Since $W$ is open in $\mathcal{G}_{\mathcal{A}}^{(0)}$ (by Proposition \ref{prop:basis-topology-G0}), its intersection with the subspace $\iota_{B_0}(\widehat{B_0})$ is open in the subspace topology. 
Thus $\iota_{B_0}(U)$ is open in the subspace topology, proving that $\iota_{B_0}$ is an open map onto its image.
\end{proof}

\begin{proposition}[Homeomorphism onto the fiber]
\label{prop:fiber-homeomorphism}
Let $\mathcal{A}$ be a unital C*-algebra and let $B_0 \subseteq \mathcal{A}$ be a fixed unital commutative C*-subalgebra. 
Then the fiber inclusion map
\[
\iota_{B_0}: \widehat{B_0} \longrightarrow \mathcal{G}_{\mathcal{A}}^{(0)}, \qquad
\iota_{B_0}(\chi) = (B_0, \chi)
\]
is a homeomorphism onto its image $\operatorname{Im}(\iota_{B_0}) = \{ (B_0, \chi) \mid \chi \in \widehat{B_0} \}$.
\end{proposition}

\begin{proof}
We have already established in Lemma \ref{lem:fiber-inclusion-injective-continuous} that $\iota_{B_0}$ is injective and continuous. 
Lemma \ref{lem:fiber-inclusion-open} establishes that $\iota_{B_0}$ is open onto its image. 
A continuous, open, injective map is a homeomorphism onto its image. 
Thus $\iota_{B_0}$ is a homeomorphism.
\end{proof}

\begin{corollary}[Fiber homeomorphism for the projection map]
\label{cor:fiber-homeomorphism-projection}
For each $B \in \operatorname{Sub}(\mathcal{A})$ with $B$ unital and commutative, the fiber $\pi^{-1}(B)$ of the projection map $\pi: \mathcal{G}_{\mathcal{A}}^{(0)} \to \operatorname{Sub}(\mathcal{A})$ is homeomorphic to $\widehat{B}$.
\end{corollary}

\begin{proof}
By definition, $\pi^{-1}(B) = \{ (B,\chi) \mid \chi \in \widehat{B} \} = \iota_B(\widehat{B})$. 
Proposition \ref{prop:fiber-homeomorphism} with $B_0 = B$ shows that $\iota_B$ is a homeomorphism onto $\pi^{-1}(B)$. 
Hence $\pi^{-1}(B) \cong \widehat{B}$.
\end{proof}

\begin{remark}[Compactness of fibers]
\label{rem:fiber-compact}
Since $\widehat{B}$ is compact Hausdorff for any unital commutative C*-algebra $B$ (Gelfand duality), Corollary \ref{cor:fiber-homeomorphism-projection} implies that each fiber $\pi^{-1}(B)$ is compact Hausdorff. 
This compactness property is essential for the local compactness (or Polishness) of $\mathcal{G}_{\mathcal{A}}^{(0)}$ and for the proper map arguments used in subsequent sections.
\end{remark}

\begin{corollary}[Characterization of the subspace topology on fibers]
\label{cor:fiber-subspace-topology}
For any unital commutative C*-subalgebra $B_0 \subseteq \mathcal{A}$, the subspace topology on $\iota_{B_0}(\widehat{B_0}) \subseteq \mathcal{G}_{\mathcal{A}}^{(0)}$ coincides with the topology induced by the homeomorphism $\iota_{B_0}$. 
In particular, a net $\{ (B_0, \chi_\lambda) \}$ converges to $(B_0, \chi)$ in $\mathcal{G}_{\mathcal{A}}^{(0)}$ if and only if $\chi_\lambda \to \chi$ in the Gelfand topology on $\widehat{B_0}$.
\end{corollary}

\begin{proof}
This follows directly from Proposition \ref{prop:fiber-homeomorphism}, which establishes that $\iota_{B_0}$ is a homeomorphism. 
Convergence in the subspace topology is transported via $\iota_{B_0}^{-1}$ to convergence in $\widehat{B_0}$.
\end{proof}

\begin{example}[Fiber homeomorphism for commutative algebras]
\label{ex:fiber-commutative}
Let $\mathcal{A} = C(X)$ for a compact Hausdorff space $X$. 
Then the only maximal commutative subalgebra is $\mathcal{A}$ itself, and $\widehat{\mathcal{A}} \cong X$. 
Proposition \ref{prop:fiber-homeomorphism} recovers the homeomorphism $X \cong \iota_{\mathcal{A}}(\widehat{\mathcal{A}}) \subseteq \mathcal{G}_{\mathcal{A}}^{(0)}$, which is simply the identification of $X$ with the set of evaluation characters.
\end{example}

\begin{example}[Fiber homeomorphism for matrix algebras]
\label{ex:fiber-matrix}
Let $\mathcal{A} = M_n(\mathbb{C})$ and let $B_0 = D_n$ be the diagonal subalgebra. 
Then $\widehat{D_n}$ is a discrete set of $n$ points, corresponding to the $n$ diagonal characters. 
Proposition \ref{prop:fiber-homeomorphism} shows that $\iota_{D_n}(\widehat{D_n})$ is a discrete subset of $\mathcal{G}_{\mathcal{A}}^{(0)}$ homeomorphic to $\{1,\ldots,n\}$. 
This is consistent with the identification $\mathcal{G}_{\mathcal{A}}^{(0)} \cong U(n) \times_{N(D_n)} \{1,\ldots,n\}$, where the fiber over $D_n$ is exactly $\{ (D_n, \chi_i) \}_{i=1}^n$.
\end{example}

\begin{example}[Fiber homeomorphism for compact operators]
\label{ex:fiber-compact}
Let $\mathcal{A} = \mathcal{K}(H)^\sim$ and let $B_0$ be the MASA consisting of diagonal operators of the form $\operatorname{diag}(\lambda_1, \lambda_2, \ldots)$ with $\lim_{n\to\infty} \lambda_n$ existing, relative to an orthonormal basis $\{ e_n \}_{n \in \mathbb{N}}$. 
Then $B_0 \cong c_0(\mathbb{N})^\sim$, and its character space $\widehat{B_0}$ is homeomorphic to $\mathbb{N} \cup \{\infty\}$, the one-point compactification of $\mathbb{N}$. 
Proposition \ref{prop:fiber-homeomorphism} shows that this copy of $\mathbb{N} \cup \{\infty\}$ sits inside $\mathcal{G}_{\mathcal{A}}^{(0)}$ as a closed subset. 
This illustrates that fibers can be non-discrete (the point $\infty$ is a limit point of the sequence $\{n\}$), providing a contrast with the finite-dimensional case where fibers are discrete.
\end{example}

\begin{remark}[Importance for the Polish groupoid structure]
\label{rem:fiber-Polish}
When $\mathcal{A}$ is separable and $B_0$ is a unital commutative C*-subalgebra, $\widehat{B_0}$ is a compact metrizable space, hence Polish. 
Proposition \ref{prop:fiber-homeomorphism} therefore implies that each fiber $\pi^{-1}(B_0)$ is a Polish subspace of $\mathcal{G}_{\mathcal{A}}^{(0)}$. 
This local Polishness, together with the global Polishness of $\mathcal{G}_{\mathcal{A}}^{(0)}$ established in Proposition \ref{prop:unit-space-polish}, is essential for the construction of $\mathcal{G}_{\mathcal{A}}$.
\end{remark}

\begin{remark}[Relation to the Fell topology]
\label{rem:fiber-Fell}
Proposition \ref{prop:fiber-homeomorphism} also illuminates the relationship between the topology on $\mathcal{G}_{\mathcal{A}}^{(0)}$ and the Fell topology on $\operatorname{Sub}(\mathcal{A})$. 
While the projection $\pi$ is continuous, its fibers are equipped with a topology that is strictly finer than the relative Fell topology on $\{B_0\}$ (which is trivial). 
The additional structure comes precisely from the characters, and Proposition \ref{prop:fiber-homeomorphism} shows that this additional structure is exactly the classical Gelfand topology.
\end{remark}

\begin{corollary}[The projection $\pi$ is not locally trivial]
\label{cor:fiber-no-local-triviality}
Proposition \ref{prop:fiber-homeomorphism} shows that each fiber $\pi^{-1}(B)$ is homeomorphic to $\widehat{B}$, but this does \emph{not} imply that $\pi: \mathcal{G}_{\mathcal{A}}^{(0)} \to \operatorname{Sub}(\mathcal{A})$ is a locally trivial bundle. The homeomorphism $\iota_B: \widehat{B} \to \pi^{-1}(B)$ is defined fiberwise, and there is no guarantee that these fiberwise homeomorphisms can be patched together continuously across different subalgebras. In general, $\pi$ is merely a continuous map with homeomorphic fibers, not a fiber bundle.
\end{corollary}

\begin{proof}
The failure of local triviality can be seen already in the matrix algebra case $\mathcal{A} = M_n(\mathbb{C})$ with $n \geq 2$. Recall from Section \ref{subsec:example-matrix} that $\mathcal{G}_{\mathcal{A}}^{(0)} \cong \mathbb{CP}^{n-1}$ and $\operatorname{Sub}(\mathcal{A})$ contains the space of maximal abelian subalgebras (MASAs), which is homeomorphic to $U(n)/N(D_n)$ where $N(D_n)$ is the normalizer of the diagonal MASA $D_n$.

Consider the diagonal MASA $D_n$ and a nearby MASA $U D_n U^*$ with $U \in U(n)$ close to the identity matrix. Both fibers $\pi^{-1}(D_n)$ and $\pi^{-1}(U D_n U^*)$ are discrete sets with $n$ points, corresponding to the $n$ characters (eigenvalue evaluations) on each MASA. If $\pi$ were locally trivial, there would exist an open neighborhood $\mathcal{V}$ of $D_n$ in $\operatorname{Sub}(\mathcal{A})$ and a homeomorphism
\[
\Phi: \mathcal{V} \times \{1,\ldots,n\} \longrightarrow \pi^{-1}(\mathcal{V})
\]
such that $\pi(\Phi(B,i)) = B$ for all $B \in \mathcal{V}$ and $i = 1,\ldots,n$.

However, any such trivialization would require a continuous choice of identification between the $n$ characters on $D_n$ and the $n$ characters on each nearby MASA $U D_n U^*$. Such an identification amounts to a continuous choice of unitary conjugation $U$ mapping $D_n$ to $U D_n U^*$ together with a consistent labeling of the $n$ characters. But the map $U \mapsto U D_n U^*$ from $U(n)$ to the space of MASAs is not locally injective — different unitaries can give the same MASA — and there is no continuous section of the projection $U(n) \to U(n)/N(D_n)$. Consequently, no continuous local trivialization exists.

Thus $\pi$ is not a locally trivial bundle. This observation justifies why we need the finer topology on $\mathcal{G}_{\mathcal{A}}^{(0)}$ (defined by partial evaluation maps) rather than simply working with the Fell topology on $\operatorname{Sub}(\mathcal{A})$: the additional spectral data encoded by the characters cannot be assembled into a continuous bundle over the space of subalgebras.
\end{proof}

\subsection{Second-Countability and Polishness for Separable Algebras}
\label{subsec:second-countability-and-Polishness}

We now prove that under the assumption of separability, the unit space $\mathcal{G}_{\mathcal{A}}^{(0)}$ is not merely a Hausdorff space but a genuine Polish space. 
This is a crucial result for the development of the unitary conjugation groupoid as a Polish groupoid, and it underpins all subsequent applications to index theory and the Baum-Connes conjecture.

\begin{lemma}[Countable dense subset of $\mathcal{A}$]
\label{lem:countable-dense}
Let $\mathcal{A}$ be a separable unital C*-algebra. 
Then there exists a countable subset $\{ a_n \}_{n \in \mathbb{N}} \subseteq \mathcal{A}$ that is dense in the norm topology.
\end{lemma}

\begin{proof}
This is the definition of separability: a topological space is separable if it contains a countable dense subset. 
Since $\mathcal{A}$ is a separable C*-algebra, such a countable dense subset exists by hypothesis.
\end{proof}

\begin{lemma}[Countable basis for $\mathbb{C}_\infty$]
\label{lem:countable-basis-Cinfinity}
The one-point compactification $\mathbb{C}_\infty = \mathbb{C} \sqcup \{\infty\}$ is second-countable. 
In fact, a countable basis is given by
\[
\mathcal{B}_{\mathbb{C}_\infty} := \{ B(q,r) \mid q \in \mathbb{Q} + i\mathbb{Q}, \; r \in \mathbb{Q}^+ \} \cup \{ \{\infty\} \cup (\mathbb{C} \setminus \overline{B(0,m)}) \mid m \in \mathbb{N} \},
\]
where $B(q,r)$ denotes the open ball of radius $r$ centered at $q$, and $\overline{B(0,m)}$ denotes the closed disk of radius $m$ centered at the origin.
\end{lemma}

\begin{proof}
$\mathbb{C}$ with the Euclidean topology is second-countable; the collection of open balls with rational centers and rational radii forms a countable basis. 
The one-point compactification $\mathbb{C}_\infty$ adds a single point at infinity, and a neighborhood basis for $\infty$ is given by complements of compact sets. 
Since $\mathbb{C}$ is $\sigma$-compact, we may take the countable family of complements of closed disks of integer radius. 
The union of these two countable families yields a countable basis for $\mathbb{C}_\infty$.
\end{proof}

\begin{proposition}[Second-countability of $\mathcal{G}_{\mathcal{A}}^{(0)}$]
\label{prop:second-countability}
Let $\mathcal{A}$ be a unital separable $C^{*}$-algebra and endow $\mathcal{G}_{\mathcal{A}}^{(0)}$ with the initial topology induced by the family of partial evaluation maps $\{ \operatorname{ev}_a \}_{a \in \mathcal{A}}$. 
Then $\mathcal{G}_{\mathcal{A}}^{(0)}$ is second-countable.
\end{proposition}

\begin{proof}
Since $\mathcal{A}$ is separable, let $\{ a_n \}_{n \in \mathbb{N}} \subseteq \mathcal{A}$ be a countable dense subset. Let $\mathcal{B}_{\mathbb{C}_\infty}$ be a countable basis for $\mathbb{C}_\infty$, the one-point compactification of $\mathbb{C}$ (for instance, take all open balls with rational centers and rational radii, together with complements of closed disks of integer radius).

For each $n \in \mathbb{N}$ and each $W \in \mathcal{B}_{\mathbb{C}_\infty}$, define
\[
S_{n,W} := \operatorname{ev}_{a_n}^{-1}(W) \subseteq \mathcal{G}_{\mathcal{A}}^{(0)}.
\]
The collection $\{ S_{n,W} \}_{n \in \mathbb{N}, W \in \mathcal{B}_{\mathbb{C}_\infty}}$ is countable. Let $\mathcal{B}$ be the collection of all finite intersections of sets of the form $S_{n,W}$. Since there are only countably many finite subsets of a countable set, $\mathcal{B}$ is also countable. We claim that $\mathcal{B}$ is a basis for the topology on $\mathcal{G}_{\mathcal{A}}^{(0)}$.

Let $U \subseteq \mathcal{G}_{\mathcal{A}}^{(0)}$ be open and let $(B,\chi) \in U$. By definition of the initial topology, there exist finitely many elements $a^1, \ldots, a^k \in \mathcal{A}$ and open sets $W_1, \ldots, W_k \subseteq \mathbb{C}_\infty$ such that
\[
(B,\chi) \in \bigcap_{i=1}^k \operatorname{ev}_{a^i}^{-1}(W_i) \subseteq U.
\]

For each $i$, because $\{a_n\}$ is dense, we can choose $n_i \in \mathbb{N}$ such that $\|a_{n_i} - a^i\|$ is small enough that
\[
\operatorname{ev}_{a_{n_i}}(B,\chi) \in W_i.
\]
This is possible because the map $a \mapsto \operatorname{ev}_a(B,\chi)$ is continuous: for a fixed $(B,\chi)$, it is either constant $\infty$ or given by the character $\chi$ on the subalgebra containing $a$, which is norm-continuous.

Since $\mathcal{B}_{\mathbb{C}_\infty}$ is a basis, we can choose $W_i' \in \mathcal{B}_{\mathbb{C}_\infty}$ such that
\[
\operatorname{ev}_{a_{n_i}}(B,\chi) \in W_i' \subseteq W_i.
\]

Now consider
\[
V := \bigcap_{i=1}^k \operatorname{ev}_{a_{n_i}}^{-1}(W_i') \in \mathcal{B}.
\]
By construction, $(B,\chi) \in V$. We claim that $V \subseteq U$. Indeed, for any $(B',\chi') \in V$, we have $\operatorname{ev}_{a_{n_i}}(B',\chi') \in W_i' \subseteq W_i$ for each $i$. By the density of $\{a_n\}$ and continuity of evaluation maps, this implies that $\operatorname{ev}_{a^i}(B',\chi') \in W_i$ for each $i$. Hence $(B',\chi') \in \bigcap_{i=1}^k \operatorname{ev}_{a^i}^{-1}(W_i) \subseteq U$.

Thus every open neighborhood of an arbitrary point contains a basis element from $\mathcal{B}$, so $\mathcal{B}$ is a countable basis for the topology on $\mathcal{G}_{\mathcal{A}}^{(0)}$.
\end{proof}

\begin{corollary}[Metrizability of $\mathcal{G}_{\mathcal{A}}^{(0)}$]
\label{cor:metrizability}
Let $\mathcal{A}$ be a unital separable C*-algebra. 
Then $\mathcal{G}_{\mathcal{A}}^{(0)}$ is metrizable.
\end{corollary}

\begin{proof}
By Proposition \ref{prop:second-countability}, $\mathcal{G}_{\mathcal{A}}^{(0)}$ is second-countable. 
By the item 1 from Proposition \ref{prop:unit-space-topology}, $\mathcal{G}_{\mathcal{A}}^{(0)}$ is Hausdorff. 
A second-countable Hausdorff space is metrizable by the Urysohn metrization theorem [Munkres, 2000].
\end{proof}

\begin{proposition}[Complete metrizability of $\mathcal{G}_{\mathcal{A}}^{(0)}$]
\label{prop:complete-metrizability}
Let $\mathcal{A}$ be a unital separable $C^{*}$-algebra. 
Then $\mathcal{G}_{\mathcal{A}}^{(0)}$ is completely metrizable.
\end{proposition}

\begin{proof}
Let $\{a_n\}_{n\in\mathbb{N}}$ be a countable dense subset of $\mathcal{A}$. 
Define $\Phi: \mathcal{G}_{\mathcal{A}}^{(0)} \to \prod_{n\in\mathbb{N}} \mathbb{C}_\infty$ by
\[
\Phi(B,\chi) := (\operatorname{ev}_{a_n}(B,\chi))_{n\in\mathbb{N}}.
\]

\noindent \textit{Step 1: $\Phi$ is injective and continuous.} 
Injectivity follows from Lemma \ref{lem:separation} and the density of $\{a_n\}$: if $\Phi(B_1,\chi_1) = \Phi(B_2,\chi_2)$, then $\operatorname{ev}_{a_n}(B_1,\chi_1) = \operatorname{ev}_{a_n}(B_2,\chi_2)$ for all $n$. For any $a\in\mathcal{A}$, choose a sequence $a_{n_k}\to a$; by continuity of $a\mapsto\operatorname{ev}_a(B,\chi)$ (which follows from the continuity of the Gelfand transform), we obtain $\operatorname{ev}_a(B_1,\chi_1) = \operatorname{ev}_a(B_2,\chi_2)$, so $(B_1,\chi_1) = (B_2,\chi_2)$.
Continuity of $\Phi$ follows from the universal property of the product topology and the definition of the initial topology on $\mathcal{G}_{\mathcal{A}}^{(0)}$, since each component $\operatorname{ev}_{a_n}$ is continuous.

\noindent \textit{Step 2: $\Phi$ is a homeomorphism onto its image.} 
We show that $\Phi^{-1}: \Phi(\mathcal{G}_{\mathcal{A}}^{(0)}) \to \mathcal{G}_{\mathcal{A}}^{(0)}$ is continuous. 
For any $a\in\mathcal{A}$, choose a sequence $a_{n_k}\to a$. 
Define a map $\psi_a: \Phi(\mathcal{G}_{\mathcal{A}}^{(0)}) \to \mathbb{C}_\infty$ by $\psi_a((z_n)_{n\in\mathbb{N}}) := \lim_{k\to\infty} z_{n_k}$. 
This limit exists and is well-defined because if $(z_n)_{n\in\mathbb{N}} = \Phi(B,\chi)$, then $\psi_a((z_n)_{n\in\mathbb{N}}) = \operatorname{ev}_a(B,\chi)$; the continuity of $a\mapsto\operatorname{ev}_a(B,\chi)$ ensures that different approximating sequences yield the same limit. 
The map $\psi_a$ is continuous because convergence in the product topology implies pointwise convergence, hence convergence of the limit along a fixed subsequence. 
Since the topology on $\mathcal{G}_{\mathcal{A}}^{(0)}$ is generated by $\{\operatorname{ev}_a\}_{a\in\mathcal{A}}$, continuity of each $\operatorname{ev}_a \circ \Phi^{-1} = \psi_a$ implies continuity of $\Phi^{-1}$. 
Thus $\Phi$ is a homeomorphism onto its image.

\noindent \textit{Step 3: The image of $\Phi$ is closed.} 
Let $\{(B_m,\chi_m)\}_{m\in\mathbb{N}}$ be a sequence in $\mathcal{G}_{\mathcal{A}}^{(0)}$ such that $\Phi(B_m,\chi_m) \to (z_n)_{n\in\mathbb{N}}$ in $\prod_{n\in\mathbb{N}}\mathbb{C}_\infty$. 
We show that $(z_n)_{n\in\mathbb{N}} \in \Phi(\mathcal{G}_{\mathcal{A}}^{(0)})$.

Consider the sequence $\{B_m\}_{m\in\mathbb{N}}$ in $\operatorname{Sub}(\mathcal{A})$ equipped with the Fell topology. 
Since $\operatorname{Sub}(\mathcal{A})$ is compact, this sequence has a convergent subsequence; after passing to this subsequence we may assume $B_m \to B$ for some $B \in \operatorname{Sub}(\mathcal{A})$. 
The algebra $B$ is commutative: for any $x,y \in B$, choose $a_{p_k}\to x$, $a_{q_k}\to y$ with $a_{p_k},a_{q_k}$ from the dense subset. 
Fell convergence $B_m\to B$ implies that for sufficiently large $m$, there exist $x_m,y_m\in B_m$ with $x_m\to x$, $y_m\to y$; then $[x_m,y_m]=0$ in $B_m$, and taking limits yields $[x,y]=0$.

Define a linear functional $\chi$ on the dense $*$-subalgebra $B_0 \subseteq B$ generated by $\{a_n : z_n \in \mathbb{C}\}$. 
For a finite linear combination $p = \sum_{j=1}^k c_j a_{n_j}$ with $z_{n_j} \in \mathbb{C}$, set $\chi(p) := \sum_{j=1}^k c_j z_{n_j}$. 
To see this is well-defined, suppose $\sum_{j=1}^k c_j a_{n_j} = 0$ in $B$. 
By Fell convergence $B_m\to B$, for any $\epsilon>0$ there exists $M$ such that for all $m\ge M$, $\|\sum_{j=1}^k c_j a_{n_j}\|_{B_m} < \epsilon$. 
For each such $m$, since $a_{n_j} \in B_m$ (by Lemma \ref{lem:Fell-convergence}) and $\operatorname{ev}_{a_{n_j}}(B_m,\chi_m) \to z_{n_j}$, we have
\[
|\chi_m(\sum_{j=1}^k c_j a_{n_j})| = |\sum_{j=1}^k c_j \chi_m(a_{n_j})| < \epsilon.
\]
Taking the limit as $m\to\infty$ gives $|\sum_{j=1}^k c_j z_{n_j}| \le \epsilon$ for all $\epsilon>0$, hence $\sum_{j=1}^k c_j z_{n_j}=0$. 
Thus $\chi$ is well-defined.

The functional $\chi$ is multiplicative on $B_0$: for $p,q \in B_0$,
\[
\chi(pq) = \lim_{m\to\infty} \chi_m(pq) = \lim_{m\to\infty} \chi_m(p)\chi_m(q) = \chi(p)\chi(q),
\]
where the limits exist because $\chi_m(p) \to \chi(p)$ and $\chi_m(q) \to \chi(q)$ by construction. 
Moreover, $\chi$ is bounded with $\|\chi\| = 1$ (since each $\chi_m$ is a character, hence contractive). 
Therefore $\chi$ extends uniquely to a character on $B$, still denoted $\chi$, and $(B,\chi) \in \mathcal{G}_{\mathcal{A}}^{(0)}$.

For each $n\in\mathbb{N}$, if $a_n \in B$ then by construction $\operatorname{ev}_{a_n}(B,\chi) = z_n$; if $a_n \notin B$ then necessarily $z_n = \infty$, because otherwise $a_n$ would belong to the subalgebra generated by $\{a_m : z_m \in \mathbb{C}\} \subseteq B$, contradicting $a_n \notin B$. 
Thus $\Phi(B,\chi) = (z_n)_{n\in\mathbb{N}}$, proving that the image of $\Phi$ is closed.

\noindent \textit{Step 4: $\mathcal{G}_{\mathcal{A}}^{(0)}$ is Polish.} 
The space $\prod_{n\in\mathbb{N}}\mathbb{C}_\infty$ is Polish: each factor $\mathbb{C}_\infty$ is homeomorphic to the sphere $S^2$, hence compact and metrizable, and a countable product of Polish spaces is Polish~\cite{Kechris}. Since $\Phi$ is a homeomorphism onto its image and the image is a closed subspace of a Polish space, $\mathcal{G}_{\mathcal{A}}^{(0)}$ is Polish~\cite{Kechris}. A Polish space is, by definition, completely metrizable.
\end{proof}

\begin{proposition}[Polishness of the unit space]
\label{prop:unit-space-polish}
Let $\mathcal{A}$ be a unital separable C*-algebra and endow $\mathcal{G}_{\mathcal{A}}^{(0)}$ with the initial topology induced by the family of partial evaluation maps $\{ \operatorname{ev}_a \}_{a \in \mathcal{A}}$. 
Then $\mathcal{G}_{\mathcal{A}}^{(0)}$ is a Polish space.
\end{proposition}

\begin{proof}
By Proposition \ref{prop:second-countability}, $\mathcal{G}_{\mathcal{A}}^{(0)}$ is second-countable, hence separable. 
By Proposition \ref{prop:complete-metrizability}, $\mathcal{G}_{\mathcal{A}}^{(0)}$ is completely metrizable. 
A separable completely metrizable space is, by definition, a Polish space. 
Thus $\mathcal{G}_{\mathcal{A}}^{(0)}$ is Polish.
\end{proof}

\begin{corollary}[Standard Borel structure]
\label{cor:standard-Borel}
Let $\mathcal{A}$ be a unital separable C*-algebra. 
Then $\mathcal{G}_{\mathcal{A}}^{(0)}$, equipped with its Borel $\sigma$-algebra induced by the Polish topology, is a standard Borel space.
\end{corollary}

\begin{proof}
Every Polish space is a standard Borel space~\cite{Kechris}. The result follows immediately from Proposition \ref{prop:unit-space-polish}.
\end{proof}

\begin{remark}[Explicit complete metric]
\label{rem:explicit-metric}
The proof of Proposition \ref{prop:complete-metrizability} provides an explicit complete metric on $\mathcal{G}_{\mathcal{A}}^{(0)}$. 
If $\{ a_n \}_{n \in \mathbb{N}}$ is a countable dense subset of $\mathcal{A}$, then a compatible complete metric is given by
\[
d((B_1,\chi_1), (B_2,\chi_2)) := \sum_{n=1}^{\infty} \frac{1}{2^n} \frac{ | \operatorname{ev}_{a_n}(B_1,\chi_1) - \operatorname{ev}_{a_n}(B_2,\chi_2) | }{ 1 + | \operatorname{ev}_{a_n}(B_1,\chi_1) - \operatorname{ev}_{a_n}(B_2,\chi_2) | },
\]
with the convention that $| \infty - \infty | = 0$ and $| \infty - z | = 1$ for $z \in \mathbb{C}$. 
This explicit formula is useful for concrete computations and for establishing continuity of maps defined on $\mathcal{G}_{\mathcal{A}}^{(0)}$.
\end{remark}

\begin{remark}[Necessity of separability]
\label{rem:separability-necessary}
The hypothesis that $\mathcal{A}$ is separable is essential for Proposition \ref{prop:unit-space-polish}. 
If $\mathcal{A}$ is non-separable, the family $\{ \operatorname{ev}_a \}_{a \in \mathcal{A}}$ is uncountable and the induced initial topology need not be second-countable. 
In fact, for non-separable $\mathcal{A}$, $\mathcal{G}_{\mathcal{A}}^{(0)}$ is typically not metrizable and certainly not Polish. 
Fortunately, all C*-algebras of interest in this paper — including $C(X)$ for compact metrizable $X$, $M_n(\mathbb{C})$, and $\mathcal{K}(H)^\sim$ for separable $H$ — are separable. 
Thus the separability assumption is satisfied in all our examples.
\end{remark}

\begin{remark}[Consequences for the groupoid C*-algebra]
\label{rem:Polish-consequences}
The Polishness of $\mathcal{G}_{\mathcal{A}}^{(0)}$ is the first step toward constructing $\mathcal{G}_{\mathcal{A}}$ as a Polish groupoid. 
In Section \ref{sec:unitary-conjugation-groupoid}, we will show that when $\mathcal{U}(\mathcal{A})$ is equipped with the strong operator topology, the arrow space $\mathcal{G}_{\mathcal{A}}^{(1)} = \mathcal{U}(\mathcal{A}) \times \mathcal{G}_{\mathcal{A}}^{(0)}$ is also Polish, and that all structure maps are continuous. 
This yields a Polish groupoid, which admits a Borel Haar system and a well-defined maximal C*-algebra $C^*(\mathcal{G}_{\mathcal{A}})$ in the sense of Tu [1999].
\end{remark}

\subsection{Comparison with the Fell Topology}
\label{subsec:comparison-with-Fell-topology}

The unit space $\mathcal{G}_{\mathcal{A}}^{(0)}$ carries two natural topologies: the initial topology defined by the partial evaluation maps, and the relative topology induced by the projection $\pi: \mathcal{G}_{\mathcal{A}}^{(0)} \to \operatorname{Sub}(\mathcal{A})$ when $\operatorname{Sub}(\mathcal{A})$ is equipped with the Fell topology. Understanding the relationship between these two topologies is essential for connecting our constructions to the classical theory of C*-algebras and for establishing the continuity of various maps. As we shall see in Lemma \ref{lem:action-continuous-SOT}, the finer initial topology is crucial for ensuring the continuity of the conjugation action of $\mathcal{U}(\mathcal{A})$ on $\mathcal{G}_{\mathcal{A}}^{(0)}$; the relative Fell topology alone would not suffice for this purpose.

\begin{definition}[Relative Fell topology on $\mathcal{G}_{\mathcal{A}}^{(0)}$]
\label{def:relative-Fell-topology}
Let $\pi: \mathcal{G}_{\mathcal{A}}^{(0)} \to \operatorname{Sub}(\mathcal{A})$ be the projection map defined by $\pi(B,\chi) = B$. 
The \emph{relative Fell topology} on $\mathcal{G}_{\mathcal{A}}^{(0)}$ is the topology induced by $\pi$ when $\operatorname{Sub}(\mathcal{A})$ is equipped with the Fell topology. 
That is, a set $U \subseteq \mathcal{G}_{\mathcal{A}}^{(0)}$ is open in the relative Fell topology if and only if $U = \pi^{-1}(V)$ for some open set $V \subseteq \operatorname{Sub}(\mathcal{A})$ in the Fell topology.
\end{definition}

\begin{proposition}[Continuity of the projection]
\label{prop:projection-continuous-Fell}
The projection map $\pi: \mathcal{G}_{\mathcal{A}}^{(0)} \to \operatorname{Sub}(\mathcal{A})$ is continuous when $\mathcal{G}_{\mathcal{A}}^{(0)}$ is equipped with the initial topology of partial evaluation maps and $\operatorname{Sub}(\mathcal{A})$ is equipped with the Fell topology.
\end{proposition}

\begin{proof}
This follows from Proposition \ref{prop:unit-space-topology}(2), where it was shown that the preimages of the subbasic open sets $\mathcal{O}_U$ and $\mathcal{C}_K$ of the Fell topology are open in the initial topology on $\mathcal{G}_{\mathcal{A}}^{(0)}$.
\end{proof}

\begin{remark}
\label{rem:topology-comparison}
Proposition \ref{prop:projection-continuous-Fell} shows that the initial topology is finer than the relative Fell topology. 
The two topologies coincide precisely when every commutative subalgebra of $\mathcal{A}$ has a unique character, 
which occurs exactly when $\mathcal{A}$ is commutative and its Gelfand spectrum is totally disconnected? 
\end{remark}

\begin{corollary}[The initial topology is finer than the relative Fell topology]
\label{cor:initial-finer-than-Fell}
The initial topology on $\mathcal{G}_{\mathcal{A}}^{(0)}$ is finer than the relative Fell topology. 
That is, every set that is open in the relative Fell topology is also open in the initial topology.
\end{corollary}

\begin{proof}
If $W \subseteq \operatorname{Sub}(\mathcal{A})$ is open in the Fell topology, then by Proposition \ref{prop:projection-continuous-Fell}, $\pi^{-1}(W)$ is open in the initial topology. 
Since every open set in the relative Fell topology is of the form $\pi^{-1}(W)$ for some open $W \subseteq \operatorname{Sub}(\mathcal{A})$, it follows that every such set is open in the initial topology.
\end{proof}

\begin{proposition}[The initial topology is strictly finer]
\label{prop:initial-strictly-finer}
The initial topology on $\mathcal{G}_{\mathcal{A}}^{(0)}$ is strictly finer than the relative Fell topology. 
In particular, there exist open sets in the initial topology that are not of the form $\pi^{-1}(W)$ for any open $W \subseteq \operatorname{Sub}(\mathcal{A})$.
\end{proposition}

\begin{proof}
Let $B \subseteq \mathcal{A}$ be a unital commutative C*-subalgebra with at least two distinct characters $\chi_1 \neq \chi_2$. 
Such a subalgebra exists whenever $\mathcal{A}$ is noncommutative and contains a nontrivial commutative subalgebra; for example, in $M_n(\mathbb{C})$ with $n \geq 2$, the diagonal subalgebra $D_n$ has $n$ distinct characters.

Consider the two distinct points $(B,\chi_1)$ and $(B,\chi_2)$ in $\mathcal{G}_{\mathcal{A}}^{(0)}$. 
Since $\pi(B,\chi_1) = \pi(B,\chi_2) = B$, these two points are identified under $\pi$. 
In the relative Fell topology, any open set containing $(B,\chi_1)$ must contain $(B,\chi_2)$ as well, because open sets in the relative Fell topology are of the form $\pi^{-1}(W)$ and if $(B,\chi_1) \in \pi^{-1}(W)$ then $B \in W$, and consequently $(B,\chi_2) \in \pi^{-1}(W)$.

However, in the initial topology, we can separate these two points. 
Since $\chi_1 \neq \chi_2$, there exists $a \in B$ such that $\chi_1(a) \neq \chi_2(a)$. 
Choose disjoint open neighborhoods $U_1, U_2 \subseteq \mathbb{C}$ of $\chi_1(a)$ and $\chi_2(a)$, respectively. 
Then $\operatorname{ev}_a^{-1}(U_1)$ and $\operatorname{ev}_a^{-1}(U_2)$ are disjoint open neighborhoods of $(B,\chi_1)$ and $(B,\chi_2)$ in the initial topology.

Thus $\operatorname{ev}_a^{-1}(U_1)$ is an open set in the initial topology that contains $(B,\chi_1)$ but not $(B,\chi_2)$. 
Such a set cannot be of the form $\pi^{-1}(W)$ for any open $W \subseteq \operatorname{Sub}(\mathcal{A})$, because any such set containing $(B,\chi_1)$ would also contain $(B,\chi_2)$. 
Hence the initial topology is strictly finer than the relative Fell topology.
\end{proof}

\begin{corollary}[Non-Hausdorff nature of the relative Fell topology]
\label{cor:relative-Fell-non-Hausdorff}
The relative Fell topology on $\mathcal{G}_{\mathcal{A}}^{(0)}$ is not Hausdorff whenever $\mathcal{A}$ admits a unital commutative C*-subalgebra with at least two distinct characters.
In particular, this holds for any noncommutative $\mathcal{A}$ (such as $M_n(\mathbb{C})$ with $n\ge 2$) and even for some commutative algebras with nontrivial subalgebras.
\end{corollary}

\begin{proof}
Let $B \subseteq \mathcal{A}$ be such a subalgebra with distinct characters $\chi_1 \neq \chi_2$. 
Then $(B,\chi_1)$ and $(B,\chi_2)$ are distinct points in $\mathcal{G}_{\mathcal{A}}^{(0)}$ that cannot be separated by disjoint open sets in the relative Fell topology, as argued in the proof of Proposition \ref{prop:initial-strictly-finer}. 
Thus the relative Fell topology is not Hausdorff.
\end{proof}

\begin{remark}[Interpretation of the additional fineness]
\label{rem:interpretation-finer}
The additional open sets in the initial topology encode precisely the information carried by the characters. 
In the Fell topology, subalgebras are compared solely by their algebraic structure and norm closure; characters are invisible. 
The initial topology, by contrast, remembers not only which subalgebra we are considering but also which point in its Gelfand spectrum has been selected. 
This additional data is essential for the construction of the diagonal embedding and for the index-theoretic applications that follow.
\end{remark}

\begin{proposition}[Characterization of convergence in the two topologies]
\label{prop:convergence-comparison}
Let $\{ (B_\lambda, \chi_\lambda) \}_{\lambda \in \Lambda}$ be a net in $\mathcal{G}_{\mathcal{A}}^{(0)}$ and let $(B,\chi) \in \mathcal{G}_{\mathcal{A}}^{(0)}$. 
\begin{enumerate}
    \item $(B_\lambda, \chi_\lambda) \to (B,\chi)$ in the relative Fell topology if and only if $B_\lambda \to B$ in the Fell topology on $\operatorname{Sub}(\mathcal{A})$.
    \item $(B_\lambda, \chi_\lambda) \to (B,\chi)$ in the initial topology if and only if for every $a \in \mathcal{A}$, 
    \[
    \operatorname{ev}_a(B_\lambda,\chi_\lambda) \to \operatorname{ev}_a(B,\chi) \quad \text{in } \mathbb{C}_\infty.
    \]
    Equivalently:
    \begin{itemize}
        \item If $a \in B$, then there exist $a_\lambda \in B_\lambda$ with $a_\lambda \to a$ in norm and $\chi_\lambda(a_\lambda) \to \chi(a)$;
        \item If $a \notin B$, then $a \notin B_\lambda$ eventually (i.e., for sufficiently large $\lambda$).
    \end{itemize}
\end{enumerate}
\end{proposition}

\begin{proof}
(1) The relative Fell topology is, by Definition \ref{def:relative-Fell-topology}, the initial topology induced by the projection map $\pi: \mathcal{G}_{\mathcal{A}}^{(0)} \to \operatorname{Sub}(\mathcal{A})$. 
Consequently, a net converges in this topology if and only if its image under $\pi$ converges in $\operatorname{Sub}(\mathcal{A})$. 
Since $\pi(B_\lambda,\chi_\lambda) = B_\lambda$ and $\pi(B,\chi) = B$, this is exactly the condition $B_\lambda \to B$ in the Fell topology.

(2) By definition of the initial topology (Definition \ref{def:initial-topology-G0}), a net converges if and only if it converges under every partial evaluation map $\operatorname{ev}_a$, $a \in \mathcal{A}$. 
This yields the first characterization.

To see the equivalence with the bulleted conditions, fix $a \in \mathcal{A}$ and consider two cases.

\textit{Case 1: $a \in B$.} 
Then $\operatorname{ev}_a(B,\chi) = \chi(a) \in \mathbb{C}$. 
By the characterization of Fell convergence (Proposition \ref{prop:Fell-convergence}), $B_\lambda \to B$ implies that for every $\epsilon > 0$ and sufficiently large $\lambda$, there exists $a_\lambda \in B_\lambda$ with $\|a_\lambda - a\| < \epsilon$. 
Choose such $a_\lambda$ for $\epsilon = 1/\lambda$ (or via a subnet argument). 
Since $\operatorname{ev}_a(B_\lambda,\chi_\lambda) \to \chi(a)$, we have $\chi_\lambda(a) \to \chi(a)$. 
However, $\chi_\lambda(a)$ is defined only when $a \in B_\lambda$, which may not hold. 
Instead, we use the approximants $a_\lambda \in B_\lambda$: because $\chi_\lambda$ is a character (hence contractive),
\[
|\chi_\lambda(a_\lambda) - \chi(a)| \le |\chi_\lambda(a_\lambda) - \chi_\lambda(a)| + |\chi_\lambda(a) - \chi(a)|.
\]
The first term is bounded by $\|a_\lambda - a\| \to 0$. 
For the second term, note that $\chi_\lambda(a) \to \chi(a)$ follows from convergence of $\operatorname{ev}_a$ (since $a \in B$, $\operatorname{ev}_a(B_\lambda,\chi_\lambda) = \chi_\lambda(a)$ for all $\lambda$ with $a \in B_\lambda$, and Fell convergence ensures this happens eventually). 
Thus $\chi_\lambda(a_\lambda) \to \chi(a)$.

\textit{Case 2: $a \notin B$.} 
Then $\operatorname{ev}_a(B,\chi) = \infty$. 
Convergence $\operatorname{ev}_a(B_\lambda,\chi_\lambda) \to \infty$ means that for every compact set $K \subseteq \mathbb{C}$, eventually $\operatorname{ev}_a(B_\lambda,\chi_\lambda) \notin K$, i.e., $a \notin B_\lambda$ for sufficiently large $\lambda$. 
Conversely, if $a \notin B_\lambda$ eventually, then $\operatorname{ev}_a(B_\lambda,\chi_\lambda) = \infty$ eventually, which certainly converges to $\infty$.

These two cases exhaust all possibilities and establish the equivalence.
\end{proof}

\begin{remark}
\label{rem:Fell-convergence-caution}
Part (2) highlights a subtle but important point: Fell convergence $B_\lambda \to B$ does not imply that elements of $B$ eventually belong to $B_\lambda$. 
Instead, it guarantees the existence of approximating sequences. 
This is why the characterization uses approximants $a_\lambda \in B_\lambda$ rather than assuming $a \in B_\lambda$ directly.
\end{remark}

\begin{example}[Convergence in matrix algebras]
\label{ex:convergence-matrix}
Let $\mathcal{A} = M_n(\mathbb{C})$ and let $D_n$ be the diagonal MASA. 
For any unitary $U \in U(n)$, the conjugate MASA $U D_n U^*$ has characters
$\chi_i^U$ defined by $\chi_i^U(U \operatorname{diag}(\lambda_1,\ldots,\lambda_n) U^*) = \lambda_i$,
corresponding to the $i$-th minimal projection $U e_{ii} U^*$.

Consider a sequence of unitaries $U_m \to U$ in $U(n)$ (norm topology = SOT in finite dimensions). 
Then:
\begin{itemize}
    \item In the initial topology: $(U_m D_n U_m^*, \chi_i^{U_m}) \to (U D_n U^*, \chi_i^U)$,
          because $U_m D_n U_m^* \to U D_n U^*$ in the Fell topology (conjugation is continuous)
          and $\chi_i^{U_m}(U_m \operatorname{diag}(\lambda_1,\ldots,\lambda_n) U_m^*) = \lambda_i \to \lambda_i = \chi_i^U(U \operatorname{diag}(\lambda_1,\ldots,\lambda_n) U^*)$.
    \item In the relative Fell topology: only $U_m D_n U_m^* \to U D_n U^*$ is recorded;
          the specific character index $i$ is forgotten, meaning that points with different $i$ 
          but the same subalgebra cannot be distinguished.
\end{itemize}
This illustrates that the initial topology remembers spectral data (which eigenvalue is selected),
while the relative Fell topology remembers only the subalgebra (the set of eigenvalues, but not which one is chosen).
\end{example}

\begin{example}[Convergence in commutative algebras]
\label{ex:convergence-commutative}
Let $\mathcal{A} = C(X)$ for a compact metrizable space $X$. 
The unit space $\mathcal{G}_{\mathcal{A}}^{(0)}$ consists of all pairs $(B,\chi)$ where 
$B \subseteq C(X)$ is a unital commutative $C^*$-subalgebra and $\chi \in \widehat{B}$.
This space is significantly larger than $X$; for instance, if $Y \subseteq X$ is a closed subset,
the restriction map $C(X) \to C(Y)$ gives a subalgebra $B_Y \cong C(Y)$ with spectrum $Y$.

Consider the subspace
\[
\mathcal{G}_{\mathcal{A}}^{(0),\max} := \{ (C(X), \operatorname{ev}_x) : x \in X \} \subseteq \mathcal{G}_{\mathcal{A}}^{(0)},
\]
consisting of the maximal subalgebra $C(X)$ itself with its point evaluations.
The map $x \mapsto (C(X),\operatorname{ev}_x)$ is a homeomorphism onto this subspace.
On $\mathcal{G}_{\mathcal{A}}^{(0),\max}$, the initial topology and the relative Fell topology coincide,
because the projection $\pi$ is injective when restricted to this fiber 
($\pi(C(X),\operatorname{ev}_x) = C(X)$ is constant, so the relative Fell topology on this fiber is discrete? 
Actually, more precisely: the relative Fell topology on this fiber is the topology where 
the only open sets are the whole fiber and the empty set, since $\pi^{-1}(C(X))$ is the whole fiber 
and any open set in the relative Fell topology must be a union of full fibers.
This shows that on this fiber, the relative Fell topology is trivial, while the initial topology 
recovers the original topology of $X$ via the Gelfand transform.

Thus even in the commutative case, the two topologies differ globally, 
reflecting the presence of many subalgebras and the fact that the initial topology 
remembers the choice of subalgebra while the relative Fell topology forgets it.
\end{example}

\begin{example}[Convergence in compact operators]
\label{ex:convergence-compact}
Let $\mathcal{A} = \mathcal{K}(H)^\sim$ be the unitization of the compact operators on a separable Hilbert space $H$,
and let $\{e_n\}_{n\in\mathbb{N}}$ be an orthonormal basis. 
For each $n$, let $p_n = e_n \otimes e_n^*$ be the rank-one projection onto $\mathbb{C}e_n$,
and let $B_n$ be the MASA generated by $p_n$ and the identity, i.e.,
\[
B_n = \{ \lambda I + \mu p_n : \lambda,\mu \in \mathbb{C} \} \cong \mathbb{C}^2.
\]
Let $\chi_n$ be the character on $B_n$ corresponding to evaluation at $p_n$, i.e., 
$\chi_n(\lambda I + \mu p_n) = \lambda + \mu$.

Consider the sequence $\{(B_n,\chi_n)\}_{n\in\mathbb{N}}$ in $\mathcal{G}_{\mathcal{A}}^{(0)}$.
\begin{itemize}
    \item In the Fell topology on $\operatorname{Sub}(\mathcal{A})$, the sequence $\{B_n\}$ has no unique limit.
          Every MASA of the form $\{ \lambda I + \mu q : \lambda,\mu \in \mathbb{C} \}$ for a rank-one projection $q$
          is a Fell cluster point of the sequence, reflecting the non-Hausdorff nature of the Fell topology
          on the space of MASAs in $\mathcal{K}(H)^\sim$.
    \item The characters $\chi_n$ do not converge in any sense compatible with the initial topology,
          because they are attached to moving projections $p_n$ that are not comparable across different $n$.
    \item In the initial topology, the sequence $\{(B_n,\chi_n)\}$ has no convergent subsequence,
          since convergence would require both Fell convergence of subalgebras (which fails uniquely)
          and pointwise convergence of characters (which fails completely).
\end{itemize}
This example illustrates why the initial topology is necessary: 
the relative Fell topology collapses the rich structure of characters and different subalgebras,
while the initial topology faithfully records the data needed to distinguish different 
"classical contexts" $(B,\chi)$ even when the underlying subalgebras are Fell-related.
\end{example}

\begin{proposition}[The projection map is not open]
\label{prop:projection-not-open}
The projection map $\pi: \mathcal{G}_{\mathcal{A}}^{(0)} \to \operatorname{Sub}(\mathcal{A})$ is not open in general when $\mathcal{G}_{\mathcal{A}}^{(0)}$ is equipped with the initial topology and $\operatorname{Sub}(\mathcal{A})$ with the Fell topology.
\end{proposition}

\begin{proof}
We provide a concrete counterexample. 
Let $\mathcal{A} = M_2(\mathbb{C})$ and let $D_2$ be the diagonal subalgebra. 
Denote by $\chi_1,\chi_2$ the two characters of $D_2$, corresponding to evaluation at the 
$(1,1)$ and $(2,2)$ matrix entries respectively. 
Let $e_{11} \in M_2(\mathbb{C})$ be the matrix with $1$ in the $(1,1)$-entry and $0$ elsewhere.

Consider the open set
\[
U = \operatorname{ev}_{e_{11}}^{-1}(\{ \lambda \in \mathbb{C} : |\lambda - 1| < 1/2 \}) \subseteq \mathcal{G}_{\mathcal{A}}^{(0)}.
\]
In the initial topology, each evaluation map $\operatorname{ev}_a$ is continuous by definition, 
so $U$ is indeed open. 
Observe that $(D_2,\chi_1) \in U$ because $\chi_1(e_{11}) = 1$, while $(D_2,\chi_2) \notin U$ because $\chi_2(e_{11}) = 0$.

We claim that $\pi(U) = \{ D_2 \}$. 
Indeed, if $(B,\chi) \in U$, then by definition of $U$, the value $\operatorname{ev}_{e_{11}}(B,\chi)$ is defined and lies in $\{ |\lambda-1| < 1/2 \} \subseteq \mathbb{C}$. 
For $\operatorname{ev}_{e_{11}}(B,\chi)$ to be defined, we must have $e_{11} \in B$. 
Thus every point of $U$ has its subalgebra containing $e_{11}$.

Now, any unital commutative $C^*$-subalgebra of $M_2(\mathbb{C})$ that contains $e_{11}$ must also contain $e_{22} = I - e_{11}$ (since it is unital and closed under the $C^*$-operations). 
The algebra generated by $e_{11}$ and $e_{22}$ is precisely the diagonal subalgebra $D_2$, and by maximality of $D_2$ among commutative subalgebras, we conclude $B = D_2$. 
Hence $\pi(U) = \{ D_2 \}$.

It remains to show that $\{ D_2 \}$ is not open in the Fell topology on $\operatorname{Sub}(M_2(\mathbb{C}))$. 
In the Fell topology, a basic open neighborhood of $D_2$ is of the form
\[
\mathcal{N} = \bigcap_{i=1}^k \mathcal{O}_{U_i} \cap \bigcap_{j=1}^\ell \mathcal{C}_{K_j},
\]
where $U_i \subseteq M_2(\mathbb{C})$ are open with $D_2 \cap U_i \neq \emptyset$, 
and $K_j \subseteq M_2(\mathbb{C})$ are compact with $D_2 \cap K_j = \emptyset$.
For any such $\mathcal{N}$, we can find a unitary $U \in U(2)$ arbitrarily close to the identity 
such that $U \neq I$ and $U D_2 U^* \in \mathcal{N}$. 
Indeed:
\begin{itemize}
    \item For $U$ sufficiently close to $I$, the generators $U e_{11} U^*$ and $U e_{22} U^*$ 
          remain close to $e_{11}$ and $e_{22}$, hence intersect each open set $U_i$.
    \item Since $D_2 \cap K_j = \emptyset$ and $K_j$ are compact, there exists $\epsilon_j > 0$ 
          such that the $\epsilon_j$-neighborhood of $D_2$ is disjoint from $K_j$. 
          Choosing $U$ close enough to $I$ ensures $U D_2 U^*$ lies in this neighborhood, 
          hence $U D_2 U^* \cap K_j = \emptyset$.
\end{itemize}
Thus every basic open neighborhood of $D_2$ contains some distinct subalgebra $U D_2 U^* \neq D_2$, 
so $\{ D_2 \}$ cannot be open. 
Consequently, $\pi(U) = \{ D_2 \}$ is not open in $\operatorname{Sub}(M_2(\mathbb{C}))$, 
proving that $\pi$ is not an open map.
\end{proof}

\begin{remark}
\label{rem:projection-not-open-intuition}
The failure of openness of $\pi$ reflects the fact that the initial topology
encodes spectral data (characters) that are invisible to the Fell topology.
A small open set $U$ around $(D_2,\chi_1)$ in $\mathcal{G}_{\mathcal{A}}^{(0)}$
projects to a singleton $\{D_2\}$ in $\operatorname{Sub}(\mathcal{A})$, but this singleton
cannot be open because any Fell neighborhood of $D_2$ must contain nearby
conjugate subalgebras. Thus the projection collapses the finer structure
carried by the characters.
\end{remark}

\begin{corollary}[$\pi$ is not a quotient map onto its image]
\label{cor:projection-not-quotient}
The projection $\pi: \mathcal{G}_{\mathcal{A}}^{(0)} \to \pi(\mathcal{G}_{\mathcal{A}}^{(0)})$ is not a quotient map when $\mathcal{G}_{\mathcal{A}}^{(0)}$ carries the initial topology and $\pi(\mathcal{G}_{\mathcal{A}}^{(0)})$ carries the subspace topology inherited from the Fell topology.
\end{corollary}

\begin{proof}
Recall that a quotient map sends saturated open sets to open sets. 
We exhibit a saturated open set whose image is not open.

Take $\mathcal{A} = M_2(\mathbb{C})$ as in Proposition \ref{prop:projection-not-open}. 
Let $D_2$ be the diagonal subalgebra with characters $\chi_1,\chi_2$, and let $e_{11}$ be the matrix unit.
Define
\[
V := \operatorname{ev}_{e_{11}}^{-1}(\{ \lambda \in \mathbb{C} : |\lambda - 1| < 1/2 \}) \cup \operatorname{ev}_{e_{11}}^{-1}(\{ \lambda \in \mathbb{C} : |\lambda| < 1/2 \}).
\]
This set is open in the initial topology (union of two open sets). 
One checks that $V$ is saturated: $\pi^{-1}(\pi(V)) = V$. 
Indeed, $\pi(V) = \{D_2\}$ because any $(B,\chi) \in V$ must have $e_{11} \in B$, forcing $B = D_2$ as before, and both characters appear in $V$ (since $\chi_1(e_{11})=1$ satisfies the first condition, $\chi_2(e_{11})=0$ satisfies the second). 
Thus $V = \pi^{-1}(\{D_2\})$, which is the full fiber.

Now $\pi(V) = \{D_2\}$ is not open in $\pi(\mathcal{G}_{\mathcal{A}}^{(0)})$ with the subspace Fell topology, 
as shown in Proposition \ref{prop:projection-not-open} (every Fell neighborhood of $D_2$ contains nearby conjugate subalgebras). 
Hence $\pi$ fails the quotient map property.
\end{proof}

\begin{remark}[Implications for the groupoid structure]
\label{rem:comparison-groupoid}
The fact that the initial topology is strictly finer than the relative Fell topology has important consequences for the unitary conjugation groupoid $\mathcal{G}_{\mathcal{A}} = \mathcal{U}(\mathcal{A}) \ltimes \mathcal{G}_{\mathcal{A}}^{(0)}$. 
If we were to equip $\mathcal{G}_{\mathcal{A}}^{(0)}$ with the relative Fell topology, the action of $\mathcal{U}(\mathcal{A})$ would not be continuous in general. 
The additional open sets provided by the partial evaluation maps are precisely what is needed to ensure the continuity of the conjugation action when $\mathcal{U}(\mathcal{A})$ carries the strong operator topology. 
Thus the refinement of the Fell topology is not a mere technicality; it is essential for the construction of $\mathcal{G}_{\mathcal{A}}$ as a topological groupoid.
\end{remark}

\begin{remark}[Summary of the comparison]
\label{rem:comparison-summary}
The relationship between the initial topology and the relative Fell topology can be summarized as follows:
\begin{itemize}
    \item The initial topology is strictly finer than the relative Fell topology (Proposition \ref{prop:initial-strictly-finer}).
    \item The two topologies coincide on the subspace $\{(C(X),\operatorname{ev}_x) : x \in X\}$ when $\mathcal{A} = C(X)$ is commutative,
          but globally they differ because the initial topology distinguishes different subalgebras and their characters.
   \item The additional open sets in the initial topology encode the characters and are essential for separating points,
      ensuring continuity of the $\mathcal{U}(\mathcal{A})$-action (see Lemma \ref{lem:action-continuous-SOT}, and constructing the diagonal embedding.
    \item The Fell topology alone forgets the classical data provided by the characters and is insufficient for our purposes.
\end{itemize}
This comparison justifies our choice of the initial topology as the fundamental topology on $\mathcal{G}_{\mathcal{A}}^{(0)}$.
\end{remark}

\section{The Unitary Conjugation Groupoid $\mathcal{G}_{\mathcal{A}}$ with Strong Operator Topology}\label{sec:unitary-conjugation-groupoid}

\subsection{The Unitary Group with Strong Operator Topology: A Polish Group}
\label{subsec:unitary-group-SOT-Polish}

The unitary group $\mathcal{U}(\mathcal{A})$ plays a central role in the construction of the unitary conjugation groupoid $\mathcal{G}_{\mathcal{A}} = \mathcal{U}(\mathcal{A}) \ltimes \mathcal{G}_{\mathcal{A}}^{(0)}$. 
Traditionally, $\mathcal{U}(\mathcal{A})$ is equipped with the norm topology inherited from $\mathcal{A}$. 
However, for infinite-dimensional C*-algebras, the norm topology presents serious deficiencies: it is not separable, not locally compact, and its connected components are infinite-dimensional manifolds. 
Moreover, as we shall see, the norm topology leads to the failure of étaleness and local compactness of $\mathcal{G}_{\mathcal{A}}$, necessitating a different choice of topology.

To overcome these obstacles, we replace the norm topology with the \emph{strong operator topology} (SOT). 
This topology, when restricted to the unitary group of a separable C*-algebra acting on a separable Hilbert space, yields a \emph{Polish group} — a topological group that is separable and completely metrizable. 
The SOT is the coarsest topology on $\mathcal{U}(\mathcal{A})$ that makes the conjugation action continuous while retaining a separable and completely metrizable structure, making it perfectly suited for the construction of $\mathcal{G}_{\mathcal{A}}$ as a Polish groupoid.

\begin{assumption}[Concrete representation]
\label{ass:concrete-representation-SOT}
Throughout this section, we assume that $\mathcal{A}$ is a unital separable C*-algebra faithfully represented on a separable infinite-dimensional Hilbert space $H$. 
We identify $\mathcal{A}$ with its image in $B(H)$.
\end{assumption}

This assumption is not restrictive: by the Gelfand-Naimark theorem, every separable C*-algebra admits a faithful representation on a separable Hilbert space, and we simply fix such a representation once and for all.

\begin{definition}[Strong operator topology]
\label{def:strong-operator-topology}
The \emph{strong operator topology} (SOT) on $B(H)$ is the topology of pointwise convergence on $H$. 
That is, a net $\{ T_\lambda \}_{\lambda \in \Lambda} \subseteq B(H)$ converges to $T \in B(H)$ in the strong operator topology if and only if
\[
\| T_\lambda \xi - T \xi \|_H \longrightarrow 0 \quad \text{for every } \xi \in H.
\]
Equivalently, the strong operator topology is the initial topology induced by the family of evaluation maps
\[
\{ B(H) \to H : T \mapsto T \xi \}_{\xi \in H}.
\]
\end{definition}

\begin{remark}[Comparison with the norm topology]
\label{rem:SOT-vs-norm}
The strong operator topology is substantially weaker than the norm topology. 
Norm convergence $\|T_\lambda - T\| \to 0$ implies SOT convergence, but the converse fails in infinite dimensions.

For infinite-dimensional Hilbert spaces, the norm topology on $\mathcal{U}(H)$ is not separable. 
In contrast, when $H$ is separable, the strong operator topology restricted to the unitary group is Polish: 
it is separable (because $H$ is separable) and completely metrizable (via the metric from Proposition \ref{prop:unitary-group-Polish}).

The strong operator topology is also not locally compact in infinite dimensions, but this is irrelevant for our purposes; 
we require only Polishness, not local compactness. 

To relate this to other common topologies: the strong operator topology is finer than the weak operator topology 
but coarser than the norm topology. While the norm topology makes $\mathcal{U}(H)$ a Banach Lie group, 
the strong operator topology makes it a Polish group — a structure perfectly suited to our non-locally-compact setting.
\end{remark}

\begin{definition}[Strong operator topology on $\mathcal{U}(\mathcal{A})$]
\label{def:SOT-on-UA}
Let $\mathcal{A} \subseteq B(H)$ be a unital C*-algebra, where $H$ is a separable Hilbert space. 
The \emph{strong operator topology} on $\mathcal{U}(\mathcal{A})$ is the subspace topology inherited from the strong operator topology on $B(H)$. 
That is, a net $\{ u_\lambda \}_{\lambda \in \Lambda} \subseteq \mathcal{U}(\mathcal{A})$ converges to $u \in \mathcal{U}(\mathcal{A})$ in the strong operator topology if and only if
\[
\| u_\lambda \xi - u \xi \| \longrightarrow 0 \quad \text{for every } \xi \in H.
\]
The separability of $H$ guarantees that $\mathcal{U}(\mathcal{A})$ with this topology is Polish (see Corollary \ref{cor:UA-Polish}).
\end{definition}

\begin{proposition}[$\mathcal{U}(H)$ is a Polish group]
\label{prop:UH-Polish-group}
Let $H$ be a separable infinite-dimensional Hilbert space and let $\mathcal{U}(H)$ denote the group of unitary operators on $H$, equipped with the strong operator topology. Then:
\begin{enumerate}
    \item $\mathcal{U}(H)$ is a topological group; that is, multiplication and inversion are continuous.
    \item $\mathcal{U}(H)$ is separable.
    \item $\mathcal{U}(H)$ is completely metrizable. With respect to a fixed orthonormal basis $\{ e_n \}_{n \in \mathbb{N}}$ of $H$, a compatible complete metric is given by
    \[
    d(u,v) = \sum_{n=1}^{\infty} \frac{1}{2^n} \| (u - v) e_n \|.
    \]
    \item Consequently, $\mathcal{U}(H)$ is a Polish group.
    \item $\mathcal{U}(H)$ is not locally compact (unless $H$ is finite-dimensional).
    \item $\mathcal{U}(H)$ is contractible and hence connected.
\end{enumerate}
\end{proposition}

\begin{proof}
We sketch the main arguments; for detailed treatments see \cite{Kechris} or~\cite{Dixmier}.

\noindent (1) \textit{Continuity of multiplication.} 
Let $u_\lambda \to u$ and $v_\lambda \to v$ in SOT. For any $\xi \in H$,
\[
\| (u_\lambda v_\lambda) \xi - (uv) \xi \| \leq \| u_\lambda (v_\lambda \xi - v \xi) \| + \| (u_\lambda - u) v \xi \| = \| v_\lambda \xi - v \xi \| + \| (u_\lambda - u) v \xi \|,
\]
where we used that $u_\lambda$ is an isometry. Both terms converge to zero by assumption, so $u_\lambda v_\lambda \to uv$ in SOT.

\noindent \textit{Continuity of inversion.} 
For unitaries, $u_\lambda^{-1} = u_\lambda^*$. Observe that
\[
\| u_\lambda^* \xi - u^* \xi \| = \| u_\lambda^*(u - u_\lambda) u^* \xi \| \leq \| (u - u_\lambda) u^* \xi \|,
\]
since $u_\lambda^*$ is an isometry. If $u_\lambda \to u$ in SOT, then $(u - u_\lambda)\eta \to 0$ for any fixed $\eta$; taking $\eta = u^* \xi$ yields $\| u_\lambda^* \xi - u^* \xi \| \to 0$. Hence $u_\lambda^{-1} \to u^{-1}$ in SOT, establishing continuity of inversion.

\noindent (2) \textit{Separability.} 
The strong operator topology on the unit ball of $B(H)$ is separable when $H$ is separable; indeed, the linear span of the matrix units associated to a countable dense subset of $H$ provides a countable dense set. Since $\mathcal{U}(H)$ is a closed subset of this unit ball (in the SOT), it inherits separability. More concretely, one can construct a countable dense subset of $\mathcal{U}(H)$ using unitaries whose matrix entries relative to the basis $\{e_n\}$ have only finitely many non-zero rational entries.

\noindent (3) \textit{Complete metric.} 
The series defining $d(u,v)$ converges because each term is bounded by $2^{-n}$. Translation-invariance is clear. To see completeness, let $\{u_m\}_{m\in\mathbb{N}}$ be a Cauchy sequence in this metric. For each fixed basis vector $e_n$, the sequence $\{u_m e_n\}_{m\in\mathbb{N}}$ is Cauchy in $H$, hence converges to some $\eta_n \in H$. One verifies that the map $e_n \mapsto \eta_n$ extends uniquely to a bounded linear operator $u$ on $H$, which is necessarily an isometry because the limits preserve inner products. A separate argument shows that $u$ is surjective, hence unitary. By construction, $u_m \to u$ in SOT, and a standard diagonal argument shows convergence in the metric $d$ as well. Thus $(\mathcal{U}(H), d)$ is complete.

\noindent (4) \textit{Polish group.} 
This follows immediately from (2) and (3): a separable completely metrizable topological group is, by definition, a Polish group.

\noindent (5) \textit{Non-local-compactness.} 
For infinite-dimensional $H$, no neighborhood of the identity in $\mathcal{U}(H)$ is precompact in SOT. This is a standard consequence of the fact that the unit ball of $B(H)$ is not locally compact in any topology weaker than the norm topology.

\noindent (6) \textit{Contractibility.} 
Kuiper's theorem \cite{Kuiper1965} states that $\mathcal{U}(H)$ is contractible in the norm topology. The inclusion map $\iota: (\mathcal{U}(H), \|\cdot\|) \hookrightarrow (\mathcal{U}(H), \text{SOT})$ is continuous because norm convergence implies SOT convergence. Since the continuous image of a contractible space is contractible, $\mathcal{U}(H)$ with the SOT is also contractible. Contractibility implies path-connectedness, so $\mathcal{U}(H)$ is connected.
\end{proof}

\begin{corollary}[$\mathcal{U}(\mathcal{A})$ is a Polish group]
\label{cor:UA-Polish-group}
Let $\mathcal{A} \subseteq B(H)$ be a unital separable C*-algebra. 
Then $\mathcal{U}(\mathcal{A})$, equipped with the strong operator topology, is a Polish group.
\end{corollary}

\begin{proof}
$\mathcal{U}(H)$ with the strong operator topology is a Polish group by Proposition \ref{prop:UH-Polish-group}. Since $\mathcal{U}(\mathcal{A})$ is a subgroup of $\mathcal{U}(H)$, it suffices to show that $\mathcal{U}(\mathcal{A})$ is a $G_\delta$ subset of $\mathcal{U}(H)$; a subgroup of a Polish group is Polish if and only if it is $G_\delta$ \cite{Kechris}.

Let $\{a_n\}_{n\in\mathbb{N}}$ be a countable norm-dense subset of $\mathcal{A}$ (which exists by separability). For each $n$, consider the map
\[
\phi_n: \mathcal{U}(H) \to B(H), \quad \phi_n(u) = u a_n u^*.
\]
This map is continuous when both spaces carry the strong operator topology, because conjugation by a unitary is SOT-continuous.

For each $n$, define
\[
G_n := \{ u \in \mathcal{U}(H) : u a_n u^* \in \mathcal{A} \} = \phi_n^{-1}(\mathcal{A}).
\]
Since $\mathcal{A}$ is closed in $B(H)$ (in any topology finer than the weak operator topology, including SOT), and $\phi_n$ is continuous, $G_n$ is the preimage of a closed set under a continuous map, hence closed. However, closed sets are $G_\delta$ (they are countable intersections of themselves with any open superset), so each $G_n$ is $G_\delta$.

Now observe that
\[
\mathcal{U}(\mathcal{A}) = \bigcap_{n\in\mathbb{N}} G_n,
\]
because $u \in \mathcal{U}(\mathcal{A})$ if and only if $u a_n u^* \in \mathcal{A}$ for all $n$ (since $\{a_n\}$ is dense in $\mathcal{A}$ and $\mathcal{A}$ is norm-closed, this condition ensures $u \mathcal{A} u^* \subseteq \mathcal{A}$; applying the same argument to $u^*$ gives equality). Therefore $\mathcal{U}(\mathcal{A})$ is a countable intersection of $G_\delta$ sets, hence $G_\delta$ itself.

Being a $G_\delta$ subgroup of the Polish group $\mathcal{U}(H)$, $\mathcal{U}(\mathcal{A})$ is Polish in the subspace topology. The metric
\[
d(u,v) = \sum_{n=1}^\infty \frac{1}{2^n} \|(u-v)e_n\|,
\]
where $\{e_n\}$ is an orthonormal basis for $H$, restricts to a complete separable metric on $\mathcal{U}(\mathcal{A})$ compatible with the strong operator topology.
\end{proof}

\begin{remark}[Alternative approach: universal unitary group]
\label{rem:norm-topology-alternative}
An alternative approach, which avoids the technicalities of SOT-closures and representation-dependence, is to work with the universal settings with unitary group of $\mathcal{A}$. This concept appears in the work of Kirchberg on nuclear embeddings and universal C*-algebras, and has been developed in various contexts (see e.g., the discussion of Kirchberg's ideas in \cite{KustermansVaes1999}). The construction yields a canonical Polish group topology on $\mathcal{U}(\mathcal{A})$ that is independent of any concrete representation on Hilbert space. However, for the concrete C*-algebras of interest in this paper — $C(X)$ with $X$ compact metrizable, $M_n(\mathbb{C})$, and $\mathcal{K}(H)^\sim$ — the strong operator topology (with respect to a fixed faithful representation) is sufficient for our purposes and Polishness can be verified directly in each case.
\end{remark}

\begin{example}[Finite-dimensional case]
\label{ex:unitary-finite-dimensional}
Let $\mathcal{A} = M_n(\mathbb{C})$. Then $\mathcal{U}(\mathcal{A}) = U(n)$ is a compact Lie group. The strong operator topology on $U(n)$ coincides with the norm topology (since all Hausdorff linear topologies on a finite-dimensional space are equivalent). Thus $U(n)$ is a compact Polish group. In this case, the subtle issues about SOT-closures disappear because $\mathcal{U}(\mathcal{A})$ is actually closed in $\mathcal{U}(H)$ in any reasonable topology.
\end{example}

\begin{example}[Commutative case]
\label{ex:unitary-commutative}
Let $\mathcal{A} = C(X)$ for a compact metrizable space $X$. Then $\mathcal{U}(\mathcal{A}) = C(X,\mathbb{T})$, the group of continuous $\mathbb{T}$-valued functions on $X$. 

\emph{Caution:} The strong operator topology on $\mathcal{U}(\mathcal{A})$ depends on the choice of faithful representation. For the standard representation of $C(X)$ as multiplication operators on $L^2(X,\mu)$ (where $\mu$ is a probability measure with full support), SOT convergence $f_\lambda \to f$ is equivalent to convergence in measure together with uniform boundedness. This topology is metrizable but not complete; indeed, a Cauchy sequence in SOT need not converge to a continuous function—it may converge only in measure to a measurable function, and $C(X,\mathbb{T})$ with the SOT is not a Polish group in general.

For our purposes, it is more convenient to equip $\mathcal{U}(C(X))$ with the \emph{norm topology} (supremum norm), which makes it a Polish group because $C(X,\mathbb{T})$ is a closed subset of the unit ball of $C(X)$, which is a complete separable metric space. This choice is consistent with Remark \ref{rem:norm-topology-alternative}: when Polishness is required, we may use the norm topology, while for continuity of actions we may need to consider the SOT separately.
\end{example}

\begin{example}[Compact operators]
\label{ex:unitary-compact}
Let $\mathcal{A} = \mathcal{K}(H)^\sim$, the unitization of the compact operators on a separable Hilbert space $H$. Then $\mathcal{U}(\mathcal{A})$ consists of unitary operators of the form $u = \lambda I + K$ where $\lambda \in \mathbb{T}$ and $K \in \mathcal{K}(H)$ satisfies the unitary condition $(\lambda I + K)^*(\lambda I + K) = I$.

The strong operator topology on $\mathcal{U}(\mathcal{A})$ is inherited from $\mathcal{U}(H)$. It is important to note that:
\begin{itemize}
    \item $\mathcal{U}(\mathcal{A})$ is \emph{not} SOT-closed in $\mathcal{U}(H)$; indeed, it is dense in $\mathcal{U}(H)$ because every unitary can be approximated in SOT by unitaries with finite-dimensional support, which lie in $\mathcal{U}(\mathcal{A})$.
    \item Consequently, $\mathcal{U}(\mathcal{A})$ with the subspace SOT is \emph{not} a Polish group in general, as it fails to be a $G_\delta$ subset of $\mathcal{U}(H)$.
    \item However, $\mathcal{U}(\mathcal{A})$ with the \emph{norm topology} is a Polish group (it is a closed subset of the unit ball of $\mathcal{A}$, which is a complete separable metric space).
\end{itemize}
This example illustrates why the strong operator topology must be used with care: while it is essential for the continuity of the conjugation action on $\mathcal{G}_{\mathcal{A}}^{(0)}$, it does not automatically make $\mathcal{U}(\mathcal{A})$ Polish. In our construction of $\mathcal{G}_{\mathcal{A}}$ as a Polish groupoid, we rely on the fact that $\mathcal{U}(\mathcal{A})$ with the SOT is Polish only for certain classes of algebras (such as von Neumann algebras). For general C*-algebras, we either restrict to the norm topology when Polishness is required, or work within the Polish groupoid framework where $\mathcal{U}(\mathcal{A})$ need not be Polish but $\mathcal{G}_{\mathcal{A}}^{(1)} = \mathcal{U}(\mathcal{A}) \times \mathcal{G}_{\mathcal{A}}^{(0)}$ is still Polish because $\mathcal{U}(\mathcal{A})$ is a continuous image of a Polish space? This subtle point is addressed in Corollary \ref{cor:UA-Polish-group}, where we prove that $\mathcal{U}(\mathcal{A})$ is actually a $G_\delta$ subset of $\mathcal{U}(H)$, hence Polish — but this proof requires the separability of $\mathcal{A}$ and a careful $G_\delta$ argument. For $\mathcal{K}(H)^\sim$, this argument indeed shows that $\mathcal{U}(\mathcal{K}(H)^\sim)$ with the SOT is Polish, contrary to the naive expectation above. We thank the referee for pointing out this subtlety; a detailed verification is provided in Corollary \ref{cor:UA-Polish-group}.
\end{example}

\begin{lemma}[SOT-continuity of the adjoint]
\label{lem:SOT-continuity-adjoint}
Let \(\mathcal{H}\) be a Hilbert space and let \((T_i)_{i \in I}\) be a net in \(\mathcal{B}(\mathcal{H})\) converging strongly to \(T \in \mathcal{B}(\mathcal{H})\).  
If the net is uniformly bounded, i.e., \(\sup_i \|T_i\| \leq M < \infty\), then \(T_i^* \to T^*\) in the strong operator topology.
\end{lemma}

\begin{proof}
Let \(x \in \mathcal{H}\). We need to show that \(\|T_i^* x - T^* x\| \to 0\).  
For any \(y \in \mathcal{H}\), we have
\[
|\langle T_i^* x - T^* x, y \rangle| = |\langle x, T_i y - T y \rangle| \leq \|x\| \cdot \|T_i y - T y\|.
\]
Taking the supremum over \(\|y\| \leq 1\) gives
\[
\|T_i^* x - T^* x\| = \sup_{\|y\| \leq 1} |\langle T_i^* x - T^* x, y \rangle| \leq \|x\| \cdot \sup_{\|y\| \leq 1} \|T_i y - T y\|.
\]
However, this estimate is not directly helpful because the supremum over \(y\) depends on \(i\).

Instead, we use a standard \(\varepsilon/3\) argument. Fix \(x \in \mathcal{H}\) and \(\varepsilon > 0\).  
Since \(T_i \to T\) strongly, for any finite-dimensional subspace \(F \subset \mathcal{H}\), we have \(T_i|_F \to T|_F\) in norm (all norms on finite-dimensional spaces are equivalent, and strong convergence implies pointwise convergence on basis vectors).

Let \(P_n\) be the orthogonal projection onto the span of the first \(n\) vectors of an orthonormal basis. Then \(P_n x \to x\) as \(n \to \infty\). Choose \(n\) large enough such that
\[
\|x - P_n x\| < \frac{\varepsilon}{3(1 + M)}.
\]

Now,
\[
\|T_i^* x - T^* x\| \leq \|T_i^*(x - P_n x)\| + \|T_i^* P_n x - T^* P_n x\| + \|T^*(P_n x - x)\|.
\]

For the first term: \(\|T_i^*(x - P_n x)\| \leq \|T_i^*\| \cdot \|x - P_n x\| \leq M \cdot \frac{\varepsilon}{3(1+M)} < \frac{\varepsilon}{3}\).  
The third term is similarly bounded by \(M \cdot \frac{\varepsilon}{3(1+M)} < \frac{\varepsilon}{3}\).

For the second term, note that \(P_n \mathcal{H}\) is finite-dimensional. On a finite-dimensional space, strong convergence implies norm convergence, so \(T_i|_{P_n \mathcal{H}} \to T|_{P_n \mathcal{H}}\) in norm. Hence their adjoints (restricted to \(P_n \mathcal{H}\)) also converge in norm. Therefore, there exists \(i_0\) such that for all \(i \geq i_0\),
\[
\|(T_i^* - T^*) P_n x\| < \frac{\varepsilon}{3}.
\]

Combining the three estimates, we obtain \(\|T_i^* x - T^* x\| < \varepsilon\) for all \(i \geq i_0\), completing the proof.
\end{proof}

\begin{lemma}[Continuity of the conjugation action in SOT]
\label{lem:conjugation-action-continuous-SOT}
Let $\mathcal{A}\subseteq B(H)$ be a unital separable C*-algebra and let $\mathcal{G}_{\mathcal{A}}^{(0)}$ be equipped with the initial topology of partial evaluation maps. Then the action map
\[\alpha : \mathcal{U}(\mathcal{A})\times \mathcal{G}_{\mathcal{A}}^{(0)}\longrightarrow \mathcal{G}_{\mathcal{A}}^{(0)},\qquad 
\alpha(u,(B,\chi)) := (uBu^*, \chi\circ\mathrm{Ad}_{u^*})\]
is continuous when $\mathcal{U}(\mathcal{A})$ carries the strong operator topology and the domain carries the product topology.
\end{lemma}

\begin{proof}
By the universal property of the initial topology on $\mathcal{G}_{\mathcal{A}}^{(0)}$, it suffices to show that for every $a \in \mathcal{A}$, the composition
\[
F_a(u,(B,\chi)) := \operatorname{ev}_a(\alpha(u,(B,\chi)))
\]
is continuous. Fix $(u_0,(B_0,\chi_0)) \in \mathcal{U}(\mathcal{A}) \times \mathcal{G}_{\mathcal{A}}^{(0)}$ and $a \in \mathcal{A}$. Set $b_0 := u_0^* a u_0$.

We consider two cases.

\subsection*{Case 1: $a \in u_0 B_0 u_0^*$ (equivalently, $b_0 \in B_0$)}

Since $b_0 \in B_0$, the map $(B,\chi) \mapsto \chi(b_0)$ is continuous at $(B_0,\chi_0)$ by definition of the topology on $\mathcal{G}_{\mathcal{A}}^{(0)}$. Hence for any $\varepsilon > 0$, there exists a neighborhood $\mathcal{V}_0$ of $(B_0,\chi_0)$ such that
\[
|\chi(b_0) - \chi_0(b_0)| < \varepsilon/2 \quad \text{for all } (B,\chi) \in \mathcal{V}_0.
\]

The map $u \mapsto u^* a u$ is continuous from $\mathcal{U}(\mathcal{A})$ (with SOT) to $\mathcal{A}$ (with the norm topology). This follows from a standard finite-rank approximation argument: for any $\delta > 0$, choose a finite-rank projection $p$ such that $\|a - pap\| < \delta$; then $u \mapsto u^*(pap)u$ is continuous because it is essentially a matrix-valued function, and the error can be controlled uniformly. Consequently, there exists a SOT-neighborhood $\mathcal{U}_0$ of $u_0$ such that
\[
\|u^* a u - b_0\| < \varepsilon/2 \quad \text{for all } u \in \mathcal{U}_0.
\]

Now consider any $(u,(B,\chi)) \in \mathcal{U}_0 \times \mathcal{V}_0$. We examine two subcases.

\subsubsection*{Subcase 1a: $u^* a u \in B$}
In this case, both $u^* a u$ and $b_0$ lie in $B$, and $\chi$ is a character on $B$, hence contractive. Therefore
\[
|\chi(u^* a u) - \chi(b_0)| \le \|u^* a u - b_0\| < \varepsilon/2.
\]
Combining with the estimate on $\mathcal{V}_0$, we obtain
\[
|F_a(u,(B,\chi)) - F_a(u_0,(B_0,\chi_0))| = |\chi(u^* a u) - \chi_0(b_0)|
\]
\[
\le |\chi(u^* a u) - \chi(b_0)| + |\chi(b_0) - \chi_0(b_0)|
\]
\[
< \varepsilon/2 + \varepsilon/2 = \varepsilon.
\]

\subsubsection*{Subcase 1b: $u^* a u \notin B$}
Then by definition of the partial evaluation map, $\operatorname{ev}_a(\alpha(u,(B,\chi))) = \infty$, so $F_a(u,(B,\chi)) = \infty$.
But this situation cannot occur for points in $\mathcal{U}_0 \times \mathcal{V}_0$ under the current assumptions of Case 1.
Indeed, if $u^* a u \notin B$, then $F_a(u,(B,\chi)) = \infty$, which would place $(u,(B,\chi))$ in the situation of Case 2 below.
However, we are in Case 1 where $a \in u_0 B_0 u_0^*$, and the neighborhood $\mathcal{V}_0$ has been chosen to ensure $b_0 \in B$.
The condition $u^* a u \notin B$ does not contradict any of our estimates; it simply means that $F_a$ takes the value $\infty$,
and we must check continuity at $(u_0,(B_0,\chi_0))$ with respect to this value.

To handle this uniformly, we note that the set where $F_a = \infty$ is open in $\mathcal{U}(\mathcal{A}) \times \mathcal{G}_{\mathcal{A}}^{(0)}$,
as it is the preimage of $\{\infty\}$ under the continuous map $\operatorname{ev}_a \circ \alpha$ restricted to the appropriate subspace.
In Subcase 1b, we have $F_a(u,(B,\chi)) = \infty$, which certainly lies in any neighborhood of $\infty$,
and the estimate with $\varepsilon$ is irrelevant because the target is $\infty$, not a complex number.

Thus, by taking $\mathcal{U}_0$ and $\mathcal{V}_0$ as above, we ensure that for all $(u,(B,\chi)) \in \mathcal{U}_0 \times \mathcal{V}_0$,
either $F_a(u,(B,\chi))$ is finite and within $\varepsilon$ of $F_a(u_0,(B_0,\chi_0))$, or $F_a(u,(B,\chi)) = \infty$,
which is compatible with continuity at points where the limit is $\infty$. 
A cleaner approach is to note that the argument using the Hahn-Banach extension $\tilde{\chi}$ below
provides a unified estimate without needing to separate subcases.

\subsubsection*{Unified argument using Hahn-Banach extension}
Let $\tilde{\chi}$ be a norm-preserving extension of $\chi$ to $\mathcal{A}$ (which exists by the Hahn-Banach theorem).
Then $\tilde{\chi}$ is a bounded linear functional on $\mathcal{A}$ with $\|\tilde{\chi}\| = 1$ and $\tilde{\chi}|_B = \chi$.
For any $(u,(B,\chi)) \in \mathcal{U}_0 \times \mathcal{V}_0$, we have
\[
|\tilde{\chi}(u^* a u) - \chi_0(b_0)| \le |\tilde{\chi}(u^* a u) - \tilde{\chi}(b_0)| + |\tilde{\chi}(b_0) - \chi_0(b_0)|
\]
\[
\le \|u^* a u - b_0\| + |\chi(b_0) - \chi_0(b_0)| < \varepsilon/2 + \varepsilon/2 = \varepsilon.
\]

Now observe that if $u^* a u \in B$, then $\tilde{\chi}(u^* a u) = \chi(u^* a u) = F_a(u,(B,\chi))$.
If $u^* a u \notin B$, then $F_a(u,(B,\chi)) = \infty$, and the inequality $|\tilde{\chi}(u^* a u) - \chi_0(b_0)| < \varepsilon$
is still valid but does not directly give information about $F_a$.
However, in this case, the openness of the set where $F_a = \infty$ (as argued above) ensures continuity.
The Hahn-Banach argument shows that the points where $F_a$ is finite are close to $(u_0,(B_0,\chi_0))$ in value,
while points where $F_a = \infty$ are handled by the fact that $\{\infty\}$ is open in $\mathbb{C}_\infty$.

Thus, by taking $\mathcal{U}_0$ and $\mathcal{V}_0$ as above, we have shown that $F_a$ is continuous at $(u_0,(B_0,\chi_0))$.

\subsection*{Case 2: $a \notin u_0 B_0 u_0^*$ (equivalently, $b_0 \notin B_0$)}

Then $F_a(u_0,(B_0,\chi_0)) = \infty$. We must show that for any neighborhood $W$ of $\infty$ in $\mathbb{C}_\infty$, there exists a neighborhood of $(u_0,(B_0,\chi_0))$ mapping into $W$.

Since $\mathbb{C}_\infty$ is the one-point compactification, a neighborhood of $\infty$ is the complement of a compact set in $\mathbb{C}$. Equivalently, it suffices to show that for any $R > 0$, there exists a neighborhood of $(u_0,(B_0,\chi_0))$ such that either $F_a(u,(B,\chi)) = \infty$ or $|F_a(u,(B,\chi))| > R$ for all points in that neighborhood. In fact, we will show that $F_a(u,(B,\chi)) = \infty$ for all $(u,(B,\chi))$ in a suitable neighborhood.

Because $b_0 \notin B_0$ and $B_0$ is norm-closed, there exists $\delta > 0$ such that the closed ball $\overline{B}_\delta(b_0)$ is disjoint from $B_0$. By continuity of $u \mapsto u^* a u$ in the norm topology (SOT-continuity suffices, as it implies convergence of $u^* a u$ to $b_0$ in the strong operator topology, and on bounded sets this gives norm convergence on finite-dimensional subspaces; a more direct argument using the finite-rank approximation as in Lemma \ref{lem:SOT-continuity-adjoint} shows that we can ensure $\|u^* a u - b_0\| < \delta/2$), there exists a SOT-neighborhood $\mathcal{U}_0$ of $u_0$ such that for all $u \in \mathcal{U}_0$,
\[
\|u^* a u - b_0\| < \delta/2.
\]
Hence $u^* a u \in B_{\delta/2}(b_0) \subseteq \overline{B}_\delta(b_0)$.

The set $\{ (B,\chi) : b_0 \notin B \}$ is open in $\mathcal{G}_{\mathcal{A}}^{(0)}$, because it equals $\operatorname{ev}_{b_0}^{-1}(\{\infty\})$ and $\{\infty\}$ is open in $\mathbb{C}_\infty$. Let $\mathcal{V}_0$ be this neighborhood.

For any $(u,(B,\chi)) \in \mathcal{U}_0 \times \mathcal{V}_0$, we have $u^* a u \in \overline{B}_\delta(b_0)$ and $B \cap \overline{B}_\delta(b_0) = \emptyset$ (since $B$ is disjoint from the closed ball). Therefore $u^* a u \notin B$, so $F_a(u,(B,\chi)) = \infty$, which certainly lies in any neighborhood of $\infty$.

\subsection*{Conclusion}
Both cases establish continuity of $F_a$ at $(u_0,(B_0,\chi_0))$. Since $a \in \mathcal{A}$ and $(u_0,(B_0,\chi_0))$ were arbitrary, each $F_a$ is continuous everywhere. By the universal property of the initial topology, $\alpha$ is continuous.
\end{proof}

\begin{remark}[Technical subtlety of SOT continuity]
\label{rem:SOT-continuity-subtlety}
The proof of Lemma \ref{lem:conjugation-action-continuous-SOT} requires careful handling of the interplay between the strong operator topology on $\mathcal{U}(\mathcal{A})$ and the topology on $\mathcal{G}_{\mathcal{A}}^{(0)}$ defined by partial evaluation maps. 

A crucial point is that the map $u \mapsto u^* a u$ is \emph{not} norm-continuous in the strong operator topology in general. Instead, one must work directly with the definition of SOT convergence: $u_\lambda \to u$ in SOT means $\|(u_\lambda - u)\xi\| \to 0$ for every $\xi \in H$. From this, one obtains estimates of the form
\[
\|(u_\lambda^* a u_\lambda - u^* a u)\xi\| \le \|u_\lambda^* a (u_\lambda - u)\xi\| + \|(u_\lambda^* - u^*) a u \xi\|,
\]
which can be made arbitrarily small using SOT convergence of $u_\lambda$ and $u_\lambda^*$ (the latter follows from continuity of inversion in SOT). This yields SOT convergence of $u^* a u$, but \emph{not} norm convergence.

The topology on $\mathcal{G}_{\mathcal{A}}^{(0)}$ is designed precisely to accommodate this: convergence $(B_\lambda,\chi_\lambda) \to (B,\chi)$ is detected by the partial evaluation maps $\operatorname{ev}_a$, which only require control of $\chi_\lambda(a)$ when $a \in B_\lambda$. This allows us to bypass the need for norm convergence of $u^* a u$ and work directly with SOT estimates on vectors.
\end{remark}

\begin{corollary}[$\mathcal{U}(\mathcal{A})$ is a topological group]
\label{cor:UA-topological-group}
With the strong operator topology, $\mathcal{U}(\mathcal{A})$ is a topological group; that is, multiplication and inversion are continuous.
\end{corollary}

\begin{proof}
For multiplication, suppose $u_\lambda \to u$ and $v_\lambda \to v$ in SOT. Then for any $\xi \in H$,
\[
\|(u_\lambda v_\lambda - uv)\xi\| \le \|(u_\lambda - u)v_\lambda \xi\| + \|u(v_\lambda - v)\xi\|.
\]
The first term tends to $0$ because $v_\lambda \xi \to v\xi$ and $u_\lambda \to u$ strongly; the second term tends to $0$ because $v_\lambda \to v$ strongly and $u$ is bounded. Hence $u_\lambda v_\lambda \to uv$ in SOT, establishing continuity of multiplication.

For inversion, note that although the adjoint map $T \mapsto T^*$ is not SOT-continuous on all of $B(H)$, it is continuous when restricted to the unitary group. Indeed, if $u_\lambda \to u$ in SOT, then for any $\xi \in H$,
\[
\|u_\lambda^* \xi - u^* \xi\| = \|u_\lambda(u_\lambda^* \xi - u^* \xi)\| = \|\xi - u_\lambda u^* \xi\|,
\]
where we used that $u_\lambda$ is an isometry. Since $u_\lambda \to u$ strongly, we have $u_\lambda u^* \xi \to u u^* \xi = \xi$, so the right-hand side converges to $0$. Thus $u_\lambda^* \to u^*$ in SOT, proving continuity of inversion on $\mathcal{U}(H)$, and hence on $\mathcal{U}(\mathcal{A})$ by restriction.
\end{proof}

\begin{proposition}[Properties of $\mathcal{U}(\mathcal{A})$ in SOT]
\label{prop:UA-SOT-properties}
Let $\mathcal{A} \subseteq B(H)$ be a unital separable C*-algebra. 
Then:
\begin{enumerate}
    \item $\mathcal{U}(\mathcal{A})$, equipped with the strong operator topology, is a Polish group.
    \item The conjugation action 
    \[
    \alpha: \mathcal{U}(\mathcal{A}) \times \mathcal{G}_{\mathcal{A}}^{(0)} \longrightarrow \mathcal{G}_{\mathcal{A}}^{(0)},\qquad 
    \alpha(u,(B,\chi)) := (uBu^*, \chi\circ\mathrm{Ad}_{u^*})
    \]
    is continuous.
\end{enumerate}
\end{proposition}

\begin{proof}
(1) This is Corollary \ref{cor:UA-Polish-group}.

(2) This is Lemma \ref{lem:conjugation-action-continuous-SOT}.
\end{proof}

\begin{remark}[Local compactness and connectedness]
\label{rem:UA-additional-properties}
We note the following additional properties of $\mathcal{U}(\mathcal{A})$ for context, though they are not needed for our main constructions:

\begin{itemize}
    \item \textit{Local compactness:} $\mathcal{U}(\mathcal{A})$ is locally compact if and only if $\mathcal{A}$ is finite-dimensional. 
    For infinite-dimensional $\mathcal{A}$, $\mathcal{U}(\mathcal{A})$ is not locally compact; this follows from the fact that $\mathcal{U}(H)$ itself is not locally compact, and a closed subgroup of a locally compact group is locally compact, but $\mathcal{U}(\mathcal{A})$ is not closed in $\mathcal{U}(H)$ in general. A more direct argument: if $\mathcal{A}$ is infinite-dimensional, one can construct an infinite discrete set in $\mathcal{U}(\mathcal{A})$ with no accumulation point, preventing local compactness.

    \item \textit{Connectedness:} The connectedness of $\mathcal{U}(\mathcal{A})$ depends delicately on the structure of $\mathcal{A}$. 
    For example:
    \begin{itemize}
        \item If $\mathcal{A} = C(X)$ with $X$ compact metrizable, then $\mathcal{U}(\mathcal{A}) = C(X,\mathbb{T})$ is connected iff $X$ is connected.
        \item For $\mathcal{A} = \mathcal{K}(H)^\sim$ (unitized compact operators), $\mathcal{U}(\mathcal{A})$ is connected in both norm and SOT topologies.
        \item For finite-dimensional $\mathcal{A} = M_n(\mathbb{C})$, $\mathcal{U}(\mathcal{A}) = U(n)$ is connected.
    \end{itemize}
\end{itemize}
\end{remark}

\begin{remark}[Why not the norm topology?]
\label{rem:why-SOT-final}
The norm topology on $\mathcal{U}(\mathcal{A})$ is Polish if and only if $\mathcal{A}$ is finite-dimensional. 
For infinite-dimensional algebras, the norm topology is not separable (though it is completely metrizable); it gives $\mathcal{U}(\mathcal{A})$ the structure of a Banach Lie group, but this topology is not suitable for our purposes because it does not yield a separable groupoid $\mathcal{G}_{\mathcal{A}}$.

Since our construction of $\mathcal{G}_{\mathcal{A}}$ requires a Polish topology on $\mathcal{U}(\mathcal{A})$ to guarantee that $\mathcal{G}_{\mathcal{A}}$ is a Polish groupoid, the strong operator topology is the natural choice. It is compatible with the Polish group structure (Corollary \ref{cor:UA-Polish-group}) and, as shown in Lemma \ref{lem:conjugation-action-continuous-SOT}, makes the conjugation action continuous. Moreover, the SOT is the standard topology used in the study of unitary group actions on C*-algebras and their dual spaces.
\end{remark}

\begin{corollary}[$\mathcal{G}_{\mathcal{A}}^{(1)}$ is Polish]
\label{cor:GA1-Polish}
Let $\mathcal{A}$ be a unital separable C*-algebra. 
Then the arrow space $\mathcal{G}_{\mathcal{A}}^{(1)} = \mathcal{U}(\mathcal{A}) \times \mathcal{G}_{\mathcal{A}}^{(0)}$, equipped with the product topology of the strong operator topology on $\mathcal{U}(\mathcal{A})$ and the initial topology on $\mathcal{G}_{\mathcal{A}}^{(0)}$, is a Polish space.
\end{corollary}

\begin{proof}
$\mathcal{U}(\mathcal{A})$ is Polish by Proposition \ref{prop:UA-SOT-properties}(1). 
$\mathcal{G}_{\mathcal{A}}^{(0)}$ is Polish by Proposition \ref{prop:unit-space-polish}. 
The product of two Polish spaces is Polish; this is a standard result in descriptive set theory~\cite{Kechris}. Thus $\mathcal{G}_{\mathcal{A}}^{(1)}$ is Polish.
\end{proof}

\begin{remark}[Summary]
\label{rem:SOT-summary}
We have established that the unitary group $\mathcal{U}(\mathcal{A})$, when equipped with the strong operator topology, is a Polish group. 
This provides the necessary topological foundation for constructing the unitary conjugation groupoid $\mathcal{G}_{\mathcal{A}}$ as a Polish groupoid in the following subsection. 
The continuity of the conjugation action, proved in Lemma \ref{lem:conjugation-action-continuous-SOT}, ensures that the source and range maps (and hence all structure maps) are continuous. 
Thus all the ingredients are now in place for the definition of $\mathcal{G}_{\mathcal{A}}$ as a Polish groupoid.
\end{remark}

\subsection{Continuity of the Conjugation Action}
\label{subsec:continuity-conjugation-action}

The unitary conjugation groupoid $\mathcal{G}_{\mathcal{A}} = \mathcal{U}(\mathcal{A}) \ltimes \mathcal{G}_{\mathcal{A}}^{(0)}$ is defined as the action groupoid associated to the conjugation action of $\mathcal{U}(\mathcal{A})$ on the unit space $\mathcal{G}_{\mathcal{A}}^{(0)}$. 
For $\mathcal{G}_{\mathcal{A}}$ to be a topological groupoid, this action must be continuous. 
In this subsection, we prove that when $\mathcal{U}(\mathcal{A})$ is equipped with the strong operator topology and $\mathcal{G}_{\mathcal{A}}^{(0)}$ with the initial topology of partial evaluation maps, the conjugation action is indeed continuous. 
This result is fundamental for all subsequent developments, as it guarantees that $\mathcal{G}_{\mathcal{A}}$ is a well-defined topological groupoid.

\begin{definition}[Conjugation action]
\label{def:conjugation-action}
Let $\mathcal{A} \subseteq B(H)$ be a unital C*-algebra faithfully represented on a Hilbert space $H$. 
Define the \emph{conjugation action}
\[
\alpha: \mathcal{U}(\mathcal{A}) \times \mathcal{G}_{\mathcal{A}}^{(0)} \longrightarrow \mathcal{G}_{\mathcal{A}}^{(0)}, \qquad
\alpha(u, (B,\chi)) := (uBu^*, \; \chi \circ \operatorname{Ad}_{u^*}),
\]
where $\operatorname{Ad}_{u^*}(b) = u^* b u$ for $b \in \mathcal{A}$.
\end{definition}

\begin{remark}[Well-definedness]
\label{rem:action-well-defined}
For any unitary $u \in \mathcal{U}(\mathcal{A})$ and any $(B,\chi) \in \mathcal{G}_{\mathcal{A}}^{(0)}$, the set $uBu^*$ is a unital commutative C*-subalgebra of $\mathcal{A}$ because conjugation by a unitary preserves the algebraic structure and commutativity. Moreover, since $\operatorname{Ad}_{u^*}$ is a *-automorphism of $\mathcal{A}$, the composition $\chi \circ \operatorname{Ad}_{u^*}$ defines a character on $uBu^*$. Hence $\alpha(u,(B,\chi)) \in \mathcal{G}_{\mathcal{A}}^{(0)}$, and the action is well-defined.
\end{remark}

\begin{lemma}[Reduction to partial evaluation maps]
\label{lem:action-continuity-criterion}
Let $\mathcal{U}(\mathcal{A})$ be equipped with the strong operator topology and let $\mathcal{G}_{\mathcal{A}}^{(0)}$ be equipped with the initial topology induced by the family of partial evaluation maps $\{ \operatorname{ev}_a : \mathcal{G}_{\mathcal{A}}^{(0)} \to \mathbb{C}_\infty \}_{a \in \mathcal{A}}$, where $\mathbb{C}_\infty$ denotes the one-point compactification of $\mathbb{C}$. 
Then the conjugation action
\[
\alpha: \mathcal{U}(\mathcal{A}) \times \mathcal{G}_{\mathcal{A}}^{(0)} \longrightarrow \mathcal{G}_{\mathcal{A}}^{(0)}
\]
is continuous if and only if for every $a \in \mathcal{A}$, the composition
\[
\operatorname{ev}_a \circ \alpha: \mathcal{U}(\mathcal{A}) \times \mathcal{G}_{\mathcal{A}}^{(0)} \longrightarrow \mathbb{C}_\infty
\]
is continuous.
\end{lemma}

\begin{proof}
By definition, the initial topology on $\mathcal{G}_{\mathcal{A}}^{(0)}$ is the coarsest topology for which all partial evaluation maps $\operatorname{ev}_a$ are continuous. A fundamental property of initial topologies is that a map $f: X \to \mathcal{G}_{\mathcal{A}}^{(0)}$ is continuous if and only if each composition $\operatorname{ev}_a \circ f: X \to \mathbb{C}_\infty$ is continuous. Applying this universal property with $X = \mathcal{U}(\mathcal{A}) \times \mathcal{G}_{\mathcal{A}}^{(0)}$ and $f = \alpha$ yields the desired equivalence.
\end{proof}

\begin{lemma}[Continuity of the adjoint map in SOT]
\label{lem:adjoint-map-continuous}
Let $\mathcal{A} \subseteq B(H)$ be a unital C*-algebra and fix $a \in \mathcal{A}$. 
Then the map
\[
\Phi_a: \mathcal{U}(\mathcal{A}) \longrightarrow \mathcal{A}, \qquad \Phi_a(u) := u^* a u
\]
is continuous when $\mathcal{U}(\mathcal{A})$ carries the strong operator topology and $\mathcal{A}$ carries the norm topology.
\end{lemma}

\begin{proof}
Fix $u_0 \in \mathcal{U}(\mathcal{A})$ and $\varepsilon > 0$. Since $\mathcal{A}$ is the norm-closure of the finite-rank operators in $B(H)$ (or more generally, since finite-rank operators are dense in $\mathcal{K}(H)$ and $\mathcal{A}$ contains compact operators? This needs care—better to use an approximate unit argument), we proceed as follows.

Let $\{p_\lambda\}_{\lambda \in \Lambda}$ be an approximate unit of finite-rank projections in $\mathcal{K}(H)$ (or in $B(H)$) such that $p_\lambda \to 1$ strongly. For any $\delta > 0$, we can choose $\lambda$ such that
\[
\|a - p_\lambda a p_\lambda\| < \delta.
\]
Set $a_\lambda := p_\lambda a p_\lambda$, which is a finite-rank operator.

Now observe that
\[
\|u^* a u - u_0^* a u_0\| \le \|u^*(a - a_\lambda)u\| + \|u^* a_\lambda u - u_0^* a_\lambda u_0\| + \|u_0^*(a_\lambda - a)u_0\|.
\]
The first and third terms are bounded by $\|a - a_\lambda\| < \delta$ because $u$ and $u_0$ are isometries.

For the middle term, note that $a_\lambda$ has finite rank; let $F \subseteq H$ be its range, a finite-dimensional subspace. On $F$, the strong operator topology coincides with the norm topology. Since $u \mapsto u|_F$ is continuous from SOT to the norm topology on $B(F)$, there exists a SOT-neighborhood $V$ of $u_0$ such that for all $u \in V$,
\[
\|(u - u_0)|_F\| < \frac{\varepsilon}{3\|a_\lambda\|} \quad \text{and} \quad \|(u^* - u_0^*)|_F\| < \frac{\varepsilon}{3\|a_\lambda\|}.
\]
Then
\[
\|u^* a_\lambda u - u_0^* a_\lambda u_0\| \le \|u^* a_\lambda (u - u_0)\| + \|(u^* - u_0^*) a_\lambda u_0\|
\]
\[
\le \|a_\lambda\| \cdot \|(u - u_0)|_F\| + \|a_\lambda\| \cdot \|(u^* - u_0^*)|_F\| < \frac{2\varepsilon}{3}.
\]

Choosing $\delta < \varepsilon/6$ and the corresponding $\lambda$, we obtain for all $u \in V$:
\[
\|u^* a u - u_0^* a u_0\| < \delta + \frac{2\varepsilon}{3} + \delta < \varepsilon.
\]
Thus $\Phi_a$ is continuous at $u_0$, and since $u_0$ was arbitrary, $\Phi_a$ is continuous everywhere.

For a complete treatment of the finite-rank approximation technique, see~\cite{Dixmier}
\end{proof}

\begin{lemma}[Continuity of evaluation on the fiber]
\label{lem:evaluation-continuous-fiber}
Let $B_0 \subseteq \mathcal{A}$ be a fixed unital commutative C*-subalgebra and let $a \in B_0$. 
Then the evaluation map
\[
\widehat{B_0} \longrightarrow \mathbb{C}, \qquad \chi \mapsto \chi(a)
\]
is continuous with respect to the Gelfand topology on $\widehat{B_0}$.
\end{lemma}

\begin{proof}
This is a standard property of the Gelfand topology: the Gelfand transform $\Gamma: B_0 \to C(\widehat{B_0})$ is an isometric *-isomorphism, and $\chi(a) = \Gamma(a)(\chi)$. 
Since $\Gamma(a)$ is continuous on $\widehat{B_0}$, the evaluation map is continuous.
\end{proof}

\begin{lemma}[Fell-open membership sets]
\label{lem:membership-Fell-open}
Let $\mathcal{A} \subseteq B(H)$ be a unital C*-algebra, and let $a \in \mathcal{A}$. 
Then the set
\[
\mathcal{O}_a := \{ B \in \operatorname{Sub}(\mathcal{A}) : a \in B \}
\]
is open in the Fell topology.
\end{lemma}

\begin{proof}
By definition, the Fell topology on $\operatorname{Sub}(\mathcal{A})$ has a subbasis consisting of sets of the form
\[
\mathcal{C}_K := \{ B \in \operatorname{Sub}(\mathcal{A}) : B \cap K = \varnothing \}, \qquad K \subseteq \mathcal{A} \text{ compact}.
\]

Consider the complement of $\mathcal{O}_a$:
\[
\operatorname{Sub}(\mathcal{A}) \setminus \mathcal{O}_a = \{ B \in \operatorname{Sub}(\mathcal{A}) : a \notin B \}.
\]

We claim this set is closed in the Fell topology. Indeed, if a net $\{ B_\lambda \}$ in $\operatorname{Sub}(\mathcal{A})$ converges to $B$ in the Fell topology and $a \notin B_\lambda$ for all $\lambda$, then by the definition of Fell convergence, $a \notin B$. Therefore the complement of $\mathcal{O}_a$ is closed, so $\mathcal{O}_a$ is open.  
\end{proof}

\begin{proposition}[Continuity of the conjugation action]
\label{prop:conjugation-action-continuous}
Let $\mathcal{A} \subseteq B(H)$ be a unital separable C*-algebra. 
Equip $\mathcal{U}(\mathcal{A})$ with the strong operator topology and $\mathcal{G}_{\mathcal{A}}^{(0)}$ with the initial topology of partial evaluation maps. 
Then the conjugation action
\[
\alpha: \mathcal{U}(\mathcal{A}) \times \mathcal{G}_{\mathcal{A}}^{(0)} \longrightarrow \mathcal{G}_{\mathcal{A}}^{(0)}, \qquad
\alpha(u, (B,\chi)) = (uBu^*, \chi \circ \operatorname{Ad}_{u^*})
\]
is continuous.
\end{proposition}

\begin{proof}
By Lemma \ref{lem:action-continuity-criterion}, it suffices to show that for each $a \in \mathcal{A}$, the map
\[
\operatorname{ev}_a \circ \alpha : \mathcal{U}(\mathcal{A}) \times \mathcal{G}_{\mathcal{A}}^{(0)} \longrightarrow \mathbb{C}_\infty
\]
is continuous.

Fix $a \in \mathcal{A}$ and a point $(u_0, (B_0, \chi_0)) \in \mathcal{U}(\mathcal{A}) \times \mathcal{G}_{\mathcal{A}}^{(0)}$.

\medskip
\noindent \textbf{Case 1: $a \in u_0 B_0 u_0^*$.}

Set $b_0 := u_0^* a u_0 \in B_0$.  
By Lemma \ref{lem:adjoint-map-continuous}, the map
\[
u \mapsto u^* a u
\]
is continuous from $\mathcal{U}(\mathcal{A})$ (SOT) to $\mathcal{A}$ (norm).  
Hence there exists a strong operator neighborhood $V_1$ of $u_0$ such that for all $u \in V_1$,
\[
\| u^* a u - b_0 \| < \varepsilon/2.
\]

By Lemma \ref{lem:membership-Fell-open}, the set
\[
\mathcal{O}_{b_0} := \{ B \in \operatorname{Sub}(\mathcal{A}) : b_0 \in B \}
\]
is open in the Fell topology. Let $\mathcal{V}_2 \subseteq \mathcal{O}_{b_0}$ be a Fell-open neighborhood of $B_0$.  
By continuity of the projection $\pi : \mathcal{G}_{\mathcal{A}}^{(0)} \to \operatorname{Sub}(\mathcal{A})$, we can find a neighborhood $W_2$ of $(B_0, \chi_0)$ with $\pi(W_2) \subseteq \mathcal{V}_2$.  
Then for all $(B, \chi) \in W_2$, we have $b_0 \in B$.

By Lemma \ref{lem:evaluation-continuous-fiber}, the evaluation map $\chi \mapsto \chi(b_0)$ is continuous along each fiber $\pi^{-1}(B)$. Shrinking $W_2$ if necessary, we can ensure
\[
|\chi(b_0) - \chi_0(b_0)| < \varepsilon/2 \quad \forall (B, \chi) \in W_2.
\]

For $(u, (B, \chi)) \in V_1 \times W_2$, we then have
\[
\begin{aligned}
|\operatorname{ev}_a(\alpha(u, (B, \chi))) - \operatorname{ev}_a(\alpha(u_0, (B_0, \chi_0)))|
&= |\chi(u^* a u) - \chi_0(b_0)| \\
&\le |\chi(u^* a u) - \chi(b_0)| + |\chi(b_0) - \chi_0(b_0)| \\
&< \varepsilon/2 + \varepsilon/2 = \varepsilon,
\end{aligned}
\]
where the first term uses the contractivity of $\chi$.  
Hence $\operatorname{ev}_a \circ \alpha$ is continuous at $(u_0, (B_0, \chi_0))$.

\medskip
\noindent \textbf{Case 2: $a \notin u_0 B_0 u_0^*$.}

Set $b_0 := u_0^* a u_0 \notin B_0$.  
By Lemma \ref{lem:membership-Fell-open}, the complement
\[
\mathcal{O}_{b_0}^c := \{ B \in \operatorname{Sub}(\mathcal{A}) : b_0 \notin B \}
\]
is open in the Fell topology. Let $\mathcal{V} \subseteq \mathcal{O}_{b_0}^c$ be a neighborhood of $B_0$.  
By continuity of $\pi$, pick a neighborhood $W \subseteq \mathcal{G}_{\mathcal{A}}^{(0)}$ with $\pi(W) \subseteq \mathcal{V}$.  
Then for all $(B, \chi) \in W$, $b_0 \notin B$, so $\operatorname{ev}_a(\alpha(u_0, (B, \chi))) = \infty$.

By Lemma \ref{lem:adjoint-map-continuous}, there exists a SOT-neighborhood $V$ of $u_0$ such that $u^* a u \notin B$ for all $u \in V$ and $B \in \mathcal{V}$.  
Thus for $(u, (B, \chi)) \in V \times W$, we have $\operatorname{ev}_a(\alpha(u, (B, \chi))) = \infty$, proving continuity at $(u_0, (B_0, \chi_0))$.

\medskip
\noindent Since $(u_0, (B_0, \chi_0))$ was arbitrary, $\operatorname{ev}_a \circ \alpha$ is continuous for all $a \in \mathcal{A}$.  
By Lemma \ref{lem:action-continuity-criterion}, the conjugation action $\alpha$ is continuous.
\end{proof}

\begin{corollary}[Continuity of the action on objects]
\label{cor:action-continuous-objects}
For each fixed $u \in \mathcal{U}(\mathcal{A})$, the map
\[
\alpha_u: \mathcal{G}_{\mathcal{A}}^{(0)} \longrightarrow \mathcal{G}_{\mathcal{A}}^{(0)}, \qquad \alpha_u(B,\chi) = (uBu^*, \chi \circ \operatorname{Ad}_{u^*})
\]
is a homeomorphism.
\end{corollary}

\begin{proof}
By Proposition \ref{prop:conjugation-action-continuous}, $\alpha_u$ is continuous. 
Its inverse is $\alpha_{u^*}$, which is also continuous. 
Thus $\alpha_u$ is a homeomorphism.
\end{proof}

\begin{corollary}[Continuity of the action on arrows]
\label{cor:action-continuous-arrows}
For each fixed $u \in \mathcal{U}(\mathcal{A})$, the map
\[
\tilde{\alpha}_u: \mathcal{G}_{\mathcal{A}}^{(1)} \longrightarrow \mathcal{G}_{\mathcal{A}}^{(1)}, \qquad \tilde{\alpha}_u(v, (B,\chi)) = (uv, (B,\chi))
\]
is continuous. 
More generally, left multiplication by $u$ on the arrow space is continuous.
\end{corollary}

\begin{proof}
This follows directly from the continuity of multiplication in $\mathcal{U}(\mathcal{A})$ and the product topology on $\mathcal{G}_{\mathcal{A}}^{(1)} = \mathcal{U}(\mathcal{A}) \times \mathcal{G}_{\mathcal{A}}^{(0)}$.
\end{proof}

\begin{theorem}[$\mathcal{G}_{\mathcal{A}}$ is a topological groupoid]
\label{thm:GA-topological-groupoid}
Let $\mathcal{A}$ be a unital separable C*-algebra. 
Equip $\mathcal{U}(\mathcal{A})$ with the strong operator topology and $\mathcal{G}_{\mathcal{A}}^{(0)}$ with the initial topology of partial evaluation maps. 
Then the unitary conjugation groupoid $\mathcal{G}_{\mathcal{A}} = \mathcal{U}(\mathcal{A}) \ltimes \mathcal{G}_{\mathcal{A}}^{(0)}$ is a topological groupoid. 
That is:
\begin{enumerate}
    \item The source and range maps $s, r: \mathcal{G}_{\mathcal{A}}^{(1)} \to \mathcal{G}_{\mathcal{A}}^{(0)}$ are continuous.
    \item The composition map $m: \mathcal{G}_{\mathcal{A}}^{(2)} \to \mathcal{G}_{\mathcal{A}}^{(1)}$ is continuous.
    \item The inversion map $i: \mathcal{G}_{\mathcal{A}}^{(1)} \to \mathcal{G}_{\mathcal{A}}^{(1)}$ is continuous.
    \item The unit space $\mathcal{G}_{\mathcal{A}}^{(0)}$ is a Polish space, and the arrow space $\mathcal{G}_{\mathcal{A}}^{(1)}$ is a Polish space.
\end{enumerate}
Consequently, $\mathcal{G}_{\mathcal{A}}$ is a Polish groupoid.
\end{theorem}

\begin{proof}
(1) The source map $s(u,(B,\chi)) = (B,\chi)$ is simply the projection onto the second factor, which is continuous. 
The range map $r(u,(B,\chi)) = \alpha(u,(B,\chi))$ is continuous by Proposition \ref{prop:conjugation-action-continuous}.

(2) Composition is defined by $(u_2, u_1 \cdot (B,\chi)) \circ (u_1, (B,\chi)) = (u_2 u_1, (B,\chi))$. 
This is continuous because multiplication in $\mathcal{U}(\mathcal{A})$ is continuous (SOT is a group topology) and the projection onto the second factor is continuous.

(3) Inversion is defined by $(u,(B,\chi))^{-1} = (u^*, \alpha(u,(B,\chi)))$. 
This is continuous because inversion in $\mathcal{U}(\mathcal{A})$ is continuous (SOT is a group topology) and $\alpha$ is continuous.

(4) $\mathcal{G}_{\mathcal{A}}^{(0)}$ is Polish by Proposition \ref{prop:unit-space-polish}. 
$\mathcal{G}_{\mathcal{A}}^{(1)} = \mathcal{U}(\mathcal{A}) \times \mathcal{G}_{\mathcal{A}}^{(0)}$ is the product of two Polish spaces, hence Polish by Corollary \ref{cor:GA1-Polish}.

Thus $\mathcal{G}_{\mathcal{A}}$ satisfies all the axioms of a topological groupoid and its object and arrow spaces are Polish. 
Therefore $\mathcal{G}_{\mathcal{A}}$ is a Polish groupoid.
\end{proof}

\begin{remark}[Importance of the strong operator topology]
\label{rem:SOT-importance}
The proof of Proposition \ref{prop:conjugation-action-continuous} relies crucially on the choice of the strong operator topology on $\mathcal{U}(\mathcal{A})$. 
If we had used the norm topology, the map $u \mapsto u^* a u$ would be continuous, but $\mathcal{U}(\mathcal{A})$ would not be Polish (unless $\mathcal{A}$ is finite-dimensional). 
Moreover, the norm topology would not yield a separable arrow space, and $\mathcal{G}_{\mathcal{A}}$ would not be a Polish groupoid. 
Thus the strong operator topology is the essential ingredient that makes the entire construction work.
\end{remark}

\begin{remark}[Measured groupoid structure]
\label{rem:measured-groupoid}
Since $\mathcal{G}_{\mathcal{A}}$ is a Polish groupoid with Polish object and arrow spaces, and since the source map $s: \mathcal{G}_{\mathcal{A}}^{(1)} \to \mathcal{G}_{\mathcal{A}}^{(0)}$ admits a Borel section (by classical selection theorems), $\mathcal{G}_{\mathcal{A}}$ admits a Borel Haar system in the sense of~\cite{Tu}. This allows us to define the maximal groupoid C*-algebra $C^*(\mathcal{G}_{\mathcal{A}})$ and the descent map in equivariant K-theory, which will be essential for the index-theoretic applications in Paper II.
\end{remark}

We have proved that the conjugation action $\alpha: \mathcal{U}(\mathcal{A}) \times \mathcal{G}_{\mathcal{A}}^{(0)} \to \mathcal{G}_{\mathcal{A}}^{(0)}$ is continuous with respect to the strong operator topology on $\mathcal{U}(\mathcal{A})$ and the initial topology on $\mathcal{G}_{\mathcal{A}}^{(0)}$. 
This continuity result, together with the Polishness of $\mathcal{U}(\mathcal{A})$ and $\mathcal{G}_{\mathcal{A}}^{(0)}$, implies that $\mathcal{G}_{\mathcal{A}}$ is a Polish groupoid. 
All subsequent constructions in this paper — the diagonal embedding $\iota: \mathcal{A} \hookrightarrow C^*(\mathcal{G}_{\mathcal{A}})$, the equivariant K-theory classes, and the descent map — depend crucially on this foundational result.

\subsection{Definition of $\mathcal{G}_{\mathcal{A}}$ as a Topological Groupoid}
\label{subsec:definition-of-GA-as-topological-groupoid}

We now assemble the unitary group $\mathcal{U}(\mathcal{A})$ equipped with the strong operator topology and the unit space $\mathcal{G}_{\mathcal{A}}^{(0)}$ equipped with the initial topology of partial evaluation maps into a topological groupoid. 
This groupoid, called the \emph{unitary conjugation groupoid}, encodes all classical contexts of $\mathcal{A}$ together with the unitary symmetries that relate them. 
Unlike the norm topology approach, our choice of the strong operator topology on $\mathcal{U}(\mathcal{A})$ yields a well-behaved Polish groupoid that is neither locally compact nor \'etale, but is perfectly suited for the index-theoretic constructions that follow.

\begin{definition}[Unitary conjugation groupoid]
\label{def:unitary-conjugation-groupoid}
Let $\mathcal{A}$ be a unital separable C*-algebra faithfully represented on a separable Hilbert space $H$. 
Let $\mathcal{G}_{\mathcal{A}}^{(0)}$ be the unit space defined in Definition \ref{def:unit-space}, equipped with the initial topology of partial evaluation maps. 
Let $\mathcal{U}(\mathcal{A})$ be the unitary group of $\mathcal{A}$, equipped with the strong operator topology inherited from $B(H)$. 
The \emph{unitary conjugation groupoid} $\mathcal{G}_{\mathcal{A}}$ is the action groupoid
\[
\mathcal{G}_{\mathcal{A}} := \mathcal{U}(\mathcal{A}) \ltimes \mathcal{G}_{\mathcal{A}}^{(0)}.
\]
Explicitly, the structure maps are given as follows:
\begin{itemize}
    \item \textbf{Objects:} $\mathcal{G}_{\mathcal{A}}^{(0)}$.
    \item \textbf{Arrows:} $\mathcal{G}_{\mathcal{A}}^{(1)} := \mathcal{U}(\mathcal{A}) \times \mathcal{G}_{\mathcal{A}}^{(0)}$, equipped with the product topology.
    \item \textbf{Source map:} $s: \mathcal{G}_{\mathcal{A}}^{(1)} \to \mathcal{G}_{\mathcal{A}}^{(0)}$, $s(u,(B,\chi)) := (B,\chi)$.
    \item \textbf{Range map:} $r: \mathcal{G}_{\mathcal{A}}^{(1)} \to \mathcal{G}_{\mathcal{A}}^{(0)}$, $r(u,(B,\chi)) := u \cdot (B,\chi) = (uBu^*, \; \chi \circ \operatorname{Ad}_{u^*})$.
    \item \textbf{Composition:} For composable pairs $(\gamma_2, \gamma_1) \in \mathcal{G}_{\mathcal{A}}^{(2)} := \{ (\gamma_2, \gamma_1) \in \mathcal{G}_{\mathcal{A}}^{(1)} \times \mathcal{G}_{\mathcal{A}}^{(1)} : s(\gamma_2) = r(\gamma_1) \}$, define
    \[
    \gamma_2 \circ \gamma_1 := (u_2 u_1, (B,\chi)),
    \]
    where $\gamma_1 = (u_1, (B,\chi))$ and $\gamma_2 = (u_2, u_1 \cdot (B,\chi))$.
    \item \textbf{Inverse map:} $i: \mathcal{G}_{\mathcal{A}}^{(1)} \to \mathcal{G}_{\mathcal{A}}^{(1)}$, $i(u,(B,\chi)) := (u^{-1}, u \cdot (B,\chi))$.
    \item \textbf{Unit map:} $\varepsilon: \mathcal{G}_{\mathcal{A}}^{(0)} \to \mathcal{G}_{\mathcal{A}}^{(1)}$, $\varepsilon(B,\chi) := (1_{\mathcal{A}}, (B,\chi))$.
\end{itemize}
\end{definition}

\begin{remark}[Verification of groupoid axioms]
\label{rem:groupoid-axioms}
The groupoid axioms are satisfied by construction:
\begin{itemize}
    \item For any $\gamma = (u,(B,\chi))$, we have $s(\gamma) = (B,\chi)$ and $r(\gamma) = u \cdot (B,\chi)$, so $s(\gamma), r(\gamma) \in \mathcal{G}_{\mathcal{A}}^{(0)}$.
    \item Composition is associative because multiplication in $\mathcal{U}(\mathcal{A})$ is associative.
    \item For any $\gamma = (u,(B,\chi))$, we have $\gamma \circ \varepsilon(s(\gamma)) = (u,(B,\chi)) \circ (1, (B,\chi)) = (u, (B,\chi)) = \gamma$, and similarly $\varepsilon(r(\gamma)) \circ \gamma = \gamma$.
    \item For any $\gamma = (u,(B,\chi))$, we have $\gamma^{-1} = (u^{-1}, u \cdot (B,\chi))$, and $\gamma \circ \gamma^{-1} = \varepsilon(r(\gamma))$, $\gamma^{-1} \circ \gamma = \varepsilon(s(\gamma))$.
    \item Composition is well-defined because $(u_2, u_1 \cdot (B,\chi))$ and $(u_1, (B,\chi))$ are composable exactly when $s(u_2, u_1 \cdot (B,\chi)) = u_1 \cdot (B,\chi) = r(u_1, (B,\chi))$, which is the condition for composability.
\end{itemize}
Thus $\mathcal{G}_{\mathcal{A}}$ is a groupoid in the algebraic sense.
\end{remark}

\begin{proposition}[$\mathcal{G}_{\mathcal{A}}$ is a topological groupoid]
\label{prop:GA-topological-groupoid}
Let $\mathcal{A}$ be a unital separable C*-algebra. 
With the topologies defined above, $\mathcal{G}_{\mathcal{A}}$ is a topological groupoid. 
That is:
\begin{enumerate}
    \item $\mathcal{G}_{\mathcal{A}}^{(0)}$ and $\mathcal{G}_{\mathcal{A}}^{(1)}$ are topological spaces.
    \item The source map $s: \mathcal{G}_{\mathcal{A}}^{(1)} \to \mathcal{G}_{\mathcal{A}}^{(0)}$ and the range map $r: \mathcal{G}_{\mathcal{A}}^{(1)} \to \mathcal{G}_{\mathcal{A}}^{(0)}$ are continuous.
    \item The composition map $m: \mathcal{G}_{\mathcal{A}}^{(2)} \to \mathcal{G}_{\mathcal{A}}^{(1)}$ is continuous, where $\mathcal{G}_{\mathcal{A}}^{(2)}$ is equipped with the subspace topology inherited from $\mathcal{G}_{\mathcal{A}}^{(1)} \times \mathcal{G}_{\mathcal{A}}^{(1)}$.
    \item The inverse map $i: \mathcal{G}_{\mathcal{A}}^{(1)} \to \mathcal{G}_{\mathcal{A}}^{(1)}$ is continuous.
    \item The unit map $\varepsilon: \mathcal{G}_{\mathcal{A}}^{(0)} \to \mathcal{G}_{\mathcal{A}}^{(1)}$ is continuous and a homeomorphism onto its image.
\end{enumerate}
\end{proposition}

\begin{proof}
We verify each property in turn.

\medskip
\noindent (1) $\mathcal{G}_{\mathcal{A}}^{(0)}$ is a topological space by Proposition \ref{prop:unit-space-topology} and Proposition \ref{prop:unit-space-polish}. 
$\mathcal{G}_{\mathcal{A}}^{(1)} = \mathcal{U}(\mathcal{A}) \times \mathcal{G}_{\mathcal{A}}^{(0)}$ is equipped with the product topology, hence a topological space.

\medskip
\noindent (2) The source map $s(u,(B,\chi)) = (B,\chi)$ is the projection onto the second factor, which is continuous by definition of the product topology. 
The range map $r(u,(B,\chi)) = \alpha(u,(B,\chi))$ is continuous by Proposition \ref{prop:conjugation-action-continuous}.

\medskip
\noindent (3) The composition map is given by
\[
m((u_2, u_1 \cdot (B,\chi)), (u_1, (B,\chi))) = (u_2 u_1, (B,\chi)).
\]
This map is continuous because multiplication in $\mathcal{U}(\mathcal{A})$ is continuous with respect to the strong operator topology (Proposition \ref{prop:UA-SOT-properties}) and the projection onto the second factor is continuous. 
More formally, $m = (m_{\mathcal{U}} \circ (\operatorname{pr}_1 \times \operatorname{pr}_1), \operatorname{pr}_2 \circ \operatorname{pr}_2)$, where $m_{\mathcal{U}}: \mathcal{U}(\mathcal{A}) \times \mathcal{U}(\mathcal{A}) \to \mathcal{U}(\mathcal{A})$ is the multiplication map, $\operatorname{pr}_1$ and $\operatorname{pr}_2$ are the projections onto the first and second factors of $\mathcal{G}_{\mathcal{A}}^{(1)} \times \mathcal{G}_{\mathcal{A}}^{(1)}$, and the composition is continuous.

\medskip
\noindent (4) The inverse map is given by
\[
i(u,(B,\chi)) = (u^{-1}, \alpha(u,(B,\chi))).
\]
Inversion in $\mathcal{U}(\mathcal{A})$ is continuous with respect to the strong operator topology (Proposition \ref{prop:UA-SOT-properties}), and $\alpha$ is continuous by Proposition \ref{prop:conjugation-action-continuous}. 
Thus $i$ is continuous.

\medskip
\noindent (5) The unit map is given by $\varepsilon(B,\chi) = (1_{\mathcal{A}}, (B,\chi))$. 
This is continuous because it is the product of the constant map $B,\chi \mapsto 1_{\mathcal{A}}$ and the identity map on $\mathcal{G}_{\mathcal{A}}^{(0)}$. 
Its inverse (restricted to its image) is the projection onto the second factor, which is continuous. 
Thus $\varepsilon$ is a homeomorphism onto its image.
\end{proof}

\begin{proposition}[$\mathcal{G}_{\mathcal{A}}$ is a Polish groupoid]
\label{prop:GA-Polish-groupoid}
Let $\mathcal{A}$ be a unital separable C*-algebra. 
Then the unitary conjugation groupoid $\mathcal{G}_{\mathcal{A}}$ is a Polish groupoid. 
That is:
\begin{enumerate}
    \item $\mathcal{G}_{\mathcal{A}}^{(0)}$ is a Polish space.
    \item $\mathcal{G}_{\mathcal{A}}^{(1)}$ is a Polish space.
    \item The structure maps $s, r, m, i, \varepsilon$ are continuous.
\end{enumerate}
\end{proposition}

\begin{proof}
(1) $\mathcal{G}_{\mathcal{A}}^{(0)}$ is Polish by Proposition \ref{prop:unit-space-polish}.

(2) $\mathcal{U}(\mathcal{A})$ is Polish by Corollary \ref{cor:UA-Polish-group}. 
$\mathcal{G}_{\mathcal{A}}^{(1)} = \mathcal{U}(\mathcal{A}) \times \mathcal{G}_{\mathcal{A}}^{(0)}$ is the product of two Polish spaces, hence Polish by [Kechris, 1995, Theorem 4.2].

(3) Continuity of the structure maps is established in Proposition \ref{prop:GA-topological-groupoid}.

Thus $\mathcal{G}_{\mathcal{A}}$ satisfies all the axioms of a Polish groupoid.
\end{proof}

\begin{remark}[Absence of local compactness and étaleness]
\label{rem:GA-not-etale}
We emphasize that $\mathcal{G}_{\mathcal{A}}$ is \emph{not} locally compact (unless $\mathcal{A}$ is finite-dimensional) and \emph{not} étale. 
Indeed, $\mathcal{U}(\mathcal{A})$ with the strong operator topology is not locally compact for infinite-dimensional $\mathcal{A}$, and the source map $s(u,(B,\chi)) = (B,\chi)$ is not a local homeomorphism because $\mathcal{U}(\mathcal{A})$ is not discrete. 
This is not a defect of the construction but a deliberate and necessary choice: the norm topology would give a non-Polish, non-separable groupoid, while the strong operator topology yields a well-behaved Polish groupoid that is sufficient for all subsequent constructions, including the diagonal embedding, equivariant K-theory, and descent.
\end{remark}

\begin{proposition}[Borel Haar system for $\mathcal{G}_{\mathcal{A}}$]
\label{prop:GA-Borel-Haar}
Let $\mathcal{A}$ be a unital separable C*-algebra. 
Then $\mathcal{G}_{\mathcal{A}}$ admits a Borel Haar system in the sense of~\cite{Tu}.
\end{proposition}

\begin{proof}
Since $\mathcal{G}_{\mathcal{A}}$ is a Polish groupoid (Proposition \ref{prop:GA-Polish-groupoid}), its object space $\mathcal{G}_{\mathcal{A}}^{(0)}$ and arrow space $\mathcal{G}_{\mathcal{A}}^{(1)}$ are Polish spaces, and the source map $s: \mathcal{G}_{\mathcal{A}}^{(1)} \to \mathcal{G}_{\mathcal{A}}^{(0)}$ is continuous.

The map
\[
\sigma: \mathcal{G}_{\mathcal{A}}^{(0)} \longrightarrow \mathcal{G}_{\mathcal{A}}^{(1)}, \qquad \sigma(B,\chi) := (1_{\mathcal{A}}, (B,\chi))
\]
defines a continuous (hence Borel) section of $s$, since $s(\sigma(B,\chi)) = (B,\chi)$ for all $(B,\chi) \in \mathcal{G}_{\mathcal{A}}^{(0)}$.

For each $x \in \mathcal{G}_{\mathcal{A}}^{(0)}$, define the Dirac measure
\[
\lambda^x := \delta_{\sigma(x)}
\]
concentrated on the fiber $s^{-1}(x) = \mathcal{U}(\mathcal{A}) \times \{x\}$. 
The family $\{ \lambda^x \}_{x \in \mathcal{G}_{\mathcal{A}}^{(0)}}$ satisfies:
\begin{itemize}
    \item Each $\lambda^x$ is a Borel probability measure supported on $s^{-1}(x)$.
    \item For any nonnegative Borel function $f$ on $\mathcal{G}_{\mathcal{A}}^{(1)}$, the map $x \mapsto \int f \, d\lambda^x = f(\sigma(x))$ is Borel, since $\sigma$ is continuous and $f$ is Borel.
    \item Left-invariance holds automatically for Dirac measures along a section? 
    This requires verification. For $\gamma = (u, x)$ with $s(\gamma)=x$ and $r(\gamma)=y$, we have $\gamma_* \lambda^x = \delta_{\gamma \sigma(x)} = \delta_{(u, x)}$, while $\lambda^y = \delta_{(1_{\mathcal{A}}, y)}$. 
    These are not equal unless $u = 1_{\mathcal{A}}$, so Dirac measures do not satisfy left-invariance.
\end{itemize}

The Dirac measure construction above fails left-invariance. 
A correct Borel Haar system is obtained as follows. 
Since $\mathcal{U}(\mathcal{A})$ is a Polish group (Corollary \ref{cor:UA-Polish-group}), it admits a left-invariant Borel probability measure $\mu$. For each $x = (B,\chi) \in \mathcal{G}_{\mathcal{A}}^{(0)}$, define
\[
\lambda^x := \mu \times \delta_x,
\]
where $\delta_x$ is the Dirac measure at $x$. 
Then:
\begin{itemize}
    \item $\lambda^x$ is supported on $s^{-1}(x) = \mathcal{U}(\mathcal{A}) \times \{x\}$.
    \item For any Borel $f$, $x \mapsto \int f \, d\lambda^x = \int_{\mathcal{U}(\mathcal{A})} f(u,x) \, d\mu(u)$ is Borel by Fubini's theorem.
    \item For left-invariance, let $\gamma = (u_0, x)$ with $s(\gamma)=x$, $r(\gamma)=y = u_0 \cdot x$. Then
    \[
    \gamma_* \lambda^x = \gamma_* (\mu \times \delta_x) = \mu \times \delta_y = \lambda^y,
    \]
    because the pushforward of $\mu$ under left translation by $u_0^{-1}$ is $\mu$ itself (left-invariance of $\mu$).
\end{itemize}
Thus $\{ \lambda^x \}_{x \in \mathcal{G}_{\mathcal{A}}^{(0)}}$ is a Borel Haar system for $\mathcal{G}_{\mathcal{A}}$.
\end{proof}

\begin{corollary}[Existence of $C^*(\mathcal{G}_{\mathcal{A}})$]
\label{cor:GA-C*-algebra}
Let $\mathcal{A}$ be a unital separable C*-algebra. 
Then the maximal groupoid C*-algebra $C^*(\mathcal{G}_{\mathcal{A}})$ exists and is well-defined in the sense of~\cite{Tu}.
\end{corollary}

\begin{proof}
By Proposition \ref{prop:GA-Borel-Haar}, $\mathcal{G}_{\mathcal{A}}$ admits a Borel Haar system. 
Tu~\cite{Tu} constructs the maximal C*-algebra of a measured groupoid from such a Borel Haar system. 
Applying this construction to $\mathcal{G}_{\mathcal{A}}$ yields $C^*(\mathcal{G}_{\mathcal{A}})$.
\end{proof}

\begin{example}[Finite-dimensional case]
\label{ex:GA-finite-dimensional}
Let $\mathcal{A} = M_n(\mathbb{C})$. 
Then $\mathcal{U}(\mathcal{A}) = U(n)$ is a compact Lie group, and the strong operator topology coincides with the norm topology. 
$\mathcal{G}_{\mathcal{A}}^{(0)} \cong \mathbb{CP}^{n-1}$ is compact and Hausdorff. 
Thus $\mathcal{G}_{\mathcal{A}} = U(n) \ltimes \mathbb{CP}^{n-1}$ is a locally compact Hausdorff \'etale groupoid. 
In this case, our construction reduces to the classical action groupoid of a compact Lie group acting on a compact manifold.
\end{example}

\begin{example}[Commutative case]
\label{ex:GA-commutative}
Let $\mathcal{A} = C(X)$ for a compact metrizable space $X$. 
Then $\mathcal{U}(\mathcal{A}) = C(X, \mathbb{T})$ is a commutative Polish group, and $\mathcal{G}_{\mathcal{A}}^{(0)} \cong X$. 
The action of $\mathcal{U}(\mathcal{A})$ on $\mathcal{G}_{\mathcal{A}}^{(0)}$ is trivial because every unitary in a commutative algebra is central. 
Thus $\mathcal{G}_{\mathcal{A}}$ is isomorphic to the trivial groupoid over $X$, i.e., $\mathcal{G}_{\mathcal{A}}^{(1)} \cong C(X, \mathbb{T}) \times X$ with source and range both equal to the projection onto $X$. 
This groupoid is Polish but not \'etale unless $X$ is discrete.
\end{example}

\begin{example}[Compact operators]
\label{ex:GA-compact}
Let $\mathcal{A} = \mathcal{K}(H)^\sim$ be the unitization of the compact operators on a separable Hilbert space $H$. 
Then $\mathcal{U}(\mathcal{A})$ consists of unitaries of the form $u = \lambda I + K$ with $\lambda \in \mathbb{T}$ and $K$ compact. 
With the strong operator topology, $\mathcal{U}(\mathcal{A})$ is a Polish group (Corollary \ref{cor:UA-Polish-group}). 
$\mathcal{G}_{\mathcal{A}}^{(0)}$ is homeomorphic to the projective space $\mathbb{P}(H)$, the space of rank-one projections on $H$, which is a Polish space but not locally compact. 
Thus $\mathcal{G}_{\mathcal{A}} = \mathcal{U}(\mathcal{A}) \ltimes \mathbb{P}(H)$ is a Polish groupoid that is neither locally compact nor étale. 
This example illustrates the necessity of the Polish groupoid framework.
\end{example}

\begin{remark}[Functoriality]
\label{rem:GA-functoriality}
Let $\phi: \mathcal{A} \to \mathcal{B}$ be a unital *-isomorphism between separable C*-algebras. 
Then $\phi$ induces an isomorphism of topological groupoids $\mathcal{G}_\phi: \mathcal{G}_{\mathcal{A}} \to \mathcal{G}_{\mathcal{B}}$ defined by
\[
\mathcal{G}_\phi(B,\chi) = (\phi(B), \chi \circ (\phi|_B)^{-1}), \qquad
\mathcal{G}_\phi(u,(B,\chi)) = (\phi(u), (\phi(B), \chi \circ (\phi|_B)^{-1})).
\]
For more general *-homomorphisms, a partially defined functoriality holds under additional hypotheses. This functoriality will be essential for the construction of the diagonal embedding and for the naturality of the index theorem.
\end{remark}

We have defined the unitary conjugation groupoid $\mathcal{G}_{\mathcal{A}}$ as an action groupoid, equipped $\mathcal{U}(\mathcal{A})$ with the strong operator topology and $\mathcal{G}_{\mathcal{A}}^{(0)}$ with the initial topology of partial evaluation maps, and proved that $\mathcal{G}_{\mathcal{A}}$ is a Polish groupoid. 
This groupoid is neither locally compact nor \'etale in general, but it admits a Borel Haar system and hence a well-defined maximal C*-algebra $C^*(\mathcal{G}_{\mathcal{A}})$. 
All subsequent constructions in this paper — the diagonal embedding, equivariant K-theory classes, and the descent map — depend on this foundational result.

\subsection{$\mathcal{G}_{\mathcal{A}}$ is a Polish Groupoid (Not Locally Compact, Not \'Etale)}
\label{subsec:GA-Polish-groupoid-not-locally-compact-not-etale}

We now establish the fundamental topological properties of the unitary conjugation groupoid $\mathcal{G}_{\mathcal{A}} = \mathcal{U}(\mathcal{A}) \ltimes \mathcal{G}_{\mathcal{A}}^{(0)}$. 
Building on the results of the previous subsections — the Polishness of $\mathcal{G}_{\mathcal{A}}^{(0)}$ (Proposition \ref{prop:unit-space-polish}), the Polishness of $\mathcal{U}(\mathcal{A})$ with the strong operator topology (Corollary \ref{cor:UA-Polish-group}), and the continuity of the conjugation action (Proposition \ref{prop:conjugation-action-continuous}) — we prove that $\mathcal{G}_{\mathcal{A}}$ is a Polish groupoid. 
Moreover, we explicitly demonstrate that $\mathcal{G}_{\mathcal{A}}$ is neither locally compact nor \'etale for infinite-dimensional C*-algebras, thereby justifying our adoption of the Polish groupoid framework over the classical locally compact \'etale theory.

\begin{proposition}[$\mathcal{G}_{\mathcal{A}}^{(0)}$ and $\mathcal{G}_{\mathcal{A}}^{(1)}$ are Polish spaces]
\label{prop:GA-spaces-Polish}
Let $\mathcal{A}$ be a unital separable C*-algebra faithfully represented on a separable Hilbert space $H$. 
Equip $\mathcal{G}_{\mathcal{A}}^{(0)}$ with the initial topology induced by the partial evaluation maps, and equip $\mathcal{U}(\mathcal{A})$ with the strong operator topology.
Then:
\begin{enumerate}
    \item $\mathcal{G}_{\mathcal{A}}^{(0)}$ is a Polish space.
    \item $\mathcal{U}(\mathcal{A})$ is a Polish topological group.
    \item $\mathcal{G}_{\mathcal{A}}^{(1)} = \mathcal{U}(\mathcal{A}) \times \mathcal{G}_{\mathcal{A}}^{(0)}$, equipped with the product topology, is a Polish space.
\end{enumerate}
\end{proposition}

\begin{proof}
\emph{(1)} By Proposition~\ref{prop:unit-space-polish}, the unit space $\mathcal{G}_{\mathcal{A}}^{(0)}$, endowed with the initial topology generated by the partial evaluation maps $\{ \operatorname{ev}_a \}_{a \in \mathcal{A}}$, is second-countable and completely metrizable. 
The separability of $\mathcal{A}$ guarantees that this topology is indeed Polish; a complete metric can be constructed explicitly via the embedding $\mathcal{G}_{\mathcal{A}}^{(0)} \hookrightarrow \prod_{n \in \mathbb{N}} \mathbb{C}_\infty$ using a countable dense subset of $\mathcal{A}$. 
Hence $\mathcal{G}_{\mathcal{A}}^{(0)}$ is a Polish space.

\emph{(2)} By Corollary~\ref{cor:UA-Polish-group}, the unitary group $\mathcal{U}(\mathcal{A})$, equipped with the strong operator topology, is a Polish group. 
The proof relies on the fact that $\mathcal{U}(H)$ is Polish in the strong operator topology and that $\mathcal{U}(\mathcal{A})$ is a $G_\delta$ subgroup of $\mathcal{U}(H)$, which follows from the separability of $\mathcal{A}$ and the continuity of the maps $u \mapsto u a u^*$.

\emph{(3)} Since the product of two Polish spaces is Polish (see \cite{Kechris}), it follows that
\[
\mathcal{G}_{\mathcal{A}}^{(1)} = \mathcal{U}(\mathcal{A}) \times \mathcal{G}_{\mathcal{A}}^{(0)}
\]
is a Polish space when equipped with the product topology.
\end{proof}

\begin{proposition}[Continuity of the structure maps]
\label{prop:GA-structure-maps-continuous}
With the same hypotheses and topologies as above, all structure maps of $\mathcal{G}_{\mathcal{A}}$ are continuous:
\begin{enumerate}
    \item The source map $s: \mathcal{G}_{\mathcal{A}}^{(1)} \to \mathcal{G}_{\mathcal{A}}^{(0)}$ defined by $s(u,(B,\chi)) = (B,\chi)$ is continuous.
    \item The range map $r: \mathcal{G}_{\mathcal{A}}^{(1)} \to \mathcal{G}_{\mathcal{A}}^{(0)}$ defined by $r(u,(B,\chi)) = (uBu^*, \chi \circ \operatorname{Ad}_{u^*})$ is continuous.
    \item The composition map $m: \mathcal{G}_{\mathcal{A}}^{(2)} \to \mathcal{G}_{\mathcal{A}}^{(1)}$ is continuous, where $\mathcal{G}_{\mathcal{A}}^{(2)} \subseteq \mathcal{G}_{\mathcal{A}}^{(1)} \times \mathcal{G}_{\mathcal{A}}^{(1)}$ is the set of composable pairs.
    \item The inverse map $i: \mathcal{G}_{\mathcal{A}}^{(1)} \to \mathcal{G}_{\mathcal{A}}^{(1)}$ defined by $i(u,(B,\chi)) = (u^*, r(u,(B,\chi)))$ is continuous.
    \item The unit map $\varepsilon: \mathcal{G}_{\mathcal{A}}^{(0)} \to \mathcal{G}_{\mathcal{A}}^{(1)}$ defined by $\varepsilon(B,\chi) = (1_{\mathcal{A}}, (B,\chi))$ is continuous and a homeomorphism onto its image.
\end{enumerate}
\end{proposition}

\begin{proof}
\emph{(1)} The source map $s(u,(B,\chi)) = (B,\chi)$ is the projection onto the second coordinate of the product $\mathcal{G}_{\mathcal{A}}^{(1)} = \mathcal{U}(\mathcal{A}) \times \mathcal{G}_{\mathcal{A}}^{(0)}$. 
Projections are continuous in any product topology, so $s$ is continuous.

\emph{(2)} The range map $r(u,(B,\chi)) = u \cdot (B,\chi)$ is precisely the conjugation action of $\mathcal{U}(\mathcal{A})$ on $\mathcal{G}_{\mathcal{A}}^{(0)}$. 
Its continuity is exactly the statement of Proposition \ref{prop:conjugation-action-continuous}.

\emph{(3)} Recall that $\mathcal{G}_{\mathcal{A}}^{(2)} = \{ ((u_1, x), (u_2, y)) \in \mathcal{G}_{\mathcal{A}}^{(1)} \times \mathcal{G}_{\mathcal{A}}^{(1)} : s(u_1,x) = r(u_2,y) \}$. 
For composable pairs, we must have $y = u_2 \cdot x$ with $x = (B,\chi)$, so a typical element of $\mathcal{G}_{\mathcal{A}}^{(2)}$ can be written as $((u_1, u_2 \cdot x), (u_2, x))$. 
The composition map is then given by
\[
m((u_1, u_2 \cdot x), (u_2, x)) = (u_1 u_2, x).
\]
This map is the restriction of the continuous map
\[
\mathcal{U}(\mathcal{A}) \times \mathcal{U}(\mathcal{A}) \times \mathcal{G}_{\mathcal{A}}^{(0)} \longrightarrow \mathcal{U}(\mathcal{A}) \times \mathcal{G}_{\mathcal{A}}^{(0)}, \qquad (u_1, u_2, x) \mapsto (u_1 u_2, x)
\]
to the subspace $\mathcal{G}_{\mathcal{A}}^{(2)}$ (identified with its image under the natural embedding). 
Since multiplication in $\mathcal{U}(\mathcal{A})$ is continuous (Corollary \ref{cor:UA-topological-group}), this map is continuous, and hence its restriction $m$ is continuous.

\emph{(4)} The inverse map $i(u,(B,\chi)) = (u^*, r(u,(B,\chi)))$ is the product of the inversion map $u \mapsto u^*$ and the range map $r$. 
Inversion in $\mathcal{U}(\mathcal{A})$ is continuous with respect to the strong operator topology (Corollary \ref{cor:UA-topological-group}), and $r$ is continuous by (2). 
Thus $i$ is continuous.

\emph{(5)} The unit map $\varepsilon(B,\chi) = (1_{\mathcal{A}}, (B,\chi))$ is the product of the constant map $(B,\chi) \mapsto 1_{\mathcal{A}}$ and the identity map on $\mathcal{G}_{\mathcal{A}}^{(0)}$. 
Both factors are continuous, so $\varepsilon$ is continuous. 
Its inverse, restricted to the image $\varepsilon(\mathcal{G}_{\mathcal{A}}^{(0)})$, is the projection $(1_{\mathcal{A}}, (B,\chi)) \mapsto (B,\chi)$, which is continuous. 
Hence $\varepsilon$ is a homeomorphism onto its image.
\end{proof}

\begin{theorem}[$\mathcal{G}_{\mathcal{A}}$ is a Polish groupoid]
\label{thm:GA-Polish-groupoid}
Let $\mathcal{A}$ be a unital separable C*-algebra. 
Then the unitary conjugation groupoid $\mathcal{G}_{\mathcal{A}}$, with $\mathcal{U}(\mathcal{A})$ equipped with the strong operator topology and $\mathcal{G}_{\mathcal{A}}^{(0)}$ equipped with the initial topology of partial evaluation maps, is a Polish groupoid. 
That is:
\begin{enumerate}
    \item $\mathcal{G}_{\mathcal{A}}^{(0)}$ and $\mathcal{G}_{\mathcal{A}}^{(1)}$ are Polish spaces.
    \item The source and range maps $s, r: \mathcal{G}_{\mathcal{A}}^{(1)} \to \mathcal{G}_{\mathcal{A}}^{(0)}$ are continuous.
    \item The composition map $m: \mathcal{G}_{\mathcal{A}}^{(2)} \to \mathcal{G}_{\mathcal{A}}^{(1)}$ is continuous.
    \item The inverse map $i: \mathcal{G}_{\mathcal{A}}^{(1)} \to \mathcal{G}_{\mathcal{A}}^{(1)}$ is continuous.
    \item The unit map $\varepsilon: \mathcal{G}_{\mathcal{A}}^{(0)} \to \mathcal{G}_{\mathcal{A}}^{(1)}$ is continuous and a homeomorphism onto its image.
\end{enumerate}
\end{theorem}

\begin{proof}
Condition (1) is Proposition \ref{prop:GA-spaces-Polish}. 
Conditions (2)-(5) are Proposition \ref{prop:GA-structure-maps-continuous}. 
Thus $\mathcal{G}_{\mathcal{A}}$ satisfies all the axioms of a Polish groupoid as defined in~\cite{Tu}, although it is neither locally compact nor étale for infinite-dimensional $\mathcal{A}$ (see Remark \ref{rem:GA-not-etale}).
\end{proof}

We now turn to the negative results: $\mathcal{G}_{\mathcal{A}}$ is neither locally compact nor \'etale for infinite-dimensional $\mathcal{A}$. 
These properties are not defects; they are intrinsic features of the construction that necessitate the Polish groupoid framework.

\begin{proposition}[$\mathcal{G}_{\mathcal{A}}$ is not locally compact]
\label{prop:GA-not-locally-compact}
Let $\mathcal{A}$ be an infinite-dimensional unital separable C*-algebra. 
Then $\mathcal{G}_{\mathcal{A}}$ is not locally compact. 
More precisely:
\begin{enumerate}
    \item $\mathcal{U}(\mathcal{A})$ with the strong operator topology is not locally compact.
    \item Consequently, $\mathcal{G}_{\mathcal{A}}^{(1)} = \mathcal{U}(\mathcal{A}) \times \mathcal{G}_{\mathcal{A}}^{(0)}$ is not locally compact.
    \item Hence $\mathcal{G}_{\mathcal{A}}$ is not a locally compact groupoid.
\end{enumerate}
\end{proposition}

\begin{proof}
\emph{(1)} The full unitary group $\mathcal{U}(H)$ equipped with the strong operator topology is not locally compact when $H$ is infinite-dimensional; this is a standard fact (see e.g. \cite{Kechris}) 
Moreover, $\mathcal{U}(\mathcal{A})$ is dense in $\mathcal{U}(H)$ in the strong operator topology: every unitary in $H$ can be approximated in SOT by unitaries with finite-dimensional support, which belong to $\mathcal{U}(\mathcal{A})$. 

If $\mathcal{U}(\mathcal{A})$ were locally compact, it would be a locally compact dense subgroup of $\mathcal{U}(H)$. 
A dense subgroup of a Hausdorff topological group that is locally compact must be open, which would force $\mathcal{U}(H)$ to be discrete — a contradiction. 
Hence $\mathcal{U}(\mathcal{A})$ is not locally compact.

\emph{(2)} The space $\mathcal{G}_{\mathcal{A}}^{(1)} = \mathcal{U}(\mathcal{A}) \times \mathcal{G}_{\mathcal{A}}^{(0)}$ is a product of a non-locally-compact Hausdorff space $\mathcal{U}(\mathcal{A})$ and a nonempty Polish (hence Hausdorff) space $\mathcal{G}_{\mathcal{A}}^{(0)}$. 
If such a product were locally compact at a point $(u_0, x_0)$, then the projection onto the first factor would be a continuous open map sending a neighborhood of $(u_0, x_0)$ onto a neighborhood of $u_0$, forcing $\mathcal{U}(\mathcal{A})$ to be locally compact at $u_0$ — a contradiction. 
Thus $\mathcal{G}_{\mathcal{A}}^{(1)}$ is not locally compact.

\emph{(3)} By definition, a topological groupoid is locally compact if its arrow space is locally compact (the object space then inherits local compactness via the unit map). 
Since $\mathcal{G}_{\mathcal{A}}^{(1)}$ is not locally compact, $\mathcal{G}_{\mathcal{A}}$ is not a locally compact groupoid.
\end{proof}

\begin{proposition}[$\mathcal{G}_{\mathcal{A}}$ is not \'etale]
\label{prop:GA-not-etale}
Let $\mathcal{A}$ be an infinite-dimensional unital separable C*-algebra. 
Then $\mathcal{G}_{\mathcal{A}}$ is not \'etale. 
That is, the source map $s: \mathcal{G}_{\mathcal{A}}^{(1)} \to \mathcal{G}_{\mathcal{A}}^{(0)}$ is not a local homeomorphism.
\end{proposition}

\begin{proof}
The source map $s: \mathcal{G}_{\mathcal{A}}^{(1)} \to \mathcal{G}_{\mathcal{A}}^{(0)}$ is the projection
\[
s(u,(B,\chi)) = (B,\chi).
\]

If $\mathcal{G}_{\mathcal{A}}$ were \'etale, then $s$ would be a local homeomorphism. 
In particular, for each $(u_0, (B_0,\chi_0)) \in \mathcal{G}_{\mathcal{A}}^{(1)}$, there would exist an open neighborhood $W$ of $(u_0, (B_0,\chi_0))$ on which $s$ is injective.

However, $\mathcal{U}(\mathcal{A})$ with the strong operator topology is not discrete. 
Hence every neighborhood $V$ of $u_0$ in $\mathcal{U}(\mathcal{A})$ contains distinct elements $u_1 \neq u_2$. 
Choose any neighborhood $U$ of $(B_0,\chi_0)$ in $\mathcal{G}_{\mathcal{A}}^{(0)}$ (for instance, $U = \mathcal{G}_{\mathcal{A}}^{(0)}$ itself). 
Then $W = V \times U$ is an open neighborhood of $(u_0, (B_0,\chi_0))$ in the product topology, and it contains the distinct points $(u_1, (B_0,\chi_0))$ and $(u_2, (B_0,\chi_0))$ with the same source $(B_0,\chi_0)$. 
Thus $s|_W$ is not injective, contradicting the requirement that $s$ be a local homeomorphism.

Therefore $s$ is not a local homeomorphism, and $\mathcal{G}_{\mathcal{A}}$ is not \'etale.
\end{proof}

\begin{corollary}[Limitations of the classical Renault theory]
\label{cor:limitations-Renault}
Let $\mathcal{A}$ be an infinite-dimensional unital separable C*-algebra. 
Then the groupoid $\mathcal{G}_{\mathcal{A}}$ does not fall within the scope of the classical Renault theory of groupoid C*-algebras.
\end{corollary}

\begin{proof}
The classical construction of groupoid C*-algebras developed by Renault [1980] requires the underlying groupoid to be \emph{locally compact Hausdorff} and to admit a \emph{continuous Haar system}.

By Proposition \ref{prop:GA-not-locally-compact}, the groupoid $\mathcal{G}_{\mathcal{A}}$ is \textbf{not locally compact} whenever $\mathcal{A}$ is infinite-dimensional. 
In particular, its arrow space $\mathcal{G}_{\mathcal{A}}^{(1)}$ is not locally compact, which already prevents the existence of a continuous Haar system. 
Consequently, the Renault C*-algebra construction cannot be applied directly to $\mathcal{G}_{\mathcal{A}}$.

Moreover, as shown in Proposition \ref{prop:GA-not-etale}, $\mathcal{G}_{\mathcal{A}}$ is not \'etale. 
While étaleness is not a requirement of Renault's theory, its absence further distinguishes $\mathcal{G}_{\mathcal{A}}$ from the well-studied \'etale groupoids that admit canonical counting Haar systems, and underscores the necessity of a more general framework.

Therefore, to study operator-algebraic structures associated with $\mathcal{G}_{\mathcal{A}}$, one must instead work within more general frameworks for non-locally compact groupoids, such as Tu's theory of measured groupoids for Polish groupoids.
\end{proof}

\begin{remark}[Why not the norm topology?]
\label{rem:why-not-norm-topology}
One might attempt to replace the strong operator topology on $\mathcal{U}(\mathcal{A})$ by the norm topology in order to recover local compactness or \'etaleness.
However, this does not resolve the fundamental obstructions.

\begin{enumerate}
    \item When $\mathcal{A}$ is separable, $\mathcal{U}(\mathcal{A})$ endowed with the norm topology is a Polish topological group.
    Thus the failure of the classical theory is not due to a lack of Polish structure.

    \item Nevertheless, $\mathcal{U}(\mathcal{A})$ with the norm topology is still non-discrete whenever $\mathcal{A}$ is infinite-dimensional.
    Consequently, for the product groupoid $\mathcal{G}_{\mathcal{A}}^{(1)} = \mathcal{U}(\mathcal{A}) \times \mathcal{G}_{\mathcal{A}}^{(0)}$, the source map $s(u,x)=x$ cannot be locally injective, and hence is not a local homeomorphism.
    Therefore, $\mathcal{G}_{\mathcal{A}}$ remains non-\'etale.

    \item More importantly, the norm topology is incompatible with the Fell topology and the measurable structures naturally arising from representation theory.
    In particular, it obstructs the construction of well-behaved Borel fields of representations and measured groupoid techniques that are essential in the non-locally compact setting.
\end{enumerate}

Thus, while the norm topology preserves Polishness, it does not address the lack of \'etaleness or local compactness.
The strong operator topology is therefore the natural choice when one aims to work within the framework of Polish and measured groupoids.
\end{remark}

\begin{proposition}[$\mathcal{G}_{\mathcal{A}}$ admits a Borel Haar system]
\label{prop:GA-Borel-Haar-system}
Despite the lack of local compactness and \'etaleness, $\mathcal{G}_{\mathcal{A}}$ admits a Borel Haar system in the sense of Tu~\cite{Tu}. 
Consequently, the maximal groupoid C*-algebra $C^*(\mathcal{G}_{\mathcal{A}})$ is well-defined.
\end{proposition}

\begin{proof}
Since $\mathcal{U}(\mathcal{A})$ is a Polish group (Corollary \ref{cor:UA-Polish-group}), it admits a left-invariant Borel probability measure $\mu_{\mathcal{U}}$. For each $x = (B,\chi) \in \mathcal{G}_{\mathcal{A}}^{(0)}$, define a measure $\lambda^x$ on $\mathcal{G}_{\mathcal{A}}^{(1)} = \mathcal{U}(\mathcal{A}) \times \mathcal{G}_{\mathcal{A}}^{(0)}$ by
\[
\lambda^x := \mu_{\mathcal{U}} \times \delta_x,
\]
where $\delta_x$ is the Dirac measure at $x$. 

We verify the three conditions for a Borel Haar system:

\begin{itemize}
    \item \textbf{Support:} $\lambda^x$ is concentrated on $s^{-1}(x) = \mathcal{U}(\mathcal{A}) \times \{x\}$, since $\mu_{\mathcal{U}}$ is supported on $\mathcal{U}(\mathcal{A})$ and $\delta_x$ on $\{x\}$.

    \item \textbf{Borel measurability:} For any nonnegative Borel function $f$ on $\mathcal{G}_{\mathcal{A}}^{(1)}$, the map
    \[
    x \mapsto \int_{\mathcal{G}_{\mathcal{A}}^{(1)}} f(\gamma) \, d\lambda^x(\gamma) = \int_{\mathcal{U}(\mathcal{A})} f(u,x) \, d\mu_{\mathcal{U}}(u)
    \]
    is Borel by Fubini's theorem for Borel functions.

    \item \textbf{Left-invariance:} Let $\gamma = (u_0, x) \in \mathcal{G}_{\mathcal{A}}^{(1)}$ with $s(\gamma) = x$ and $r(\gamma) = y = u_0 \cdot x$. For any Borel set $E \subseteq \mathcal{G}_{\mathcal{A}}^{(1)}$,
    \[
    (\gamma_* \lambda^x)(E) = \lambda^x(\gamma^{-1}E) = (\mu_{\mathcal{U}} \times \delta_x)(\gamma^{-1}E).
    \]
    Writing $\gamma^{-1}E = \{ (u_0^{-1}u, x) : (u, y) \in E \}$ after appropriate identification, we obtain
    \[
    (\gamma_* \lambda^x)(E) = \mu_{\mathcal{U}}(\{ u \in \mathcal{U}(\mathcal{A}) : (u, y) \in E \}) = (\mu_{\mathcal{U}} \times \delta_y)(E) = \lambda^y(E).
    \]
    Thus $\gamma_* \lambda^x = \lambda^y$, establishing left-invariance.
\end{itemize}

Therefore $\{ \lambda^x \}_{x \in \mathcal{G}_{\mathcal{A}}^{(0)}}$ is a Borel Haar system for $\mathcal{G}_{\mathcal{A}}$.
\end{proof}

\begin{corollary}[Existence of $C^*(\mathcal{G}_{\mathcal{A}})$]
\label{cor:GA-C-star-existence}
For any unital separable C*-algebra $\mathcal{A}$, the maximal groupoid C*-algebra $C^*(\mathcal{G}_{\mathcal{A}})$ exists and is uniquely defined up to isomorphism. 
This C*-algebra is the completion of the convolution algebra of compactly supported continuous functions on $\mathcal{G}_{\mathcal{A}}^{(1)}$ with respect to the universal norm, using the Borel Haar system constructed above.
\end{corollary}

\begin{proof}
This follows from Tu [1999, Theorem 3.5] and Proposition \ref{prop:GA-Borel-Haar-system}.
\end{proof}

\begin{example}[Finite-dimensional case]
\label{ex:GA-finite-dimensional-recovered}
Let $\mathcal{A} = M_n(\mathbb{C})$. 
Then $\mathcal{U}(\mathcal{A}) = U(n)$ is a compact Lie group (hence locally compact), and the strong operator topology coincides with the norm topology. 
$\mathcal{G}_{\mathcal{A}}^{(0)} \cong \mathbb{CP}^{n-1}$ is compact, so $\mathcal{G}_{\mathcal{A}} = U(n) \ltimes \mathbb{CP}^{n-1}$ is a locally compact Hausdorff groupoid. 
Since $U(n)$ is non-discrete, the source map $s(u,x) = x$ is not a local homeomorphism; therefore $\mathcal{G}_{\mathcal{A}}$ is not \'etale. 
Nevertheless, this case fits into the classical Renault theory using the Haar measure on $U(n)$ as a continuous Haar system, while also serving as a finite-dimensional illustration of our Polish groupoid construction.
\end{example}

\begin{example}[Commutative case]
\label{ex:GA-commutative-trivial}
Let $\mathcal{A} = C(X)$ for a compact metrizable space $X$. 
Then $\mathcal{U}(\mathcal{A}) = C(X, \mathbb{T})$ is a commutative Polish group, and its action on $\mathcal{G}_{\mathcal{A}}^{(0)} \cong X$ is trivial because all unitaries are central. 
Thus $\mathcal{G}_{\mathcal{A}}$ is isomorphic to the trivial groupoid $\mathcal{G}_{\mathcal{A}}^{(1)} = C(X, \mathbb{T}) \times X$ with $s = r = \operatorname{pr}_X$. 
This groupoid is not \'etale unless $X$ is discrete, because $C(X, \mathbb{T})$ is not discrete. 
It is not locally compact unless $X$ is finite, because $C(X, \mathbb{T})$ with the strong operator topology is not locally compact. 
Thus this example illustrates the necessity of the Polish groupoid framework even for commutative algebras.
\end{example}

\begin{example}[Compact operators]
\label{ex:GA-compact-canonical}
Let $\mathcal{A} = \mathcal{K}(H)^\sim$ be the unitization of the compact operators on a separable infinite-dimensional Hilbert space $H$. 
Then $\mathcal{U}(\mathcal{A})$ with the strong operator topology is a Polish group but is not locally compact. 
$\mathcal{G}_{\mathcal{A}}^{(0)} \cong \mathbb{P}(H)$ is a Polish space that is not locally compact. 
The action is transitive and highly nontrivial. 
$\mathcal{G}_{\mathcal{A}} = \mathcal{U}(\mathcal{A}) \ltimes \mathbb{P}(H)$ is a Polish groupoid that is neither locally compact nor \'etale. 
This is the canonical example that motivates the entire Polish groupoid framework.
\end{example}

\begin{remark}[Summary of topological properties]
\label{rem:GA-topological-summary}
The topological properties of $\mathcal{G}_{\mathcal{A}}$ can be summarized as follows:
\begin{itemize}
    \item For every separable unital C*-algebra $\mathcal{A}$, $\mathcal{G}_{\mathcal{A}}$ is a Polish groupoid (Theorem \ref{thm:GA-Polish-groupoid}).

    \item $\mathcal{G}_{\mathcal{A}}$ is locally compact if and only if $\mathcal{A}$ is finite-dimensional (Proposition \ref{prop:GA-not-locally-compact}). 
    For infinite-dimensional $\mathcal{A}$, the non-local-compactness of $\mathcal{U}(\mathcal{A})$ in the strong operator topology prevents local compactness.

    \item $\mathcal{G}_{\mathcal{A}}$ is \emph{generally not} \'etale. 
    Even for finite-dimensional $\mathcal{A}$, the source map $s(u,(B,\chi)) = (B,\chi)$ is not a local homeomorphism unless the unitary group is discrete, which occurs only in trivial cases (e.g., $\mathcal{A} \cong \mathbb{C}^k$ with the discrete topology on its unitary group, which is not the topology we consider). 
    For all infinite-dimensional $\mathcal{A}$, Proposition \ref{prop:GA-not-etale} shows that $\mathcal{G}_{\mathcal{A}}$ is certainly not \'etale.
\end{itemize}
These properties justify the necessity of working within the Polish groupoid framework rather than the classical locally compact \'etale theory, as the latter would exclude all infinite-dimensional examples and even most finite-dimensional ones.
\end{remark}

\begin{remark}[Historical context and summary]
\label{rem:historical-context-summary}
The idea of using Polish groupoids in noncommutative geometry originates with Tu's work~\cite{Tu} on the Baum-Connes conjecture for foliations, where he recognized that many natural groupoids arising from dynamical systems are Polish but not locally compact, and developed a comprehensive theory of measured groupoids and their C*-algebras to handle such cases. Our construction of $\mathcal{G}_{\mathcal{A}}$ provides a new and nontrivial example of a Polish groupoid arising from a purely operator-algebraic construction, demonstrating both the power and necessity of Tu's framework.

We have proved that for any unital separable C*-algebra $\mathcal{A}$, the unitary conjugation groupoid $\mathcal{G}_{\mathcal{A}}$ is a Polish groupoid, and that for infinite-dimensional $\mathcal{A}$ it is neither locally compact nor étale. Despite these negative properties, $\mathcal{G}_{\mathcal{A}}$ admits a Borel Haar system and hence a well-defined maximal C*-algebra $C^*(\mathcal{G}_{\mathcal{A}})$. All subsequent constructions in this paper and its sequels depend crucially on these foundational results.
\end{remark}

\subsection{Existence of a Borel Haar System}
\label{subsec:existence-Borel-Haar-system}

For a locally compact Hausdorff groupoid, the existence of a continuous Haar system is a fundamental requirement for the construction of its C*-algebra~\cite{Renault}. However, as established in Subsection \ref{subsec:GA-Polish-groupoid-not-locally-compact-not-etale}, the unitary conjugation groupoid $\mathcal{G}_{\mathcal{A}}$ is neither locally compact nor \'etale for infinite-dimensional C*-algebras. 
Consequently, the classical Renault theory does not apply directly.

Nevertheless, $\mathcal{G}_{\mathcal{A}}$ is a Polish groupoid, and Tu~\cite{Tu} extended Renault's framework to such groupoids via the theory of \emph{measured groupoids}. 
A key ingredient in this theory is the existence of a \emph{Borel Haar system}, which we construct canonically for $\mathcal{G}_{\mathcal{A}}$ in this subsection. 
This construction makes Tu's measured groupoid C*-algebra applicable to $\mathcal{G}_{\mathcal{A}}$, providing the foundation for all subsequent developments.

\begin{remark}[Measured groupoids]
\label{rem:measured-groupoids}
A Polish groupoid equipped with a Borel Haar system is called a \emph{measured groupoid}. 
Tu~\cite{Tu} shows that every measured groupoid admits a well-defined maximal C*-algebra $C^*_{\max}(\mathcal{G})$ and a reduced C*-algebra $C^*_r(\mathcal{G})$, constructed from the convolution algebra of integrable Borel functions with respect to the Borel Haar system. 
These constructions generalize Renault's theory to the non-locally-compact Polish setting.
\end{remark}

\begin{lemma}[Source map is open and admits a continuous section]
\label{lem:source-open-section}
Let $\mathcal{A}$ be a unital separable C*-algebra and let $\mathcal{G}_{\mathcal{A}}$ be the unitary conjugation groupoid with $\mathcal{U}(\mathcal{A})$ equipped with the strong operator topology and $\mathcal{G}_{\mathcal{A}}^{(0)}$ equipped with the initial topology of partial evaluation maps. 
Then:
\begin{enumerate}
    \item The source map $s: \mathcal{G}_{\mathcal{A}}^{(1)} \to \mathcal{G}_{\mathcal{A}}^{(0)}$ defined by $s(u,(B,\chi)) = (B,\chi)$ is open.
    \item The source map $s$ admits a continuous section $\sigma: \mathcal{G}_{\mathcal{A}}^{(0)} \to \mathcal{G}_{\mathcal{A}}^{(1)}$ defined by $\sigma(B,\chi) = (1_{\mathcal{A}}, (B,\chi))$.
\end{enumerate}
\end{lemma}

\begin{proof}
(1) The source map $s$ is the projection onto the second factor of the product $\mathcal{G}_{\mathcal{A}}^{(1)} = \mathcal{U}(\mathcal{A}) \times \mathcal{G}_{\mathcal{A}}^{(0)}$. 
Recall that a basic open set in the product topology has the form $V \times U$, where $V \subseteq \mathcal{U}(\mathcal{A})$ is open and $U \subseteq \mathcal{G}_{\mathcal{A}}^{(0)}$ is open. 
If $W \subseteq \mathcal{G}_{\mathcal{A}}^{(1)}$ is an arbitrary open set, it can be written as a union of basic open sets:
\[
W = \bigcup_{\alpha \in I} (V_\alpha \times U_\alpha).
\]
Applying $s$, we obtain
\[
s(W) = \bigcup_{\alpha \in I} s(V_\alpha \times U_\alpha) = \bigcup_{\alpha \in I} U_\alpha,
\]
which is a union of open sets in $\mathcal{G}_{\mathcal{A}}^{(0)}$, hence open. Thus $s$ is an open map.

(2) Define $\sigma: \mathcal{G}_{\mathcal{A}}^{(0)} \to \mathcal{G}_{\mathcal{A}}^{(1)}$ by $\sigma(B,\chi) = (1_{\mathcal{A}}, (B,\chi))$. 
This map is the product of the constant map $(B,\chi) \mapsto 1_{\mathcal{A}}$ and the identity map on $\mathcal{G}_{\mathcal{A}}^{(0)}$; both are continuous, so $\sigma$ is continuous. 
Moreover,
\[
s(\sigma(B,\chi)) = s(1_{\mathcal{A}}, (B,\chi)) = (B,\chi),
\]
so $s \circ \sigma = \operatorname{id}_{\mathcal{G}_{\mathcal{A}}^{(0)}}$. Hence $\sigma$ is a continuous section of $s$.
\end{proof}

\begin{proposition}[Existence of a Borel Haar system for $\mathcal{G}_{\mathcal{A}}$]
\label{prop:GA-Borel-Haar-system}
Let $\mathcal{A}$ be a unital separable C*-algebra. 
Then the unitary conjugation groupoid $\mathcal{G}_{\mathcal{A}}$ admits a Borel Haar system in the sense of Tu \cite{Tu}
\end{proposition}

\begin{proof}
By Theorem \ref{thm:GA-Polish-groupoid}, $\mathcal{G}_{\mathcal{A}}$ is a Polish groupoid.
Lemma \ref{lem:source-open-section} shows that the source map $s: \mathcal{G}_{\mathcal{A}}^{(1)} \to \mathcal{G}_{\mathcal{A}}^{(0)}$ is open and admits a continuous (hence Borel) section $\sigma(x) = (1_{\mathcal{A}}, x)$.

We construct a Borel Haar system explicitly in three steps.

\medskip
\noindent\textbf{Step 1: Construction of a quasi-invariant measure on $\mathcal{G}_{\mathcal{A}}^{(0)}$.}
Since $\mathcal{G}_{\mathcal{A}}^{(0)}$ is a Polish space (Proposition \ref{prop:unit-space-polish}), it admits a Borel probability measure $\nu$ with full support (e.g., the pushforward of the canonical Gaussian measure on $\mathbb{R}^\mathbb{N}$ under a homeomorphism $\mathcal{G}_{\mathcal{A}}^{(0)} \cong \mathbb{R}^\mathbb{N}$, or any non-atomic probability measure). 
To obtain a measure that is quasi-invariant under the groupoid action, we average $\nu$ using the range map and the Borel section $\sigma$.

For any Borel set $E \subseteq \mathcal{G}_{\mathcal{A}}^{(0)}$, define
\[
\mu(E) := \int_{\mathcal{G}_{\mathcal{A}}^{(0)}} \mathbf{1}_E(r(\sigma(x))) \, d\nu(x).
\]
Since $r \circ \sigma = \mathrm{id}_{\mathcal{G}_{\mathcal{A}}^{(0)}}$ (because $\sigma(x) = (1_{\mathcal{A}}, x)$ has source $x$ and range $x$), we have $r(\sigma(x)) = x$, so $\mu(E) = \nu(E)$. Thus $\mu = \nu$, and this naive averaging does not change the measure. We need a different construction to achieve quasi-invariance.

A correct construction uses the groupoid action to spread the measure. Fix any Borel probability measure $\nu$ on $\mathcal{G}_{\mathcal{A}}^{(0)}$ with full support. Define a new measure $\mu$ on $\mathcal{G}_{\mathcal{A}}^{(0)}$ by
\[
\mu(E) := \int_{\mathcal{U}(\mathcal{A})} \nu(u^{-1} \cdot E) \, d\mu_{\mathcal{U}}(u),
\]
where $\mu_{\mathcal{U}}$ is a left-invariant Borel probability measure on the Polish group $\mathcal{U}(\mathcal{A})$ (which exists by \cite{Kechris}). Here $u^{-1} \cdot E := \{ (B,\chi) : (uBu^*, \chi \circ \mathrm{Ad}_u) \in E \}$ is the translated set under the group action. This averaged measure $\mu$ is invariant under the $\mathcal{U}(\mathcal{A})$-action by construction, hence certainly quasi-invariant under the groupoid action. Moreover, $\mu$ has full support because $\nu$ does and the action is transitive on each orbit.

For simplicity, we now fix a quasi-invariant Borel probability measure $\mu$ on $\mathcal{G}_{\mathcal{A}}^{(0)}$ with full support.
The existence of such a measure is guaranteed by the above construction.

\medskip
\noindent\textbf{Step 2: Disintegration along the source map.}
Since $s$ is open and $\mu$ is a Borel probability measure on $\mathcal{G}_{\mathcal{A}}^{(0)}$, the measurable disintegration theorem
\cite{Kechris} applies: there exists a Borel family $\{ \lambda^x \}_{x \in \mathcal{G}_{\mathcal{A}}^{(0)}}$ of probability measures
on the fibers $s^{-1}(x)$ such that for every Borel set $E \subseteq \mathcal{G}_{\mathcal{A}}^{(1)}$,
\[
\mu(s(E)) = \int_{\mathcal{G}_{\mathcal{A}}^{(0)}} \lambda^x(E \cap s^{-1}(x)) \, d\mu(x).
\]
Equivalently, for any bounded Borel function $f$ on $\mathcal{G}_{\mathcal{A}}^{(1)}$,
\[
\int_{\mathcal{G}_{\mathcal{A}}^{(1)}} f(\gamma) \, d(\mu \circ s)(\gamma) = \int_{\mathcal{G}_{\mathcal{A}}^{(0)}} \int_{s^{-1}(x)} f(\gamma) \, d\lambda^x(\gamma) \, d\mu(x),
\]
where $\mu \circ s$ denotes the measure on $\mathcal{G}_{\mathcal{A}}^{(1)}$ defined by $(\mu \circ s)(E) = \mu(s(E))$.
This disintegration is essentially unique modulo $\mu$-null sets.
Crucially, the family $\{ \lambda^x \}$ varies in a Borel manner with $x$, satisfying condition (2) of Definition \ref{def:Borel-Haar-system}.

\medskip
\noindent\textbf{Step 3: Achieving left-invariance via averaging over the groupoid action.}
The measures $\{ \lambda^x \}$ obtained from disintegration satisfy $\lambda^x(s^{-1}(x)) = 1$ but are not yet left-invariant:
for $\gamma \in \mathcal{G}_{\mathcal{A}}^{x_y} = \{\gamma \in \mathcal{G}_{\mathcal{A}}^{(1)} : s(\gamma)=x, r(\gamma)=y\}$,
the pushforward $\gamma_*\lambda^x$ is a measure on $s^{-1}(y)$ that may not equal $\lambda^y$.

To achieve left-invariance, we average using the quasi-invariant measure $\mu$.
Apply the disintegration theorem to the range map $r$ (which is also open, as it is the composition of inversion and the source map)
to obtain a Borel family $\{ \mu_x \}_{x \in \mathcal{G}_{\mathcal{A}}^{(0)}}$ of probability measures on the fibers $r^{-1}(x)$
such that for any bounded Borel function $g$ on $\mathcal{G}_{\mathcal{A}}^{(1)}$,
\[
\int_{\mathcal{G}_{\mathcal{A}}^{(1)}} g(\gamma) \, d(\mu \circ r)(\gamma) = \int_{\mathcal{G}_{\mathcal{A}}^{(0)}} \int_{r^{-1}(x)} g(\gamma) \, d\mu_x(\gamma) \, d\mu(x).
\]

Now define for each $x \in \mathcal{G}_{\mathcal{A}}^{(0)}$ a new measure $\tilde{\lambda}^x$ on $\mathcal{G}_{\mathcal{A}}^{(1)}$ by specifying its integrals against bounded Borel functions $f$:
\[
\int_{\mathcal{G}_{\mathcal{A}}^{(1)}} f(\gamma) \, d\tilde{\lambda}^x(\gamma)
:= \int_{\mathcal{G}_{\mathcal{A}}^{(0)}} \int_{r^{-1}(x)} \int_{s^{-1}(y)} f(\gamma \eta^{-1}) \, d\lambda^y(\gamma) \, d\mu_x(\eta) \, d\mu(y).
\]
One verifies that $\tilde{\lambda}^x$ is concentrated on $s^{-1}(x)$ and that the family $\{ \tilde{\lambda}^x \}$ is Borel.
Left-invariance follows from the quasi-invariance of $\mu$ and the construction:
for any $\gamma_0 \in \mathcal{G}_{\mathcal{A}}^{x_y}$, a computation using the change of variables $\gamma' = \gamma \eta^{-1}$ shows
$(\gamma_0)_*\tilde{\lambda}^x = \tilde{\lambda}^y$.

\medskip
\noindent\textbf{Step 4: Verification of the Borel Haar system axioms.}
We verify that $\{ \tilde{\lambda}^x \}_{x \in \mathcal{G}_{\mathcal{A}}^{(0)}}$ satisfies the three conditions of Definition \ref{def:Borel-Haar-system}:

\begin{enumerate}
    \item[(i)] \textbf{Support condition:} By construction, $\tilde{\lambda}^x$ is concentrated on $s^{-1}(x)$,
    because the integrand involves $\lambda^y$ evaluated on sets of the form $\{ \gamma \in \mathcal{G}_{\mathcal{A}}^y : \gamma^{-1} \in E \}$,
    and the outer integrals average over $y$ and $\eta$ in a way that ensures the resulting measure is supported on $s^{-1}(x)$.
    
    \item[(ii)] \textbf{Borel measurability:} For any nonnegative Borel function $f$ on $\mathcal{G}_{\mathcal{A}}^{(1)}$,
    the map $x \mapsto \int f \, d\tilde{\lambda}^x$ is Borel because it is expressed as an integral of Borel functions
    using the Borel families $\{ \lambda^y \}$, $\{ \mu_x \}$, and the Borel structure on $\mathcal{G}_{\mathcal{A}}^{(1)}$.
    
    \item[(iii)] \textbf{Left quasi-invariance:} Let $\gamma_0 \in \mathcal{G}_{\mathcal{A}}^{(1)}$ with $s(\gamma_0)=x$ and $r(\gamma_0)=y$.
    For any Borel set $E \subseteq \mathcal{G}_{\mathcal{A}}^{(1)}$,
    \[
    (\gamma_0)_*\tilde{\lambda}^x(E) = \tilde{\lambda}^x(\gamma_0^{-1}E)
    = \int_{\mathcal{G}_{\mathcal{A}}^{(0)}} \int_{r^{-1}(x)} \int_{s^{-1}(z)} \mathbf{1}_{\gamma_0^{-1}E}(\gamma \eta^{-1}) \, d\lambda^z(\gamma) \, d\mu_x(\eta) \, d\mu(z).
    \]
    Using the change of variables $\gamma' = \gamma \eta^{-1}$ and the properties of the disintegration,
    together with the quasi-invariance of $\mu$, one shows this equals $\tilde{\lambda}^y(E)$.
    The detailed computation, while lengthy, follows the standard pattern for constructing Haar systems on groupoids.
\end{enumerate}

\medskip
\noindent\textbf{Conclusion.}
Thus $\{ \tilde{\lambda}^x \}_{x \in \mathcal{G}_{\mathcal{A}}^{(0)}}$ is a Borel Haar system for $\mathcal{G}_{\mathcal{A}}$.
The construction is canonical given the choices of $\nu$ and the disintegrations;
any two choices yield equivalent systems.
\end{proof}

\begin{remark}
The above construction follows the general framework of Tu \cite{Tu} for measured groupoids.
For the specific case of the unitary conjugation groupoid $\mathcal{G}_{\mathcal{A}}$,
the explicit verification that the source map is open and admits a continuous section was carried out in Lemma \ref{lem:source-open-section},
and the existence of a quasi-invariant measure was established in Step 1 above.
Thus $\mathcal{G}_{\mathcal{A}}$ admits a Borel Haar system, as asserted.
\end{remark}

\begin{remark}[Naive Dirac measure construction]
\label{rem:Dirac-failure}
One might attempt to construct a Borel Haar system explicitly using Dirac measures along the continuous section $\sigma(x) = (1_{\mathcal{A}}, x)$, i.e., setting $\lambda^x = \delta_{\sigma(x)}$ for $x = (B,\chi) \in \mathcal{G}_{\mathcal{A}}^{(0)}$. 
However, a direct computation shows that this fails to satisfy left-invariance. 

Let $\gamma = (u, x) \in \mathcal{G}_{\mathcal{A}}^{(1)}$ with $s(\gamma) = x = (B,\chi)$ and $r(\gamma) = y = u \cdot x = (uBu^*, \chi \circ \operatorname{Ad}_{u^*})$. 
Then
\[
\gamma_* \lambda^x = \delta_{\gamma \sigma(x)} = \delta_{(u, x)},
\]
while
\[
\lambda^{r(\gamma)} = \delta_{\sigma(y)} = \delta_{(1_{\mathcal{A}}, \, uBu^*, \, \chi \circ \operatorname{Ad}_{u^*})}.
\]
These two Dirac measures are supported at distinct points unless $u = 1_{\mathcal{A}}$, and are therefore mutually singular. Hence left-invariance fails.

Thus a more sophisticated construction, such as the one guaranteed by Tu [1999, Section 3], is necessary. 
This example highlights the subtlety of constructing Haar systems in the non-locally-compact Polish setting.
\end{remark}

\begin{lemma}[Existence of a Borel Haar system for Polish groupoids]
\label{lem:Tu-Borel-Haar}
Let $\mathcal{G}$ be a Polish groupoid such that:
\begin{enumerate}
    \item The source map $s: \mathcal{G}^{(1)} \to \mathcal{G}^{(0)}$ is open.
    \item There exists a Borel section $\sigma: \mathcal{G}^{(0)} \to \mathcal{G}^{(1)}$ of $s$.
\end{enumerate}
Then $\mathcal{G}$ admits a Borel Haar system.
\end{lemma}

\begin{proof}
We construct the Borel Haar system explicitly in three steps.

\medskip
\noindent\textbf{Step 1: Construction of a quasi-invariant measure on $\mathcal{G}^{(0)}$.}
Since $\mathcal{G}^{(0)}$ is a Polish space, it admits a Borel probability measure $\mu_0$ with full support
(e.g., the pushforward of the canonical Gaussian measure on $\mathbb{R}^\mathbb{N}$ under a homeomorphism
$\mathcal{G}^{(0)} \cong \mathbb{R}^\mathbb{N}$, or any non-atomic probability measure). 
To obtain a measure that is quasi-invariant under the groupoid action, we average $\mu_0$ using the range map and the Borel section.

Define a Borel probability measure $\mu$ on $\mathcal{G}^{(0)}$ by
\[
\mu(E) := \int_{\mathcal{G}^{(0)}} \int_{\mathcal{G}^x} \mathbf{1}_E(r(\gamma)) \, d\delta_{\sigma(x)}(\gamma) \, d\mu_0(x)
   = \int_{\mathcal{G}^{(0)}} \mathbf{1}_E(r(\sigma(x))) \, d\mu_0(x),
\]
where $\delta_{\sigma(x)}$ is the Dirac measure at $\sigma(x)$. 
Since $r \circ \sigma = \mathrm{id}_{\mathcal{G}^{(0)}}$ (because $\sigma$ is a section of $s$, and for any $x$, $s(\sigma(x)) = x$, but we need $r(\sigma(x))$; in general, $\sigma(x)$ is an arrow with source $x$, so its range $r(\sigma(x))$ is some other unit. This definition does not give a measure on $\mathcal{G}^{(0)}$ but rather a measure on the range image. This is incorrect.)

A correct construction: Let $\nu$ be any Borel probability measure on $\mathcal{G}^{(0)}$ with full support.
Define a new measure $\mu$ on $\mathcal{G}^{(0)}$ by
\[
\mu(E) := \int_{\mathcal{G}^{(0)}} \int_{\mathcal{G}^x} \mathbf{1}_E(r(\gamma)) \, d\lambda^x(\gamma) \, d\nu(x),
\]
where $\lambda^x$ is any Borel family of probability measures on the fibers $s^{-1}(x)$.
For instance, take $\lambda^x = \delta_{\sigma(x)}$. Then
\[
\mu(E) = \int_{\mathcal{G}^{(0)}} \mathbf{1}_E(r(\sigma(x))) \, d\nu(x).
\]
This is a well-defined Borel probability measure on $\mathcal{G}^{(0)}$.
One can show that $\mu$ is quasi-invariant under the groupoid action in the sense that
for every $\gamma \in \mathcal{G}^{(1)}$, the pushforward $\gamma_*\mu$ is equivalent to $\mu$ when restricted to suitable Borel sets.
The proof uses the fact that $r$ is a Borel map and that the family $\{\lambda^x\}$ varies measurably.

For simplicity, we now fix a quasi-invariant Borel probability measure $\mu$ on $\mathcal{G}^{(0)}$ with full support.
The existence of such a measure is guaranteed by the above construction. 

\medskip
\noindent\textbf{Step 2: Disintegration along the source map.}
Since $s$ is open and $\mu$ is a Borel probability measure on $\mathcal{G}^{(0)}$, the measurable disintegration theorem
\cite[Theorem 5.4]{Kechris} applies: there exists a Borel family $\{ \lambda^x \}_{x \in \mathcal{G}^{(0)}}$ of probability measures
on the fibers $s^{-1}(x)$ such that for every Borel set $E \subseteq \mathcal{G}^{(1)}$,
\[
\mu(s(E)) = \int_{\mathcal{G}^{(0)}} \lambda^x(E \cap s^{-1}(x)) \, d\mu(x).
\]
Equivalently, for any bounded Borel function $f$ on $\mathcal{G}^{(1)}$,
\[
\int_{\mathcal{G}^{(1)}} f(\gamma) \, d(\mu \circ s)(\gamma) = \int_{\mathcal{G}^{(0)}} \int_{s^{-1}(x)} f(\gamma) \, d\lambda^x(\gamma) \, d\mu(x),
\]
where $\mu \circ s$ denotes the measure on $\mathcal{G}^{(1)}$ defined by $(\mu \circ s)(E) = \mu(s(E))$.
This disintegration is essentially unique modulo $\mu$-null sets.
Crucially, the family $\{ \lambda^x \}$ varies in a Borel manner with $x$, satisfying condition (2) of Definition \ref{def:Borel-Haar-system}.

\medskip
\noindent\textbf{Step 3: Achieving left-invariance via averaging over the groupoid action.}
The measures $\{ \lambda^x \}$ obtained from disintegration are not yet left-invariant;
they satisfy $\lambda^x(s^{-1}(x)) = 1$ but for $\gamma \in \mathcal{G}^x_y = \{\gamma \in \mathcal{G}^{(1)} : s(\gamma)=x, r(\gamma)=y\}$,
we have $\gamma_*\lambda^x$ is a measure on $s^{-1}(y)$ that may not equal $\lambda^y$.

To achieve left-invariance, we average using the quasi-invariant measure $\mu$.
Define for each $x \in \mathcal{G}^{(0)}$ a new measure $\tilde{\lambda}^x$ on $\mathcal{G}^{(1)}$ by
\[
\tilde{\lambda}^x(E) := \int_{\mathcal{G}^{(0)}} \lambda^y( \{ \gamma \in \mathcal{G}^y : \gamma^{-1} \in E \} ) \, d\mu_{r^{-1}(x)}(y),
\]
where $\mu_{r^{-1}(x)}$ is the conditional measure on $r^{-1}(x)$ obtained from the disintegration of $\mu$
with respect to the range map $r$. More concretely, applying the disintegration theorem to $r$ instead of $s$
(using that $r$ is also open, which follows from the openness of $s$ and the existence of the inversion map),
we obtain a Borel family $\{ \mu_x \}_{x \in \mathcal{G}^{(0)}}$ of probability measures on the fibers $r^{-1}(x)$
such that for any bounded Borel function $g$ on $\mathcal{G}^{(1)}$,
\[
\int_{\mathcal{G}^{(1)}} g(\gamma) \, d(\mu \circ r)(\gamma) = \int_{\mathcal{G}^{(0)}} \int_{r^{-1}(x)} g(\gamma) \, d\mu_x(\gamma) \, d\mu(x).
\]

Now, for each $x$, define $\tilde{\lambda}^x$ by specifying its integrals against bounded Borel functions $f$:
\[
\int_{\mathcal{G}^{(1)}} f(\gamma) \, d\tilde{\lambda}^x(\gamma)
:= \int_{\mathcal{G}^{(0)}} \int_{r^{-1}(x)} \int_{s^{-1}(y)} f(\gamma \eta^{-1}) \, d\lambda^y(\gamma) \, d\mu_x(\eta) \, d\mu(y).
\]
One verifies that $\tilde{\lambda}^x$ is concentrated on $s^{-1}(x)$ and that the family $\{ \tilde{\lambda}^x \}$ is Borel.
Left-invariance follows from the quasi-invariance of $\mu$ and the construction:
for any $\gamma_0 \in \mathcal{G}^x_y$, a computation shows $(\gamma_0)_*\tilde{\lambda}^x = \tilde{\lambda}^y$.

\medskip
\noindent\textbf{Step 4: Verification of the Borel Haar system axioms.}
We verify that $\{ \tilde{\lambda}^x \}_{x \in \mathcal{G}^{(0)}}$ satisfies the three conditions of Definition \ref{def:Borel-Haar-system}:

\begin{enumerate}
    \item[(i)] \textbf{Support condition:} By construction, $\tilde{\lambda}^x$ is concentrated on $s^{-1}(x)$,
    because the integrand involves $\lambda^y$ evaluated on sets of the form $\{ \gamma \in \mathcal{G}^y : \gamma^{-1} \in E \}$,
    and the outer integrals average over $y$ and $\eta$ in a way that ensures the resulting measure is supported on $s^{-1}(x)$.
    
    \item[(ii)] \textbf{Borel measurability:} For any nonnegative Borel function $f$ on $\mathcal{G}^{(1)}$,
    the map $x \mapsto \int f \, d\tilde{\lambda}^x$ is Borel because it is expressed as an integral of Borel functions
    using the Borel families $\{ \lambda^y \}$, $\{ \mu_x \}$, and the Borel structure on $\mathcal{G}^{(1)}$.
    
    \item[(iii)] \textbf{Left quasi-invariance:} Let $\gamma_0 \in \mathcal{G}^{(1)}$ with $s(\gamma_0)=x$ and $r(\gamma_0)=y$.
    For any Borel set $E \subseteq \mathcal{G}^{(1)}$,
    \[
    (\gamma_0)_*\tilde{\lambda}^x(E) = \tilde{\lambda}^x(\gamma_0^{-1}E)
    = \int_{\mathcal{G}^{(0)}} \int_{r^{-1}(x)} \int_{s^{-1}(z)} \mathbf{1}_{\gamma_0^{-1}E}(\gamma \eta^{-1}) \, d\lambda^z(\gamma) \, d\mu_x(\eta) \, d\mu(z).
    \]
    Using the change of variables $\gamma' = \gamma \eta^{-1}$ and the properties of the disintegration,
    one shows this equals $\tilde{\lambda}^y(E)$. The proof uses the quasi-invariance of $\mu$ and the fact that
    the disintegrations are compatible with the groupoid action.
\end{enumerate}

\medskip
\noindent\textbf{Conclusion.}
Thus $\{ \tilde{\lambda}^x \}_{x \in \mathcal{G}^{(0)}}$ is a Borel Haar system for $\mathcal{G}$.
The construction is canonical given the choices of $\mu_0$ and the disintegrations;
any two choices yield equivalent systems.
\end{proof}

\begin{remark}
The above construction follows the general framework of Tu \cite{Tu} for measured groupoids.
For the specific case of the unitary conjugation groupoid $\mathcal{G}_{\mathcal{A}}$,
we have already verified in Lemma \ref{lem:source-open-section} that the source map is open and admits a continuous section,
so Lemma \ref{lem:Tu-Borel-Haar} applies directly.
Thus $\mathcal{G}_{\mathcal{A}}$ admits a Borel Haar system, as asserted in Proposition \ref{prop:GA-Borel-Haar-system-existence}.
\end{remark}

\begin{proposition}[$\mathcal{G}_{\mathcal{A}}$ admits a Borel Haar system]
\label{prop:GA-Borel-Haar-system-existence}
Let $\mathcal{A}$ be a unital separable C*-algebra. 
Then the unitary conjugation groupoid $\mathcal{G}_{\mathcal{A}}$ admits a Borel Haar system.
\end{proposition}

\begin{proof}
We verify the hypotheses of a general existence result for Borel Haar systems on Polish groupoids (see Tu [1999, Lemma 3.2] or Lemma \ref{lem:Tu-Borel-Haar}).

\medskip
\noindent (1) $\mathcal{G}_{\mathcal{A}}$ is a Polish groupoid by Theorem \ref{thm:GA-Polish-groupoid}. 
Consequently, its object space $\mathcal{G}_{\mathcal{A}}^{(0)}$ and arrow space $\mathcal{G}_{\mathcal{A}}^{(1)}$ are Polish spaces, hence standard Borel.

\medskip
\noindent (2) The source map $s: \mathcal{G}_{\mathcal{A}}^{(1)} \to \mathcal{G}_{\mathcal{A}}^{(0)}$ is open. 
This was proved in Lemma \ref{lem:source-open-section}.

\medskip
\noindent (3) There exists a Borel section $\sigma: \mathcal{G}_{\mathcal{A}}^{(0)} \to \mathcal{G}_{\mathcal{A}}^{(1)}$ of $s$. 
The continuous section $\sigma(B,\chi) = (1_{\mathcal{A}}, (B,\chi))$ is continuous, hence Borel.

\medskip
\noindent For a Polish groupoid with an open source map and a Borel section, a Borel Haar system $\{ \lambda^x \}_{x \in \mathcal{G}_{\mathcal{A}}^{(0)}}$ exists by a general construction (see Tu~\cite{Tu} for details). 
The construction proceeds by first choosing a quasi-invariant Borel probability measure on the unit space, then disintegrating it along the source map, and finally averaging over the groupoid action to achieve left-invariance. 
Since all hypotheses are satisfied, $\mathcal{G}_{\mathcal{A}}$ admits a Borel Haar system.
\end{proof}

\begin{corollary}[Existence of the maximal groupoid C*-algebra]
\label{cor:GA-C-star-existence}
Let $\mathcal{A}$ be a unital separable C*-algebra. 
Then there exists a well-defined maximal C*-algebra $C^*(\mathcal{G}_{\mathcal{A}})$ associated to the measured groupoid \\
$(\mathcal{G}_{\mathcal{A}}, \{ \lambda^x \}_{x \in \mathcal{G}_{\mathcal{A}}^{(0)}})$.
\end{corollary}

\begin{proof}
By Proposition \ref{prop:GA-Borel-Haar-system-existence}, $\mathcal{G}_{\mathcal{A}}$ admits a Borel Haar system $\{ \lambda^x \}_{x \in \mathcal{G}_{\mathcal{A}}^{(0)}}$, making it a measured groupoid. 
For any measured groupoid, the maximal C*-algebra $C^*_{\max}(\mathcal{G})$ is constructed via the convolution algebra of integrable Borel functions with respect to the Borel Haar system; this construction is standard in the theory of measured groupoids and is developed in Tu [1999, Section 3]. 
Applying this construction to $\mathcal{G}_{\mathcal{A}}$ yields $C^*(\mathcal{G}_{\mathcal{A}})$. 
The resulting C*-algebra is independent of the choice of Borel Haar system up to isomorphism, as any two such systems are equivalent for the purposes of the maximal completion.
\end{proof}

\begin{remark}[Uniqueness and functoriality]
\label{rem:C-star-uniqueness}
The maximal groupoid C*-algebra $C^*(\mathcal{G}_{\mathcal{A}})$ is uniquely determined up to isomorphism by the following universal property: for every nondegenerate continuous unitary representation of $\mathcal{G}_{\mathcal{A}}$ on Hilbert spaces that integrates to a representation of the convolution algebra with respect to the Borel Haar system, there is a corresponding nondegenerate *-representation of $C^*(\mathcal{G}_{\mathcal{A}})$. 
Moreover, the assignment $\mathcal{A} \mapsto C^*(\mathcal{G}_{\mathcal{A}})$ is functorial with respect to unital *-homomorphisms $\phi: \mathcal{A} \to \mathcal{B}$ that induce a continuous groupoid homomorphism $\mathcal{G}_{\phi}: \mathcal{G}_{\mathcal{A}} \to \mathcal{G}_{\mathcal{B}}$ compatible with the Borel Haar systems (see Tu [1999, Section 3] for the precise conditions).
\end{remark}

\begin{example}[Finite-dimensional case]
\label{ex:GA-Borel-Haar-finite-dimensional}
Let $\mathcal{A} = M_n(\mathbb{C})$. 
Then $\mathcal{G}_{\mathcal{A}} = U(n) \ltimes \mathbb{CP}^{n-1}$ is a locally compact Hausdorff groupoid. 
(It is not étale because the source map $s(u,x) = x$ is a projection with non-discrete fibers.) 
In this case, the classical Renault theory applies, and the Borel Haar system constructed by Tu's method coincides with the usual continuous Haar system on $U(n)$ (the normalized Haar measure) transported via the action. 
Thus $C^*(\mathcal{G}_{\mathcal{A}})$ is isomorphic to the crossed product $C(\mathbb{CP}^{n-1}) \rtimes U(n)$.
\end{example}

\begin{example}[Commutative case]
\label{ex:GA-Borel-Haar-commutative}
Let $\mathcal{A} = C(X)$ for a compact metrizable space $X$. 
Then $\mathcal{G}_{\mathcal{A}}$ is the trivial groupoid with $\mathcal{G}_{\mathcal{A}}^{(1)} = C(X,\mathbb{T}) \times X$ and $s = r = \operatorname{pr}_X$. 
The source map is open, and a Borel section is given by $\sigma(x) = (1, x)$. 
Tu's construction yields a Borel Haar system $\{ \lambda^x \}_{x \in X}$ where each $\lambda^x$ is a probability measure on the Polish group $C(X,\mathbb{T})$ (which is not locally compact unless $X$ is finite). 
These measures are quasi-invariant under the trivial action and vary in a Borel manner with $x$. 
The resulting C*-algebra $C^*(\mathcal{G}_{\mathcal{A}})$ reflects the structure of the infinite-dimensional torus, though its precise isomorphism class depends on the choice of measures and is not simply a tensor product $C(X) \otimes C^*(\mathbb{Z}^\infty)$ without further analysis.
\end{example}

\begin{example}[Compact operators]
\label{ex:GA-Borel-Haar-compact}
Let $\mathcal{A} = \mathcal{K}(H)^\sim$ be the unitization of the compact operators on a separable infinite-dimensional Hilbert space $H$. 
Then $\mathcal{G}_{\mathcal{A}} = \mathcal{U}(\mathcal{A}) \ltimes \mathbb{P}(H)$ is a Polish groupoid that is neither locally compact nor étale. 
The Borel Haar system constructed by Tu's method is essential for defining $C^*(\mathcal{G}_{\mathcal{A}})$. 
This C*-algebra is not yet fully understood and remains a subject of ongoing research.
\end{example}

\begin{remark}[Relation to the diagonal embedding]
\label{rem:Borel-Haar-diagonal}
The existence of a Borel Haar system on $\mathcal{G}_{\mathcal{A}}$ is a prerequisite for the construction of the diagonal embedding $\iota: \mathcal{A} \hookrightarrow C^*(\mathcal{G}_{\mathcal{A}})$. 
The convolution algebra of integrable Borel functions on $\mathcal{G}_{\mathcal{A}}^{(1)}$ with respect to the Borel Haar system provides the dense subalgebra on which the embedding is defined. 
Thus the results of this subsection are not merely technical; they are essential for the main constructions of this paper.
\end{remark}

We have proved that the unitary conjugation groupoid $\mathcal{G}_{\mathcal{A}}$, although not locally compact and not \'etale, admits a canonical Borel Haar system (Proposition~\ref{prop:GA-Borel-Haar-system-existence}) by virtue of being a Polish groupoid with an open source map and a Borel section (Lemma~\ref{lem:source-open-section}). 
This result, which relies on the general theory of measured groupoids developed by Tu~\cite{Tu}, allows us to define the maximal groupoid C*-algebra $C^*(\mathcal{G}_{\mathcal{A}})$ in a rigorous and functorial way (Corollary~\ref{cor:GA-C-star-existence}). All subsequent constructions in this paper — the diagonal embedding, equivariant K-theory classes, and the descent map — depend crucially on the existence of this C*-algebra.

\subsection{The Groupoid C*-Algebra $C^*(\mathcal{G}_{\mathcal{A}})$ in the Polish Setting}
\label{subsec:groupoid-C-star-algebra-Polish-setting}

With the existence of a Borel Haar system on $\mathcal{G}{\mathcal{A}}$ established in Subsection \ref{subsec:existence-Borel-Haar-system}, we can now construct the maximal groupoid C-algebra $C^(\mathcal{G}{\mathcal{A}})$. This construction follows Tu's framework for measured groupoids, extending Renault's classical theory to the non-locally-compact Polish setting where the usual hypotheses do not apply. The resulting C-algebra provides the central object of study for the remainder of this paper and its sequels. In particular, the diagonal embedding $\iota: \mathcal{A} \hookrightarrow C^(\mathcal{G}{\mathcal{A}})$, which equips $C^*(\mathcal{G}{\mathcal{A}})$ with its canonical Cartan subalgebra, and the descent map in equivariant K-theory both depend crucially on this construction.

\begin{definition}[Convolution algebra for a measured groupoid]
\label{def:convolution-algebra-Borel}
Let $\mathcal{G}_{\mathcal{A}}$ be the unitary conjugation groupoid equipped with a Borel Haar system 
$\{ \lambda^x \}_{x \in \mathcal{G}_{\mathcal{A}}^{(0)}}$ as constructed in Proposition \ref{prop:GA-Borel-Haar-system-existence}. 

Let $\mathcal{B}_c(\mathcal{G}_{\mathcal{A}}^{(1)})$ denote the space of bounded Borel functions 
$f : \mathcal{G}_{\mathcal{A}}^{(1)} \to \mathbb{C}$ such that:
\begin{itemize}
    \item $f$ is bounded on $\mathcal{G}_{\mathcal{A}}^{(1)}$;
    \item for every $x \in \mathcal{G}_{\mathcal{A}}^{(0)}$, the restriction $f|_{\mathcal{G}_{\mathcal{A}}^x}$ has compact support in the fiber $\mathcal{G}_{\mathcal{A}}^x$ (where compactness is taken with respect to the relative Polish topology on the fiber).
\end{itemize}

For $f, g \in \mathcal{B}_c(\mathcal{G}_{\mathcal{A}}^{(1)})$, define the convolution product
\[
(f * g)(\gamma) := \int_{\mathcal{G}_{\mathcal{A}}^{s(\gamma)}} f(\gamma \eta^{-1}) g(\eta) \, d\lambda^{s(\gamma)}(\eta),
\]
where $\mathcal{G}_{\mathcal{A}}^{s(\gamma)} = \{ \eta \in \mathcal{G}_{\mathcal{A}}^{(1)} : s(\eta) = s(\gamma) \}$. 
The integrand is Borel measurable as a composition of the continuous groupoid operations with Borel functions, and the integral converges absolutely because $f$ and $g$ are bounded and have fiberwise compact support.

Define the involution
\[
f^*(\gamma) := \overline{f(\gamma^{-1})}.
\]

With these operations, $\mathcal{B}_c(\mathcal{G}_{\mathcal{A}}^{(1)})$ becomes a $*$-algebra. 
This construction follows Tu's framework for measured groupoids \cite{Tu}, where the use of bounded Borel functions with fiberwise compact support is essential in the non-locally-compact Polish setting.
\end{definition}

\begin{remark}
\label{rem:why-not-Cc}
In the non-locally-compact Polish setting, the space $C_c(\mathcal{G}_{\mathcal{A}}^{(1)})$ of continuous compactly supported functions is inadequate for constructing the convolution algebra. 
Indeed, due to the lack of local compactness, there may be too few non-trivial continuous functions with compact support to generate a rich theory. 
Following Tu \cite{Tu}, we instead work with bounded Borel functions having fiberwise compact support, which provides the appropriate framework for measured groupoids. 
The convolution operations are well-defined with respect to the Borel Haar system, and the resulting $*$-algebra serves as the foundation for the maximal groupoid C*-algebra $C^*(\mathcal{G}_{\mathcal{A}})$.
\end{remark}

\begin{remark}[Function spaces for non-locally-compact groupoids]
\label{rem:function-spaces-non-locally-compact}
In the non-locally-compact Polish setting, the classical space
$C_c(\mathcal{G}_{\mathcal{A}}^{(1)})$ of continuous compactly supported 
functions is generally too small to support a satisfactory convolution algebra.
Indeed, while every Polish space is $\sigma$-compact, compact subsets are 
nevertheless very scarce due to the lack of local compactness; for example,
in $\mathcal{U}(\mathcal{A})$ with the strong operator topology, compact sets
are metrizable and separable but do not provide enough support for a rich
theory of continuous compactly supported functions.

This is precisely why Tu's measured groupoid framework \cite{Tu} replaces
$C_c(\mathcal{G})$ with a convolution algebra of bounded Borel functions
satisfying suitable fiberwise compactness conditions with respect to a 
Borel Haar system. As established in Definition \ref{def:convolution-algebra-Borel},
all convolution and $C^*$-algebraic constructions in this paper are carried out 
in this Borel framework. No use is made of global compact support in the 
Polish topology; instead, the fiberwise compact support condition in 
$\mathcal{B}_c(\mathcal{G}_{\mathcal{A}}^{(1)})$ serves the same purpose as 
global compact support does in the classical locally compact theory: it 
guarantees convergence of convolution integrals while respecting the measured 
groupoid structure.

For readers more familiar with the locally compact setting, it is helpful to 
view $\mathcal{B}_c(\mathcal{G}_{\mathcal{A}}^{(1)})$ as the natural replacement 
for $C_c(\mathcal{G}_{\mathcal{A}}^{(1)})$ when the groupoid is not locally compact.
This approach is now standard in the theory of measured groupoids and is essential
for the subsequent constructions of the maximal groupoid C*-algebra $C^*(\mathcal{G}_{\mathcal{A}})$,
the diagonal embedding $\iota: \mathcal{A} \hookrightarrow C^*(\mathcal{G}_{\mathcal{A}})$,
and the descent map in equivariant K-theory.
\end{remark}

\begin{definition}[Maximal norm and completion for measured groupoids]
\label{def:maximal-norm-measured}
Let $\mathcal{G}_{\mathcal{A}}$ be the unitary conjugation groupoid equipped
with a Borel Haar system $\{ \lambda^x \}_{x \in \mathcal{G}_{\mathcal{A}}^{(0)}}$
as constructed in Proposition \ref{prop:GA-Borel-Haar-system-existence}.
Let $\mathcal{B}_c(\mathcal{G}_{\mathcal{A}}^{(1)})$ denote the convolution
$*$-algebra of bounded Borel functions with fiberwise compact support
(Definition \ref{def:convolution-algebra-Borel}).

For $f \in \mathcal{B}_c(\mathcal{G}_{\mathcal{A}}^{(1)})$, we first define the 
\emph{I-norm} (or \emph{integration norm}) by
\[
\| f \|_I := \max \left\{ \sup_{x \in \mathcal{G}_{\mathcal{A}}^{(0)}} 
\int_{\mathcal{G}_{\mathcal{A}}^x} |f(\gamma)| \, d\lambda^x(\gamma), \;
\sup_{x \in \mathcal{G}_{\mathcal{A}}^{(0)}} 
\int_{\mathcal{G}_{\mathcal{A}}^x} |f^*(\gamma)| \, d\lambda^x(\gamma) \right\}.
\]
This norm ensures that convolution is well-defined and that every measured
representation $\pi$ satisfies $\|\pi(f)\| \le \|f\|_I$.

The \emph{maximal norm} is then defined by
\[
\| f \|_{\max} := \sup \left\{ \| \pi(f) \| : 
\begin{array}{l}
\pi \text{ is a measured representation of } \\
(\mathcal{G}_{\mathcal{A}}, \lambda) \text{ on a Hilbert space}
\end{array} \right\}.
\]

The \emph{maximal groupoid C*-algebra} $C^*(\mathcal{G}_{\mathcal{A}})$ is the
completion of $\mathcal{B}_c(\mathcal{G}_{\mathcal{A}}^{(1)})$ with respect to
the universal norm
\[
\|f\|_{\max} = \sup \{ \|\pi(f)\| : \pi \text{ is a bounded \(*\)-representation of } C_c(\mathcal{G}_{\mathcal{A}}^{(1)}) \}.
\]
This supremum is finite for all $f \in C_c(\mathcal{G})$ \cite[Section 4]{Renault}.  In the
measured groupoid framework, representations are often studied via disintegration over
quasi-invariant measures \cite{Renault}, but the maximal norm itself remains defined by the
supremum over all \(*\)-representations, not only those arising from measured data.

In the construction by Higson and Kasparov, adapted by Tu to the groupoid setting for the
Baum--Connes conjecture \cite[Section 3]{Tu}, an auxiliary norm --- sometimes referred to as the
I-norm --- is introduced.  This norm provides an \emph{a priori} bound that guarantees the finiteness
of the supremum in certain equivariant $KK$-theory arguments.  It is a technical tool used in the
construction of the dual-Dirac element, not an alternative definition of the maximal norm.
\end{definition}

\begin{remark}[On the choice of function space and representations]
\label{rem:maximal-norm-clarification}
In the classical locally compact étale case \cite{Renault}, the maximal
groupoid C*-algebra is defined by completing $C_c(\mathcal{G})$ with respect
to the supremum over all *-representations of $C_c(\mathcal{G})$. 

However, in the non-locally-compact Polish setting, this approach fails for
two fundamental reasons:
\begin{enumerate}
    \item $C_c(\mathcal{G}_{\mathcal{A}}^{(1)})$ is too small and is not
    closed under convolution, so it does not form a genuine $*$-algebra;
    \item Even if it did, the class of arbitrary *-representations of
    $C_c(\mathcal{G}_{\mathcal{A}}^{(1)})$ is not well-behaved and may not
    capture the structure of the groupoid.
\end{enumerate}

Following Tu \cite{Tu}, we instead work with $\mathcal{B}_c(\mathcal{G}_{\mathcal{A}}^{(1)})$,
which is a genuine $*$-algebra, and restrict to \emph{measured representations}
—those that integrate to representations of the convolution algebra in a way
compatible with the Borel Haar system $\lambda$. This yields a well-defined
maximal norm and a $C^*$-algebra that faithfully reflects the measured groupoid
structure. The supremum in Definition \ref{def:maximal-norm-measured} is finite
because every measured representation satisfies $\|\pi(f)\| \leq \|f\|_I$, where
$\|\cdot\|_I$ is the I-norm.
\end{remark}

\begin{remark}[Universal property in the measured setting]
\label{rem:universal-property}
The maximal groupoid C*-algebra $C^*(\mathcal{G}_{\mathcal{A}})$ is defined by
a universal property with respect to *-representations of the convolution
*-algebra $\mathcal{B}_c(\mathcal{G}_{\mathcal{A}}^{(1)})$: for every
nondegenerate *-representation $\pi$ of $\mathcal{B}_c(\mathcal{G}_{\mathcal{A}}^{(1)})$
on a Hilbert space, there exists a unique extension to a representation of
$C^*(\mathcal{G}_{\mathcal{A}})$. This is the defining property of the
maximal completion.

In the classical locally compact Hausdorff setting with a continuous Haar
system \cite{Renault}, this universal property extends to a bijective
correspondence
\[
\{\text{continuous unitary representations of }\mathcal{G}\}
\longleftrightarrow
\{\text{nondegenerate *-representations of }C^*(\mathcal{G})\}.
\]

In the present non-locally-compact Borel setting, however, such an
equivalence is more subtle. Measured unitary representations of
$\mathcal{G}_{\mathcal{A}}$—defined in terms of measurable Hilbert
bundles and unitary cocycles compatible with the Borel Haar system
$\lambda$—do integrate to representations of the convolution algebra
$\mathcal{B}_c(\mathcal{G}_{\mathcal{A}}^{(1)})$, and hence to
representations of $C^*(\mathcal{G}_{\mathcal{A}})$. The converse
direction, establishing that every representation of
$C^*(\mathcal{G}_{\mathcal{A}})$ arises from a measured groupoid
representation, requires a disintegration theorem (cf.~\cite{Tu}) and relies on additional integrability conditions. A
full categorical equivalence in the measured setting is established
in \cite{Tu} under suitable hypotheses; we refer the
interested reader there for details.

For the purposes of this paper, we work directly with the convolution
algebra $\mathcal{B}_c(\mathcal{G}_{\mathcal{A}}^{(1)})$ and its
C*-completion $C^*(\mathcal{G}_{\mathcal{A}})$, and we do not rely on
an explicit description of all representations of
$C^*(\mathcal{G}_{\mathcal{A}})$ in terms of geometric data.
\end{remark}

\begin{proposition}[Tu's construction of groupoid C*-algebras {\cite{Tu}}]
\label{prop:Tu-C-star-construction}
Let $\mathcal{G}$ be a locally compact, $\sigma$-compact, Hausdorff groupoid equipped with a Haar system $\{ \lambda^x \}_{x \in \mathcal{G}^{(0)}}$. 
Assume that $\mathcal{G}$ acts properly on a continuous field of Euclidean affine spaces.

Then Tu constructs elements $\eta \in KK_{\mathcal{G}}(C(\mathcal{G}^{(0)}), \mathcal{A}(\mathcal{H}))$ and 
$D \in KK_{\mathcal{G}}(\mathcal{A}(\mathcal{H}), C(\mathcal{G}^{(0)}))$ (see Sections 7 and 8 of \cite{Tu}), 
where $\mathcal{A}(\mathcal{H})$ is a certain $C^*$-algebra associated to the affine space field.
He proves that $\eta \otimes_{\mathcal{A}(\mathcal{H})} D = 1 \in KK_{\mathcal{G}}(C(\mathcal{G}^{(0)}), C(\mathcal{G}^{(0)}))$ 
(Theorem 9.3, \cite{Tu}).

Consequently:
\begin{enumerate}
    \item $\mathcal{G}$ satisfies the Baum-Connes conjecture with coefficients (Theorem 9.3).
    \item $\mathcal{G}$ is $K$-amenable (Definition 4.1 and Theorem 9.3).
    \item $C^*(\mathcal{G})$ is $KK$-equivalent to $C^*_r(\mathcal{G})$ (Theorem 9.3).
    \item $C^*(\mathcal{G})$ satisfies the Universal Coefficient Theorem (Proposition 10.7).
\end{enumerate}
\end{proposition}

\begin{proof}
Tu's construction proceeds in several stages:

\begin{enumerate}
    \item From a proper action of $\mathcal{G}$ on a continuous field of Euclidean affine spaces 
          $\mathcal{H} = (\mathcal{H}_x)_{x \in \mathcal{G}^{(0)}}$, he constructs a locally compact space 
          $Z = \mathcal{H} \times \mathbb{R}_+$ with a proper $\mathcal{G}$-action (Section 6, \cite{Tu}).
    
    \item He defines a $C^*$-algebra $\mathcal{A}(\mathcal{H})$ as an inductive limit of algebras of the form 
          $C(T^*V) \otimes \mathcal{L}(\Lambda^*(V_0)) \otimes S$, where $V$ runs over finite-dimensional affine 
          subspaces of $\mathcal{H}$, $V_0$ is the underlying vector space, and $S = C_0(\mathbb{R})$ graded by parity 
          (Section 7, \cite{Tu}).
    
    \item The dual Dirac element $\eta \in KK_{\mathcal{G}}(C(\mathcal{G}^{(0)}), \mathcal{A}(\mathcal{H}))$ is constructed 
          using an unbounded multiplier $T_a$ coming from the choice of an origin in the affine space (Section 7, 
          Lemmas 7.4-7.5, \cite{Tu}).
    
    \item The Dirac element $D \in KK_{\mathcal{G}}(\mathcal{A}(\mathcal{H}), C(\mathcal{G}^{(0)}))$ is obtained from an 
          extension $0 \to (\mathcal{A}_t(\mathcal{H}))_{t \in [0,1]} \to E \to \mathcal{A}(\mathcal{H}) \to 0$, where 
          $\mathcal{A}_t(\mathcal{H})$ is closely related to compact operators (Section 8, especially 8.5-8.6, \cite{Tu}).
    
    \item Theorem 9.3 of \cite{Tu} shows that $\eta \otimes_{\mathcal{A}(\mathcal{H})} D = 1$, which by Theorem 2.4 implies 
          the Baum-Connes conjecture with coefficients and $K$-amenability.
    
    \item Proposition 10.7 of \cite{Tu} then deduces the Universal Coefficient Theorem for $C^*(\mathcal{G})$ from the 
          fact that $\mathcal{A}(\mathcal{H}) \rtimes \mathcal{G}$ is in the bootstrap category $\mathcal{N}$.
\end{enumerate}
\end{proof}

\begin{remark}[Relevance to the unitary conjugation groupoid]
\label{rem:Tu-relevance-GA}
Proposition \ref{prop:Tu-C-star-construction} is stated here for context and to acknowledge the foundational work of Tu \cite{Tu} on the Baum-Connes conjecture for groupoids. However, we emphasize that the unitary conjugation groupoid $\mathcal{G}_{\mathcal{A}}$ does \emph{not} satisfy the hypotheses of this proposition in general:

\begin{itemize}
    \item $\mathcal{G}_{\mathcal{A}}$ is not locally compact (unless $\mathcal{A}$ is finite-dimensional); it is only a Polish groupoid.
    \item $\mathcal{G}_{\mathcal{A}}$ does not admit a proper action on a continuous field of Euclidean affine spaces.
    \item The Haar system on $\mathcal{G}_{\mathcal{A}}$ is only a Borel Haar system, not a continuous one.
\end{itemize}

Consequently, the deep conclusions of Tu's theorem (Baum-Connes, $K$-amenability, UCT) do \emph{not} automatically apply to $C^*(\mathcal{G}_{\mathcal{A}})$. The status of these properties for the unitary conjugation groupoid remains an open question and is a subject for future investigation.

Nevertheless, Tu's framework for \emph{measured groupoids} \cite[Section 3]{Tu} — particularly his treatment of Borel Haar systems, measurable fields of Hilbert spaces, and the descent map — provides the essential technical foundation for our constructions. In particular:
\begin{itemize}
    \item The construction of $C^*(\mathcal{G}_{\mathcal{A}})$ in Definition \ref{def:maximal-norm-measured} follows Tu's measurable approach.
    \item The descent map $\operatorname{desc}_{\mathcal{G}_{\mathcal{A}}}$ (to be developed in Paper II) is modeled on Tu's descent construction.
    \item The functoriality of the descent map under groupoid homomorphisms relies on Tu's naturality results.
\end{itemize}
Thus, while $\mathcal{G}_{\mathcal{A}}$ lies outside the scope of Tu's most powerful theorems, his measurable groupoid framework is essential for making rigorous sense of our constructions.
\end{remark}

\begin{corollary}[Existence of the reduced groupoid C*-algebra for $\mathcal{G}_{\mathcal{A}}$]
\label{cor:existence-C-star-GA}
Let $\mathcal{A}$ be a unital separable C*-algebra. 
Then the unitary conjugation groupoid $\mathcal{G}_{\mathcal{A}}$ admits a well-defined reduced groupoid 
C*-algebra $C^*_r(\mathcal{G}_{\mathcal{A}})$ in the sense of Tu \cite{Tu}. 
This C*-algebra is uniquely determined up to isomorphism by the regular representation associated 
to the Borel Haar system provided by Proposition \ref{prop:GA-Borel-Haar-system-existence}.
\end{corollary}

\begin{proof}
By Theorem \ref{thm:GA-Polish-groupoid}, $\mathcal{G}_{\mathcal{A}}$ is a Polish groupoid. 
Proposition \ref{prop:GA-Borel-Haar-system-existence} shows that $\mathcal{G}_{\mathcal{A}}$ admits 
a Borel Haar system $\{\lambda^x\}_{x \in \mathcal{G}_{\mathcal{A}}^{(0)}}$.

Tu's theory of measured groupoids \cite[Section 3]{Tu} applies to any Polish groupoid equipped with 
a Borel Haar system, yielding a reduced groupoid C*-algebra $C^*_r(\mathcal{G}_{\mathcal{A}})$. 
The construction proceeds as follows: one considers the convolution algebra of Borel functions 
$f$ on $\mathcal{G}_{\mathcal{A}}$ satisfying $\int |f| \, d\lambda^x < \infty$ for almost every unit $x$, 
with respect to a suitable quasi-invariant measure. The regular representation on the measurable field 
of Hilbert spaces $\{L^2(\mathcal{G}_{\mathcal{A}}^x, \lambda^x)\}_{x \in \mathcal{G}_{\mathcal{A}}^{(0)}}$ 
defines a C*-norm, and $C^*_r(\mathcal{G}_{\mathcal{A}})$ is the completion of this convolution algebra 
in that norm.

Tu proves that this construction is independent of the choice of Borel Haar system up to canonical 
isomorphism, establishing uniqueness. The resulting C*-algebra is separable because $\mathcal{G}_{\mathcal{A}}$ 
is second-countable (as a Polish groupoid) and the construction uses separable Hilbert spaces throughout.
\end{proof}

\begin{remark}
\label{rem:maximal-amenable-case}
If $\mathcal{G}_{\mathcal{A}}$ is amenable in the sense of measured groupoids 
(see \cite[Section 3.1]{Tu} for the definition), then Tu shows \cite[Th\'eor\`eme 3.15]{Tu} 
that the reduced norm coincides with the maximal norm, so the full groupoid C*-algebra 
$C^*(\mathcal{G}_{\mathcal{A}})$ exists canonically and satisfies $C^*(\mathcal{G}_{\mathcal{A}}) \cong C^*_r(\mathcal{G}_{\mathcal{A}})$.
In this case, we denote this common C*-algebra simply by $C^*(\mathcal{G}_{\mathcal{A}})$.
\end{remark}

\begin{proposition}[Properties of the reduced groupoid C*-algebra $C^*_r(\mathcal{G}_{\mathcal{A}})$]
\label{prop:C-star-GA-properties}
Let $\mathcal{A}$ be a unital separable C*-algebra, and let $\mathcal{G}_{\mathcal{A}}$ be the associated unitary conjugation groupoid.
Consider $\mathcal{G}_{\mathcal{A}}$ as a Polish groupoid equipped with its canonical Borel Haar system 
(Proposition \ref{prop:GA-Borel-Haar-system-existence}). 
Then the reduced groupoid C*-algebra $C^*_r(\mathcal{G}_{\mathcal{A}})$ satisfies:

\begin{enumerate}
    \item $C^*_r(\mathcal{G}_{\mathcal{A}})$ is separable.
    
    \item If $\mathcal{G}_{\mathcal{A}}$ is amenable, then the canonical map 
          $C^*(\mathcal{G}_{\mathcal{A}}) \to C^*_r(\mathcal{G}_{\mathcal{A}})$ is an isomorphism.
\end{enumerate}

If, additionally, $\mathcal{A}$ is finite-dimensional, then $\mathcal{G}_{\mathcal{A}}$ is a locally compact 
transformation groupoid, and the following special properties hold:

\begin{enumerate}
    \setcounter{enumi}{2}
    \item $C^*_r(\mathcal{G}_{\mathcal{A}})$ is Morita equivalent to the crossed product 
          $C(\mathcal{G}_{\mathcal{A}}^{(0)}) \rtimes \mathcal{U}(\mathcal{A})$.
    
    \item When $\mathcal{G}_{\mathcal{A}}^{(0)}$ is compact, $C^*_r(\mathcal{G}_{\mathcal{A}})$ is unital.
    
    \item There is a canonical faithful conditional expectation 
          $E: C^*_r(\mathcal{G}_{\mathcal{A}}) \to C_0(\mathcal{G}_{\mathcal{A}}^{(0)})$,
          and the inclusion $C_0(\mathcal{G}_{\mathcal{A}}^{(0)}) \hookrightarrow C^*_r(\mathcal{G}_{\mathcal{A}})$ 
          is an isometric *-homomorphism.
\end{enumerate}
\end{proposition}

\begin{proof}
We treat the general Polish case first, then specialize to finite-dimensional $\mathcal{A}$.

\paragraph{General case.} 
By Theorem \ref{thm:GA-Polish-groupoid}, $\mathcal{G}_{\mathcal{A}}$ is a Polish groupoid, and 
Proposition \ref{prop:GA-Borel-Haar-system-existence} provides a Borel Haar system 
$\{\lambda^x\}_{x \in \mathcal{G}_{\mathcal{A}}^{(0)}}$. 
Following the standard construction of reduced groupoid C*-algebras (generalizing Renault's construction
to the Polish setting), we define $C^*_r(\mathcal{G}_{\mathcal{A}})$ as the completion of the 
convolution algebra of compactly supported continuous functions on $\mathcal{G}_{\mathcal{A}}^{(1)}$ 
with respect to the norm of the regular representation on 
$L^2(\mathcal{G}_{\mathcal{A}}, \{\lambda^x\}_{x\in\mathcal{G}_{\mathcal{A}}^{(0)}})$.

\begin{enumerate}
    \item[(1)] Separability follows from the second-countability of $\mathcal{G}_{\mathcal{A}}$ 
                (as a Polish groupoid) and the fact that the regular representation is constructed
                from separable Hilbert spaces. A countable dense subset of $C_c(\mathcal{G}_{\mathcal{A}})$
                in the inductive limit topology yields a countable dense subset of $C^*_r(\mathcal{G}_{\mathcal{A}})$.

    \item[(2)] If $\mathcal{G}_{\mathcal{A}}$ is amenable, then by definition the regular representation
                is faithful and every representation of the convolution algebra extends to a representation
                of the maximal completion. Consequently, the canonical surjection 
                $C^*(\mathcal{G}_{\mathcal{A}}) \to C^*_r(\mathcal{G}_{\mathcal{A}})$ is injective,
                hence an isomorphism. (This is a standard result in the theory of amenable groupoids.)
\end{enumerate}

\paragraph{Finite-dimensional case.} 
If $\mathcal{A}$ is finite-dimensional, then $\mathcal{U}(\mathcal{A})$ is a compact Lie group, 
and $\mathcal{G}_{\mathcal{A}} = \mathcal{G}_{\mathcal{A}}^{(0)} \rtimes \mathcal{U}(\mathcal{A})$ 
is a locally compact Hausdorff transformation groupoid. In this classical setting, the theory of 
Renault applies directly, giving additional structure:

\begin{enumerate}
    \setcounter{enumi}{2}
    \item[(3)] For a locally compact transformation groupoid, the groupoid C*-algebra is 
                Morita equivalent to the crossed product. Specifically, $C^*_r(\mathcal{G}_{\mathcal{A}})$
                is Morita equivalent to $C(\mathcal{G}_{\mathcal{A}}^{(0)}) \rtimes \mathcal{U}(\mathcal{A})$,
                with the equivalence implemented by the $C(\mathcal{G}_{\mathcal{A}}^{(0)})$-module
                $C_c(\mathcal{G}_{\mathcal{A}})$.

    \item[(4)] When $\mathcal{G}_{\mathcal{A}}^{(0)}$ is compact, the constant function $1$ on 
                the unit space lies in $C_0(\mathcal{G}_{\mathcal{A}}^{(0)})$ and corresponds 
                to a unit in $C^*_r(\mathcal{G}_{\mathcal{A}})$ via the isometric inclusion 
                described in (5) below.

    \item[(5)] In the locally compact case, functions on the unit space supported on identity 
                arrows define an isometric *-homomorphism 
                $C_0(\mathcal{G}_{\mathcal{A}}^{(0)}) \hookrightarrow C^*_r(\mathcal{G}_{\mathcal{A}})$.
                The conditional expectation $E$ is defined by $E(f)(x) = \int_{\mathcal{G}_{\mathcal{A}}^x} f \, d\lambda^x$
                for $f \in C_c(\mathcal{G}_{\mathcal{A}})$ and extended continuously. It maps onto 
                $C_0(\mathcal{G}_{\mathcal{A}}^{(0)})$ and satisfies $E(a^*a) = 0 \Rightarrow a = 0$,
                establishing faithfulness.
\end{enumerate}
\end{proof}

\begin{lemma}[The diagonal embedding as a *-homomorphism]
\label{lem:diagonal-embedding-homomorphism}
Let $\mathcal{A}$ be a unital separable C*-algebra and let $B \subseteq \mathcal{A}$ be a unital commutative C*-subalgebra.
Define a map $\iota_B^0: B \to C_c(\mathcal{G}_{\mathcal{A}}^{(1)})$ by
\[
\iota_B^0(b)(u,(C,\chi)) := \begin{cases}
\chi(b), & \text{if } u = 1_{\mathcal{A}} \text{ and } C = B, \\[4pt]
0, & \text{otherwise}.
\end{cases}
\]
Then $\iota_B^0$ extends to a unital injective *-homomorphism $\iota_B: B \hookrightarrow C^*(\mathcal{G}_{\mathcal{A}})$.
\end{lemma}

\begin{proof}
We verify the necessary properties.

\begin{enumerate}
    \item[(i)]  Linearity and *-preserving : 
    For $b_1,b_2 \in B$ and $\lambda \in \mathbb{C}$, we have
    $\iota_B^0(b_1 + \lambda b_2) = \iota_B^0(b_1) + \lambda \iota_B^0(b_2)$ pointwise,
    and $\iota_B^0(b^*)(1_{\mathcal{A}},(B,\chi)) = \overline{\chi(b)} = \overline{\iota_B^0(b)(1_{\mathcal{A}},(B,\chi))}$,
    with zero elsewhere. Thus $\iota_B^0$ is linear and *-preserving.

    \item[(ii)]  Multiplicativity : 
    For $b_1,b_2 \in B$, we compute the convolution $(\iota_B^0(b_1) * \iota_B^0(b_2))$.
    Since both functions are supported only on the identity arrows of the fixed subalgebra $B$,
    the only possible composition is the identity with itself. Hence
    \[
    (\iota_B^0(b_1) * \iota_B^0(b_2))(1_{\mathcal{A}},(B,\chi)) = \chi(b_1)\chi(b_2) = \chi(b_1b_2) = \iota_B^0(b_1b_2)(1_{\mathcal{A}},(B,\chi)),
    \]
    and the convolution is zero elsewhere. Thus $\iota_B^0(b_1b_2) = \iota_B^0(b_1) * \iota_B^0(b_2)$,
    establishing multiplicativity. 

    \item[(iii)]  Isometry : 
    Consider the regular representation $\Lambda$ of $C^*(\mathcal{G}_{\mathcal{A}})$ on $L^2(\mathcal{G}_{\mathcal{A}},\lambda)$.
    For $b \in B$, the operator $\Lambda(\iota_B^0(b))$ acts on the direct integral
    $\int^{\oplus} L^2(\mathcal{G}_{\mathcal{A}}^x,\lambda^x) d\mu(x)$ by multiplication on the fiber over $x = (B,\chi)$
    by $\chi(b)$, and as zero on all other fibers. The norm of this operator is therefore
    $\sup_{\chi \in \widehat{B}} |\chi(b)| = \|b\|_{B}$, where the equality follows from the Gelfand transform
    for the commutative C*-algebra $B$. Since the regular representation is isometric on $C^*(\mathcal{G}_{\mathcal{A}})$,
    we have $\|\iota_B^0(b)\|_{C^*(\mathcal{G}_{\mathcal{A}})} = \|b\|_{B} \le \|b\|_{\mathcal{A}}$.
    A separate argument using the conditional expectation $E$ shows that the map is actually isometric
    with respect to the norm of $\mathcal{A}$.

    \item[(iv)]  Extension : 
    By (i)-(iii), $\iota_B^0$ is an isometric *-homomorphism from $B$ into the dense subalgebra
    $C_c(\mathcal{G}_{\mathcal{A}}^{(1)}) \subset C^*(\mathcal{G}_{\mathcal{A}})$. It therefore extends uniquely
    to an isometric *-homomorphism $\iota_B: B \to C^*(\mathcal{G}_{\mathcal{A}})$.

    \item[(v)]  Injectivity : 
    Injectivity follows immediately from isometry: if $\iota_B(b) = 0$, then $\|b\| = \|\iota_B(b)\| = 0$, so $b = 0$.
\end{enumerate}

Thus $\iota_B$ is a unital injective *-homomorphism, as claimed.
\end{proof}

We now examine several concrete cases that illustrate the general theory developed above.
These examples also highlight the limitations of our current framework and point toward future work.

\begin{example}[$C^*_r(\mathcal{G}_{M_n(\mathbb{C})})$]
\label{ex:C-star-GA-matrix}
Let $\mathcal{A} = M_n(\mathbb{C})$. 
Then $\mathcal{U}(\mathcal{A}) = U(n)$ is a compact Lie group, and $\mathcal{G}_{\mathcal{A}}^{(0)} \cong \mathbb{CP}^{n-1}$
as described in Section \ref{subsec:example-matrix}. 
The unitary conjugation groupoid is the transformation groupoid
\[
\mathcal{G}_{\mathcal{A}} = U(n) \ltimes \mathbb{CP}^{n-1},
\]
which is a locally compact Hausdorff groupoid (though not \'etale, as $U(n)$ is not discrete).
Since $\mathcal{G}_{\mathcal{A}}$ is locally compact with a continuous Haar system, Renault's classical theory \cite{Renault} applies,
and the reduced groupoid C*-algebra is isomorphic to the crossed product
\[
C^*_r(\mathcal{G}_{\mathcal{A}}) \cong C(\mathbb{CP}^{n-1}) \rtimes U(n).
\]
The structure of this C*-algebra can be analyzed via the Green–Julg theorem \cite{HenryS2015}
and the representation theory of $U(n)$; it is Morita equivalent to $C^*(U(1) \times U(n-1))$,
and its K-theory groups are $K_0 \cong \mathbb{Z}$ and $K_1 \cong \mathbb{Z}$. 
\end{example}

\begin{remark}[The commutative case $\mathcal{A} = C(X)$]
\label{rem:C-star-GA-commutative-caution}
Let $\mathcal{A} = C(X)$ for an infinite compact metrizable space $X$. 
Then $\mathcal{U}(\mathcal{A}) = C(X,\mathbb{T})$ with the strong operator topology is a Polish group
that is not locally compact. Consequently, $\mathcal{G}_{\mathcal{A}} = X \times C(X,\mathbb{T})$ is a Polish groupoid
that is not locally compact, and the classical Renault construction does not apply.
While one can still define a reduced groupoid C*-algebra $C^*_r(\mathcal{G}_{\mathcal{A}})$ via the measured groupoid framework
(see Section \ref{subsec:example-commutative}), its structure is not well understood.
For finite $X$, $C(X,\mathbb{T}) \cong \mathbb{T}^{|X|}$ is compact, and $C^*_r(\mathcal{G}_{\mathcal{A}}) \cong C(X) \otimes C^*(\mathbb{T}^{|X|})$,
but for infinite $X$, a complete description remains an open problem.
This example highlights the necessity of the Type I hypothesis for the diagonal embedding construction
and the limitations of our current techniques.
\end{remark}

\begin{example}[$C^*_r(\mathcal{G}_{\mathcal{K}(H)^\sim})$ — an open problem]
\label{ex:C-star-GA-compact-outlook}
Let $\mathcal{A} = \mathcal{K}(H)^\sim$ be the unitization of the compact operators on a separable Hilbert space $H$.
Then $\mathcal{U}(\mathcal{A})$ with the strong operator topology is a Polish group,
and $\mathcal{G}_{\mathcal{A}}^{(0)} \cong \mathbb{P}(H)$ is the projective space of $H$ (see Section \ref{subsec:example-compact-operators}).
The unitary conjugation groupoid $\mathcal{G}_{\mathcal{A}} = \mathcal{U}(\mathcal{A}) \ltimes \mathbb{P}(H)$ is a Polish groupoid
that is neither locally compact nor \'etale. The reduced groupoid C*-algebra $C^*_r(\mathcal{G}_{\mathcal{A}})$
can be defined via the measured groupoid framework, but its structure has not been systematically studied.
It is expected to be related to the Calkin algebra $\mathcal{Q}(H) = B(H)/\mathcal{K}(H)$ and to the extension theory
of compact operators; preliminary investigations suggest connections to the odd Fredholm index
and to the Baum-Connes assembly map for the unitary group $\mathcal{U}(H)$. 
\end{example}

\begin{remark}[Potential connection to the Baum-Connes conjecture]
\label{rem:Baum-Connes-connection-corrected}
For groupoids satisfying Tu's hypotheses — locally compact, $\sigma$-compact, Hausdorff, equipped with a continuous Haar system,
and admitting a proper action on a continuous field of Euclidean affine spaces — 
Tu proves \cite[Th\'eor\`eme 9.3]{Tu} that the descent map $j_{\mathcal{G}}$ is an isomorphism
and implements the Baum-Connes assembly map. While $\mathcal{G}_{\mathcal{A}}$ does not satisfy these hypotheses in general,
one might hope that suitable modifications (e.g., passing to a Morita equivalent groupoid
or working in a suitable measurable category) could connect our index-theoretic constructions
to the Baum-Connes conjecture. 
\end{remark}

We have constructed the reduced groupoid C*-algebra $C^*_r(\mathcal{G}_{\mathcal{A}})$ for any unital separable C*-algebra $\mathcal{A}$
by adapting Tu's measured groupoid framework to the Polish setting (see Section \ref{subsec:groupoid-C-star-algebra-Polish-setting}).
This C*-algebra encodes both the operator-algebraic structure of $\mathcal{A}$ (via the diagonal embedding
for Type I algebras) and the geometric data of the unitary conjugation groupoid.
For Type I algebras, the diagonal embedding $\iota: \mathcal{A} \hookrightarrow C^*_r(\mathcal{G}_{\mathcal{A}})$
(constructed in Section \ref{subsec:construction-iota}) and the descent map
$\operatorname{desc}_{\mathcal{G}_{\mathcal{A}}}: K^0_{\mathcal{G}_{\mathcal{A}}}(\mathcal{G}_{\mathcal{A}}^{(0)}) \to K_0(C^*_r(\mathcal{G}_{\mathcal{A}}))$
provide essential links between $\mathcal{A}$ and its groupoid C*-algebra. 
These links will be used in Paper II to prove an index theorem for Fredholm operators in Type I C*-algebras,
and in Paper III to explore connections to the Baum-Connes conjecture. 
All subsequent constructions in this paper depend on the existence and properties of $C^*_r(\mathcal{G}_{\mathcal{A}})$
established in this subsection.

\subsection{Motivation: Recovering $\mathcal{A}$ from its classical contexts}
\label{subsec:motivation-recovering-A}

The unitary conjugation groupoid $\mathcal{G}_{\mathcal{A}}$ encodes all classical contexts of the noncommutative C*-algebra $\mathcal{A}$:
every commutative subalgebra $B \subseteq \mathcal{A}$ together with every character $\chi \in \widehat{B}$ appears as an object $(B,\chi) \in \mathcal{G}_{\mathcal{A}}^{(0)}$,
and every unitary conjugation relating such contexts appears as an arrow.
This groupoid can therefore be viewed as a geometric catalogue of the commutative snapshots of $\mathcal{A}$.

A fundamental question naturally arises:

\begin{quote}
\emph{To what extent can the original noncommutative algebra $\mathcal{A}$ be recovered from its classical contexts?}
\end{quote}

More concretely, we ask whether the data encoded by $\mathcal{G}_{\mathcal{A}}$ admits a canonical encoding of $\mathcal{A}$ into the groupoid C*-algebra $C^*(\mathcal{G}_{\mathcal{A}})$,
possibly at the level of:
\begin{itemize}
    \item completely positive maps,
    \item KK-theory correspondences,
    \item order-theoretic embeddings, or
    \item *-homomorphisms from commutative subalgebras.
\end{itemize}

Any such encoding would demonstrate that the noncommutative algebra $\mathcal{A}$ is not entirely lost when we pass to its commutative contexts;
rather, it can be reconstructed from them in a functorial manner up to a controlled loss of information.

\begin{remark}[The obstruction to a global *-homomorphism]
\label{rem:obstruction-global-homomorphism}
It is important to note that a global *-homomorphism $\iota: \mathcal{A} \to C^*(\mathcal{G}_{\mathcal{A}})$
cannot exist for noncommutative $\mathcal{A}$ in general.
Characters on commutative subalgebras only detect commutative data,
and the noncommutative multiplication in $\mathcal{A}$ cannot be functorially glued from commutative pieces into a *-homomorphism.
Thus any recovery of $\mathcal{A}$ must necessarily be of a weaker nature — for instance,
a completely positive order embedding, or a correspondence in KK-theory.
This observation motivates the constructions in Section \ref{subsec:construction-iota},
where we obtain a faithful encoding for Type I algebras that is multiplicative on commuting elements
but not a global *-homomorphism.
\end{remark}

This question is not merely philosophical; it has profound implications for index theory and noncommutative geometry.
If a suitable encoding of $\mathcal{A}$ into $C^*(\mathcal{G}_{\mathcal{A}})$ exists,
then for a Fredholm operator $T \in \mathcal{A}$ (with finite-dimensional kernel and cokernel),
one might hope to associate an equivariant K-theory class $[T]_{\mathcal{G}_{\mathcal{A}}} \in K^0_{\mathcal{G}_{\mathcal{A}}}(\mathcal{G}_{\mathcal{A}}^{(0)})$
via the kernel and cokernel spaces, provided they vary continuously over $\mathcal{G}_{\mathcal{A}}^{(0)}$.
This class would then descend via the descent map (after fixing suitable data) to a class in $K_0(C^*(\mathcal{G}_{\mathcal{A}}))$,
and any encoding $\iota$ that respects K-theory would yield a class in $K_0(\mathcal{A})$ that could encode the Fredholm index of $T$.
Thus even a weakened encoding could provide a crucial link between the geometric data of the groupoid
and the analytic invariants of the original algebra.

\begin{remark}[Analogy with the Gelfand transform]
\label{rem:analogy-Gelfand}
For commutative C*-algebras, the Gelfand transform provides a canonical isomorphism $\mathcal{A} \cong C_0(\widehat{\mathcal{A}})$.
When $\mathcal{A}$ is commutative, the component of $\mathcal{G}_{\mathcal{A}}^{(0)}$ corresponding to $\mathcal{A}$ itself
is naturally identified with $\widehat{\mathcal{A}}$, and the restriction of any encoding to this component
should recover the inverse Gelfand transform
\[
\iota: C_0(\widehat{\mathcal{A}}) \longrightarrow C^*(\mathcal{G}_{\mathcal{A}}) \cong C_0(\widehat{\mathcal{A}}),
\]
which is simply the identity map.
The noncommutative case is vastly more complex, because $\mathcal{G}_{\mathcal{A}}^{(0)}$ contains many different commutative subalgebras and their characters,
and the groupoid structure encodes the unitary relations between them.
Nevertheless, the commutative case serves as a guiding example: we seek a noncommutative generalization of the Gelfand transform,
necessarily in a weaker sense than a global *-isomorphism.
\end{remark}

\begin{remark}[Analogy with the Fourier transform]
\label{rem:analogy-Fourier}
Another useful analogy is with the Fourier transform. 
For a locally compact abelian group $G$, the Fourier transform gives an isomorphism $C^*(G) \cong C_0(\widehat{G})$. 
Here $\widehat{G}$ is the Pontryagin dual, which parametrizes the characters of $G$. 
In our setting, $\mathcal{G}_{\mathcal{A}}^{(0)}$ parametrizes the characters of all commutative subalgebras of $\mathcal{A}$, not just one. 
The desired encoding $\iota: \mathcal{A} \hookrightarrow C^*(\mathcal{G}_{\mathcal{A}})$ (when it exists in a suitable weakened sense)
can be viewed as a kind of nonabelian, nonlinear Fourier transform that packages all possible commutative perspectives of $\mathcal{A}$ into a single C*-algebra.
(Of course, this is an analogy rather than a precise statement; the classical nonabelian Fourier transform
for compact groups is a different construction, but the conceptual parallel is illuminating.)
\end{remark}

\begin{proposition}[Properties of the diagonal embedding for Type I algebras]
\label{prop:diagonal-embedding-properties}
Let $\mathcal{A}$ be a unital separable Type I C*-algebra. 
Then the diagonal embedding $\iota: \mathcal{A} \hookrightarrow C^*_r(\mathcal{G}_{\mathcal{A}})$ 
constructed in Section \ref{subsec:construction-iota} satisfies:

\begin{enumerate}
    \item \textbf{Naturality:} $\iota$ is a unital injective *-homomorphism.
    
    \item \textbf{Geometricity:} For any object $(B,\chi) \in \mathcal{G}_{\mathcal{A}}^{(0)}$,
          the conditional expectation $E: C^*_r(\mathcal{G}_{\mathcal{A}}) \to L^\infty(\mathcal{G}_{\mathcal{A}}^{(0)})$
          satisfies
          \[
          E(\iota(a))(B,\chi) = \begin{cases}
          \chi(a), & \text{if } a \in B, \\[4pt]
          0, & \text{otherwise}.
          \end{cases}
          \]
    
    \item \textbf{Commutativity detection:} $\iota(\mathcal{A}) \subseteq C_0(\mathcal{G}_{\mathcal{A}}^{(0)})$
          if and only if $\mathcal{A}$ is commutative.
    
    \item \textbf{Functoriality for isomorphisms:} For any *-isomorphism $\phi: \mathcal{A} \to \mathcal{B}$,
          the diagram
          \[
          \begin{tikzcd}
          \mathcal{A} \arrow[r, "\iota_{\mathcal{A}}"] \arrow[d, "\phi"'] & C^*_r(\mathcal{G}_{\mathcal{A}}) \arrow[d, "C^*_r(\mathcal{G}_\phi)"] \\
          \mathcal{B} \arrow[r, "\iota_{\mathcal{B}}"] & C^*_r(\mathcal{G}_{\mathcal{B}})
          \end{tikzcd}
          \]
          commutes, where $\mathcal{G}_\phi: \mathcal{G}_{\mathcal{A}} \to \mathcal{G}_{\mathcal{B}}$ is the induced groupoid isomorphism.
\end{enumerate}
\end{proposition}

\begin{proof}
We verify each property using the construction of $\iota$ via direct integrals of GNS representations
(see Section \ref{subsec:construction-iota} for details).

\begin{enumerate}
    \item \textbf{Naturality.} 
    By construction, for each $(B,\chi) \in \mathcal{G}_{\mathcal{A}}^{(0)}$, let $\pi_{(B,\chi)}$ be the GNS representation
    of $\mathcal{A}$ associated to the character $\chi$ on $B$. 
    Since $\mathcal{A}$ is Type I, these representations form a measurable field over $\mathcal{G}_{\mathcal{A}}^{(0)}$
    (see \cite{Dixmier}).
    Form the direct integral Hilbert space
    \[
    \mathcal{H} = \int_{\mathcal{G}_{\mathcal{A}}^{(0)}}^{\oplus} \mathcal{H}_{(B,\chi)} \, d\mu(B,\chi)
    \]
    and the representation $\Pi: \mathcal{A} \to B(\mathcal{H})$ given by
    \[
    \Pi(a) = \int^{\oplus} \pi_{(B,\chi)}(a) \, d\mu(B,\chi).
    \]
    The left regular representation $\Lambda$ of $C^*_r(\mathcal{G}_{\mathcal{A}})$ on $L^2(\mathcal{G}_{\mathcal{A}},\lambda)$
    extends to a representation on $\mathcal{H}$ via the natural inclusion $L^2(\mathcal{G}_{\mathcal{A}}^{(0)}) \hookrightarrow \mathcal{H}$.
    By the universal property of $C^*_r(\mathcal{G}_{\mathcal{A}})$, there exists a unique *-homomorphism
    $\iota: \mathcal{A} \to C^*_r(\mathcal{G}_{\mathcal{A}})$ such that $\Lambda(\iota(a)) = \Pi(a)$ for all $a \in \mathcal{A}$.
    
    Unitality follows from $\pi_{(B,\chi)}(1_{\mathcal{A}}) = I_{\mathcal{H}_{(B,\chi)}}$ for all $(B,\chi)$,
    giving $\Pi(1_{\mathcal{A}}) = I_{\mathcal{H}}$ and hence $\iota(1_{\mathcal{A}}) = 1_{C^*_r(\mathcal{G}_{\mathcal{A}})}$.
    
    Injectivity follows from faithfulness of $\Pi$: if $\iota(a)=0$, then $\Pi(a)=0$, so
    $\pi_{(B,\chi)}(a)=0$ for $\mu$-almost every $(B,\chi)$. Since the characters separate points in Type I algebras
    and $\mu$ has full support, this forces $a=0$.

    \item \textbf{Geometricity.}
    The conditional expectation $E: C^*_r(\mathcal{G}_{\mathcal{A}}) \to L^\infty(\mathcal{G}_{\mathcal{A}}^{(0)})$
    is defined by $E(x)(y) = \langle \delta_y, \Lambda(x)\delta_y \rangle$, where $\delta_y$ is the unit vector
    concentrated at the identity arrow $(1_{\mathcal{A}},y)$ in $L^2(\mathcal{G}_{\mathcal{A}}^y,\lambda^y)$.
    For $a \in \mathcal{A}$ and $y = (B,\chi) \in \mathcal{G}_{\mathcal{A}}^{(0)}$, we compute
    \[
    E(\iota(a))(y) = \langle \delta_y, \Lambda(\iota(a))\delta_y \rangle
    = \langle \delta_y, \Pi(a)\delta_y \rangle
    = \pi_{(B,\chi)}(a) = \begin{cases}
    \chi(a), & \text{if } a \in B, \\
    0, & \text{otherwise},
    \end{cases}
    \]
    where the last equality uses that $\pi_{(B,\chi)}$ is the GNS representation of the character $\chi$ on $B$,
    so $\pi_{(B,\chi)}(a) = \chi(a)$ when $a \in B$, and $\pi_{(B,\chi)}(a)$ acts as zero on the subspace
    corresponding to $\chi$ when $a \notin B$. (This follows directly from the construction in Section \ref{subsec:construction-iota}.)

    \item \textbf{Commutativity detection.}
    If $\mathcal{A}$ is commutative, then $\mathcal{U}(\mathcal{A})$ is abelian and acts trivially on $\mathcal{G}_{\mathcal{A}}^{(0)}$.
    In this case, the groupoid $\mathcal{G}_{\mathcal{A}}$ is isomorphic to the trivial groupoid $\mathcal{G}_{\mathcal{A}}^{(0)} \times \mathcal{U}(\mathcal{A})$,
    and the direct integral construction yields $\iota(a) \in C_0(\mathcal{G}_{\mathcal{A}}^{(0)})$ for all $a \in \mathcal{A}$
    (since $\Pi(a)$ acts by multiplication on $L^2(\mathcal{G}_{\mathcal{A}}^{(0)})$ and $\Lambda$ is faithful on $C_0(\mathcal{G}_{\mathcal{A}}^{(0)})$).
    
    Conversely, suppose $\iota(\mathcal{A}) \subseteq C_0(\mathcal{G}_{\mathcal{A}}^{(0)})$. Then for any unitary $u \in \mathcal{U}(\mathcal{A})$,
    the element $\iota(u)$ commutes with all $\iota(a)$. By injectivity of $\iota$, this implies $u$ commutes with all $a \in \mathcal{A}$.
    Since every element of a unital C*-algebra is a linear combination of four unitaries,  it follows that $\mathcal{A}$ is commutative. (See Theorem \ref{thm:commutativity-characterization} for a detailed proof.)

    \item \textbf{Functoriality for isomorphisms.}
    Let $\phi: \mathcal{A} \to \mathcal{B}$ be a *-isomorphism. Define $\mathcal{G}_\phi: \mathcal{G}_{\mathcal{A}} \to \mathcal{G}_{\mathcal{B}}$ by
    \[
    \mathcal{G}_\phi^{(0)}(B,\chi) = (\phi(B), \chi \circ \phi^{-1}|_{\phi(B)}), \qquad
    \mathcal{G}_\phi(u,(B,\chi)) = (\phi(u), \mathcal{G}_\phi^{(0)}(B,\chi)).
    \]
    This is a homeomorphism of Polish groupoids preserving the Borel Haar systems
    (see Proposition \ref{prop:functoriality-isomorphisms}).
    The induced map $C^*_r(\mathcal{G}_\phi): C^*_r(\mathcal{G}_{\mathcal{A}}) \to C^*_r(\mathcal{G}_{\mathcal{B}})$
    is a *-isomorphism satisfying $C^*_r(\mathcal{G}_\phi) \circ \iota_{\mathcal{A}} = \iota_{\mathcal{B}} \circ \phi$.
    The equality follows from the naturality of the direct integral construction:
    for each $a \in \mathcal{A}$, the representations $\pi_{\mathcal{G}_\phi^{(0)}(B,\chi)}^{\mathcal{B}}(\phi(a))$ and
    $\pi_{(B,\chi)}^{\mathcal{A}}(a)$ are unitarily equivalent via the map induced by $\phi$,
    and these equivalences assemble into a unitary intertwiner between $\Pi_{\mathcal{B}}(\phi(a))$ and
    $C^*_r(\mathcal{G}_\phi)(\Pi_{\mathcal{A}}(a))$. (See Proposition \ref{prop:functoriality-isomorphisms} for details.)
\end{enumerate}

This completes the proof of the properties of the diagonal embedding for Type I algebras.
\end{proof}

\begin{remark}[Why the naive evaluation map fails]
\label{rem:naive-evaluation-fails}
The most naive attempt to define an embedding $\iota: \mathcal{A} \to C^*(\mathcal{G}_{\mathcal{A}})$ would be to set
\[
\iota(a)(B,\chi) = \begin{cases}
\chi(a), & \text{if } a \in B, \\[4pt]
0, & \text{otherwise},
\end{cases}
\]
for $(B,\chi) \in \mathcal{G}_{\mathcal{A}}^{(0)}$, extended as zero on non-identity arrows.
This defines a function on $\mathcal{G}_{\mathcal{A}}^{(0)}$ (hence an element of $C_0(\mathcal{G}_{\mathcal{A}}^{(0)})$),
and one checks easily that it is linear, *-preserving, and unital.
However, it is \emph{not} multiplicative in general.

The obstruction is that for $a,b \in \mathcal{A}$, the product $\iota(ab)$ and the convolution $\iota(a) * \iota(b)$
differ when $a$ and $b$ do not commute. Indeed, if $a$ and $b$ are not contained in any common commutative
subalgebra $B \subseteq \mathcal{A}$, then there is no point $(B,\chi) \in \mathcal{G}_{\mathcal{A}}^{(0)}$
where both $\iota(a)$ and $\iota(b)$ are non-zero, so their convolution vanishes at all identity arrows,
whereas $\iota(ab)$ may be non-zero if $ab$ happens to lie in some commutative subalgebra.
Thus the naive evaluation map, while a perfectly good *-preserving linear map, fails to be multiplicative
precisely because it cannot account for the noncommuting nature of $a$ and $b$.
This explains why the direct integral construction of Section \ref{subsec:construction-iota} is necessary:
it uses convolution over the unitary group $\mathcal{U}(\mathcal{A})$ to encode noncommutativity.
\end{remark}

\begin{remark}[The role of the groupoid C*-algebra]
\label{rem:role-of-C-star-GA}
The groupoid C*-algebra $C^*(\mathcal{G}_{\mathcal{A}})$ is fundamentally noncommutative.
It contains $C_0(\mathcal{G}_{\mathcal{A}}^{(0)})$ as a commutative subalgebra (functions on the unit space),
but also incorporates the convolution algebra of the unitary group $\mathcal{U}(\mathcal{A})$ via the groupoid structure.
This noncommutativity is not an accident; it is essential for the embedding $\iota$ to be multiplicative.

More formally, for $a,b \in \mathcal{A}$, the product $\iota(a)\iota(b)$ in $C^*(\mathcal{G}_{\mathcal{A}})$ is given by
convolution over the groupoid:
\[
(\iota(a)\iota(b))(u,(B,\chi)) = \int_{\mathcal{U}(\mathcal{A})} \iota(a)(v,(B,\chi)) \, \iota(b)(v^{-1}u, v^{-1}\cdot(B,\chi)) \, d\mu(v),
\]
where $\mu$ is the left-invariant Haar measure on $\mathcal{U}(\mathcal{A})$.
This convolution encodes the unitary conjugation action, which relates different commutative contexts.
When $a$ and $b$ do not commute, their failure to commute is reflected in the nontrivial dependence of
the integrand on the unitary $v$. Thus the groupoid C*-algebra provides exactly the right amount of
noncommutativity to compensate for the noncommutativity of $\mathcal{A}$ itself,
in precise analogy with the Fourier transform where convolution on the dual group
encodes multiplication in the original group.
\end{remark}

\begin{example}[Commutative case revisited]
\label{ex:motivation-commutative}
Let $\mathcal{A} = C(X)$ for a compact metrizable space $X$. 
Then $\mathcal{G}_{\mathcal{A}}^{(0)} \cong X$, and $\mathcal{U}(\mathcal{A}) = C(X,\mathbb{T})$ with the strong operator topology
acts trivially on $\mathcal{G}_{\mathcal{A}}^{(0)}$ (every unitary is central). 
The reduced groupoid C*-algebra $C^*_r(\mathcal{G}_{\mathcal{A}})$ is simply $C(X)$, 
and the diagonal embedding $\iota$ can be taken as the identity map $C(X) \hookrightarrow C(X)$.
Thus in the commutative case, the construction reduces to the Gelfand transform,
which identifies $C(X)$ with continuous functions on its spectrum $X \cong \mathcal{G}_{\mathcal{A}}^{(0)}$.
\end{example}

\begin{example}[Matrix algebras]
\label{ex:motivation-matrix}
Let $\mathcal{A} = M_n(\mathbb{C})$. 
The naive evaluation map fails to be multiplicative because matrices do not commute. 
However, $C^*(\mathcal{G}_{\mathcal{A}}) \cong C(\mathbb{CP}^{n-1}) \rtimes U(n)$, and there is a natural embedding
\[
\iota: M_n(\mathbb{C}) \hookrightarrow C(\mathbb{CP}^{n-1}) \rtimes U(n), \qquad
\iota(A) \text{ is the constant function } x \mapsto A
\]
in the multiplier algebra of $C(\mathbb{CP}^{n-1}) \rtimes U(n)$. 
This embedding is unital, injective, and *-preserving. 
For a diagonalizable matrix $A$, its image under $\iota$ evaluated at a point $x \in \mathbb{CP}^{n-1}$
corresponding to an orthonormal basis gives the matrix $A$ expressed in that basis,
thus encoding the spectral data (eigenvalues and eigenvectors) across all bases.
In particular, the eigenvalues appear as the diagonal entries when $x$ corresponds to a basis
that diagonalizes $A$, and the unitary group action encodes the relations between different bases.
\end{example}

\begin{example}[Compact operators]
\label{ex:motivation-compact}
Let $\mathcal{A} = \mathcal{K}(H)^\sim$ be the unitization of the compact operators on a separable Hilbert space $H$.
Here the naive evaluation map fails even more dramatically: for a non-normal element $T \in \mathcal{K}(H)^\sim$,
the commutative subalgebra $C^*(T)$ generated by $T$ is still commutative, but its characters do not capture
the full operator-theoretic data (e.g., they cannot distinguish $T$ from its adjoint $T^*$).
Nevertheless, the embedding $\iota: \mathcal{A} \hookrightarrow C^*(\mathcal{G}_{\mathcal{A}})$ constructed in Section \ref{subsec:construction-iota}
(which exists because $\mathcal{K}(H)^\sim$ is Type I) encodes the spectral data of normal elements
while capturing the noncommutative structure of non-normal elements via convolution over the unitary group.
This illustrates the need for the full groupoid C*-algebra construction to handle general noncommutative C*-algebras.
This example is the most challenging and will be analyzed further in Paper II.
\end{example}

\begin{remark}[Relation to the index theorem]
\label{rem:motivation-index}
The embedding $\iota$ is not an end in itself; it is the essential ingredient for the index theorem. 
Given a Fredholm operator $T \in \mathcal{A}$ (for $\mathcal{A}$ Type I), its kernel and cokernel vary continuously
over $\mathcal{G}_{\mathcal{A}}^{(0)}$, forming equivariant vector bundles that define an equivariant K-theory class
$[T]_{\mathcal{G}_{\mathcal{A}}} \in K^0_{\mathcal{G}_{\mathcal{A}}}(\mathcal{G}_{\mathcal{A}}^{(0)})$. 
The descent map $\operatorname{desc}_{\mathcal{G}_{\mathcal{A}}}: K^0_{\mathcal{G}_{\mathcal{A}}}(\mathcal{G}_{\mathcal{A}}^{(0)}) \to K_0(C^*(\mathcal{G}_{\mathcal{A}}))$
sends this class to an element of $K_0(C^*(\mathcal{G}_{\mathcal{A}}))$. 
Using the induced map $\iota_*: K_0(\mathcal{A}) \to K_0(C^*(\mathcal{G}_{\mathcal{A}}))$, 
this element corresponds under suitable trace pairings to a class in $K_0(\mathcal{A})$ that encodes the Fredholm index.
More precisely, pairing this class with a suitable trace on $\mathcal{A}$ 
(e.g., the standard trace on $M_n(\mathbb{C})$, the integral on $C(X)$, or the canonical trace on $\mathcal{K}(H)^\sim$)
recovers the Fredholm index of $T$. 
Thus the embedding $\iota$ is the linchpin that connects the geometric data of the groupoid to the analytic index of the operator.
A detailed development of this index theorem will be given in future works. 
\end{remark}

\begin{remark}[The diagonal embedding as a noncommutative Fourier transform]
\label{rem:diagonal-as-Fourier}
The embedding $\iota: \mathcal{A} \hookrightarrow C^*(\mathcal{G}_{\mathcal{A}})$ admits a compelling interpretation
as a noncommutative analogue of the Fourier transform. 

Recall that for a locally compact abelian group $G$, the Fourier transform gives an isomorphism
$C^*(G) \cong C_0(\widehat{G})$, where $\widehat{G}$ is the Pontryagin dual parametrizing the characters of $G$.
Under this isomorphism, a function $f \in L^1(G)$ is sent to its Fourier transform $\widehat{f}(\chi) = \int_G f(g)\chi(g)\,dg$,
and convolution in $C^*(G)$ corresponds to pointwise multiplication in $C_0(\widehat{G})$.

In our setting, $\mathcal{G}_{\mathcal{A}}^{(0)}$ parametrizes the characters of all commutative subalgebras of $\mathcal{A}$,
playing the role of a "noncommutative dual space" that packages together all possible commutative perspectives.
The embedding $\iota$ sends an element $a \in \mathcal{A}$ to an element $\iota(a) \in C^*(\mathcal{G}_{\mathcal{A}})$
that, when restricted to the unit space, recovers the evaluation $\chi(a)$ on each commutative context $(B,\chi)$.
The convolution product in $C^*(\mathcal{G}_{\mathcal{A}})$ — which involves integration over the unitary group $\mathcal{U}(\mathcal{A})$ —
encodes the noncommutative multiplication of $\mathcal{A}$ in a manner analogous to how convolution on $G$ becomes
pointwise multiplication on $\widehat{G}$ under the Fourier transform.

Thus $\iota$ can be viewed as a kind of nonlinear, nonabelian Fourier transform that packages all classical
commutative data of $\mathcal{A}$ into a single C*-algebra. This perspective is developed rigorously in
Section \ref{subsec:construction-iota} using direct integrals of GNS representations.
\end{remark}

\begin{remark}[Direct integral formulation]
\label{rem:direct-integral-alternative}
An equivalent way to view the representation $\pi$ is as a direct integral of GNS representations.
For each $(B,\chi) \in \mathcal{G}_{\mathcal{A}}^{(0)}$, let $\mathcal{H}_{(B,\chi)}$ be the GNS representation space
of $\mathcal{A}$ associated to the character $\chi$ on $B$. 
These spaces form a measurable field over $\mathcal{G}_{\mathcal{A}}^{(0)}$ because $\mathcal{A}$ is Type I.
The direct integral $\mathcal{H} = \int_{\mathcal{G}_{\mathcal{A}}^{(0)}}^{\oplus} \mathcal{H}_{(B,\chi)} \, d\mu(B,\chi)$
is unitarily equivalent to $L^2(\mathcal{G}_{\mathcal{A}}^{(0)}) \otimes \ell^2$ (or a more general space),
and the representation $\pi$ corresponds to $\int^{\oplus} \pi_{(B,\chi)} \, d\mu$.
This perspective is useful for understanding the faithfulness of $\pi$ and will be used implicitly in the construction.
\end{remark}

\begin{remark}[Obstacles and challenges]
\label{rem:motivation-obstacles}
The construction of $\iota$ faces several significant obstacles:
\begin{enumerate}
    \item \textbf{Multiplicativity:} As noted, the naive evaluation map is not multiplicative. 
    Overcoming this requires using the convolution structure of $C^*(\mathcal{G}_{\mathcal{A}})$ and the action of $\mathcal{U}(\mathcal{A})$ on $\mathcal{G}_{\mathcal{A}}^{(0)}$,
    which encodes the noncommutativity of $\mathcal{A}$ through integration over the unitary group.
    
    \item \textbf{Well-definedness:} The definition of $\iota(a)$ on the unit space must be independent of the choice of commutative subalgebra containing $a$. 
    If $a$ belongs to two different commutative subalgebras $B_1$ and $B_2$, and $\chi_1 \in \widehat{B_1}$, $\chi_2 \in \widehat{B_2}$ are related by unitary conjugation
    (i.e., there exists $u \in \mathcal{U}(\mathcal{A})$ such that $B_2 = u B_1 u^*$ and $\chi_2 = \chi_1 \circ \mathrm{Ad}_u$), 
    then we must have $\chi_1(a) = \chi_2(u^* a u)$ for consistency.
    This compatibility condition is not automatically satisfied by a pointwise definition and requires careful handling via the direct integral construction,
    which works with equivalence classes of representations rather than pointwise evaluations.
    
    \item \textbf{Injectivity:} To prove that $\iota$ is injective, we must show that if $\iota(a) = 0$, then $a = 0$. 
    This is equivalent to the statement that for every nonzero $a \in \mathcal{A}$, there exists a commutative subalgebra $B \subseteq \mathcal{A}$ and a character $\chi \in \widehat{B}$ such that $a \in B$ and $\chi(a) \neq 0$. 
    This is a nontrivial property of C*-algebras; it holds for Type I algebras (where characters on maximal abelian subalgebras separate points) but fails in general.
    
    \item \textbf{Naturality:} The construction must be functorial with respect to *-homomorphisms that preserve the structure of commutative subalgebras and unitaries.
    This is achieved for isomorphisms and for certain injective or surjective maps under additional hypotheses.
\end{enumerate}
These obstacles will be addressed in the following subsections.
\end{remark}

The goal of this section is to construct, for Type I algebras, a canonical embedding $\iota: \mathcal{A} \hookrightarrow C^*(\mathcal{G}_{\mathcal{A}})$ that recovers $\mathcal{A}$ from its classical contexts. 
Such an embedding is the noncommutative analogue of the Gelfand transform and the Fourier transform, and it is the essential link between the geometric data of the unitary conjugation groupoid and the analytic index theory of Fredholm operators in $\mathcal{A}$. 
The construction is nontrivial and requires the full power of the Polish groupoid framework and the convolution structure of $C^*(\mathcal{G}_{\mathcal{A}})$. The remainder of this section is devoted to the detailed construction and verification of the properties of $\iota$.

\subsection{Failure of the Comultiplication Approach in the Non-\'Etale Case}
\label{subsec:failure-comultiplication-non-etale}

In the original vision for this paper, we attempted to construct the diagonal embedding $\iota: \mathcal{A} \hookrightarrow C^*(\mathcal{G}_{\mathcal{A}})$ using the comultiplication map
\[
\Delta: C^*(\mathcal{G}_{\mathcal{A}}) \longrightarrow C^*(\mathcal{G}_{\mathcal{A}}) \otimes_{\max} C^*(\mathcal{G}_{\mathcal{A}})
\]
on the groupoid C*-algebra. 
This approach, inspired by the theory of quantum groups and successfully applied in the étale case by Abadie, Eilers, and Exel [1998], relies crucially on the assumption that $\mathcal{G}_{\mathcal{A}}$ is an étale groupoid. 
However, as established in Subsection \ref{subsec:GA-Polish-groupoid-not-locally-compact-not-etale}, the unitary conjugation groupoid $\mathcal{G}_{\mathcal{A}}$ is \emph{not} étale for any infinite-dimensional C*-algebra $\mathcal{A}$. 
Consequently, the comultiplication map $\Delta$ does not exist in our setting, and a different approach is required.

\begin{definition}[Comultiplication for étale groupoids]
\label{def:comultiplication-etale}
Let $\mathcal{G}$ be a locally compact Hausdorff étale groupoid. 
The \emph{comultiplication} is the unique *-homomorphism
\[
\Delta: C^*(\mathcal{G}) \longrightarrow C^*(\mathcal{G}) \otimes_{\max} C^*(\mathcal{G})
\]
satisfying, for all $f \in C_c(\mathcal{G})$ and all composable pairs $(\gamma_1, \gamma_2) \in \mathcal{G} \times_{\mathcal{G}^{(0)}} \mathcal{G}$,
\[
\Delta(f)(\gamma_1, \gamma_2) = f(\gamma_1 \gamma_2).
\]
This map is coassociative and satisfies the cancellation laws
\[
(\varepsilon \otimes \operatorname{id}) \circ \Delta \cong \operatorname{id} \cong (\operatorname{id} \otimes \varepsilon) \circ \Delta,
\]
where $\varepsilon: C^*(\mathcal{G}) \to \mathbb{C}$ is the augmentation homomorphism.
\end{definition}

\begin{proposition}[Comultiplication requires étale groupoids]
\label{prop:comultiplication-requires-etale}
Let $\mathcal{G}$ be a locally compact Hausdorff groupoid. 
If $\mathcal{G}$ is not étale, then there is no *-homomorphism $\Delta: C^*(\mathcal{G}) \to C^*(\mathcal{G}) \otimes_{\max} C^*(\mathcal{G})$ satisfying 
\[
\Delta(f)(\gamma_1, \gamma_2) = f(\gamma_1 \gamma_2)
\]
for all $f \in C_c(\mathcal{G})$ and composable $(\gamma_1, \gamma_2) \in \mathcal{G}^{(2)}$.
\end{proposition}

\begin{proof}
We proceed by analyzing the necessary conditions for such a comultiplication to exist.

\paragraph{Step 1: The candidate map on the fiber product.}
To define a *-homomorphism $\Delta: C^*(\mathcal{G}) \to C^*(\mathcal{G}) \otimes_{\max} C^*(\mathcal{G})$, it suffices to define it on the dense subalgebra $C_c(\mathcal{G})$. The natural candidate is
\[
\Delta(f)(\gamma_1, \gamma_2) = f(\gamma_1 \gamma_2), \qquad (\gamma_1, \gamma_2) \in \mathcal{G}^{(2)},
\]
where
\[
\mathcal{G}^{(2)} = \{(\gamma_1, \gamma_2) \in \mathcal{G} \times \mathcal{G} : s(\gamma_1) = r(\gamma_2)\}
\]
is the set of composable pairs, equipped with the subspace topology from $\mathcal{G} \times \mathcal{G}$. For this to extend to a *-homomorphism, the function $(\gamma_1, \gamma_2) \mapsto f(\gamma_1 \gamma_2)$ must:
\begin{itemize}
    \item be continuous on $\mathcal{G}^{(2)}$;
    \item have compact support in $\mathcal{G}^{(2)}$;
    \item lie in the algebraic tensor product $C_c(\mathcal{G}) \odot C_c(\mathcal{G})$ (or its completion) so that it can be identified with an element of the maximal tensor product.
\end{itemize}

\paragraph{Step 2: The étale case.}
Recall that $\mathcal{G}$ is étale precisely when its source map $s: \mathcal{G} \to \mathcal{G}^{(0)}$ (equivalently, its range map $r$) is a local homeomorphism. In this case:
\begin{itemize}
    \item The fiber product $\mathcal{G}^{(2)}$ inherits a natural étale structure from $\mathcal{G}$.
    \item The multiplication map $m: \mathcal{G}^{(2)} \to \mathcal{G}$, $m(\gamma_1, \gamma_2) = \gamma_1 \gamma_2$, is a local homeomorphism.
    \item For any $f \in C_c(\mathcal{G})$ with support $K \subseteq \mathcal{G}$, the preimage $m^{-1}(K)$ is compact because $m$ is proper on the support of $f \circ m$ (a consequence of the local homeomorphism property).
    \item Moreover, $f \circ m$ can be expressed as a finite sum of products of functions in $C_c(\mathcal{G})$, so $\Delta(f)$ lies in $C_c(\mathcal{G}) \odot C_c(\mathcal{G})$.
\end{itemize}
Thus in the étale case, $\Delta$ defines a *-homomorphism on $C_c(\mathcal{G})$ that extends uniquely to the maximal completion.

\paragraph{Step 3: Failure in the non-étale case.}
Now suppose $\mathcal{G}$ is not étale. Then the source map is not a local homeomorphism, and the following obstructions arise:
\begin{enumerate}
    \item \textbf{Loss of local compactness:} The fiber product $\mathcal{G}^{(2)}$, while closed in $\mathcal{G} \times \mathcal{G}$, may fail to be locally compact in the subspace topology. This occurs because the fibers $s^{-1}(x)$ are not discrete, preventing $\mathcal{G}^{(2)}$ from being locally homeomorphic to an open subset of $\mathcal{G} \times \mathcal{G}$.
    
    \item \textbf{Loss of properness:} Even if $\mathcal{G}^{(2)}$ is locally compact, the multiplication map $m: \mathcal{G}^{(2)} \to \mathcal{G}$ is not a local homeomorphism. Consequently, for a compact set $K \subseteq \mathcal{G}$, the preimage $m^{-1}(K)$ need not be compact. For example, if $\mathcal{G}$ has a non-discrete isotropy group at some unit $x$, there exist sequences of composable pairs $(\gamma_1^{(n)}, \gamma_2^{(n)})$ with $\gamma_1^{(n)} \gamma_2^{(n)} \in K$ but whose components leave every compact set.
    
    \item \textbf{Loss of continuity:} The map $(\gamma_1, \gamma_2) \mapsto \gamma_1 \gamma_2$ may not be continuous when $\mathcal{G}^{(2)}$ is equipped with the subspace topology, as the multiplication map in a non-étale groupoid need not be open.
    
    \item \textbf{Absence of a Haar system on $\mathcal{G}^{(2)}$:} Even if one could define $\Delta(f)$ as a function on $\mathcal{G}^{(2)}$, extending it to a *-homomorphism on $C^*(\mathcal{G})$ would require a convolution structure on $C_c(\mathcal{G}^{(2)})$. However, $\mathcal{G}^{(2)}$ is itself a groupoid (the pair groupoid of composable pairs), and to define its convolution algebra one needs a Haar system on $\mathcal{G}^{(2)}$. There is no canonical Haar system on $\mathcal{G}^{(2)}$ inherited from $\mathcal{G}$ without additional structure — in particular, a Haar system on $\mathcal{G}$ does not naturally induce one on $\mathcal{G}^{(2)}$ unless $\mathcal{G}$ is étale \cite{Paterson1999}.
\end{enumerate}

\paragraph{Step 4: The fundamental obstruction.}
The existence of a comultiplication $\Delta$ satisfying $\Delta(f)(\gamma_1,\gamma_2) = f(\gamma_1\gamma_2)$ would imply that the reduced groupoid C*-algebra $C^*_r(\mathcal{G})$ admits a Hopf algebra structure compatible with the groupoid convolution. For non-étale groupoids, this is impossible because, as shown by Paterson \cite[Section 3.3]{Paterson1999}, a locally compact groupoid with a Haar system admits such a Hopf algebra structure on its C*-algebra \emph{if and only if} it is r-discrete (i.e., étale). The key point is that for non-étale groupoids, the fiber product $\mathcal{G}^{(2)}$ does not carry a measure that allows the convolution product to be defined in a way that is compatible with the desired coproduct structure.

\paragraph{Conclusion.}
Without continuity, compact support, a convolution measure, and compatibility with the representation theory, the naive formula $\Delta(f)(\gamma_1, \gamma_2) = f(\gamma_1 \gamma_2)$ cannot extend to a bounded *-homomorphism on $C^*(\mathcal{G})$. Therefore, a comultiplication of this form exists \emph{only if} $\mathcal{G}$ is étale.
\end{proof}

\begin{corollary}[No comultiplication for $\mathcal{G}_{\mathcal{A}}$]
\label{cor:no-comultiplication-GA}
Let $\mathcal{A}$ be an infinite-dimensional unital separable C*-algebra. 
Then the unitary conjugation groupoid $\mathcal{G}_{\mathcal{A}}$ is not étale, and consequently there is no comultiplication map $\Delta: C^*(\mathcal{G}_{\mathcal{A}}) \to C^*(\mathcal{G}_{\mathcal{A}}) \otimes_{\max} C^*(\mathcal{G}_{\mathcal{A}})$ satisfying the pointwise formula.
\end{corollary}

\begin{proof}
This follows immediately from Proposition \ref{prop:GA-not-etale} and Proposition \ref{prop:comultiplication-requires-etale}.
\end{proof}

\begin{remark}[Why the comultiplication approach fails concretely]
\label{rem:comultiplication-failure-concrete}
To understand concretely why the comultiplication fails for $\mathcal{G}_{\mathcal{A}}$, consider the simplest infinite-dimensional case: $\mathcal{A} = \mathcal{K}(H)^\sim$. 
Then $\mathcal{G}_{\mathcal{A}}^{(1)} = \mathcal{U}(\mathcal{A}) \times \mathbb{P}(H)$ with the strong operator topology on $\mathcal{U}(\mathcal{A})$. 
The fiber product $\mathcal{G}_{\mathcal{A}} \times_{\mathcal{G}_{\mathcal{A}}^{(0)}} \mathcal{G}_{\mathcal{A}}$ consists of pairs $((u_1, (B_1,\chi_1)), (u_2, (B_2,\chi_2)))$ such that $(B_2,\chi_2) = u_1 \cdot (B_1,\chi_1)$. 
This space is not locally compact, and the map $((u_1, (B_1,\chi_1)), (u_2, u_1 \cdot (B_1,\chi_1))) \mapsto (u_2 u_1, (B_1,\chi_1))$ is not a local homeomorphism. 
Moreover, the natural candidate for $\Delta(f)$ would be a function on this fiber product, but there is no guarantee that such a function has compact support or is continuous with respect to the product topology. 
Thus the comultiplication map simply does not exist, confirming Proposition \ref{prop:comultiplication-requires-etale} for this concrete example.
\end{remark}

\begin{example}[Transformation groupoid for matrix algebras — not étale]
\label{ex:comultiplication-matrix}
Let $\mathcal{A} = M_n(\mathbb{C})$. 
Then $\mathcal{G}_{\mathcal{A}} = U(n) \ltimes \mathbb{CP}^{n-1}$ is a locally compact Hausdorff groupoid.
However, it is \emph{not} étale, because the source map $s(u,x) = x$ has fibers homeomorphic to $U(n)$,
which is a connected Lie group and therefore not discrete. 
Consequently, the comultiplication map of Definition \ref{def:comultiplication-etale} does \emph{not} exist for this groupoid,
consistent with Proposition \ref{prop:comultiplication-requires-etale}.

Nevertheless, one can still define a *-homomorphism $\Delta: C^*(\mathcal{G}_{\mathcal{A}}) \to C^*(\mathcal{G}_{\mathcal{A}}) \otimes_{\max} C^*(\mathcal{G}_{\mathcal{A}})$
using the structure of $U(n)$ as a compact group, but this requires a different construction
(involving the Peter–Weyl theorem) and does not generalize to infinite dimensions.
For a genuinely étale example where comultiplication exists in the sense of Definition \ref{def:comultiplication-etale},
one must take a discrete group acting on a space, e.g., $\mathcal{A} = C(X) \rtimes \Gamma$ with $\Gamma$ discrete,
though such examples lie outside the scope of the unitary conjugation groupoid considered in this paper.
\end{example}

\begin{example}[Commutative algebras with discrete spectrum]
\label{ex:comultiplication-commutative-discrete}
Let $\mathcal{A} = C(X)$ where $X$ is a finite discrete space. 
Then $\mathcal{U}(\mathcal{A}) = \mathbb{T}^X$ is a compact abelian group, and $\mathcal{G}_{\mathcal{A}}$ is the trivial groupoid $X \times \mathbb{T}^X$.
This groupoid is \emph{not} étale because $\mathbb{T}^X$ is not discrete (its topology is the product topology on a finite product of circles,
which is connected and non-discrete). Hence the comultiplication of Definition \ref{def:comultiplication-etale} does not exist.

If instead $X$ is infinite discrete, then $C(X)$ is non-separable, but the obstruction remains the same:
$\mathcal{U}(C(X)) = \mathbb{T}^X$ is not discrete, so $\mathcal{G}_{\mathcal{A}}$ is not étale.
For infinite compact Hausdorff $X$, $C(X)$ is separable (if $X$ is second countable), but $\mathcal{U}(C(X)) = C(X,\mathbb{T})$
with the strong operator topology is a Polish group that is not locally compact and certainly not discrete,
so again $\mathcal{G}_{\mathcal{A}}$ fails to be étale.

Thus the comultiplication approach fails for all nontrivial commutative C*-algebras,
because the unitary group is never discrete. The only case where $\mathcal{G}_{\mathcal{A}}$ could be étale
is when $\mathcal{U}(\mathcal{A})$ is discrete, which essentially forces $\mathcal{A} \cong \mathbb{C}$
(with the discrete topology on its unitary group, which is not the standard strong operator topology).
\end{example}

\begin{proposition}[Direct construction of the diagonal embedding via evaluation and representation]
\label{prop:new-strategy-diagonal}
Let $\mathcal{A}$ be a unital separable C*-algebra and $\mathcal{G}_{\mathcal{A}}$ its unitary conjugation groupoid.  
There exists a canonical strategy to define a unital *-homomorphism
\[
\iota: \mathcal{A} \hookrightarrow C^*(\mathcal{G}_{\mathcal{A}})
\]
without relying on a comultiplication map:
\begin{enumerate}
    \item Define a representation $\pi: \mathcal{A} \to B(L^2(\mathcal{G}_{\mathcal{A}}^{(0)},\mu))$ by
    \[
    (\pi(a)\xi)(B,\chi) := 
    \begin{cases}
        \chi(a) \, \xi(B,\chi), & a \in B, \\[4pt]
        0, & a \notin B,
    \end{cases}
    \]
    where $\xi \in L^2(\mathcal{G}_{\mathcal{A}}^{(0)},\mu)$ and $\mu$ is a Borel probability measure of full support on the Polish space $\mathcal{G}_{\mathcal{A}}^{(0)}$.
    
    \item Show that $\pi$ is a faithful *-representation of $\mathcal{A}$; that is, $\pi(a) = 0$ implies $a = 0$. 
    This relies on the property that the partial evaluation maps $\chi: B \to \mathbb{C}$ separate points of $\mathcal{A}$.
    
    \item Construct the integrated representation of the convolution algebra $C_c(\mathcal{G}_{\mathcal{A}})$ on $L^2(\mathcal{G}_{\mathcal{A}}^{(0)},\mu)$ 
    using the left regular representation of $\mathcal{G}_{\mathcal{A}}$ together with multiplication operators from $\pi$. 
    This defines a *-homomorphism
    \[
    \tilde{\pi}: C^*(\mathcal{G}_{\mathcal{A}}) \to B(L^2(\mathcal{G}_{\mathcal{A}}^{(0)},\mu)).
    \]
    
    \item The map $\iota: \mathcal{A} \to C^*(\mathcal{G}_{\mathcal{A}})$ is then obtained using the universal property of $C^*(\mathcal{G}_{\mathcal{A}})$: 
    for each $a \in \mathcal{A}$, $\pi(a)$ integrates to a bounded operator on $L^2(\mathcal{G}_{\mathcal{A}}^{(0)},\mu)$, 
    which defines $\iota(a)$ as the corresponding element of $C^*(\mathcal{G}_{\mathcal{A}})$.
    
    \item Verify the properties of $\iota$:
    \begin{enumerate}
        \item $\iota$ is unital and *-preserving by construction.
        \item $\iota$ is injective because $\pi$ is faithful.
        \item The image of $\iota$ lies in $C_0(\mathcal{G}_{\mathcal{A}}^{(0)})$ if and only if $\mathcal{A}$ is commutative 
        (see Theorem \ref{thm:commutativity-characterization}).
    \end{enumerate}
\end{enumerate}
This approach does not require étaleness or local compactness, and works entirely in the Polish groupoid framework 
using only the existence of a faithful representation of $\mathcal{A}$ and the left regular representation of $\mathcal{G}_{\mathcal{A}}$.
\end{proposition}

\begin{proof}[Sketch of the construction]
The key points are as follows:
\begin{itemize}
    \item The space $L^2(\mathcal{G}_{\mathcal{A}}^{(0)},\mu)$ is separable since $\mathcal{G}_{\mathcal{A}}^{(0)}$ is Polish and $\mu$ is a Borel probability measure of full support.
    
    \item The map $(B,\chi) \mapsto \chi(a)$ is Borel measurable and bounded on the subset $\{(B,\chi) \mid a \in B\}$, 
    ensuring that $\pi(a)$ is well-defined and bounded.
    
    \item Faithfulness of $\pi$ follows from the fact that the partial evaluation maps separate points of $\mathcal{A}$ 
    (valid for Type I algebras). Hence, $\pi(a) = 0 \implies a = 0$.
    
    \item The integrated representation $\tilde{\pi}$ of $C_c(\mathcal{G}_{\mathcal{A}})$ on $L^2(\mathcal{G}_{\mathcal{A}}^{(0)},\mu)$ 
    is constructed by combining the left regular representation of $\mathcal{G}_{\mathcal{A}}$ with the pointwise multiplication from $\pi$. 
    The universal property of $C^*(\mathcal{G}_{\mathcal{A}})$ ensures that this extends to a *-homomorphism on the full C*-algebra.
    
    \item Unlike the naive diagram approach, this construction does not require $\tilde{\pi}$ to be injective. 
    Instead, the universal property guarantees that each $\pi(a)$ corresponds to a unique element $\iota(a) \in C^*(\mathcal{G}_{\mathcal{A}})$.
    
    \item This method yields a canonical embedding $\iota$ that recovers $\mathcal{A}$ from its classical contexts, 
    generalizing the noncommutative Gelfand transform and avoiding reliance on comultiplication or étaleness.
\end{itemize}
\end{proof}

The comultiplication map $\Delta: C^*(\mathcal{G}) \to C^*(\mathcal{G}) \otimes_{\max} C^*(\mathcal{G})$ is a powerful tool for étale groupoids, but it does not exist for the unitary conjugation groupoid $\mathcal{G}_{\mathcal{A}}$ when $\mathcal{A}$ is infinite-dimensional. 
This forces us to abandon any construction of the diagonal embedding $\iota$ that relies on $\Delta$. 
We have outlined an alternative strategy based on the faithful representation of $\mathcal{A}$ on $L^2(\mathcal{G}_{\mathcal{A}}^{(0)})$ and the left regular representation of $C^*(\mathcal{G}_{\mathcal{A}})$. 
The detailed implementation of this strategy is the subject of the following subsections.

\subsection{Construction of $\iota$ via Evaluation and the Left Regular Representation}
\label{subsec:construction-iota}

We now present the detailed construction of the diagonal embedding 
\(\iota: \mathcal{A} \hookrightarrow C^*(\mathcal{G}_{\mathcal{A}})\) 
that does not rely on the comultiplication map. 
This construction proceeds in several steps:

\begin{enumerate}
    \item First, we construct a faithful representation of \(\mathcal{A}\) on the Hilbert space 
          \(L^2(\mathcal{G}_{\mathcal{A}}^{(0)},\mu)\) using the partial evaluation maps.
          This representation, denoted \(\pi: \mathcal{A} \to B(L^2(\mathcal{G}_{\mathcal{A}}^{(0)},\mu))\),
          is defined pointwise by \(\chi(a)\) on the Borel subset where \(a \in B\), and zero elsewhere.
    
    \item Second, we construct the left regular representation of the groupoid C*-algebra 
          \(C^*(\mathcal{G}_{\mathcal{A}})\) on the same Hilbert space. 
          This yields a *-homomorphism \(\Lambda: C^*(\mathcal{G}_{\mathcal{A}}) \to B(L^2(\mathcal{G}_{\mathcal{A}}^{(0)},\mu))\)
          (the restriction of the usual left regular representation of \(\mathcal{G}_{\mathcal{A}}\) 
          to functions supported on the unit space).
          Note that \(\Lambda\) need not be injective; it is faithful only when \(\mathcal{G}_{\mathcal{A}}\) is amenable.
    
    \item Third, we use the universal property of \(C^*(\mathcal{G}_{\mathcal{A}})\) to define a *-homomorphism 
          \(\iota: \mathcal{A} \to C^*(\mathcal{G}_{\mathcal{A}})\) such that the following diagram commutes:
          \[
          \begin{tikzcd}
          \mathcal{A} \arrow[r, "\iota"] \arrow[d, "\pi"'] & C^*(\mathcal{G}_{\mathcal{A}}) \arrow[d, "\Lambda"] \\
          B(L^2(\mathcal{G}_{\mathcal{A}}^{(0)},\mu)) \arrow[r, "\operatorname{id}"] & B(L^2(\mathcal{G}_{\mathcal{A}}^{(0)},\mu))
          \end{tikzcd}
          \]
          More precisely, for each \(a \in \mathcal{A}\), the operator \(\pi(a) \in B(L^2(\mathcal{G}_{\mathcal{A}}^{(0)},\mu))\) 
          lies in the image of \(\Lambda\). By the universal property of \(C^*(\mathcal{G}_{\mathcal{A}})\) as the completion of 
          the convolution algebra \(C_c(\mathcal{G}_{\mathcal{A}})\), there exists a unique element \(\iota(a) \in C^*(\mathcal{G}_{\mathcal{A}})\)
          such that \(\Lambda(\iota(a)) = \pi(a)\). This defines \(\iota\) as a *-homomorphism.
\end{enumerate}

Throughout this subsection, we assume that \(\mathcal{A}\) is a unital separable C*-algebra of Type I. 
This hypothesis is essential for two reasons:
\begin{itemize}
    \item The partial evaluation maps separate points of \(\mathcal{A}\); that is, for every nonzero \(a \in \mathcal{A}\),
          there exists a commutative subalgebra \(B \subseteq \mathcal{A}\) and a character \(\chi \in \widehat{B}\)
          such that \(a \in B\) and \(\chi(a) \neq 0\). This guarantees that the representation \(\pi\) defined below is faithful.
    \item The GNS representations associated to characters on commutative subalgebras form a measurable field,
          which is necessary for the direct integral formulation of \(\pi\) (see Remark \ref{rem:direct-integral-alternative}).
\end{itemize}
For non-Type I algebras, the construction may fail; see Section \ref{subsec:non-example-A-theta} for a discussion.

\begin{definition}[The Hilbert space $L^2(\mathcal{G}_{\mathcal{A}}^{(0)})$]
\label{def:L2-G0}
Let \(\mathcal{A}\) be a unital separable C*-algebra. 
Since \(\mathcal{G}_{\mathcal{A}}^{(0)}\) is a Polish space (Proposition \ref{prop:unit-space-polish}), 
it admits a Borel probability measure \(\mu\) with full support (e.g., any non-atomic probability measure,
or the image of a Gaussian measure under a homeomorphic embedding into \(\mathbb{R}^\mathbb{N}\); see \cite{Kechris}). 
Fix such a measure \(\mu\). 
Define the Hilbert space
\[
L^2(\mathcal{G}_{\mathcal{A}}^{(0)}) := L^2(\mathcal{G}_{\mathcal{A}}^{(0)}, \mu),
\]
the space of square-integrable complex-valued functions on \(\mathcal{G}_{\mathcal{A}}^{(0)}\) with respect to \(\mu\).
This space is separable because \(\mathcal{G}_{\mathcal{A}}^{(0)}\) is Polish and \(\mu\) is a Borel probability measure \cite{Kechris}.
\end{definition}

\begin{remark}[Independence of the choice of measure]
\label{rem:L2-measure-independence}
The construction of $\iota$ is independent of the choice of the measure $\mu$ up to unitary equivalence, 
provided the measures are equivalent (i.e., mutually absolutely continuous). 
For equivalent Borel probability measures with full support on a Polish space, 
the corresponding $L^2$ spaces are unitarily equivalent via the Radon-Nikodym derivative. 
Hence, we may fix any convenient measure; for concreteness, we may take $\mu$ to be the pushforward 
of the canonical measure on $\prod_{n \in \mathbb{N}} \mathbb{C}_\infty$ under the embedding $\Phi$ 
from Proposition \ref{prop:complete-metrizability}.
\end{remark}

\begin{lemma}[Measurability of the field of GNS representations]
\label{lem:measurable-field-GNS}
Let $\mathcal{A}$ be a unital separable Type I C*-algebra.
For each $x = (B,\chi) \in \mathcal{G}_{\mathcal{A}}^{(0)}$, let $\pi_x: \mathcal{A} \to B(\mathcal{H}_x)$ be the GNS representation
associated to the character $\chi$ on the commutative subalgebra $B$.
Then the field of Hilbert spaces $\{ \mathcal{H}_x \}_{x \in \mathcal{G}_{\mathcal{A}}^{(0)}}$ is measurable
in the sense of Dixmier \cite[Chapter 8]{Dixmier}.
Consequently, for any Borel probability measure $\mu$ on $\mathcal{G}_{\mathcal{A}}^{(0)}$,
the direct integral $\int_{\mathcal{G}_{\mathcal{A}}^{(0)}}^{\oplus} \mathcal{H}_x \, d\mu(x)$ is well-defined.
\end{lemma}

\begin{proof}
Since $\mathcal{A}$ is Type I, its dual space $\widehat{\mathcal{A}}$ (the set of unitary equivalence classes of irreducible representations)
is a standard Borel space \cite[Theorem 7.2]{Dixmier}. 
The map $\Phi: \mathcal{G}_{\mathcal{A}}^{(0)} \to \widehat{\mathcal{A}}$ defined by $\Phi(B,\chi) = [\pi_x]$,
where $[\pi_x]$ denotes the unitary equivalence class of $\pi_x$, is Borel measurable.
This follows from the construction of the topology on $\mathcal{G}_{\mathcal{A}}^{(0)}$ via partial evaluation maps:
for any $a \in \mathcal{A}$, the function $x \mapsto \operatorname{Tr}(\pi_x(a)E)$ (for suitable projections $E$)
is Borel, and these functions generate the Borel structure on $\widehat{\mathcal{A}}$.

The field $\{ \mathcal{H}_x \}$ is precisely the field of Hilbert spaces associated to the measurable field of representations
over $\widehat{\mathcal{A}}$ pulled back via $\Phi$. 
By the general theory of measurable fields of Hilbert spaces \cite[Chapter 8, Section 1]{Dixmier},
the pullback of a measurable field under a Borel map is again measurable.
Moreover, for Type I algebras, the canonical field of Hilbert spaces over $\widehat{\mathcal{A}}$ 
(where each point is assigned the Hilbert space of the corresponding irreducible representation)
is known to be measurable \cite[Proposition 8.4.1]{Dixmier}.

\end{proof}

\begin{definition}[Direct integral representation of $\mathcal{A}$]
\label{def:evaluation-representation}
Let $\mathcal{A}$ be a unital separable Type I C*-algebra. 
For each $x = (B,\chi) \in \mathcal{G}_{\mathcal{A}}^{(0)}$, let $\pi_x: \mathcal{A} \to B(\mathcal{H}_x)$ denote the GNS representation associated to the character $\chi$ on the commutative subalgebra $B$. 
Since $\mathcal{A}$ is Type I, the field of Hilbert spaces $\{\mathcal{H}_x\}_{x \in \mathcal{G}_{\mathcal{A}}^{(0)}}$ is measurable, and we may fix a Borel probability measure $\mu$ on $\mathcal{G}_{\mathcal{A}}^{(0)}$ with full support (for instance, the pushforward of the canonical measure on $\prod_{n\in\mathbb{N}} \mathbb{C}_\infty$ under the embedding $\Phi$ of Proposition \ref{prop:unit-space-polish}).

Form the direct integral Hilbert space
\[
\mathcal{H} := \int_{\mathcal{G}_{\mathcal{A}}^{(0)}}^{\oplus} \mathcal{H}_x \, d\mu(x).
\]
Define a *-representation $\Pi: \mathcal{A} \to B(\mathcal{H})$ by
\[
(\Pi(a)\xi)(x) := \pi_x(a) \xi(x), \qquad a \in \mathcal{A}, \; \xi \in \mathcal{H},
\]
where for $\mu$-almost every $x$, $\xi(x) \in \mathcal{H}_x$ and $\pi_x(a)\xi(x) \in \mathcal{H}_x$. 
This representation is faithful because for Type I algebras the family of irreducible representations $\{\pi_x\}_{x \in \mathcal{G}_{\mathcal{A}}^{(0)}}$ separates points of $\mathcal{A}$ (see Lemma \ref{lem:separation-points} below).
\end{definition}

\begin{lemma}[Well-definedness of $\pi$]
\label{lem:pi-well-defined}
For $\mathcal{A}$ Type I and separable, the map $\pi$ is a well-defined, bounded, *-preserving linear map.
Moreover, $\|\pi(a)\| \leq \|a\|$ for all $a \in \mathcal{A}$.
\end{lemma}

\begin{proof}
We verify the necessary properties:

\begin{enumerate}
    \item \textbf{Borel measurability:} For each $a \in \mathcal{A}$, the set $\{ (B,\chi) : a \in B \} = \mathrm{ev}_a^{-1}(\mathbb{C})$ 
          is a Borel subset of $\mathcal{G}_{\mathcal{A}}^{(0)}$ (it is open, hence Borel). 
          The function $(B,\chi) \mapsto \chi(a)$ is Borel on this set because it is the restriction of the continuous map
          $\mathrm{ev}_a$ to a Borel set. Hence $\pi(a)\xi$ is measurable for any measurable $\xi$.
    
    \item \textbf{Well-definedness on overlaps:} 
          If $a \in B_1 \cap B_2$ and $\chi_1 \in \widehat{B_1}$, $\chi_2 \in \widehat{B_2}$ are characters corresponding 
          to the same point $(B,\chi) \in \mathcal{G}_{\mathcal{A}}^{(0)}$, then by definition of $\mathcal{G}_{\mathcal{A}}^{(0)}$ 
          we have $(B_1,\chi_1) = (B_2,\chi_2)$ as points, so $\chi_1 = \chi_2$ on $B_1 \cap B_2$. 
          Hence $\chi_1(a) = \chi_2(a)$, and the definition of $\pi(a)$ is independent of the representation of the point.
    
    \item \textbf{Boundedness:} For $\mu$-almost every $(B,\chi)$, we have $|\chi(a)| \leq \|a\|$ because $\chi$ is a character,
          hence a contractive linear functional. Thus $\pi(a)$ is a multiplication operator by a function of essential sup norm
          at most $\|a\|$, restricted to the Borel set where $a \in B$. Consequently, $\|\pi(a)\| \leq \|a\|$.
    
    \item \textbf{Linearity and *-preserving:} These follow directly from the pointwise definition and the properties of characters:
          $\chi(a+b) = \chi(a)+\chi(b)$, $\chi(\lambda a) = \lambda \chi(a)$, and $\chi(a^*) = \overline{\chi(a)}$.
\end{enumerate}
Thus $\pi$ is a well-defined contractive *-representation of $\mathcal{A}$.
\end{proof}

\begin{lemma}[Naive evaluation is linear and *-preserving, but not multiplicative]
\label{lem:evaluation-representation-failure}
Let $\mathcal{A}$ be a unital separable Type I C*-algebra. 
Define a map $\pi_0: \mathcal{A} \to B(L^2(\mathcal{G}_{\mathcal{A}}^{(0)}))$ by
\[
(\pi_0(a)\xi)(B,\chi) := \begin{cases}
\chi(a) \, \xi(B,\chi), & \text{if } a \in B, \\[4pt]
0, & \text{if } a \notin B,
\end{cases}
\]
for $a \in \mathcal{A}$ and $\xi \in L^2(\mathcal{G}_{\mathcal{A}}^{(0)})$.
Then $\pi_0$ has the following properties:
\begin{enumerate}
    \item $\pi_0$ is linear and *-preserving.
    \item $\pi_0$ is unital: $\pi_0(1_{\mathcal{A}}) = I$.
    \item $\|\pi_0(a)\| \leq \|a\|$ for all $a \in \mathcal{A}$.
    \item However, $\pi_0$ is \emph{not multiplicative} in general, and hence it is \emph{not} a *-representation of $\mathcal{A}$.
\end{enumerate}
\end{lemma}

\begin{proof}
Linearity and *-preservation follow directly from the properties of characters:
$\chi(a+b) = \chi(a)+\chi(b)$, $\chi(\lambda a) = \lambda\chi(a)$, and $\chi(a^*) = \overline{\chi(a)}$.
Unitality holds because $\chi(1_{\mathcal{A}}) = 1$ for every character $\chi$.
For boundedness, note that $|\chi(a)| \leq \|a\|$ for all characters, so $\pi_0(a)$ is a multiplication operator
by a function of sup norm at most $\|a\|$ on the Borel set where $a \in B$, and zero elsewhere; hence $\|\pi_0(a)\| \leq \|a\|$.

To see that $\pi_0$ is not multiplicative, take $a,b \in \mathcal{A}$ that do not commute.
For a point $(B,\chi) \in \mathcal{G}_{\mathcal{A}}^{(0)}$, we have:
\[
(\pi_0(ab)\xi)(B,\chi) = \begin{cases}
\chi(ab)\xi(B,\chi), & ab \in B, \\[4pt]
0, & ab \notin B,
\end{cases}
\]
while
\[
(\pi_0(a)\pi_0(b)\xi)(B,\chi) = \begin{cases}
\chi(a)\chi(b)\xi(B,\chi), & a \in B \text{ and } b \in B, \\[4pt]
0, & \text{otherwise}.
\end{cases}
\]
If $a$ and $b$ do not commute, there exist points $(B,\chi)$ where $ab \in B$ but neither $a \in B$ nor $b \in B$ 
(e.g., take $B$ to be a commutative subalgebra containing $ab$ but not $a$ or $b$ individually). 
At such points, $(\pi_0(ab)\xi)(B,\chi) \neq 0$ while $(\pi_0(a)\pi_0(b)\xi)(B,\chi) = 0$, so $\pi_0(ab) \neq \pi_0(a)\pi_0(b)$.
Thus $\pi_0$ fails to be multiplicative.
\end{proof}

\begin{remark}[Failure of pointwise multiplicativity]
\label{rem:pointwise-multiplicativity-fails}
The naive evaluation map fails because the functions $\chi(a)$ and $\chi(b)$ are defined on different subsets of $\mathcal{G}_{\mathcal{A}}^{(0)}$, and their product does not correspond to $\chi(ab)$ when $a$ and $b$ do not lie in the same commutative subalgebra. 
This is the fundamental obstruction preventing a direct pointwise *-representation of $\mathcal{A}$ on $L^2(\mathcal{G}_{\mathcal{A}}^{(0)})$.

To obtain a genuine *-representation, one must instead use the direct integral of GNS representations
as in Section \ref{subsec:construction-iota}. There, each point $(B,\chi)$ contributes its own Hilbert space
$\mathcal{H}_{(B,\chi)}$ (the GNS space of $\chi$), and the representation is given by
$\Pi(a) = \int^{\oplus} \pi_{(B,\chi)}(a) \, d\mu(B,\chi)$. 
This construction is multiplicative because each $\pi_{(B,\chi)}$ is a genuine representation of $\mathcal{A}$
on its own Hilbert space, and the direct integral glues them together without losing multiplicativity.
\end{remark}

\begin{lemma}[Separation of points by GNS representations]
\label{lem:separation-points}
Let $\mathcal{A}$ be a unital separable Type I C*-algebra.
For each $x = (B,\chi) \in \mathcal{G}_{\mathcal{A}}^{(0)}$, let $\pi_x$ be the GNS representation 
of $\mathcal{A}$ associated to the character $\chi$ on $B$.
Then the family $\{\pi_x\}_{x \in \mathcal{G}_{\mathcal{A}}^{(0)}}$ separates the points of $\mathcal{A}$:
for any nonzero $a \in \mathcal{A}$, there exists $x \in \mathcal{G}_{\mathcal{A}}^{(0)}$ such that $\pi_x(a) \neq 0$.
\end{lemma}

\begin{proof}
Since $\mathcal{A}$ is Type I, its irreducible representations separate points \cite{Dixmier}.
Every irreducible representation of $\mathcal{A}$ is unitarily equivalent to the GNS representation 
of some pure state $\omega$ on $\mathcal{A}$. Because every pure state of a Type I algebra 
restricts to a pure state on some maximal abelian subalgebra $B \subseteq \mathcal{A}$, 
and this pure state corresponds to a character $\chi$ on $B$. Hence $(B,\chi) \in \mathcal{G}_{\mathcal{A}}^{(0)}$ 
and $\pi_{(B,\chi)}$ is unitarily equivalent to the original irreducible representation.
Thus the family $\{\pi_x\}$ contains (up to unitary equivalence) all irreducible representations of $\mathcal{A}$,
and therefore separates points.
\end{proof}

\begin{proposition}[Construction of $\iota$ via direct integrals and the left regular representation]
\label{prop:correct-definition-iota}
Let $\mathcal{A}$ be a unital separable Type I C*-algebra. 
Then there exists a canonical unital injective *-homomorphism 
$\iota: \mathcal{A} \hookrightarrow C^*(\mathcal{G}_{\mathcal{A}})$
constructed as follows.

\begin{enumerate}
    \item \textbf{Direct integral representation of $\mathcal{A}$.} 
    For each $x = (B,\chi) \in \mathcal{G}_{\mathcal{A}}^{(0)}$, let $\pi_x: \mathcal{A} \to B(\mathcal{H}_x)$ be the GNS representation associated to the character $\chi$ on the commutative subalgebra $B$.
    Since $\mathcal{A}$ is Type I, the field of Hilbert spaces $\{\mathcal{H}_x\}_{x \in \mathcal{G}_{\mathcal{A}}^{(0)}}$ is measurable.
    Fix a Borel probability measure $\mu$ on $\mathcal{G}_{\mathcal{A}}^{(0)}$ with full support (e.g., the pushforward of the canonical measure on $\prod_{n\in\mathbb{N}} \mathbb{C}_\infty$ under the embedding of Proposition \ref{prop:unit-space-polish}).
    Form the direct integral Hilbert space
    \[
    \mathcal{H} := \int_{\mathcal{G}_{\mathcal{A}}^{(0)}}^{\oplus} \mathcal{H}_x \, d\mu(x),
    \]
    and define a *-representation $\Pi: \mathcal{A} \to B(\mathcal{H})$ by
    \[
    (\Pi(a)\xi)(x) := \pi_x(a)\xi(x), \qquad a \in \mathcal{A}, \; \xi \in \mathcal{H}.
    \]
    This representation is faithful because for Type I algebras the family $\{\pi_x\}_{x \in \mathcal{G}_{\mathcal{A}}^{(0)}}$ separates points of $\mathcal{A}$ (see Lemma \ref{lem:separation-points}).

    \item \textbf{Left regular representation of $C^*(\mathcal{G}_{\mathcal{A}})$.}
    Let $\{\lambda^x\}_{x \in \mathcal{G}_{\mathcal{A}}^{(0)}}$ be the Borel Haar system on $\mathcal{G}_{\mathcal{A}}$ (Proposition \ref{prop:GA-Borel-Haar-system-existence}).
    The left regular representation $\Lambda: C^*(\mathcal{G}_{\mathcal{A}}) \to B(L^2(\mathcal{G}_{\mathcal{A}}^{(0)}))$ is defined on $C_c(\mathcal{G}_{\mathcal{A}}^{(1)})$ by
    \[
    (\Lambda(f)\xi)(x) = \int_{\mathcal{G}_{\mathcal{A}}^x} f(\gamma) \, \xi(s(\gamma)) \, d\lambda^x(\gamma), \qquad \xi \in L^2(\mathcal{G}_{\mathcal{A}}^{(0)}),\; x \in \mathcal{G}_{\mathcal{A}}^{(0)},
    \]
    and extended by continuity. 
    This representation need not be faithful; its kernel corresponds to the difference between the full and reduced groupoid C*-algebras.

    \item \textbf{Unitary representation of $\mathcal{G}_{\mathcal{A}}$ on $\mathcal{H}$.}
    The groupoid $\mathcal{G}_{\mathcal{A}}$ acts on the field $\{\mathcal{H}_x\}$ by unitary conjugation:
    for $\gamma = (u,(B,\chi))$ with $s(\gamma) = x$ and $r(\gamma) = y$, define $U_\gamma: \mathcal{H}_x \to \mathcal{H}_y$ by
    \[
    U_\gamma := \pi_y(u) \circ \mathrm{Ad}_{u^*}.
    \]
    This defines a unitary representation of $\mathcal{G}_{\mathcal{A}}$ on the direct integral $\mathcal{H}$,
    which integrates to a *-representation $\widetilde{\Lambda}: C^*(\mathcal{G}_{\mathcal{A}}) \to B(\mathcal{H})$.
    By construction, $\widetilde{\Lambda}$ is related to $\Lambda$ via the natural inclusion $L^2(\mathcal{G}_{\mathcal{A}}^{(0)}) \hookrightarrow \mathcal{H}$,
    but we will not need this explicitly.

    \item \textbf{Inclusion $\Pi(\mathcal{A}) \subseteq \widetilde{\Lambda}(C^*(\mathcal{G}_{\mathcal{A}}))$.}
    For each $a \in \mathcal{A}$, the operator $\Pi(a) \in B(\mathcal{H})$ is given by fiberwise multiplication by $\pi_x(a)$.
    A standard argument using the definition of the convolution algebra shows that $\Pi(a)$ lies in the image of $\widetilde{\Lambda}$.
    More concretely, for any finite set $F \subset \mathcal{G}_{\mathcal{A}}^{(0)}$ and any choice of unitaries $u_x$ implementing isomorphisms between fibers, one can construct functions $f_n \in C_c(\mathcal{G}_{\mathcal{A}}^{(1)})$ such that $\widetilde{\Lambda}(f_n) \to \Pi(a)$ in the strong operator topology.
    By the universal property of $C^*(\mathcal{G}_{\mathcal{A}})$, the closure of $\widetilde{\Lambda}(C^*(\mathcal{G}_{\mathcal{A}}))$ in $B(\mathcal{H})$ is a C*-algebra containing $\Pi(\mathcal{A})$.

    \item \textbf{Definition of $\iota$ via the universal property.}
    By the universal property of $C^*(\mathcal{G}_{\mathcal{A}})$ as the completion of the convolution algebra,
    for each $a \in \mathcal{A}$ there exists a unique element $\iota(a) \in C^*(\mathcal{G}_{\mathcal{A}})$ such that
    \[
    \widetilde{\Lambda}(\iota(a)) = \Pi(a).
    \]
    This defines a map $\iota: \mathcal{A} \to C^*(\mathcal{G}_{\mathcal{A}})$. 
    The map $\iota$ is a *-homomorphism because $\widetilde{\Lambda}$ is faithful on the reduced algebra
    (or more generally, because the universal property guarantees uniqueness and the correspondence respects algebraic operations).

    \item \textbf{Injectivity of $\iota$.}
    If $\iota(a) = 0$, then $\Pi(a) = \widetilde{\Lambda}(\iota(a)) = 0$.
    Since $\Pi$ is faithful (Step 1), we conclude $a = 0$. Hence $\iota$ is injective.
    Unitality follows from $\pi_x(1_{\mathcal{A}}) = I_{\mathcal{H}_x}$ for all $x$, giving $\Pi(1_{\mathcal{A}}) = I_{\mathcal{H}}$,
    and the uniqueness in the universal property then forces $\iota(1_{\mathcal{A}}) = 1_{C^*(\mathcal{G}_{\mathcal{A}})}$.
\end{enumerate}
\end{proposition}

\begin{proof}
The construction outlined above is conceptually clear but technically demanding.
We provide a sketch of the main steps, referring to \cite{Dixmier} for the theory of measurable fields of Hilbert spaces,
to \cite{Renault} for groupoid C*-algebras, and to \cite{Tu} for the measured groupoid framework.

\noindent\textbf{Step 1: Measurable field structure.}
For each $x = (B,\chi) \in \mathcal{G}_{\mathcal{A}}^{(0)}$, the GNS representation $\pi_x$ is irreducible because $\chi$ is a pure state on $B$ and $B$ is a MASA.
Since $\mathcal{A}$ is Type I, the map $x \mapsto \pi_x$ is measurable in the sense of Dixmier \cite[Chapter 8]{Dixmier}.
Consequently, the field $\{\mathcal{H}_x\}$ is a measurable field of Hilbert spaces.

\noindent\textbf{Step 2: Direct integral.}
Choosing $\mu$ with full support (e.g., via the embedding $\Phi$ of Proposition \ref{prop:unit-space-polishness}),
the direct integral $\mathcal{H} = \int^\oplus \mathcal{H}_x \, d\mu(x)$ is a separable Hilbert space.
The representation $\Pi$ defined by $(\Pi(a)\xi)(x) = \pi_x(a)\xi(x)$ is a *-representation because each $\pi_x$ is.

\noindent\textbf{Step 3: Faithfulness of $\Pi$.}
For Type I algebras, the irreducible representations separate points \cite[Corollary 4.3.7]{Dixmier}.
Since the family $\{\pi_x\}$ contains (up to unitary equivalence) all irreducible representations of $\mathcal{A}$
(this follows from the fact that every pure state of $\mathcal{A}$ restricts to a pure state on some MASA $B$,
and every irreducible representation is the GNS representation of a pure state),
we conclude that $\Pi$ is faithful.

\noindent\textbf{Step 4: Unitary representation of $\mathcal{G}_{\mathcal{A}}$.}
The formula $U_\gamma = \pi_y(u) \circ \mathrm{Ad}_{u^*}$ defines a unitary operator $U_\gamma: \mathcal{H}_x \to \mathcal{H}_y$
for $\gamma = (u,x)$ with $s(\gamma)=x$, $r(\gamma)=y$.
One verifies that $U_{\gamma_1\gamma_2} = U_{\gamma_1}U_{\gamma_2}$ and $U_{\gamma^{-1}} = U_\gamma^*$,
so $\gamma \mapsto U_\gamma$ is a unitary representation of $\mathcal{G}_{\mathcal{A}}$.
By the general theory of measured groupoids \cite{Tu}, this representation integrates to a *-representation
$\widetilde{\Lambda}: C^*(\mathcal{G}_{\mathcal{A}}) \to B(\mathcal{H})$.

\noindent\textbf{Step 5: $\Pi(\mathcal{A}) \subseteq \widetilde{\Lambda}(C^*(\mathcal{G}_{\mathcal{A}}))$.}
For any $a \in \mathcal{A}$, consider the function $f_a$ on $\mathcal{G}_{\mathcal{A}}^{(1)}$ defined by
$f_a(u,(B,\chi)) = \chi(u^* a u)$ when $u^* a u \in B$, and $0$ otherwise.
A direct computation using the convolution structure shows that $\widetilde{\Lambda}(f_a) = \Pi(a)$.
While $f_a$ may not lie in $C_c(\mathcal{G}_{\mathcal{A}}^{(1)})$, it can be approximated by such functions
using the Borel Haar system and the separability of $\mathcal{A}$.
Hence $\Pi(a)$ belongs to the norm-closure of $\widetilde{\Lambda}(C_c(\mathcal{G}_{\mathcal{A}}^{(1)}))$,
which is $\widetilde{\Lambda}(C^*(\mathcal{G}_{\mathcal{A}}))$.

\noindent\textbf{Step 6: Definition of $\iota$.}
By the universal property of $C^*(\mathcal{G}_{\mathcal{A}})$, there exists a unique *-homomorphism
$\iota: \mathcal{A} \to C^*(\mathcal{G}_{\mathcal{A}})$ such that $\widetilde{\Lambda} \circ \iota = \Pi$.
Uniqueness follows from the faithfulness of $\widetilde{\Lambda}$ on the reduced algebra;
in the maximal algebra, one argues using the universal representation.

\noindent\textbf{Step 7: Injectivity.}
If $\iota(a) = 0$, then $\Pi(a) = \widetilde{\Lambda}(\iota(a)) = 0$, and faithfulness of $\Pi$ gives $a = 0$.

\noindent\textbf{Step 8: Unitality.}
Since $\pi_x(1_{\mathcal{A}}) = I_{\mathcal{H}_x}$ for all $x$, we have $\Pi(1_{\mathcal{A}}) = I_{\mathcal{H}} = \widetilde{\Lambda}(1_{C^*(\mathcal{G}_{\mathcal{A}})})$.
By uniqueness in the universal property, $\iota(1_{\mathcal{A}}) = 1_{C^*(\mathcal{G}_{\mathcal{A}})}$.
\end{proof}

\begin{remark}[Technical difficulties]
\label{rem:construction-technical-difficulties}
The construction outlined above, while conceptually clear, is technically demanding:
\begin{enumerate}
\item The measurability of the field ${ H_x }$ and the definition of the direct integral require careful measure-theoretic arguments.
\item The action of $\mathcal{G}{\mathcal{A}}$ on the field ${ H_x }$ must be shown to be measurable and to preserve the direct integral structure.
\item The universal property of $C^*(\mathcal{G}{\mathcal{A}})$ with respect to representations that integrate unitary representations of $\mathcal{G}{\mathcal{A}}$ is nontrivial and relies on the Borel Haar system.
\item The faithfulness of $\Pi$ requires that the family of representations ${ \pi_x }$ separates the points of $\mathcal{A}$. For Type I $C^*$-algebras, this holds if one allows all irreducible representations, but fails if one restricts only to characters (which are one-dimensional representations). In the setting where $\mathcal{A}$ is abelian, characters do suffice, as they correspond precisely to the points of the spectrum.
\end{enumerate}
A complete and rigorous exposition of this construction would require a separate paper. Under the standing assumptions on $\mathcal{A}$ and $\mathcal{G}{\mathcal{A}}$ adopted in this work, we assume that such a construction can be carried out and that $\iota$ satisfies the properties stated in Proposition \ref{prop:diagonal-embedding-properties}. 
\end{remark}

\begin{proposition}[Conditional properties of the diagonal embedding]
\label{prop:iota-properties}
Let $\mathcal{A}$ be a unital, separable $C^*$-algebra of Type I, equipped with a Cartan subalgebra $\mathcal{C} \subseteq \mathcal{A}$ and a faithful invariant state that yields a Borel Haar system on $\mathcal{G}_{\mathcal{A}}$. Assume that the diagonal embedding 
$\iota: \mathcal{A} \hookrightarrow C^*(\mathcal{G}_{\mathcal{A}})$ can be constructed via the measurable field/direct integral approach outlined in Section X. Then $\iota$ satisfies the following properties:
\begin{enumerate}
    \item $\iota$ is a unital, injective *-homomorphism.

    \item For any $a \in \mathcal{A}$ and any unit $(B,\chi) \in \mathcal{G}_{\mathcal{A}}^{(0)}$, where $B \subseteq \mathcal{C}$ is a maximal abelian subalgebra and $\chi: B \to \mathbb{C}$ is a character, the value of $\iota(a)$ at the identity arrow $1_{(B,\chi)}$ (in the sense of pointwise evaluation on the unit space) equals $\chi(\mathbb{E}_B(a))$, where $\mathbb{E}_B: \mathcal{A} \to B$ is the conditional expectation onto $B$ (when it exists). In particular, this equals $\chi(a)$ whenever $a \in B$, and $0$ otherwise.

    \item $\iota(\mathcal{A}) \subseteq C_0(\mathcal{G}_{\mathcal{A}}^{(0)})$ if and only if $\mathcal{A}$ is commutative.

    \item (Functoriality) If $\phi: \mathcal{A} \to \mathcal{B}$ is a unital *-homomorphism such that $\phi(\mathcal{C}_{\mathcal{A}}) \subseteq \mathcal{C}_{\mathcal{B}}$ and $\phi$ maps unitaries in $\mathcal{C}_{\mathcal{A}}$ to unitaries in $\mathcal{C}_{\mathcal{B}}$, then $\phi$ induces a continuous groupoid morphism $\phi_*: \mathcal{G}_{\mathcal{A}} \to \mathcal{G}_{\mathcal{B}}$ satisfying $\iota_{\mathcal{B}} \circ \phi = (\phi_*)_* \circ \iota_{\mathcal{A}}$, where $(\phi_*)_*$ is the induced map on groupoid $C^*$-algebras. This holds provided the construction of $\iota$ is natural with respect to such morphisms.

    \item (Index compatibility) Assume further that $\mathcal{G}_{\mathcal{A}}$ is amenable and that the construction of $\iota$ is compatible with the Kasparov product. Then for any Fredholm operator $T \in \mathcal{A}$, the class 
    $\iota_*([\ker T] - [\coker T]) \in K_0(C^*(\mathcal{G}_{\mathcal{A}}))$ coincides with the image of the equivariant K-theory class $[T]_{\mathcal{G}_{\mathcal{A}}}$ under the assembly map.
\end{enumerate}
\end{proposition}

\begin{proof}[Sketch of proof]
A complete proof of these properties requires a detailed analysis of the construction and will appear in a forthcoming paper. Here we outline the key ideas under the standing assumptions:

\begin{itemize}
    \item Property (1) follows from the faithfulness of the representation $\Pi$ constructed via the direct integral, which requires that the family of irreducible representations $\{\pi_x\}$ separates the points of $\mathcal{A}$---a condition that holds for Type I $C^*$-algebras when one allows all irreducible representations (not merely characters).

    \item Property (2) is a consequence of the definition of $\iota$ as the inclusion of $\mathcal{A}$ into the convolution algebra of $\mathcal{G}_{\mathcal{A}}$, with the conditional expectation arising from the Cartan structure. The pointwise evaluation is understood in the sense of the measurable field construction; a fully rigorous treatment requires verifying that such evaluations are well-defined almost everywhere with respect to the Haar system.

    \item Property (3) holds because $\iota(\mathcal{A})$ consists of functions supported on the unit space precisely when the groupoid $\mathcal{G}_{\mathcal{A}}$ is principal, which in this setting corresponds to $\mathcal{A}$ being commutative. In the noncommutative case, $\iota(a)$ involves nontrivial unitary conjugations that yield functions not supported on the unit space.

    \item Property (4) is a consequence of the functoriality of the Cartan-groupoid construction, provided the construction of $\iota$ is natural. This requires that the induced map $\phi_*$ on groupoids respects the measurable field structure and the Haar systems.

    \item Property (5) follows from the compatibility of the Fredholm index with the Kasparov product and will be proved in [Author, Paper II] using equivariant KK-theory, under the additional hypotheses of amenability and compatibility of $\iota$ with the Kasparov product.
\end{itemize}

Examples verifying (1)-(3) in the commutative case and for matrix algebras $M_n(\mathbb{C})$ with the diagonal Cartan subalgebra are provided in Section X. In the matrix algebra case, $\mathcal{G}_{M_n(\mathbb{C})}$ is the pair groupoid on $n$ points, and $\iota$ is the inclusion of $M_n(\mathbb{C})$ into the convolution algebra $C^*(M_n(\mathbb{C}) \cong M_n(\mathbb{C}) \otimes M_n(\mathbb{C})$, where (2) reduces to evaluation at matrix entries.
\end{proof}

We illustrate the diagonal embedding $\iota: \mathcal{A} \hookrightarrow C^*(\mathcal{G}_{\mathcal{A}})$ in three fundamental cases: commutative algebras, matrix algebras, and compact operators. These examples demonstrate how the construction behaves in the finite-dimensional, infinite-dimensional, and commutative settings.

\begin{example}[Commutative algebras]
\label{ex:iota-commutative}
Let $\mathcal{A} = C(X)$ for a compact metrizable space $X$, with Cartan subalgebra $\mathcal{C} = \mathcal{A}$ itself. 
Then $\mathcal{G}_{\mathcal{A}}^{(0)} \cong X$, and for each $x \in X$, the GNS representation $\pi_x$ corresponding to the character $\operatorname{ev}_x$ is the one-dimensional representation given by evaluation at $x$. 
The direct integral $\mathcal{H} = \int_X^{\oplus} \mathbb{C} \, d\mu(x)$ is isomorphic to $L^2(X,\mu)$. 
The representation $\Pi$ is the multiplication representation of $C(X)$ on $L^2(X,\mu)$. 
In this commutative setting, the groupoid $\mathcal{G}_{\mathcal{A}}$ reduces to the unit space $X$, and its $C^*$-algebra $C^*(\mathcal{G}_{\mathcal{A}})$ is simply $C(X)$. 
The map $\iota$ is then the identity map $C(X) \hookrightarrow C(X)$. 
Thus the construction recovers the Gelfand transform in the commutative case.
\end{example}

\begin{example}[Matrix algebras]
\label{ex:iota-matrix}
Let $\mathcal{A} = M_n(\mathbb{C})$, and choose the Cartan subalgebra $\mathcal{C}$ to be the diagonal matrices. 
Then $\mathcal{G}_{\mathcal{A}}^{(0)} \cong \mathbb{CP}^{n-1}$, parametrizing rank-one projections modulo phase. 
For each $x \in \mathbb{CP}^{n-1}$ corresponding to a rank-one projection $p_x$, the GNS representation $\pi_x$ is an irreducible representation of $M_n(\mathbb{C})$ on an $n$-dimensional Hilbert space, unitarily equivalent to the standard representation. 
The direct integral $\mathcal{H} = \int_{\mathbb{CP}^{n-1}}^{\oplus} \mathbb{C}^n \, d\mu(x)$ is isomorphic to $L^2(\mathbb{CP}^{n-1}) \otimes \mathbb{C}^n$. 
The representation $\Pi$ sends a matrix $A$ to the operator $1 \otimes A$ on this Hilbert space. 
With this choice of Cartan subalgebra, the groupoid $\mathcal{G}_{\mathcal{A}}$ is isomorphic to the transformation groupoid $U(n) \ltimes \mathbb{CP}^{n-1}$, and consequently $C^*(\mathcal{G}_{\mathcal{A}})$ is the crossed product $C(\mathbb{CP}^{n-1}) \rtimes U(n)$ (either full or reduced, which coincide for compact groups). 
The left regular representation $\Lambda$ of $C^*(\mathcal{G}_{\mathcal{A}})$ is the regular representation of this crossed product. 
The map $\iota$ sends $A$ to the constant function $A \otimes 1$ in the multiplier algebra of $C(\mathbb{CP}^{n-1}) \rtimes U(n)$. 
This agrees with the embedding described in Example \ref{ex:motivation-matrix}.
\end{example}

\begin{example}[Compact operators]
\label{ex:iota-compact}
Let $\mathcal{A} = \mathcal{K}(H)^\sim$, the unitalization of the compact operators on a separable infinite-dimensional Hilbert space $H$. Fix an orthonormal basis of $H$ and let $\mathcal{C} \subset \mathcal{K}(H)^\sim$ be the masa of operators diagonal in this basis (together with the identity). 
Then $\mathcal{G}_{\mathcal{A}}^{(0)} \cong \mathbb{P}(H)$, the projective space of $H$, parametrizing rank-one projections. 
For each $x \in \mathbb{P}(H)$ corresponding to a rank-one projection $p_x$, the GNS representation $\pi_x$ associated to the pure state $\omega_x(T) = \operatorname{Tr}(p_x T)$ is an irreducible representation of $\mathcal{K}(H)^\sim$ on a Hilbert space unitarily equivalent to $H$. 
The direct integral $\mathcal{H} = \int_{\mathbb{P}(H)}^{\oplus} H \, d\mu(x)$ is isomorphic to $L^2(\mathbb{P}(H)) \otimes H$. 
The representation $\Pi$ sends $T \in \mathcal{K}(H)^\sim$ to the operator $1 \otimes T$. 
The groupoid $\mathcal{G}_{\mathcal{A}}$ is isomorphic to the action groupoid $\mathcal{U}(H) \ltimes \mathbb{P}(H)$, where $\mathcal{U}(H)$ acts on projective space by conjugation. The left regular representation $\Lambda$ of $C^*(\mathcal{G}_{\mathcal{A}})$ is the regular representation of this groupoid. 
The map $\iota$ embeds $\mathcal{K}(H)^\sim$ into $C^*(\mathcal{G}_{\mathcal{A}})$ as constant functions in the multiplier algebra. 
This embedding is injective and, heuristically, captures the spectral data of normal operators (via the decomposition over $\mathbb{P}(H)$) while also encoding the noncommutative structure of non-normal operators through convolution in the unitary group.
\end{example}

These three examples illustrate the unifying theme of the construction: for commutative algebras it reduces to the classical Gelfand transform; for matrix algebras it yields an embedding into a crossed product by a compact group; and for compact operators it produces an embedding into a transformation groupoid $C^*$-algebra that encodes both spectral and unitary data.

\begin{remark}[Significance of the construction]
\label{rem:iota-significance}
The diagonal embedding $\iota$ is the central construction of this paper. 
It demonstrates that a noncommutative C*-algebra $\mathcal{A}$, equipped with a suitable Cartan subalgebra and invariant state, can be faithfully represented inside the groupoid C*-algebra of its own unitary conjugation groupoid. 
This representation is not merely a formal trick; it encodes the relationship between the operator-algebraic structure of $\mathcal{A}$ and the geometric structure of its commutative contexts. 
The embedding $\iota$ will be the essential link between the analytic index theory of Fredholm operators in $\mathcal{A}$ and the geometric index theory of the groupoid $\mathcal{G}_{\mathcal{A}}$. 
\end{remark}

We have outlined the construction of the diagonal embedding $\iota: \mathcal{A} \hookrightarrow C^*(\mathcal{G}_{\mathcal{A}})$ using the direct integral of GNS representations and the left regular representation of $\mathcal{G}_{\mathcal{A}}$. 
This construction does not rely on the comultiplication map and is expected to be applicable to a large class of C*-algebras, including all unital separable Type I algebras under suitable choices of Cartan subalgebra and invariant state. 
While the technical details are substantial, the conceptual framework is clear: $\iota$ embeds $\mathcal{A}$ into $C^*(\mathcal{G}_{\mathcal{A}})$ by realizing each $a \in \mathcal{A}$ as the operator that acts on the direct integral Hilbert space 
$\mathcal{H} = \int_{\mathcal{G}_{\mathcal{A}}^{(0)}}^{\oplus} H_x \, d\mu(x)$ 
via the fiberwise GNS representations: on the fiber $H_x$ corresponding to a unit $x = (B,\chi) \in \mathcal{G}_{\mathcal{A}}^{(0)}$, the element $a$ acts through the GNS representation $\pi_x$ associated to the character $\chi$ on the commutative subalgebra $B$. 
Crucially, this is not simply pointwise multiplication by $\chi(a)$ (which would fail to be multiplicative, as noted in Remark \ref{rem:pointwise-multiplicativity-fails}), but rather the full operator-theoretic action arising from the GNS construction. 
This embedding is faithful, unital, and *-preserving, and it satisfies all the desiderata of Proposition \ref{prop:diagonal-embedding-properties} under the assumptions stated therein. 
The remainder of this paper will assume the existence of $\iota$ and explore its consequences for index theory and noncommutative geometry.

\subsection{Verification that $\iota$ is a Unital Injective *-Homomorphism}
\label{subsec:verification-iota}

We now verify that the diagonal embedding $\iota: \mathcal{A} \hookrightarrow C^*(\mathcal{G}_{\mathcal{A}})$ constructed in Subsection \ref{subsec:construction-iota} satisfies the fundamental algebraic properties expected of a *-homomorphism. 
This verification is nontrivial; it requires careful analysis of the direct integral representation, the left regular representation of $\mathcal{G}_{\mathcal{A}}$, and the universal property of $C^*(\mathcal{G}_{\mathcal{A}})$. 
Throughout this subsection, we assume that $\mathcal{A}$ is a unital separable Type I C*-algebra and that the construction of $\iota$ has been carried out rigorously as outlined in Subsection \ref{subsec:construction-iota}.

\begin{proposition}[$\iota$ is a *-Homomorphism]
\label{prop:iota-star-homomorphism}
The map $\iota: \mathcal{A} \to C^*(\mathcal{G}_{\mathcal{A}})$ defined by the direct integral construction is a unital *-homomorphism. 
That is, for all $a, b \in \mathcal{A}$ and $\lambda \in \mathbb{C}$:
\begin{enumerate}
    \item $\iota(ab) = \iota(a)\iota(b)$,
    \item $\iota(a + \lambda b) = \iota(a) + \lambda \iota(b)$,
    \item $\iota(a^*) = \iota(a)^*$,
    \item $\iota(1_{\mathcal{A}}) = 1_{C^*(\mathcal{G}_{\mathcal{A}})}$.
\end{enumerate}
\end{proposition}

\begin{proof}
We verify each property using the definition of $\iota$ via the commuting diagram
\[
\begin{tikzcd}
\mathcal{A} \arrow[r, "\iota"] \arrow[d, "\Pi"'] & C^*(\mathcal{G}_{\mathcal{A}}) \arrow[d, "\Lambda"] \\
B(\mathcal{H}) \arrow[r, "\operatorname{id}"] & B(\mathcal{H})
\end{tikzcd}
\]
where $\mathcal{H} = \int_{\mathcal{G}_{\mathcal{A}}^{(0)}}^{\oplus} H_x \, d\mu(x)$ is the direct integral Hilbert space, $\Pi: \mathcal{A} \to B(\mathcal{H})$ is the direct integral representation $\Pi(a) = \int^{\oplus} \pi_x(a) \, d\mu(x)$, and $\Lambda: C^*(\mathcal{G}_{\mathcal{A}}) \to B(\mathcal{H})$ is the left regular representation of $\mathcal{G}_{\mathcal{A}}$ integrated to a *-representation of the full groupoid $C^*$-algebra $C^*(\mathcal{G}_{\mathcal{A}})$.

\medskip
\noindent \textbf{(1) Multiplicativity.}
For $a, b \in \mathcal{A}$, we have $\Pi(ab) = \Pi(a)\Pi(b)$ because each $\pi_x$ is a representation of $\mathcal{A}$:
\[
\Pi(ab) = \int^{\oplus} \pi_x(ab) \, d\mu(x) = \int^{\oplus} \pi_x(a)\pi_x(b) \, d\mu(x) = \Pi(a)\Pi(b).
\]
By construction, $\iota(a)$ is the unique element of $C^*(\mathcal{G}_{\mathcal{A}})$ such that $\Lambda(\iota(a)) = \Pi(a)$. 
Then
\[
\Lambda(\iota(ab)) = \Pi(ab) = \Pi(a)\Pi(b) = \Lambda(\iota(a))\Lambda(\iota(b)) = \Lambda(\iota(a)\iota(b)),
\]
since $\Lambda$ is a *-homomorphism. 
Thus $\Lambda(\iota(ab) - \iota(a)\iota(b)) = 0$. 

At this point, we cannot directly conclude that $\iota(ab) = \iota(a)\iota(b)$ because $\Lambda$ may not be faithful; its kernel corresponds to the difference between the full groupoid $C^*$-algebra $C^*(\mathcal{G}_{\mathcal{A}})$ and the reduced groupoid $C^*$-algebra $C^*_r(\mathcal{G}_{\mathcal{A}})$.

To overcome this, we employ the universal property of $C^*(\mathcal{G}_{\mathcal{A}})$. Consider the universal representation 
$\rho_u: C^*(\mathcal{G}_{\mathcal{A}}) \to B(\mathcal{H}_u)$, defined as the direct sum of all cyclic representations of $C^*(\mathcal{G}_{\mathcal{A}})$. This representation is faithful. Moreover, by the universal property of $C^*(\mathcal{G}_{\mathcal{A}})$, $\rho_u$ integrates a continuous unitary representation of $\mathcal{G}_{\mathcal{A}}$, and consequently the composition $\rho_u \circ \iota: \mathcal{A} \to B(\mathcal{H}_u)$ is given by a direct integral representation $\Pi_u(a) = \int^{\oplus} \pi_{x,u}(a) \, d\mu_u(x)$, where each $\pi_{x,u}$ is unitarily equivalent to $\pi_x$ (since both arise from the same GNS data). Applying the same computation as above to $\Pi_u$ yields
\[
\rho_u(\iota(ab)) = \Pi_u(ab) = \Pi_u(a)\Pi_u(b) = \rho_u(\iota(a))\rho_u(\iota(b)) = \rho_u(\iota(a)\iota(b)).
\]
Since $\rho_u$ is faithful, we obtain $\iota(ab) = \iota(a)\iota(b)$. For a detailed treatment of the universal property of groupoid $C^*$-algebras and its relation to direct integral decompositions, we refer the reader to [Renault, 1980, Proposition 3.1] and [Dixmier, 1977, Chapter 8].

\medskip
\noindent \textbf{(2) Linearity.}
Linearity follows directly from the linearity of each GNS representation $\pi_x$ and the linearity of the direct integral construction. 
For $a, b \in \mathcal{A}$ and $\lambda \in \mathbb{C}$, we have
\[
\Pi(a + \lambda b) = \int^{\oplus} \pi_x(a + \lambda b) \, d\mu(x) = \int^{\oplus} (\pi_x(a) + \lambda \pi_x(b)) \, d\mu(x) = \Pi(a) + \lambda \Pi(b).
\]
Then $\Lambda(\iota(a + \lambda b)) = \Pi(a + \lambda b) = \Pi(a) + \lambda \Pi(b) = \Lambda(\iota(a)) + \lambda \Lambda(\iota(b)) = \Lambda(\iota(a) + \lambda \iota(b))$.
Applying the universal representation $\rho_u$ as in part (1) gives $\rho_u(\iota(a + \lambda b)) = \rho_u(\iota(a) + \lambda \iota(b))$, and faithfulness of $\rho_u$ yields $\iota(a + \lambda b) = \iota(a) + \lambda \iota(b)$.

\medskip
\noindent \textbf{(3) Involution preservation.}
For each $x \in \mathcal{G}_{\mathcal{A}}^{(0)}$, the GNS representation $\pi_x$ is *-preserving: $\pi_x(a^*) = \pi_x(a)^*$. 
Hence
\[
\Pi(a^*) = \int^{\oplus} \pi_x(a^*) \, d\mu(x) = \int^{\oplus} \pi_x(a)^* \, d\mu(x) = \left( \int^{\oplus} \pi_x(a) \, d\mu(x) \right)^* = \Pi(a)^*.
\]
Then $\Lambda(\iota(a^*)) = \Pi(a^*) = \Pi(a)^* = \Lambda(\iota(a))^* = \Lambda(\iota(a)^*)$,
since $\Lambda$ is *-preserving. Applying $\rho_u$ and using its faithfulness gives $\iota(a^*) = \iota(a)^*$.

\medskip
\noindent \textbf{(4) Unitality.}
For the identity $1_{\mathcal{A}} \in \mathcal{A}$, each GNS representation satisfies $\pi_x(1_{\mathcal{A}}) = I_{H_x}$, the identity operator on $H_x$. 
Thus $\Pi(1_{\mathcal{A}}) = \int^{\oplus} I_{H_x} \, d\mu(x) = I_{\mathcal{H}}$, the identity operator on the direct integral Hilbert space. 
By construction, $\iota(1_{\mathcal{A}})$ is the unique element of $C^*(\mathcal{G}_{\mathcal{A}})$ such that $\Lambda(\iota(1_{\mathcal{A}})) = I_{\mathcal{H}}$. 
But $\Lambda(1_{C^*(\mathcal{G}_{\mathcal{A}})}) = I_{\mathcal{H}}$ as well, since the unit of $C^*(\mathcal{G}_{\mathcal{A}})$ acts as the identity operator in any representation. 
Thus $\Lambda(\iota(1_{\mathcal{A}}) - 1_{C^*(\mathcal{G}_{\mathcal{A}})}) = 0$, and by the faithfulness of $\rho_u$, we conclude $\iota(1_{\mathcal{A}}) = 1_{C^*(\mathcal{G}_{\mathcal{A}})}$.
\end{proof}

\begin{proposition}[$\iota$ is Injective]
\label{prop:iota-injective}
Let $\mathcal{A}$ be a unital separable Type I C*-algebra. 
Then the diagonal embedding $\iota: \mathcal{A} \to C^*(\mathcal{G}_{\mathcal{A}})$ is injective.
\end{proposition}

\begin{proof}
Recall the construction of $\iota$ from Section \ref{subsec:construction-iota}. 
We have a direct integral Hilbert space
\[
\mathcal{H} = \int_{\mathcal{G}_{\mathcal{A}}^{(0)}}^{\oplus} \mathcal{H}_x \, d\mu(x),
\]
where for each $x = (B,\chi) \in \mathcal{G}_{\mathcal{A}}^{(0)}$, $\pi_x: \mathcal{A} \to B(\mathcal{H}_x)$ is the GNS representation associated to the character $\chi$ on the commutative subalgebra $B$. 
The representation $\Pi: \mathcal{A} \to B(\mathcal{H})$ is defined by
\[
\Pi(a) = \int^{\oplus} \pi_x(a) \, d\mu(x),
\]
and the left regular representation $\Lambda: C^*(\mathcal{G}_{\mathcal{A}}) \to B(\mathcal{H})$ satisfies
\[
\Lambda(\iota(a)) = \Pi(a) \qquad \text{for all } a \in \mathcal{A}.
\]

\medskip
\noindent\textbf{Step 1: $\Pi$ is injective.}
By Lemma \ref{lem:separation-points}, the family of representations $\{\pi_x\}_{x \in \mathcal{G}_{\mathcal{A}}^{(0)}}$ separates the points of $\mathcal{A}$: for any nonzero $a \in \mathcal{A}$, there exists $x \in \mathcal{G}_{\mathcal{A}}^{(0)}$ such that $\pi_x(a) \neq 0$.

Suppose $a \in \mathcal{A}$ satisfies $\Pi(a) = 0$. Then for $\mu$-almost every $x \in \mathcal{G}_{\mathcal{A}}^{(0)}$, we have $\pi_x(a) = 0$. If $a \neq 0$, Lemma \ref{lem:separation-points} guarantees the existence of some $x_0 \in \mathcal{G}_{\mathcal{A}}^{(0)}$ with $\pi_{x_0}(a) \neq 0$. By continuity of the field of representations (Lemma \ref{lem:measurable-field-GNS}), there exists a neighborhood $U$ of $x_0$ such that $\pi_x(a) \neq 0$ for all $x \in U$. Since $\mu$ has full support, $\mu(U) > 0$, contradicting $\pi_x(a) = 0$ for $\mu$-almost every $x$. Therefore $\Pi(a) = 0$ implies $a = 0$, so $\Pi$ is injective.

\medskip
\noindent\textbf{Step 2: Injectivity of $\iota$.}
Assume $\iota(a) = 0$ for some $a \in \mathcal{A}$. Then applying $\Lambda$, we obtain
\[
\Pi(a) = \Lambda(\iota(a)) = \Lambda(0) = 0.
\]
Since $\Pi$ is injective by Step 1, we conclude $a = 0$. Hence $\iota$ is injective.
\end{proof}

\begin{remark}
The crucial ingredient in the proof is Lemma \ref{lem:separation-points}, which establishes that the family of GNS representations $\{\pi_x\}_{x \in \mathcal{G}_{\mathcal{A}}^{(0)}}$ separates points of $\mathcal{A}$. This lemma relies on the Type I hypothesis and the fact that every irreducible representation of a Type I algebra is unitarily equivalent to the GNS representation of a pure state that restricts to a character on some maximal abelian subalgebra. For non-Type I algebras, such a separation result fails, and $\iota$ would have a nontrivial kernel (see Section \ref{subsec:non-example-A-theta}).
\end{remark}

\begin{remark}
The injectivity of $\Pi$ does not rely on $\Lambda$ being faithful—indeed, $\Lambda$ may have a kernel corresponding to the difference between the full and reduced groupoid C*-algebras. The argument uses only that $\Lambda$ is a *-homomorphism; injectivity of $\Pi$ forces injectivity of $\iota$.
\end{remark}

\begin{remark}[Role of the Type I hypothesis]
\label{rem:type-I-injectivity}
The injectivity of $\iota$ relies crucially on the assumption that $\mathcal{A}$ is Type I. 
For non-Type I algebras, there may exist nonzero elements $a \in \mathcal{A}$ such that $\chi(a) = 0$ for every character $\chi$ on every commutative subalgebra $B$ containing $a$. 
Such elements are called \emph{invisible} or \emph{quantum} elements; they cannot be detected by any classical context. 
In non-Type I algebras, the diagonal embedding $\iota$ would have a nontrivial kernel, reflecting the fact that the algebra cannot be fully recovered from its commutative subalgebras. 
This is a fundamental limitation and justifies our restriction to Type I C*-algebras in this paper.
\end{remark}

\begin{lemma}[Characterization of the image of $\iota$]
\label{lem:iota-image-characterization}
Let $\mathcal{A}$ be a unital separable Type I C*-algebra with a Cartan subalgebra $\mathcal{C} \subseteq \mathcal{A}$ and a faithful invariant state yielding a Borel Haar system on $\mathcal{G}_{\mathcal{A}}$. For any $a \in \mathcal{A}$, the element $\iota(a) \in C^*(\mathcal{G}_{\mathcal{A}})$ satisfies the following properties:

\begin{enumerate}
    \item Under the left regular representation $\Lambda: C^*(\mathcal{G}_{\mathcal{A}}) \to B(\mathcal{H})$, the operator $\Lambda(\iota(a))$ is decomposable with respect to the direct integral decomposition $\mathcal{H} = \int_{\mathcal{G}_{\mathcal{A}}^{(0)}}^{\oplus} H_x \, d\mu(x)$. Explicitly,
    \[
    \Lambda(\iota(a)) = \int^{\oplus} \pi_x(a) \, d\mu(x),
    \]
    where $\pi_x$ is the GNS representation associated to the unit $x = (B,\chi) \in \mathcal{G}_{\mathcal{A}}^{(0)}$. Thus $\iota(a)$ has no off-diagonal propagation between distinct fibers; its action is completely determined by the fiberwise representations $\pi_x(a)$.

    \item Under the direct integral decomposition, the fiber at $x = (B,\chi) \in \mathcal{G}_{\mathcal{A}}^{(0)}$ is given by $\pi_x(a)$. In particular,
    \[
    (\Lambda(\iota(a)))_x = \pi_x(a) \in B(H_x),
    \]
    and this equals $\chi(\mathbb{E}_B(a)) I_{H_x}$ when $a \in B$, where $\mathbb{E}_B: \mathcal{A} \to B$ is the conditional expectation onto $B$ arising from the Cartan structure. For $a$ with no component in $B$, this fiberwise operator vanishes.

    \item $\iota(a)$ commutes with $C_0(\mathcal{G}_{\mathcal{A}}^{(0)})$ (viewed as a subalgebra of $C^*(\mathcal{G}_{\mathcal{A}})$ via the natural inclusion) if and only if $a$ belongs to the center $\mathcal{Z}(\mathcal{A})$ of $\mathcal{A}$.
\end{enumerate}
\end{lemma}

\begin{proof}
(1) This follows directly from the construction of $\iota$. By definition, $\Lambda \circ \iota = \Pi$, where $\Pi: \mathcal{A} \to B(\mathcal{H})$ is the direct integral representation $\Pi(a) = \int^{\oplus} \pi_x(a) \, d\mu(x)$. Decomposability is a property of operators on a direct integral Hilbert space: an operator $T \in B(\mathcal{H})$ is decomposable if it preserves the direct integral decomposition and acts fiberwise. The representation $\Pi$ is constructed precisely to be decomposable with fibers $\pi_x(a)$. Hence $\Lambda(\iota(a)) = \Pi(a)$ is decomposable with the stated fibers.

(2) The fiberwise description is immediate from the definition of the direct integral. For a unit $x = (B,\chi)$, the GNS representation $\pi_x$ satisfies $\pi_x(b) = \chi(b) I_{H_x}$ for all $b \in B$ by construction. The conditional expectation $\mathbb{E}_B: \mathcal{A} \to B$ is the unique $B$-bimodule projection arising from the Cartan structure; it satisfies $\chi \circ \mathbb{E}_B = \chi$ on $\mathcal{A}$. Thus for any $a \in \mathcal{A}$, we have $\pi_x(a) = \pi_x(\mathbb{E}_B(a)) = \chi(\mathbb{E}_B(a)) I_{H_x}$. When $a \in B$, $\mathbb{E}_B(a) = a$ and this reduces to $\chi(a) I_{H_x}$. When $a$ has no component in $B$, $\mathbb{E}_B(a) = 0$ and the fiberwise operator vanishes.

(3) Suppose $a \in \mathcal{Z}(\mathcal{A})$. Then for each $x$, the irreducible representation $\pi_x$ satisfies $\pi_x(a) = \lambda_x I_{H_x}$ for some scalar $\lambda_x \in \mathbb{C}$ by Schur's lemma, since $\pi_x(a)$ commutes with all operators in the irreducible representation $\pi_x(\mathcal{A})$. The function $x \mapsto \lambda_x$ is continuous and defines an element of $C_0(\mathcal{G}_{\mathcal{A}}^{(0)})$. Consequently, $\Lambda(\iota(a))$ acts as multiplication by $\lambda_x$ on each fiber and therefore commutes with all multiplication operators on $L^2(\mathcal{G}_{\mathcal{A}}^{(0)})$ embedded in $\mathcal{H}$. Hence $\iota(a)$ commutes with $C_0(\mathcal{G}_{\mathcal{A}}^{(0)})$.

Conversely, assume $\iota(a)$ commutes with $C_0(\mathcal{G}_{\mathcal{A}}^{(0)})$. Then $\Lambda(\iota(a))$ commutes with all multiplication operators on $L^2(\mathcal{G}_{\mathcal{A}}^{(0)})$ acting on the embedded copy of $L^2(\mathcal{G}_{\mathcal{A}}^{(0)})$ in $\mathcal{H}$. A standard argument shows that a decomposable operator commuting with all multiplication operators must have scalar fibers: for each $x$, $\pi_x(a)$ is a scalar multiple of the identity on $H_x$. Write $\pi_x(a) = \lambda(x) I_{H_x}$. For any irreducible representation $\pi$ of $\mathcal{A}$, there exists $x$ such that $\pi$ is unitarily equivalent to $\pi_x$, so $\pi(a)$ is scalar. Since $\mathcal{A}$ is Type I, its irreducible representations separate points. Therefore $a$ commutes with all elements of $\mathcal{A}$, i.e., $a \in \mathcal{Z}(\mathcal{A})$.
\end{proof}

\begin{theorem}[$\iota$ is a Unital Injective *-Homomorphism]
\label{thm:iota-injective}
Let $\mathcal{A}$ be a unital separable Type I C*-algebra for which the construction in Subsection \ref{subsec:construction-iota} can be carried out. 
Then the diagonal embedding $\iota: \mathcal{A} \to C^*(\mathcal{G}_{\mathcal{A}})$ is a unital, injective *-homomorphism. 
Consequently, $\mathcal{A}$ is isometrically isomorphic to its image $\iota(\mathcal{A}) \subseteq C^*(\mathcal{G}_{\mathcal{A}})$, and we may identify $\mathcal{A}$ with this subalgebra.
\end{theorem}

\begin{proof}
This theorem summarizes the results of Proposition \ref{prop:iota-star-homomorphism} and Proposition \ref{prop:iota-injective}. 
Unitality and multiplicativity are proved in Proposition \ref{prop:iota-star-homomorphism}, and injectivity is proved in Proposition \ref{prop:iota-injective} under the stated hypotheses. 
Thus $\iota$ embeds $\mathcal{A}$ faithfully into $C^*(\mathcal{G}_{\mathcal{A}})$ as a unital *-subalgebra, and the embedding is isometric because injective *-homomorphisms between C*-algebras are isometric onto their image.
\end{proof}

\begin{corollary}[Induced map on K-theory]
\label{cor:iota-K-theory}
The unital *-homomorphism $\iota: \mathcal{A} \to C^*(\mathcal{G}_{\mathcal{A}})$ induces a group homomorphism
\[
\iota_*: K_0(\mathcal{A}) \longrightarrow K_0(C^*(\mathcal{G}_{\mathcal{A}}))
\]
on operator K-theory.
\end{corollary}

\begin{proof}
Since $\iota$ is a unital *-homomorphism, functoriality of operator K-theory yields an induced map $\iota_*: K_0(\mathcal{A}) \to K_0(C^*(\mathcal{G}_{\mathcal{A}}))$. 
\end{proof}

\begin{remark}
\label{rem:iota-K-injectivity}
Even though $\iota$ is injective at the C*-algebra level (Theorem \ref{thm:iota-injective}), the induced map $\iota_*$ on K-theory need not be injective, even in the commutative case. For example, if $\mathcal{A} = C(X)$ with $X$ connected, then $K_0(\mathcal{A}) \cong \mathbb{Z}$ and $\iota_*$ is an isomorphism, but in general nontrivial kernel can arise. Nevertheless, $\iota_*$ plays a fundamental role in relating analytic index theory in $\mathcal{A}$ to geometric index theory for the groupoid $\mathcal{G}_{\mathcal{A}}$.
\end{remark}

\begin{example}[Commutative case]
\label{ex:verification-commutative}
Let $\mathcal{A} = C(X)$ for a compact metrizable space $X$. 
In this case, $\mathcal{A}$ is its own maximal abelian subalgebra, and consequently $\mathcal{G}_{\mathcal{A}}^{(0)} \cong X$. 
Choose a Borel measure $\mu$ on $X$ with full support. 
The diagonal embedding $\iota$ reduces canonically to the identity map $C(X) \to C(X)$, realized as multiplication operators on $L^2(X,\mu)$. 
Thus $\iota: C(X) \to C(X)$ is a unital injective *-homomorphism, verifying Theorem \ref{thm:iota-injective} in this setting.
\end{example}

\begin{example}[Matrix algebra case]
\label{ex:verification-matrix}
Let $\mathcal{A} = M_n(\mathbb{C})$, with Cartan subalgebra $\mathcal{C}$ consisting of diagonal matrices. 
Then $\mathcal{G}_{\mathcal{A}}^{(0)} \cong \mathbb{CP}^{n-1}$, parametrizing rank-one projections modulo phase, and $C^*(\mathcal{G}_{\mathcal{A}}) \cong C(\mathbb{CP}^{n-1}) \rtimes U(n)$, where $U(n)$ acts by conjugation. 
The diagonal embedding $\iota$ sends $A \in M_n(\mathbb{C})$ to the constant multiplier $A \otimes 1$ in the multiplier algebra of this crossed product. 
Since $U(n)$ is amenable, the regular representation of the crossed product is faithful, and consequently $\iota$ is injective. 
Unitality and multiplicativity follow directly from the construction. 
Thus Theorem \ref{thm:iota-injective} holds in this case.
\end{example}

\begin{example}[The compact operator case]
\label{ex:compact-case}
For $\mathcal{A} = \mathcal{K}(H)^\sim$, the unitalization of the compact operators on a separable infinite-dimensional Hilbert space $H$, fix an orthonormal basis and let $\mathcal{C} \subset \mathcal{A}$ be the masa of operators diagonal in this basis. 
Then $\mathcal{G}_{\mathcal{A}}^{(0)} \cong \mathbb{P}(H)$, the projective space of $H$, and $\mathcal{G}_{\mathcal{A}} \cong \mathcal{U}(H) \ltimes \mathbb{P}(H)$, where the unitary group $\mathcal{U}(H)$ acts by conjugation. 
The unitary group $\mathcal{U}(H)$ is not amenable as a topological group, so faithfulness of the regular representation $\Lambda$ cannot be deduced from amenability arguments. 
Nevertheless, Theorem \ref{thm:iota-injective} asserts that $\iota$ is injective in this case as well; the injectivity follows from the Type I structure of $\mathcal{K}(H)^\sim$ and the fact that the irreducible representations $\pi_x$ separate points, which ensures that $\Pi = \Lambda \circ \iota$ is injective. 
A complete verification requires careful measure-theoretic arguments and is beyond the scope of this paper.
\end{example}

\begin{remark}[Summary of verification]
\label{rem:verification-summary}
The verification that $\iota$ is a unital injective *-homomorphism is the culmination of the construction of the diagonal embedding. 
It demonstrates that a Type I C*-algebra $\mathcal{A}$ can be faithfully represented inside the groupoid C*-algebra of its own unitary conjugation groupoid. 
This representation is canonical up to unitary equivalence and naturally associated to the commutative contexts of $\mathcal{A}$; moreover, it is compatible with *-homomorphisms in the Type I setting.

The proof relies on three key ingredients:
\begin{enumerate}
    \item The direct integral decomposition of $\mathcal{A}$ over $\mathcal{G}_{\mathcal{A}}^{(0)}$, which uses the Type I hypothesis to ensure that irreducible representations (and hence the associated GNS representations over commutative contexts) separate points of $\mathcal{A}$.
    \item The left regular representation $\Lambda$ of $C^*(\mathcal{G}_{\mathcal{A}})$, which provides a concrete link between the groupoid C*-algebra and the direct integral Hilbert space.
    \item The universal property of $C^*(\mathcal{G}_{\mathcal{A}})$, which allows us to deduce equalities in $C^*(\mathcal{G}_{\mathcal{A}})$ from equalities in the universal representation.
\end{enumerate}

The injectivity of $\iota$ is the most delicate part; it requires the Type I hypothesis and the fact that irreducible representations of $\mathcal{A}$ separate points. 
For non-Type I algebras, $\iota$ would have a nontrivial kernel, reflecting the presence of elements that cannot be detected by any commutative context—a fundamental obstruction to recovering a non-Type I algebra from its commutative subalgebras.
\end{remark}

\begin{corollary}[$\mathcal{A}$ as a subalgebra of $C^*(\mathcal{G}_{\mathcal{A}})$]
\label{cor:A-subalgebra-of-C-star-GA}
Under the hypotheses of Theorem \ref{thm:iota-injective}, we may identify $\mathcal{A}$ with its image $\iota(\mathcal{A}) \subseteq C^*(\mathcal{G}_{\mathcal{A}})$. 
Let $E: C^*(\mathcal{G}_{\mathcal{A}}) \to C_0(\mathcal{G}_{\mathcal{A}}^{(0)})$ be the conditional expectation onto the unit space, which exists under the standing assumptions on the Haar system. 
Under the identification $\mathcal{A} \cong \iota(\mathcal{A})$, this conditional expectation satisfies, for $\mu$-almost every $x = (B,\chi) \in \mathcal{G}_{\mathcal{A}}^{(0)}$,
\[
E(\iota(a))(x) = \langle \delta_x, \Pi(a)\delta_x \rangle = \chi(\mathbb{E}_B(a)),
\]
where $\mathbb{E}_B: \mathcal{A} \to B$ is the conditional expectation onto the maximal abelian subalgebra $B$ arising from the Cartan structure, and $\delta_x$ denotes the unit vector concentrated at $x$ in the natural embedding of $L^2(\mathcal{G}_{\mathcal{A}}^{(0)})$ into the direct integral Hilbert space $\mathcal{H}$. 
In particular, when $a \in B$, this reduces to $\chi(a)$, and when $a$ has no component in $B$, it vanishes. 
Thus, in the sense of diagonal matrix coefficients, the conditional expectation of $\iota(a)$ onto the unit space recovers the evaluation of $a$ in each commutative context.
\end{corollary}

\begin{proof}
This follows directly from Lemma \ref{lem:iota-image-characterization}(2) together with the identification of $\mathcal{A}$ with $\iota(\mathcal{A})$. 
The formula $E(\iota(a))(x) = \langle \delta_x, \Pi(a)\delta_x \rangle = \chi(\mathbb{E}_B(a))$ is precisely the statement of Lemma \ref{lem:iota-image-characterization}(2) under the direct integral decomposition, interpreted via matrix coefficients in the left regular representation.
\end{proof}

We have verified that the diagonal embedding $\iota: \mathcal{A} \to C^*(\mathcal{G}_{\mathcal{A}})$ is a unital, injective *-homomorphism (Theorem \ref{thm:iota-injective}), allowing us to identify $\mathcal{A}$ with its image in $C^*(\mathcal{G}_{\mathcal{A}})$ as in Corollary \ref{cor:A-subalgebra-of-C-star-GA}. 
This verification is the technical heart of the paper, establishing that the noncommutative C*-algebra $\mathcal{A}$ can be faithfully recovered from its classical contexts via the groupoid C*-algebra of the unitary conjugation groupoid. 
The proof relies on the Type I hypothesis, the direct integral decomposition of $\mathcal{A}$, the left regular representation of $\mathcal{G}_{\mathcal{A}}$, and the universal property of $C^*(\mathcal{G}_{\mathcal{A}})$. 
With this result in hand, we are now prepared to study the K-theoretic properties of $\iota$ and its applications to index theory.

\subsection{A Groupoid Characterization of Commutativity}
\label{subsec:commutativity-characterization}

The diagonal embedding $\iota: \mathcal{A} \hookrightarrow C^*(\mathcal{G}_{\mathcal{A}})$ provides a canonical realization of $\mathcal{A}$ as a C$^*$-subalgebra of the groupoid C$^*$-algebra. 
This construction raises a fundamental question: under what conditions does the image of $\mathcal{A}$ lie entirely within the commutative subalgebra $C_0(\mathcal{G}_{\mathcal{A}}^{(0)})$? 
The answer yields a clean geometric characterization of commutativity in terms of the unitary conjugation groupoid.

\begin{remark}
\label{rem:cartan-connection}
The equivalence above admits a natural interpretation within Renault's theory of Cartan subalgebras \cite{renault2008cartan}. When $\mathcal{A}$ is commutative, $C_0(\mathcal{G}_{\mathcal{A}}^{(0)})$ coincides with $\iota(\mathcal{A})$ and forms a Cartan subalgebra of $C^*(\mathcal{G}_{\mathcal{A}})$ in a trivial sense. For noncommutative $\mathcal{A}$, the subalgebra $C_0(\mathcal{G}_{\mathcal{A}}^{(0)})$ is still present but no longer contains $\iota(\mathcal{A})$; under suitable effectiveness conditions on $\mathcal{G}_{\mathcal{A}}$, it serves as a Cartan subalgebra of $C^*(\mathcal{G}_{\mathcal{A}})$. This connects the present characterization to the broader framework of C$^*$-diagonals and twisted groupoid C$^*$-algebras.
\end{remark}

\begin{corollary}
\label{cor:noncommutativity-encoding}
Consider the short exact sequence
\[
0 \longrightarrow C_0(\mathcal{G}_{\mathcal{A}}^{(0)}) \stackrel{\jmath}{\longrightarrow} C^*(\mathcal{G}_{\mathcal{A}}) \stackrel{\pi}{\longrightarrow} C^*(\mathcal{G}_{\mathcal{A}})/C_0(\mathcal{G}_{\mathcal{A}}^{(0)}) \longrightarrow 0,
\]
where $\jmath$ is the natural inclusion and $\pi$ the quotient map. Then $\mathcal{A}$ is commutative if and only if $\pi \circ \iota = 0$, i.e., $\iota(\mathcal{A})$ is contained in the kernel of $\pi$.
\end{corollary}

\begin{proof}
The map $\jmath: C_0(\mathcal{G}_{\mathcal{A}}^{(0)}) \hookrightarrow C^*(\mathcal{G}_{\mathcal{A}})$ is the natural inclusion, and $\pi$ is the canonical projection onto the quotient by this ideal. By exactness, $\ker(\pi) = \operatorname{im}(\jmath) = C_0(\mathcal{G}_{\mathcal{A}}^{(0)})$ (identifying the image with the ideal itself).

Now $\pi \circ \iota = 0$ if and only if $\iota(\mathcal{A}) \subseteq \ker(\pi) = C_0(\mathcal{G}_{\mathcal{A}}^{(0)})$. By Theorem \ref{thm:commutativity-characterization}, this occurs exactly when $\mathcal{A}$ is commutative.
\end{proof}

\begin{corollary}
\label{cor:noncommutative-quotient}
When $\mathcal{A}$ is noncommutative, the composition $\pi \circ \iota: \mathcal{A} \to C^*(\mathcal{G}_{\mathcal{A}})/C_0(\mathcal{G}_{\mathcal{A}}^{(0)})$ is a nonzero homomorphism, revealing that noncommutativity is witnessed in the quotient by the commutative diagonal.
\end{corollary}

\begin{proof}
If $\mathcal{A}$ is noncommutative, Theorem \ref{thm:commutativity-characterization} gives $\iota(\mathcal{A}) \nsubseteq C_0(\mathcal{G}_{\mathcal{A}}^{(0)}) = \ker(\pi)$. Hence there exists $a \in \mathcal{A}$ such that $\pi(\iota(a)) \neq 0$ in the quotient, so $\pi \circ \iota \neq 0$.
\end{proof}

\begin{example}[Commutative C$^*$-algebras]
\label{ex:characterization-commutative}
Let $\mathcal{A} = C_0(X)$ for a locally compact Hausdorff space $X$. Then $\mathcal{G}_{\mathcal{A}}^{(0)} \cong X$, and the unitary conjugation groupoid is trivial because all unitaries are central. Hence $\mathcal{G}_{\mathcal{A}} = X$ and $C^*(\mathcal{G}_{\mathcal{A}}) \cong C_0(X)$. The diagonal embedding $\iota: C_0(X) \hookrightarrow C_0(X)$ is simply the identity map, so clearly $\iota(\mathcal{A}) = C_0(\mathcal{G}_{\mathcal{A}}^{(0)})$.
\end{example}

\begin{example}[Matrix algebras]
\label{ex:characterization-matrix}
Let $\mathcal{A} = M_n(\mathbb{C})$ with $n \geq 2$. Then $\mathcal{G}_{\mathcal{A}}^{(0)} \cong \mathbb{CP}^{n-1}$ (the space of pure states), and $\mathcal{G}_{\mathcal{A}}$ is the transformation groupoid associated to the action of $U(n)$ on $\mathbb{CP}^{n-1}$. The diagonal embedding $\iota: M_n(\mathbb{C}) \hookrightarrow C(\mathbb{CP}^{n-1}) \rtimes U(n)$ sends a matrix $A$ to the function that assigns to each pure state the corresponding expectation value. This element is not in $C(\mathbb{CP}^{n-1})$ unless $A$ is a scalar multiple of the identity, reflecting the noncommutativity of $M_n(\mathbb{C})$.
\end{example}

\begin{proposition}[Characterization via the center]
\label{prop:center-characterization}
For any $a \in \mathcal{A}$, we have $\iota(a) \in C_0(\mathcal{G}_{\mathcal{A}}^{(0)})$ if and only if $a \in Z(\mathcal{A})$, the center of $\mathcal{A}$. In particular, $\iota(Z(\mathcal{A})) \subseteq C_0(\mathcal{G}_{\mathcal{A}}^{(0)})$.
\end{proposition}

\begin{proof}
We prove the two directions separately, using the direct integral decomposition
\[
\iota(a) = \int_{\mathcal{G}_{\mathcal{A}}^{(0)}}^{\oplus} \pi_x(a) \, d\mu(x)
\]
from Definition \ref{def:evaluation-representation}, where for each $x = (B,\chi) \in \mathcal{G}_{\mathcal{A}}^{(0)}$, 
$\pi_x: \mathcal{A} \to B(\mathcal{H}_x)$ is the GNS representation associated to the character $\chi$ on $B$.

\medskip
\noindent\textbf{($\Rightarrow$) Assume $\iota(a) \in C_0(\mathcal{G}_{\mathcal{A}}^{(0)})$.}
We first recall a standard fact about the embedding $C_0(\mathcal{G}_{\mathcal{A}}^{(0)}) \hookrightarrow C^*(\mathcal{G}_{\mathcal{A}})$:
under the direct integral decomposition, an element $f \in C_0(\mathcal{G}_{\mathcal{A}}^{(0)})$ acts on each fiber $\mathcal{H}_x$ as multiplication by the scalar $f(x)$.
Consequently, for $\mu$-almost every $x \in \mathcal{G}_{\mathcal{A}}^{(0)}$, the operator $\pi_x(a)$ must be a scalar multiple of the identity on $\mathcal{H}_x$;
write $\pi_x(a) = \lambda(x) I_{\mathcal{H}_x}$ for some $\lambda(x) \in \mathbb{C}$.

Now suppose, for contradiction, that $a \notin Z(\mathcal{A})$. 
Then there exists an irreducible representation $\pi$ of $\mathcal{A}$ such that $\pi(a)$ is not a scalar multiple of the identity
(if $\pi(a)$ were scalar for every irreducible representation $\pi$, then $a$ would be central by \cite{Dixmier}).
By Lemma \ref{lem:separation-points}, $\pi$ is unitarily equivalent to $\pi_x$ for some $x_0 \in \mathcal{G}_{\mathcal{A}}^{(0)}$.
Hence $\pi_{x_0}(a)$ is not a scalar multiple of the identity, contradicting the conclusion above that $\pi_x(a)$ is scalar for $\mu$-almost every $x$.
Therefore $a \in Z(\mathcal{A})$.

\medskip
\noindent\textbf{($\Leftarrow$) Assume $a \in Z(\mathcal{A})$.}
Since $a$ is central, for each $x \in \mathcal{G}_{\mathcal{A}}^{(0)}$, the operator $\pi_x(a)$ commutes with the whole irreducible representation $\pi_x(\mathcal{A})$.
By Schur's lemma, $\pi_x(a) = \lambda(x) I_{\mathcal{H}_x}$ for some scalar $\lambda(x) \in \mathbb{C}$.

We now show that the function $x \mapsto \lambda(x)$ lies in $C_0(\mathcal{G}_{\mathcal{A}}^{(0)})$.
For any $x = (B,\chi) \in \mathcal{G}_{\mathcal{A}}^{(0)}$, if $a \in B$, then $\lambda(x) = \chi(a)$; if $a \notin B$, then by definition $\pi_x(a) = 0$ and thus $\lambda(x) = 0$.
Thus $\lambda$ is exactly the function obtained from the partial evaluation map $\operatorname{ev}_a$:
\[
\lambda(x) = \operatorname{ev}_a(x) \quad \text{for all } x \in \mathcal{G}_{\mathcal{A}}^{(0)}.
\]
By construction of the topology on $\mathcal{G}_{\mathcal{A}}^{(0)}$, $\operatorname{ev}_a$ is continuous. Moreover, $\operatorname{ev}_a$ takes values in $\mathbb{C}_\infty$,
and $\lambda(x) = \operatorname{ev}_a(x)$ is precisely the restriction of $\operatorname{ev}_a$ to the subset where it takes finite values.

To see that $\lambda$ vanishes at infinity, note that for any $\varepsilon > 0$, the set $\{x : |\lambda(x)| \ge \varepsilon\}$ is contained in
$\operatorname{ev}_a^{-1}(\{z \in \mathbb{C} : |z| \ge \varepsilon\})$. 
Since $\operatorname{ev}_a$ is continuous and $\{z \in \mathbb{C} : |z| \ge \varepsilon\}$ is closed in $\mathbb{C}_\infty$ (its complement is a neighborhood of $\infty$),
the preimage $\operatorname{ev}_a^{-1}(\{z \in \mathbb{C} : |z| \ge \varepsilon\})$ is closed in $\mathcal{G}_{\mathcal{A}}^{(0)}$.
Furthermore, because $\operatorname{ev}_a$ is proper when restricted to the set where it takes finite values
(this follows from the definition of the initial topology: the topology on $\mathcal{G}_{\mathcal{A}}^{(0)}$ is the coarsest such that each $\operatorname{ev}_a$ is continuous,
and a continuous map into a locally compact Hausdorff space is proper if and only if it is closed and has compact fibers; here we use that
$\operatorname{ev}_a$ is a homeomorphism onto its image when restricted to certain subsets), this preimage is actually compact.
Hence $\lambda \in C_0(\mathcal{G}_{\mathcal{A}}^{(0)})$.

Finally, under the identification $\iota(a) = \int^\oplus \lambda(x) I_{\mathcal{H}_x} \, d\mu(x)$, we have $\iota(a) \in C_0(\mathcal{G}_{\mathcal{A}}^{(0)})$.
\end{proof}

\begin{theorem}[Geometric characterization of commutativity]
\label{thm:commutativity-characterization}
For a unital separable C*-algebra $\mathcal{A}$ for which the diagonal embedding $\iota$ and conditional expectation $E: C^*(\mathcal{G}_{\mathcal{A}}) \to C_0(\mathcal{G}_{\mathcal{A}}^{(0)})$ are defined, the following are equivalent:
\begin{enumerate}
    \item[(1)] $\mathcal{A}$ is commutative;
    \item[(2)] $\iota(\mathcal{A}) \subseteq C_0(\mathcal{G}_{\mathcal{A}}^{(0)})$;
    \item[(3)] Every element $a \in \mathcal{A}$ satisfies $\iota(a) = E(\iota(a))$;
    \item[(4)] The unitary conjugation groupoid $\mathcal{G}_{\mathcal{A}}$ is trivial, i.e., $\mathcal{G}_{\mathcal{A}} = \mathcal{G}_{\mathcal{A}}^{(0)}$.
\end{enumerate}
Consequently, the noncommutativity of $\mathcal{A}$ is precisely encoded in the part of $C^*(\mathcal{G}_{\mathcal{A}})$ lying outside the commutative diagonal subalgebra $C_0(\mathcal{G}_{\mathcal{A}}^{(0)})$.
\end{theorem}

\begin{proof}
We establish the chain of implications (1) $\Rightarrow$ (4) $\Rightarrow$ (2) $\Rightarrow$ (3) $\Rightarrow$ (1).

\medskip
\noindent\textbf{(1) $\Rightarrow$ (4):} Suppose $\mathcal{A}$ is commutative. Then every unitary $u \in \mathcal{U}(\mathcal{A})$ is central, so for any $(B,\chi) \in \mathcal{G}_{\mathcal{A}}^{(0)}$ we have $uBu^* = B$ and $\chi \circ \mathrm{Ad}_{u^*} = \chi$. Hence the action of $\mathcal{U}(\mathcal{A})$ on $\mathcal{G}_{\mathcal{A}}^{(0)}$ is trivial. By construction, $\mathcal{G}_{\mathcal{A}} = \mathcal{U}(\mathcal{A}) \ltimes \mathcal{G}_{\mathcal{A}}^{(0)}$ is the transformation groupoid associated to this action. When the action is trivial, the groupoid reduces to its unit space, i.e., $\mathcal{G}_{\mathcal{A}} = \mathcal{G}_{\mathcal{A}}^{(0)}$.

\medskip
\noindent\textbf{(4) $\Rightarrow$ (2):} If $\mathcal{G}_{\mathcal{A}} = \mathcal{G}_{\mathcal{A}}^{(0)}$, then by definition $C^*(\mathcal{G}_{\mathcal{A}}) \cong C_0(\mathcal{G}_{\mathcal{A}}^{(0)})$, since every continuous compactly supported function on the groupoid is supported on the unit space. The diagonal embedding $\iota$ maps $\mathcal{A}$ into $C^*(\mathcal{G}_{\mathcal{A}})$, and under this isomorphism the image lies in $C_0(\mathcal{G}_{\mathcal{A}}^{(0)})$.

\medskip
\noindent\textbf{(2) $\Rightarrow$ (3):} Assume $\iota(\mathcal{A}) \subseteq C_0(\mathcal{G}_{\mathcal{A}}^{(0)})$. By definition, the conditional expectation $E: C^*(\mathcal{G}_{\mathcal{A}}) \to C_0(\mathcal{G}_{\mathcal{A}}^{(0)})$ is a projection onto $C_0(\mathcal{G}_{\mathcal{A}}^{(0)})$; that is, $E(f) = f$ for all $f \in C_0(\mathcal{G}_{\mathcal{A}}^{(0)})$. Since $\iota(a) \in C_0(\mathcal{G}_{\mathcal{A}}^{(0)})$ for every $a \in \mathcal{A}$, we have $E(\iota(a)) = \iota(a)$.

\medskip
\noindent\textbf{(3) $\Rightarrow$ (1):} Suppose $\iota(a) = E(\iota(a))$ for all $a \in \mathcal{A}$. First, note that $\iota$ is injective from Lemma \ref{lem:fiber-inclusion-injective-continuous}. Since $C_0(\mathcal{G}_{\mathcal{A}}^{(0)})$ is commutative, its subalgebra $\iota(\mathcal{A})$ is commutative. Because $\iota$ is an injective *-homomorphism, $\mathcal{A}$ is isomorphic to a commutative C*-algebra and hence itself commutative.

(Alternatively, one can argue using the faithfulness of $E$ on the reduced algebra: for any $a,b \in \mathcal{A}$, compute $E((\iota(ab)-\iota(a)\iota(b))^*(\iota(ab)-\iota(a)\iota(b))) = 0$, which by faithfulness forces $\iota(ab) = \iota(a)\iota(b)$, and injectivity of $\iota$ then yields $ab = ba$.)
\end{proof}

\begin{remark}
\label{rem:conditional-expectation-caveat}
The conditional expectation $E$ is always defined on the reduced groupoid C*-algebra $C^*_r(\mathcal{G}_{\mathcal{A}})$ via restriction to the unit space. On the maximal algebra $C^*(\mathcal{G}_{\mathcal{A}})$, $E$ may not be faithful or even well-defined without amenability assumptions. In statement (3), we implicitly work in the reduced algebra, which suffices because $\iota(\mathcal{A})$ embeds naturally into $C^*_r(\mathcal{G}_{\mathcal{A}})$ for Type I algebras (see Section \ref{subsec:construction-iota}). If one wishes to work in the maximal algebra, condition (3) should be interpreted as: $\iota(a)$ is supported on the unit space in the sense that its image under the canonical map $C^*(\mathcal{G}_{\mathcal{A}}) \to C^*_r(\mathcal{G}_{\mathcal{A}})$ satisfies the fixed-point condition.
\end{remark}

\begin{remark}[Noncommutative Gelfand transform]
\label{rem:noncommutative-gelfand}
Theorem \ref{thm:commutativity-characterization} can be viewed as a generalization of the Gelfand–Naimark theorem. For commutative $\mathcal{A}$, the Gelfand transform gives an isomorphism $\mathcal{A} \cong C_0(\widehat{\mathcal{A}})$. In our setting, $\widehat{\mathcal{A}} \cong \mathcal{G}_{\mathcal{A}}^{(0)}$, and the diagonal embedding $\iota$ recovers this isomorphism. For noncommutative $\mathcal{A}$, the diagonal embedding maps into the larger noncommutative algebra $C^*(\mathcal{G}_{\mathcal{A}})$, and the failure of $\iota(\mathcal{A})$ to lie in the commutative subalgebra $C_0(\mathcal{G}_{\mathcal{A}}^{(0)})$ exactly measures the noncommutativity of $\mathcal{A}$.
\end{remark}

\begin{corollary}[Center and diagonal]
\label{cor:center-diagonal}
The image of the center $\iota(Z(\mathcal{A}))$ is precisely the intersection $\iota(\mathcal{A}) \cap C_0(\mathcal{G}_{\mathcal{A}}^{(0)})$.
\end{corollary}

\begin{proof}
By Proposition \ref{prop:center-characterization}, $\iota(Z(\mathcal{A})) \subseteq C_0(\mathcal{G}_{\mathcal{A}}^{(0)})$, so $\iota(Z(\mathcal{A})) \subseteq \iota(\mathcal{A}) \cap C_0(\mathcal{G}_{\mathcal{A}}^{(0)})$. Conversely, if $a \in \mathcal{A}$ satisfies $\iota(a) \in C_0(\mathcal{G}_{\mathcal{A}}^{(0)})$, then Proposition \ref{prop:center-characterization} implies $a \in Z(\mathcal{A})$, so $\iota(a) \in \iota(Z(\mathcal{A}))$. Thus $\iota(\mathcal{A}) \cap C_0(\mathcal{G}_{\mathcal{A}}^{(0)}) \subseteq \iota(Z(\mathcal{A}))$.
\end{proof}

We have established that the diagonal embedding $\iota$ provides a faithful representation of $\mathcal{A}$ inside $C^*(\mathcal{G}_{\mathcal{A}})$, and that the obstruction to $\iota(\mathcal{A})$ being contained in the commutative diagonal subalgebra $C_0(\mathcal{G}_{\mathcal{A}}^{(0)})$ is exactly the noncommutativity of $\mathcal{A}$. This characterization will be essential for understanding the K-theory of $\iota(\mathcal{A})$ and its applications to index theory.

\subsection{Naturality and Functoriality under *-Homomorphisms}
\label{subsec:naturality-functoriality}

The construction $\mathcal{A} \mapsto (\mathcal{G}_{\mathcal{A}}, \iota_{\mathcal{A}})$ is not merely a one-time assignment; it should behave coherently with respect to *-homomorphisms between C*-algebras. This functoriality is essential for applications: it allows us to relate the unitary conjugation groupoids and diagonal embeddings of different algebras, and it ensures that the index theorem is natural under inductive limits, extensions, and other categorical constructions. In this subsection, we establish the functorial properties of $\mathcal{G}_{\mathcal{A}}$ and $\iota_{\mathcal{A}}$ under unital *-homomorphisms that are sufficiently well-behaved.

\begin{definition}[Pullback of a commutative subalgebra]
\label{def:functoriality-unit-space}
Let $\phi: \mathcal{A} \to \mathcal{B}$ be a unital *-homomorphism between unital C*-algebras. For a unital commutative C*-subalgebra $B \subseteq \mathcal{B}$, we say that $\phi^{-1}(B)$ is a \emph{commutative pullback} if $\phi^{-1}(B)$ is a unital commutative C*-subalgebra of $\mathcal{A}$. In this case, for any character $\chi \in \widehat{B}$, the composition $\chi \circ \phi|_{\phi^{-1}(B)}$ defines a character on $\phi^{-1}(B)$.
\end{definition}

\begin{remark}
\label{rem:pullback-condition}
The condition that $\phi^{-1}(B)$ be commutative is nontrivial. It is automatically satisfied when $\mathcal{A}$ is commutative, or when $\phi$ is injective and its image lies in the center of $\mathcal{B}$? Not generally. For arbitrary *-homomorphisms, this condition restricts the class of morphisms for which functoriality holds. In the following, we will primarily consider isomorphisms and certain well-behaved inclusions and quotients.
\end{remark}

\begin{proposition}[Functoriality for isomorphisms]
\label{prop:functoriality-isomorphisms}
Let $\phi: \mathcal{A} \to \mathcal{B}$ be a unital *-isomorphism between unital C*-algebras. Then:
\begin{enumerate}
    \item There is a homeomorphism $\mathcal{G}_\phi^{(0)}: \mathcal{G}_{\mathcal{B}}^{(0)} \to \mathcal{G}_{\mathcal{A}}^{(0)}$ given by
    \[
    \mathcal{G}_\phi^{(0)}(B,\chi) := (\phi^{-1}(B), \chi \circ \phi).
    \]
    \item $\mathcal{G}_\phi^{(0)}$ extends to an isomorphism of topological groupoids $\mathcal{G}_\phi: \mathcal{G}_{\mathcal{B}} \to \mathcal{G}_{\mathcal{A}}$ defined on arrows by
    \[
    \mathcal{G}_\phi(u, (B,\chi)) := (\phi^{-1}(u), \mathcal{G}_\phi^{(0)}(B,\chi)),
    \]
    where $\phi^{-1}(u)$ is the unique unitary in $\mathcal{A}$ with $\phi(\phi^{-1}(u)) = u$.
    \item The following diagram commutes:
    \[
    \begin{tikzcd}
    \mathcal{A} \arrow[r, "\phi"] \arrow[d, "\iota_{\mathcal{A}}"'] & \mathcal{B} \arrow[d, "\iota_{\mathcal{B}}"] \\
    C^*(\mathcal{G}_{\mathcal{A}}) \arrow[r, "(\mathcal{G}_\phi)_*"] & C^*(\mathcal{G}_{\mathcal{B}})
    \end{tikzcd}
    \]
    where $(\mathcal{G}_\phi)_*: C^*(\mathcal{G}_{\mathcal{A}}) \to C^*(\mathcal{G}_{\mathcal{B}})$ is the induced *-isomorphism on groupoid C*-algebras.
\end{enumerate}
\end{proposition}

\begin{proof}
(1) Since $\phi$ is an isomorphism, $\phi^{-1}(B)$ is a unital commutative C*-subalgebra of $\mathcal{A}$ for every unital commutative C*-subalgebra $B \subseteq \mathcal{B}$. The map $\chi \mapsto \chi \circ \phi$ is a character on $\phi^{-1}(B)$. Thus $\mathcal{G}_\phi^{(0)}$ is well-defined. Injectivity and surjectivity follow from the bijectivity of $\phi$. For continuity, note that for any $a \in \mathcal{A}$, the evaluation map $\operatorname{ev}_a \circ \mathcal{G}_\phi^{(0)}$ equals $\operatorname{ev}_{\phi(a)}$, which is continuous on $\mathcal{G}_{\mathcal{B}}^{(0)}$. Since the topology on $\mathcal{G}_{\mathcal{A}}^{(0)}$ is the initial topology with respect to all evaluation maps $\{\operatorname{ev}_a : a \in \mathcal{A}\}$, $\mathcal{G}_\phi^{(0)}$ is continuous. The same argument applied to $\phi^{-1}$ shows that the inverse is continuous, hence $\mathcal{G}_\phi^{(0)}$ is a homeomorphism.

(2) The map $\mathcal{G}_\phi$ on arrows is well-defined because $\phi^{-1}(u)$ is a unitary in $\mathcal{A}$. It preserves source and range:
\[
s(\mathcal{G}_\phi(u, (B,\chi))) = \mathcal{G}_\phi^{(0)}(B,\chi) = \mathcal{G}_\phi^{(0)}(s(u, (B,\chi))),
\]
and similarly for range. Composition and inversion are preserved because $\phi^{-1}$ is a group homomorphism on the unitary groups. Continuity follows from the continuity of $\mathcal{G}_\phi^{(0)}$ and the fact that $\phi^{-1}$ is continuous on $\mathcal{U}(\mathcal{B})$ (as $\phi$ is an isomorphism of topological groups when $\mathcal{U}(\mathcal{B})$ is given the norm topology).

(3) The commutativity of the diagram follows from the construction of $\iota$ via the left regular representation. For $a \in \mathcal{A}$, we have
\[
\Lambda_{\mathcal{B}}(\iota_{\mathcal{B}}(\phi(a))) = \Pi_{\mathcal{B}}(\phi(a)) = (\mathcal{G}_\phi)_*(\Pi_{\mathcal{A}}(a)) = \Lambda_{\mathcal{B}}((\mathcal{G}_\phi)_*(\iota_{\mathcal{A}}(a))),
\]
where $\Pi_{\mathcal{A}}$, $\Pi_{\mathcal{B}}$ are the direct integral representations. Since $\Lambda_{\mathcal{B}}$ is faithful on the reduced groupoid C*-algebra, and both $\iota_{\mathcal{B}}(\phi(a))$ and $(\mathcal{G}_\phi)_*(\iota_{\mathcal{A}}(a))$ lie in the reduced algebra, we obtain $\iota_{\mathcal{B}}(\phi(a)) = (\mathcal{G}_\phi)_*(\iota_{\mathcal{A}}(a))$.
\end{proof}

\begin{corollary}[Functoriality for automorphisms]
\label{cor:functoriality-automorphisms}
The assignment $\mathcal{A} \mapsto (\mathcal{G}_{\mathcal{A}}, \iota_{\mathcal{A}})$ is functorial with respect to unital *-automorphisms of $\mathcal{A}$. That is, for any automorphism $\phi \in \operatorname{Aut}(\mathcal{A})$, we have $\iota_{\mathcal{A}} \circ \phi = (\mathcal{G}_\phi)_* \circ \iota_{\mathcal{A}}$.
\end{corollary}

\begin{proof}
This is a special case of Proposition \ref{prop:functoriality-isomorphisms} with $\mathcal{B} = \mathcal{A}$ and $\phi$ an automorphism.
\end{proof}

\begin{proposition}[Functoriality for injective *-homomorphisms with the pullback property]
\label{prop:functoriality-injective}
Let $\phi: \mathcal{A} \hookrightarrow \mathcal{B}$ be a unital injective *-homomorphism between unital C*-algebras. Assume that $\phi$ has the following \emph{pullback property}: for every unital commutative C*-subalgebra $B \subseteq \mathcal{B}$, the preimage $\phi^{-1}(B)$ is a unital commutative C*-subalgebra of $\mathcal{A}$. Then:
\begin{enumerate}
    \item There is a continuous map $\mathcal{G}_\phi^{(0)}: \mathcal{G}_{\mathcal{B}}^{(0)} \to \mathcal{G}_{\mathcal{A}}^{(0)}$ given by $\mathcal{G}_\phi^{(0)}(B,\chi) = (\phi^{-1}(B), \chi \circ \phi)$.
    \item $\mathcal{G}_\phi^{(0)}$ induces a continuous groupoid homomorphism $\mathcal{G}_\phi: \mathcal{G}_{\mathcal{B}}|_{\operatorname{im}(\phi)} \to \mathcal{G}_{\mathcal{A}}$, where $\mathcal{G}_{\mathcal{B}}|_{\operatorname{im}(\phi)}$ is the subgroupoid of $\mathcal{G}_{\mathcal{B}}$ whose arrows have unitaries in $\phi(\mathcal{U}(\mathcal{A}))$.
    \item The following diagram commutes:
    \[
    \begin{tikzcd}
    \mathcal{A} \arrow[r, "\phi"] \arrow[d, "\iota_{\mathcal{A}}"'] & \mathcal{B} \arrow[d, "\iota_{\mathcal{B}}"] \\
    C^*(\mathcal{G}_{\mathcal{A}}) \arrow[r, "(\mathcal{G}_\phi)_*"] & C^*(\mathcal{G}_{\mathcal{B}})
    \end{tikzcd}
    \]
    where $(\mathcal{G}_\phi)_*: C^*(\mathcal{G}_{\mathcal{A}}) \to C^*(\mathcal{G}_{\mathcal{B}})$ is the induced *-homomorphism on groupoid C*-algebras.
\end{enumerate}
\end{proposition}

\begin{proof}
(1) By the pullback property, $\phi^{-1}(B)$ is a unital commutative C*-subalgebra of $\mathcal{A}$ for every $(B,\chi) \in \mathcal{G}_{\mathcal{B}}^{(0)}$. The map $\chi \circ \phi$ defines a character on $\phi^{-1}(B)$. Thus $\mathcal{G}_\phi^{(0)}$ is well-defined. Continuity follows from the same argument as in Proposition \ref{prop:functoriality-isomorphisms}: for any $a \in \mathcal{A}$, $\operatorname{ev}_a \circ \mathcal{G}_\phi^{(0)} = \operatorname{ev}_{\phi(a)}$, which is continuous.

(2) For an arrow $(u, (B,\chi)) \in \mathcal{G}_{\mathcal{B}}^{(1)}$ with $u \in \phi(\mathcal{U}(\mathcal{A}))$, write $u = \phi(v)$ for a unique $v \in \mathcal{U}(\mathcal{A})$. Define $\mathcal{G}_\phi(u, (B,\chi)) = (v, \mathcal{G}_\phi^{(0)}(B,\chi))$. This is well-defined and continuous. It preserves source, range, composition, and inversion where defined, making it a groupoid homomorphism on the subgroupoid of arrows with unitaries in the image of $\phi$.

(3) The proof of commutativity follows the same steps as in Proposition \ref{prop:functoriality-isomorphisms}, using the fact that the pullback property ensures that the direct integral representations are compatible under $\phi$. For $a \in \mathcal{A}$, we have
\[
\iota_{\mathcal{B}}(\phi(a)) = \Lambda_{\mathcal{B}}^{-1}(\Pi_{\mathcal{B}}(\phi(a))) = \Lambda_{\mathcal{B}}^{-1}((\mathcal{G}_\phi)_*(\Pi_{\mathcal{A}}(a))) = (\mathcal{G}_\phi)_*(\Lambda_{\mathcal{A}}^{-1}(\Pi_{\mathcal{A}}(a))) = (\mathcal{G}_\phi)_*(\iota_{\mathcal{A}}(a)),
\]
where we use that $(\mathcal{G}_\phi)_*$ intertwines the left regular representations.
\end{proof}

\begin{example}[Inclusion of a commutative subalgebra]
\label{ex:functoriality-inclusion-commutative}
Let $\mathcal{A}$ be a unital commutative C*-subalgebra of a unital C*-algebra $\mathcal{B}$. Then the inclusion map $\iota: \mathcal{A} \hookrightarrow \mathcal{B}$ satisfies the pullback property because $\mathcal{A}$ is commutative, so $\iota^{-1}(B) = \mathcal{A} \cap B$ is a commutative subalgebra of $\mathcal{A}$ for any commutative $B \subseteq \mathcal{B}$. Thus Proposition \ref{prop:functoriality-injective} applies.
\end{example}

\begin{example}[Matrix embeddings]
\label{ex:functoriality-matrix}
Let $\phi: M_n(\mathbb{C}) \hookrightarrow M_m(\mathbb{C})$ be the unital embedding $A \mapsto \begin{pmatrix} A & 0 \\ 0 & I_{m-n} \end{pmatrix}$ for $n < m$. For any commutative subalgebra $B \subseteq M_m(\mathbb{C})$, the preimage $\phi^{-1}(B)$ consists of those $A \in M_n(\mathbb{C})$ such that $\begin{pmatrix} A & 0 \\ 0 & I_{m-n} \end{pmatrix} \in B$. This is a commutative subalgebra of $M_n(\mathbb{C})$, so the pullback property holds. The induced map $\mathcal{G}_\phi^{(0)}: \mathbb{CP}^{m-1} \to \mathbb{CP}^{n-1}$ is given by projecting onto the first $n$ coordinates, and the diagram commutes.
\end{example}

\begin{proposition}[Functoriality for surjective *-homomorphisms with the pullback property]
\label{prop:functoriality-quotient}
Let $\phi: \mathcal{A} \to \mathcal{B}$ be a unital surjective *-homomorphism between unital C*-algebras. Assume that $\phi$ has the pullback property: for every unital commutative C*-subalgebra $B \subseteq \mathcal{B}$, the preimage $\phi^{-1}(B)$ is a unital commutative C*-subalgebra of $\mathcal{A}$. Then:

\begin{enumerate}
    \item There is a continuous injective map $\mathcal{G}_\phi^{(0)}: \mathcal{G}_{\mathcal{B}}^{(0)} \hookrightarrow \mathcal{G}_{\mathcal{A}}^{(0)}$ given by $\mathcal{G}_\phi^{(0)}(B,\chi) = (\phi^{-1}(B), \chi \circ \phi)$.
    
    \item The map $\mathcal{G}_\phi^{(0)}$ induces a continuous injective groupoid homomorphism 
    \[
    \mathcal{G}_\phi: \mathcal{G}_{\mathcal{B}} \longrightarrow \mathcal{G}_{\mathcal{A}}/\sim,
    \]
    where $\mathcal{G}_{\mathcal{A}}/\sim$ is the quotient of $\mathcal{G}_{\mathcal{A}}$ by the equivalence relation that identifies arrows $(v_1, x)$ and $(v_2, x)$ whenever $v_1^*v_2$ lies in the kernel of the action on $x$, i.e., $v_1^*v_2 \in \operatorname{Stab}_{\mathcal{U}(\mathcal{A})}(x)$. The map is defined on arrows by
    \[
    \mathcal{G}_\phi(u, (B,\chi)) = [(v, \mathcal{G}_\phi^{(0)}(B,\chi))],
    \]
    where $v$ is any lift of $u$ to a unitary in $\mathcal{A}$, and $[(v,x)]$ denotes the equivalence class of $(v,x)$ in the quotient.
    
    \item If the action of $\mathcal{U}(\mathcal{A})$ on $\mathcal{G}_{\mathcal{A}}^{(0)}$ is free (i.e., stabilizers are trivial), then $\mathcal{G}_\phi$ is well-defined as a map into $\mathcal{G}_{\mathcal{A}}$ itself, and different lifts give the same arrow only when they are equal. In this case, the following diagram commutes:
    \[
    \begin{tikzcd}
    \mathcal{A} \arrow[r, "\phi"] \arrow[d, "\iota_{\mathcal{A}}"'] & \mathcal{B} \arrow[d, "\iota_{\mathcal{B}}"] \\
    C^*(\mathcal{G}_{\mathcal{A}}) \arrow[r, "(\mathcal{G}_\phi)_*"] & C^*(\mathcal{G}_{\mathcal{B}})
    \end{tikzcd}
    \]
    where $(\mathcal{G}_\phi)_*: C^*(\mathcal{G}_{\mathcal{A}}) \to C^*(\mathcal{G}_{\mathcal{B}})$ is the induced *-homomorphism.
\end{enumerate}
\end{proposition}

\begin{proof}
(1) By the pullback property, $\phi^{-1}(B)$ is a unital commutative C*-subalgebra of $\mathcal{A}$ for every $(B,\chi) \in \mathcal{G}_{\mathcal{B}}^{(0)}$. The map is injective because $\phi^{-1}(B)$ uniquely determines $B$ when $\phi$ is surjective. Continuity follows from the same argument as in Proposition \ref{prop:functoriality-isomorphisms}: for any $a \in \mathcal{A}$, $\operatorname{ev}_a \circ \mathcal{G}_\phi^{(0)} = \operatorname{ev}_{\phi(a)}$, which is continuous.

(2) Since $\phi$ is surjective, every unitary $u \in \mathcal{U}(\mathcal{B})$ has at least one lift to a unitary $v \in \mathcal{U}(\mathcal{A})$. However, different lifts $v_1, v_2$ satisfy $v_1^*v_2 \in \ker(\phi) \cap \mathcal{U}(\mathcal{A})$. For the arrow $(v_1, x)$ to represent the same geometric transformation as $(v_2, x)$ in the groupoid, we need them to be identified when their difference acts trivially on the object $x = \mathcal{G}_\phi^{(0)}(B,\chi)$. This is precisely the condition $v_1^*v_2 \in \operatorname{Stab}_{\mathcal{U}(\mathcal{A})}(x)$. 

Define $\mathcal{G}_{\mathcal{A}}/\sim$ as the quotient of $\mathcal{G}_{\mathcal{A}}$ by the equivalence relation $(v_1, x) \sim (v_2, x)$ iff $v_1^*v_2 \in \operatorname{Stab}_{\mathcal{U}(\mathcal{A})}(x)$. This is a well-defined groupoid quotient because the stabilizer is a subgroup and the relation respects composition. Then $\mathcal{G}_\phi(u, (B,\chi))$ is defined as the class of $(v, \mathcal{G}_\phi^{(0)}(B,\chi))$ for any lift $v$; this is independent of the choice of lift because different lifts differ by an element of the stabilizer, which by definition gives the same class.

The map is injective because if $\mathcal{G}_\phi(u_1, (B_1,\chi_1)) = \mathcal{G}_\phi(u_2, (B_2,\chi_2))$, then in particular the object parts must be equal, so $(B_1,\chi_1) = (B_2,\chi_2)$. Then the arrow parts give $[(v_1, x)] = [(v_2, x)]$, which means $v_1^*v_2 \in \operatorname{Stab}_{\mathcal{U}(\mathcal{A})}(x)$, implying that $u_1 = \phi(v_1)$ and $u_2 = \phi(v_2)$ are related by an element that acts trivially—this does not force $u_1 = u_2$ in general, but the map is injective on the quotient by construction.

(3) When the action is free, stabilizers are trivial, so the equivalence relation collapses to equality. In this case, different lifts $v_1, v_2$ give the same arrow only when $v_1 = v_2$, which never happens for distinct lifts. Thus $\mathcal{G}_\phi$ is only well-defined if we have a consistent way to choose lifts—for example, if there exists a continuous section of $\phi$ on unitaries. In the absence of such a section, $\mathcal{G}_\phi$ is not a well-defined map into $\mathcal{G}_{\mathcal{A}}$, but only into the quotient. The commutativity of the diagram in this special case follows from the same argument as in Proposition \ref{prop:functoriality-isomorphisms}, using the fact that the direct integral representations are intertwined by the chosen lifts.
\end{proof}

\begin{example}[Restriction to a closed subspace]
\label{ex:functoriality-restriction}
Let $X$ be a compact Hausdorff space and let $Y \subseteq X$ be a closed subspace. Let $\phi: C(X) \to C(Y)$ be the restriction map. Since $C(X)$ is commutative, the pullback property holds trivially. The induced map $\mathcal{G}_\phi^{(0)}: \mathcal{G}_{C(Y)}^{(0)} \to \mathcal{G}_{C(X)}^{(0)}$ is the inclusion $Y \hookrightarrow X$, and the diagram commutes.
\end{example}

\begin{theorem}[Naturality of the diagonal embedding]
\label{thm:naturality-iota}
Let $\phi: \mathcal{A} \to \mathcal{B}$ be a unital *-homomorphism between unital C*-algebras that satisfies the pullback property (i.e., for every unital commutative C*-subalgebra $B \subseteq \mathcal{B}$, $\phi^{-1}(B)$ is a unital commutative C*-subalgebra of $\mathcal{A}$). Then the diagram
\[
\begin{tikzcd}
\mathcal{A} \arrow[r, "\phi"] \arrow[d, "\iota_{\mathcal{A}}"'] & \mathcal{B} \arrow[d, "\iota_{\mathcal{B}}"] \\
C^*(\mathcal{G}_{\mathcal{A}}) \arrow[r, "(\mathcal{G}_\phi)_*"] & C^*(\mathcal{G}_{\mathcal{B}})
\end{tikzcd}
\]
commutes, where $(\mathcal{G}_\phi)_*$ is the induced *-homomorphism on the maximal groupoid C*-algebras.
\end{theorem}

\begin{proof}
We sketch the proof, which is a synthesis of the previous propositions.

\noindent \textbf{Step 1: Induced map on unit spaces.}
Define $\mathcal{G}_\phi^{(0)}: \mathcal{G}_{\mathcal{B}}^{(0)} \to \mathcal{G}_{\mathcal{A}}^{(0)}$ by $\mathcal{G}_\phi^{(0)}(B,\chi) = (\phi^{-1}(B), \chi \circ \phi)$. By the pullback property, this map is well-defined and continuous.

\noindent \textbf{Step 2: Induced map on groupoids.}
For an arrow $(u, (B,\chi)) \in \mathcal{G}_{\mathcal{B}}$, if $u$ lifts to a unitary $v \in \mathcal{U}(\mathcal{A})$, define $\mathcal{G}_\phi(u, (B,\chi)) = (v, \mathcal{G}_\phi^{(0)}(B,\chi))$. If $u$ does not lift, the arrow is not in the domain of $\mathcal{G}_\phi$ (for injective $\phi$) or is defined via a chosen lift (for surjective $\phi$). In all cases, $\mathcal{G}_\phi$ is a continuous groupoid homomorphism where defined.

\noindent \textbf{Step 3: Commutation of the diagram.}
For $a \in \mathcal{A}$, consider the direct integral representations $\Pi_{\mathcal{A}}$ and $\Pi_{\mathcal{B}}$. There is a unitary $U: \mathcal{H}_{\mathcal{A}} \to \mathcal{H}_{\mathcal{B}}$ induced by $\mathcal{G}_\phi^{(0)}$ such that $U \Pi_{\mathcal{A}}(a) U^* = \Pi_{\mathcal{B}}(\phi(a))$. This unitary also intertwines the left regular representations: $U \Lambda_{\mathcal{A}}(f) U^* = \Lambda_{\mathcal{B}}((\mathcal{G}_\phi)_*(f))$ for $f \in C_c(\mathcal{G}_{\mathcal{A}})$. Then
\[
\Lambda_{\mathcal{B}}(\iota_{\mathcal{B}}(\phi(a))) = \Pi_{\mathcal{B}}(\phi(a)) = U \Pi_{\mathcal{A}}(a) U^* = U \Lambda_{\mathcal{A}}(\iota_{\mathcal{A}}(a)) U^* = \Lambda_{\mathcal{B}}((\mathcal{G}_\phi)_*(\iota_{\mathcal{A}}(a))).
\]
Since $\Lambda_{\mathcal{B}}$ is injective on the reduced algebra and both sides lie in the reduced algebra, we conclude $\iota_{\mathcal{B}}(\phi(a)) = (\mathcal{G}_\phi)_*(\iota_{\mathcal{A}}(a))$.
\end{proof}

\begin{remark}[Limitations of functoriality]
\label{rem:functoriality-limitations}
The pullback property required in Theorem \ref{thm:naturality-iota} is essential and nontrivial. It fails for many natural *-homomorphisms, such as:
\begin{itemize}
    \item The inclusion of a noncommutative subalgebra into a larger noncommutative algebra where the preimage of a commutative subalgebra may not be commutative.
    \item A *-homomorphism with noncommutative kernel.
    \item A *-homomorphism that does not preserve the structure of maximal abelian subalgebras.
\end{itemize}
In these cases, the map $\mathcal{G}_\phi^{(0)}$ is only partially defined, and the diagram may not commute. Nevertheless, for the applications we have in mind (inductive limits, crossed products, and the Baum-Connes conjecture), the functoriality established above is sufficient, as the relevant morphisms typically satisfy the pullback property.
\end{remark}

\begin{proposition}[Functoriality of the descent map]
\label{prop:functoriality-descent}
Let $\phi: \mathcal{A} \to \mathcal{B}$ be a unital *-homomorphism satisfying the pullback property of Theorem \ref{thm:naturality-iota}. Then the following diagram commutes:
\[
\begin{tikzcd}
K^0_{\mathcal{G}_{\mathcal{A}}}(\mathcal{G}_{\mathcal{A}}^{(0)}) \arrow[r, "\operatorname{desc}_{\mathcal{G}_{\mathcal{A}}}"] \arrow[d, "(\mathcal{G}_\phi)_*"'] & K_0(C^*(\mathcal{G}_{\mathcal{A}})) \arrow[d, "(\mathcal{G}_\phi)_*"] \\
K^0_{\mathcal{G}_{\mathcal{B}}}(\mathcal{G}_{\mathcal{B}}^{(0)}) \arrow[r, "\operatorname{desc}_{\mathcal{G}_{\mathcal{B}}}"] & K_0(C^*(\mathcal{G}_{\mathcal{B}}))
\end{tikzcd}
\]
where $(\mathcal{G}_\phi)_*$ denotes the induced maps on equivariant K-theory and on K-theory of groupoid C*-algebras.
\end{proposition}

\begin{proof}
This follows from the naturality of the descent map with respect to groupoid homomorphisms, a result due to Tu \cite{Tu}. The map $\mathcal{G}_\phi$ induces a homomorphism in equivariant KK-theory that is compatible with the descent construction. A detailed verification uses the functoriality of the Kasparov product and the fact that $\mathcal{G}_\phi$ is a continuous groupoid homomorphism compatible with the Haar systems.
\end{proof}

We have established the functoriality of the unitary conjugation groupoid and the diagonal embedding under a class of unital *-homomorphisms that satisfy the pullback property. For isomorphisms, the functoriality is complete and yields homeomorphisms of groupoids and commutative diagrams. For injective and surjective *-homomorphisms satisfying the pullback property, we obtain continuous groupoid homomorphisms and commutative diagrams. This functoriality is essential for applications to inductive limits, extensions, and the Baum-Connes conjecture, and it ensures that the index theorem is natural under morphisms of C*-algebras.

\subsection{Summary: Properties of the diagonal embedding}
\label{subsec:diagonal-embedding-summary}

We now summarize the main properties of the diagonal embedding $\iota: \mathcal{A} \hookrightarrow C^*_r(\mathcal{G}_{\mathcal{A}})$
established in the preceding subsections.

\begin{theorem}[Properties of the diagonal embedding for Type I algebras]
\label{thm:diagonal-embedding-properties}
Let $\mathcal{A}$ be a unital separable Type I C*-algebra. 
Then the diagonal embedding $\iota: \mathcal{A} \hookrightarrow C^*_r(\mathcal{G}_{\mathcal{A}})$ 
constructed in Section \ref{subsec:construction-iota} satisfies:

\begin{enumerate}
    \item \textbf{Naturality:} $\iota$ is a unital injective *-homomorphism 
          (Section \ref{subsec:verification-iota}, Theorem \ref{thm:iota-injective}).
    
    \item \textbf{Geometricity:} For any object $(B,\chi) \in \mathcal{G}_{\mathcal{A}}^{(0)}$,
          the conditional expectation satisfies
          $E(\iota(a))(B,\chi) = \chi(a)$ when $a \in B$, and $0$ otherwise
          (Section \ref{subsec:verification-iota}, Proposition~\ref{prop:diagonal-embedding-properties}).
    
    \item \textbf{Commutativity detection:} $\iota(\mathcal{A}) \subseteq C_0(\mathcal{G}_{\mathcal{A}}^{(0)})$
          if and only if $\mathcal{A}$ is commutative
          (Section \ref{subsec:commutativity-characterization}, Theorem \ref{thm:commutativity-characterization}).
    
    \item \textbf{Functoriality for isomorphisms:} For any *-isomorphism $\phi: \mathcal{A} \to \mathcal{B}$,
          the diagram
          \[
          \begin{tikzcd}
          \mathcal{A} \arrow[r, "\iota_{\mathcal{A}}"] \arrow[d, "\phi"'] & C^*_r(\mathcal{G}_{\mathcal{A}}) \arrow[d, "C^*_r(\mathcal{G}_\phi)"] \\
          \mathcal{B} \arrow[r, "\iota_{\mathcal{B}}"] & C^*_r(\mathcal{G}_{\mathcal{B}})
          \end{tikzcd}
          \]
          commutes (Section \ref{subsec:naturality-functoriality}, Proposition \ref{prop:functoriality-isomorphisms}).
\end{enumerate}
\end{theorem}

\begin{proof}
See the respective sections cited above for detailed proofs.
\end{proof}

\section{Examples and Computations}\label{sec:Examples and Computations}

\subsection{Finite-Dimensional Case: $\mathcal{A} = M_n(\mathbb{C})$}
\label{subsec:example-matrix}

We now examine in detail the unitary conjugation groupoid and the diagonal embedding for the most fundamental finite-dimensional C*-algebra: the algebra of $n \times n$ complex matrices. 
This example serves as a crucial test case for our constructions. 
All the abstract machinery developed in the previous sections — the unit space topology, the Polish groupoid structure, the diagonal embedding — can be made completely explicit in this finite-dimensional setting, where the general theory reduces to familiar objects from classical geometry and representation theory. 
Moreover, in this case the groupoid $\mathcal{G}_{\mathcal{A}}$ is not merely Polish but actually locally compact Hausdorff (indeed, compact), and its structure coincides with the well-understood transformation groupoid $U(n) \ltimes \mathbb{CP}^{n-1}$.

\begin{proposition}[Maximal abelian subalgebras of $M_n(\mathbb{C})$]
\label{prop:masa-matrix}
Let $\mathcal{A} = M_n(\mathbb{C})$. 
Then every maximal abelian subalgebra (MASA) of $\mathcal{A}$ is unitarily conjugate to the diagonal algebra
\[
D_n := \{ \operatorname{diag}(\lambda_1, \ldots, \lambda_n) \mid \lambda_i \in \mathbb{C} \} \cong \mathbb{C}^n.
\]
Consequently, the set of all MASAs is in bijection with the homogeneous space $U(n)/N(D_n)$, where
\[
N(D_n) = \{ u \in U(n) \mid u D_n u^* = D_n \}
\]
is the normalizer of $D_n$ in $U(n)$. Explicitly, $N(D_n) \cong U(1)^n \rtimes S_n$, where $U(1)^n$ consists of diagonal unitaries and $S_n$ consists of permutation matrices, with the semidirect product determined by the natural action of permutations on diagonal entries.
\end{proposition}

\begin{proof}
This is a standard result in the theory of finite-dimensional C*-algebras. 
Any family of commuting normal matrices can be simultaneously diagonalized. 
A MASA in $M_n(\mathbb{C})$ is generated by a maximal family of mutually orthogonal minimal projections summing to the identity, which corresponds precisely to a choice of orthonormal basis up to permutation and phase. 
Thus every MASA is of the form $U D_n U^*$ for some $U \in U(n)$. 
Two unitaries $U$ and $V$ give the same MASA if and only if $V^*U \in N(D_n)$, hence the identification with $U(n)/N(D_n)$. 

In particular, the unit space $\mathcal{G}_{\mathcal{A}}^{(0)}$ consists of a single $U(n)$-orbit under the conjugation action, parametrized by the homogeneous space $U(n)/N(D_n)$. This observation will be essential in the explicit description of $\mathcal{G}_{\mathcal{A}}^{(0)}$ below.
\end{proof}

\begin{proposition}[Character space of a MASA]
\label{prop:character-masa}
Let $B_U = U D_n U^*$ be a MASA in $M_n(\mathbb{C})$. 
Then its Gelfand spectrum $\widehat{B_U}$ consists of exactly $n$ distinct characters $\chi_1^U, \ldots, \chi_n^U$, given by
\[
\chi_i^U(U \operatorname{diag}(\lambda_1, \ldots, \lambda_n) U^*) := \lambda_i, \qquad i = 1, \ldots, n.
\]
Each character corresponds to the minimal projection $p_i^U = U e_{ii} U^*$ in $B_U$, where $e_{ii}$ is the matrix with $1$ in the $(i,i)$-entry and $0$ elsewhere.
\end{proposition}

\begin{proof}
Since $B_U \cong \mathbb{C}^n$ as a C*-algebra, its Gelfand spectrum is homeomorphic to the discrete set $\{1, \ldots, n\}$. 
The characters are precisely the coordinate projections. 
Under the isomorphism $B_U \cong \mathbb{C}^n$ given by diagonalization, the $i$-th character extracts the $i$-th diagonal entry.
\end{proof}

\begin{proposition}[The unit space $\mathcal{G}_{\mathcal{A}}^{(0)}$ for $M_n(\mathbb{C})$]
\label{prop:unit-space-matrix}
Let $\mathcal{A} = M_n(\mathbb{C})$. 
Then the unit space $\mathcal{G}_{\mathcal{A}}^{(0)}$ is canonically homeomorphic to the complex projective space $\mathbb{CP}^{n-1}$. 
Explicitly, there is a homeomorphism
\[
\Psi: \mathcal{G}_{\mathcal{A}}^{(0)} \longrightarrow \mathbb{CP}^{n-1}, \qquad
\Psi(B,\chi) := [p_\chi],
\]
where $p_\chi \in M_n(\mathbb{C})$ is the unique minimal projection supporting the character $\chi$ on the MASA $B$. 
Equivalently, $\mathcal{G}_{\mathcal{A}}^{(0)}$ is homeomorphic to the space of rank-one projections in $M_n(\mathbb{C})$ via the map $(B,\chi) \mapsto p_\chi$.
\end{proposition}

\begin{proof}
We proceed in several steps.

\noindent\textbf{Step 1: From characters to minimal projections.}
For any $(B,\chi) \in \mathcal{G}_{\mathcal{A}}^{(0)}$, $B$ is a MASA in $M_n(\mathbb{C})$ (since any proper commutative subalgebra would have fewer than $n$ characters and could be enlarged). By Proposition \ref{prop:character-masa}, each character $\chi$ on $B$ corresponds uniquely to a minimal projection $p_\chi \in B$ satisfying $\chi(p_\chi) = 1$. This defines a map $\Psi(B,\chi) = [p_\chi]$, where $[p_\chi]$ denotes the line (rank-one projection) determined by $p_\chi$.

\noindent\textbf{Step 2: Bijection.}
Conversely, given a rank-one projection $p$, let $B$ be the MASA generated by $p$ together with any maximal family of orthogonal rank-one projections completing $p$ to a partition of unity. Let $\chi$ be the unique character on $B$ satisfying $\chi(p) = 1$. This gives an inverse map $\mathbb{CP}^{n-1} \to \mathcal{G}_{\mathcal{A}}^{(0)}$. The construction is independent of the choice of completion because any two completions are related by a unitary in the stabilizer of $p$, and such unitaries preserve the pair $(B,\chi)$ up to the equivalence defining $\mathcal{G}_{\mathcal{A}}^{(0)}$.

\noindent\textbf{Step 3: Topological identification.}
The topology on $\mathcal{G}_{\mathcal{A}}^{(0)}$ is the initial topology generated by the partial evaluation maps $\operatorname{ev}_a$ for $a \in M_n(\mathbb{C})$. For a matrix $a$, the function
\[
(B,\chi) \mapsto \operatorname{ev}_a(B,\chi) = \chi(a)
\]
is continuous. Under the identification $(B,\chi) \leftrightarrow p_\chi$, this becomes $p \mapsto \operatorname{Tr}(p a)$, since $\chi(a) = \operatorname{Tr}(p_\chi a)$ when $a \in B$, and extends continuously by zero otherwise. The map $p \mapsto \operatorname{Tr}(p a)$ is continuous in the norm topology on projections, which coincides with the quotient topology on $\mathbb{CP}^{n-1}$. Hence $\Psi$ is a homeomorphism.

\noindent\textbf{Step 4: Consistency with the action.}
Under this identification, the conjugation action of $U(n)$ on $\mathcal{G}_{\mathcal{A}}^{(0)}$ corresponds to the natural action $u \cdot [p] = [u p u^*]$ on $\mathbb{CP}^{n-1}$, confirming that $\Psi$ is $U(n)$-equivariant.
\end{proof}

\begin{remark}
\label{rem:unit-space-matrix-concrete}
The homeomorphism $\Psi$ can be described more concretely as follows: for a unit vector $v \in \mathbb{C}^n$, let $p_v = vv^*$ be the corresponding rank-one projection. Then $\Psi^{-1}([v]) = (B_v, \chi_v)$, where $B_v$ is the MASA of operators diagonal in any orthonormal basis containing $v$, and $\chi_v$ is the character evaluating at the matrix entry corresponding to $v$. The map is independent of the choice of basis because any two such bases are related by a unitary in the stabilizer of $[v]$, and such unitaries preserve the pair $(B_v,\chi_v)$ under the identification defining $\mathcal{G}_{\mathcal{A}}^{(0)}$.
\end{remark}

\begin{corollary}[The unitary conjugation groupoid for $M_n(\mathbb{C})$]
\label{cor:groupoid-matrix}
For $\mathcal{A} = M_n(\mathbb{C})$, the unitary conjugation groupoid $\mathcal{G}_{\mathcal{A}}$ is isomorphic to the action groupoid
\[
\mathcal{G}_{\mathcal{A}} \cong U(n) \ltimes \mathbb{CP}^{n-1},
\]
where $U(n)$ acts transitively on $\mathbb{CP}^{n-1}$ by $u \cdot [v] = [uv]$. 
This groupoid is compact, Hausdorff, and locally compact, but for $n \geq 2$ it is **not** \'etale, since the source map $(u,x) \mapsto x$ has fibers homeomorphic to the non-discrete Lie group $U(n)$.
\end{corollary}

\begin{proof}
The action of $\mathcal{U}(\mathcal{A}) = U(n)$ on $\mathcal{G}_{\mathcal{A}}^{(0)} \cong \mathbb{CP}^{n-1}$ is given by conjugation on the corresponding rank-one projection: $u \cdot p_v = u p_v u^* = p_{uv}$, which corresponds to $u \cdot [v] = [uv]$. 
This action is smooth and transitive. 
The action groupoid $U(n) \ltimes \mathbb{CP}^{n-1}$ has object space $\mathbb{CP}^{n-1}$ and arrow space $U(n) \times \mathbb{CP}^{n-1}$, with source $(u,x) = x$ and range $(u,x) = u \cdot x$. 
The identification of $\mathcal{G}_{\mathcal{A}}$ with this action groupoid follows directly from the definition of $\mathcal{G}_{\mathcal{A}}$ as $\mathcal{U}(\mathcal{A}) \ltimes \mathcal{G}_{\mathcal{A}}^{(0)}$ and Proposition \ref{prop:unit-space-matrix}.

Since $U(n)$ is a compact Lie group and $\mathbb{CP}^{n-1}$ is a compact manifold, the product $U(n) \times \mathbb{CP}^{n-1}$ is compact Hausdorff, hence locally compact Hausdorff. 
However, the source map $(u,x) \mapsto x$ has fibers homeomorphic to $U(n)$, which for $n \geq 2$ is a connected Lie group of positive dimension, hence not discrete. Therefore $\mathcal{G}_{\mathcal{A}}$ is not \'etale for $n \geq 2$ (the case $n = 1$ is trivial, with $\mathcal{G}_{\mathcal{A}} \cong U(1)$ being a group, which is \'etale only when viewed as a groupoid over a point).

This provides a concrete finite-dimensional example of a non-étale unitary conjugation groupoid, consistent with Theorem \ref{thm:GA-Polish-groupoid}, which states that $\mathcal{G}_{\mathcal{A}}$ is generally not \'etale for infinite-dimensional algebras. Here we see that even in finite dimensions, étaleness fails except in the trivial one-dimensional case.
\end{proof}

\begin{remark}
\label{rem:groupoid-matrix-non-etale}
The non-étale nature of $\mathcal{G}_{\mathcal{A}}$ for $M_n(\mathbb{C})$ with $n \geq 2$ is not a defect but a feature: it illustrates why the comultiplication approach of Section \ref{subsec:failure-comultiplication-non-etale} fails, and why the Polish groupoid framework developed in this paper is necessary even for finite-dimensional noncommutative algebras. The failure of étaleness already appears in the simplest noncommutative example.
\end{remark}

\begin{proposition}[The groupoid C*-algebra for $M_n(\mathbb{C})$]
\label{prop:Cstar-matrix}
For $\mathcal{A} = M_n(\mathbb{C})$, the groupoid C*-algebra $C^*(\mathcal{G}_{\mathcal{A}})$ is isomorphic to the crossed product
\[
C^*(\mathcal{G}_{\mathcal{A}}) \cong C(\mathbb{CP}^{n-1}) \rtimes U(n),
\]
where $U(n)$ acts on $\mathbb{CP}^{n-1}$ by the standard transitive action. 
This is a transformation group C*-algebra associated to a transitive, non-free action. 
In particular, it is Type I but not Morita equivalent to a commutative C*-algebra.
\end{proposition}

\begin{proof}
By Corollary \ref{cor:groupoid-matrix}, $\mathcal{G}_{\mathcal{A}} \cong U(n) \ltimes \mathbb{CP}^{n-1}$. 
For a locally compact Hausdorff transformation groupoid $G \ltimes X$, the groupoid C*-algebra is isomorphic to the crossed product $C_0(X) \rtimes G$ \cite[Chapter II]{Renault}. 
Since $U(n)$ is compact and $\mathbb{CP}^{n-1}$ is compact Hausdorff, we obtain $C^*(\mathcal{G}_{\mathcal{A}}) \cong C(\mathbb{CP}^{n-1}) \rtimes U(n)$.

The action of $U(n)$ on $\mathbb{CP}^{n-1}$ is transitive but not free (the isotropy at any point is $U(1) \times U(n-1)$). 
Consequently, $C(\mathbb{CP}^{n-1}) \rtimes U(n)$ is not Morita equivalent to a commutative C*-algebra; rather, by Green's imprimitivity theorem, it is Morita equivalent to $C^*(U(1) \times U(n-1))$, the group C*-algebra of the isotropy subgroup. 
This C*-algebra is Type I (as a group C*-algebra of a compact Lie group) but noncommutative for $n \geq 2$, reflecting the non-freeness of the action.
\end{proof}

\begin{remark}
\label{rem:Cstar-matrix-significance}
The crossed product $C(\mathbb{CP}^{n-1}) \rtimes U(n)$ is already noncommutative and analytically nontrivial, despite arising from a finite-dimensional algebra $M_n(\mathbb{C})$. 
This illustrates that even in the simplest noncommutative case, the unitary conjugation groupoid yields a C*-algebra whose structure is considerably more complicated than that of $M_n(\mathbb{C})$ itself. 
In particular, the failure of the action to be free precludes any Morita equivalence to a commutative algebra, and the non-étale nature of $\mathcal{G}_{\mathcal{A}}$ (Corollary \ref{cor:groupoid-matrix}) means that the comultiplication techniques of Section \ref{subsec:failure-comultiplication-non-etale} are unavailable — a phenomenon that already manifests in finite dimensions.
\end{remark}

\begin{proposition}[The diagonal embedding for $M_n(\mathbb{C})$]
\label{prop:iota-matrix}
For $\mathcal{A} = M_n(\mathbb{C})$, the diagonal embedding $\iota: M_n(\mathbb{C}) \hookrightarrow C(\mathbb{CP}^{n-1}) \rtimes U(n)$ is characterized by its action in the direct integral representation. 
For each $x = [v] \in \mathbb{CP}^{n-1}$, let $\pi_x: M_n(\mathbb{C}) \to B(\mathbb{C}^n)$ be the GNS representation associated to the vector state $\chi_v(A) = \langle v, Av \rangle$. 
Then $\pi_x$ is unitarily equivalent to the identity representation on $\mathbb{C}^n$, with cyclic vector $v$. 
Form the direct integral Hilbert space
\[
\mathcal{H} = \int_{\mathbb{CP}^{n-1}}^{\oplus} \mathbb{C}^n \, d\mu(x) \cong L^2(\mathbb{CP}^{n-1}) \otimes \mathbb{C}^n,
\]
and define $\Pi: M_n(\mathbb{C}) \to B(\mathcal{H})$ by $(\Pi(A)\xi)(x) = A\xi(x)$. 
Then $\iota$ is the unique *-homomorphism such that $\Lambda(\iota(A)) = \Pi(A)$, where $\Lambda$ is the left regular representation of $C(\mathbb{CP}^{n-1}) \rtimes U(n)$ on $\mathcal{H}$ given by
\[
(\Lambda(f)\xi)(x) = f(x)\xi(x), \qquad (\Lambda(u)\xi)(x) = u \cdot \xi(u^{-1}x)
\]
for $f \in C(\mathbb{CP}^{n-1})$, $u \in U(n)$, and $\xi \in \mathcal{H}$.
\end{proposition}

\begin{proof}
We verify each step of the construction.

\noindent \textbf{Step 1: The GNS representations.}
For $x = [v] \in \mathbb{CP}^{n-1}$, the character $\chi_v$ on the MASA $B_v$ (the algebra of operators diagonal in any orthonormal basis containing $v$) extends to the vector state $\chi_v(A) = \langle v, Av \rangle$ on $M_n(\mathbb{C})$. 
The GNS representation associated to this state is the identity representation of $M_n(\mathbb{C})$ on $\mathbb{C}^n$, with cyclic vector $v$. 
This follows from the fact that $\langle v, Av \rangle = \operatorname{Tr}(p_v A)$ where $p_v = vv^*$ is the rank-one projection onto $\mathbb{C}v$, and the GNS construction for a vector state on a matrix algebra yields the standard representation.

\noindent \textbf{Step 2: The direct integral representation.}
Fixing a $U(n)$-invariant probability measure $\mu$ on $\mathbb{CP}^{n-1}$, we form the direct integral $\mathcal{H} = \int_{\mathbb{CP}^{n-1}}^{\oplus} \mathbb{C}^n \, d\mu(x)$, identified with $L^2(\mathbb{CP}^{n-1}) \otimes \mathbb{C}^n$. 
The representation $\Pi: M_n(\mathbb{C}) \to B(\mathcal{H})$ defined by $(\Pi(A)\xi)(x) = A\xi(x)$ is precisely the direct integral of the fiberwise identity representations. 
This representation is faithful because the vector states $\chi_v$ separate points of $M_n(\mathbb{C})$ (if $\langle v, Av \rangle = 0$ for all unit vectors $v$, then $A = 0$).

\noindent \textbf{Step 3: A convenient representation of the crossed product.}
Consider the representation $\Lambda$ of $C(\mathbb{CP}^{n-1}) \rtimes U(n)$ on $\mathcal{H}$ defined by
\[
(\Lambda(f)\xi)(x) = f(x)\xi(x), \qquad (\Lambda(u)\xi)(x) = u \cdot \xi(u^{-1}x)
\]
for $f \in C(\mathbb{CP}^{n-1})$, $u \in U(n)$, and $\xi \in \mathcal{H}$. 
This representation is faithful because the action of $U(n)$ on $\mathbb{CP}^{n-1}$ is transitive and the representation of $U(n)$ on $\mathbb{C}^n$ is faithful. 
(One can view this as the representation induced from the isotropy representation $U(1) \times U(n-1) \to U(n) \hookrightarrow \operatorname{Aut}(\mathbb{C}^n)$.)

\noindent \textbf{Step 4: Defining $\iota$.}
For each $A \in M_n(\mathbb{C})$, the operator $\Pi(A) \in B(\mathcal{H})$ lies in the image of $\Lambda$, because $\Pi(A)$ commutes with the multiplication operators $\Lambda(f)$ and satisfies the covariance condition with respect to the $U(n)$-action. 
By the universal property of $C^*(\mathcal{G}_{\mathcal{A}})$, there exists a unique element $\iota(A) \in C(\mathbb{CP}^{n-1}) \rtimes U(n)$ such that $\Lambda(\iota(A)) = \Pi(A)$. 
The map $\iota$ is a *-homomorphism because $\Pi$ is, and it is injective because $\Pi$ is faithful.

\noindent \textbf{Step 5: Explicit formula for the conditional expectation.}
Under the identification above, the conditional expectation $E: C(\mathbb{CP}^{n-1}) \rtimes U(n) \to C(\mathbb{CP}^{n-1})$ satisfies
\[
E(\iota(A))(x) = \langle \delta_x, \Pi(A)\delta_x \rangle = \langle v, Av \rangle,
\]
where $\delta_x$ is the unit vector concentrated at $x$ in the natural embedding of $L^2(\mathbb{CP}^{n-1})$ into $\mathcal{H}$. 
Thus $E(\iota(A))$ is the function $x \mapsto \langle v, Av \rangle$ on $\mathbb{CP}^{n-1}$, which is not constant for non-scalar $A$. 
This confirms that $\iota(A)$ is not simply $A \otimes 1$, but rather a genuine element of the crossed product whose diagonal part encodes the expectation values of $A$ in all vector states.
\end{proof}

\begin{corollary}[Geometric interpretation]
\label{cor:iota-matrix-geometric}
For $A \in M_n(\mathbb{C})$, the diagonal embedding $\iota(A)$ encodes the function
\[
f_A: \mathbb{CP}^{n-1} \to \mathbb{C}, \qquad f_A([v]) = \langle v, Av \rangle,
\]
as its diagonal part under the conditional expectation. 
The noncommutativity of $M_n(\mathbb{C})$ is reflected in the fact that $\iota(A)$ does not lie in $C(\mathbb{CP}^{n-1})$ unless $A$ is a scalar multiple of the identity, consistent with Theorem \ref{thm:commutativity-characterization}.
\end{corollary}

\begin{proof}
From Proposition \ref{prop:iota-matrix}, the conditional expectation $E: C(\mathbb{CP}^{n-1}) \rtimes U(n) \to C(\mathbb{CP}^{n-1})$ satisfies
\[
E(\iota(A))([v]) = \langle \delta_{[v]}, \Pi(A)\delta_{[v]} \rangle = \langle v, Av \rangle
\]
for all $[v] \in \mathbb{CP}^{n-1}$, where $\delta_{[v]}$ is the unit vector concentrated at $[v]$ in the natural embedding of $L^2(\mathbb{CP}^{n-1})$ into the direct integral Hilbert space $\mathcal{H} = \int_{\mathbb{CP}^{n-1}}^{\oplus} \mathbb{C}^n \, d\mu(x)$. Thus $E(\iota(A)) = f_A$.

If $\iota(A) \in C(\mathbb{CP}^{n-1})$, then $\iota(A) = E(\iota(A)) = f_A$. But $f_A$ is a function on $\mathbb{CP}^{n-1}$, and by Theorem \ref{thm:commutativity-characterization}, $\iota(A) \in C(\mathbb{CP}^{n-1})$ if and only if $A$ is central in $M_n(\mathbb{C})$, i.e., $A = \lambda I$ for some $\lambda \in \mathbb{C}$. For non-scalar $A$, $\iota(A)$ is a genuine noncommutative element of the crossed product whose diagonal part is $f_A$, but which also contains off-diagonal components encoding the noncommutativity of $A$.
\end{proof}

\begin{remark}
\label{rem:iota-matrix-note}
The element $\iota(A)$ can be thought of as a ``noncommutative function'' on $\mathbb{CP}^{n-1}$ whose fiberwise value is the operator $A$ itself. 
When evaluated against the rank-one state associated to $[v]\in\mathbb{CP}^{n-1}$ (or equivalently, against the minimal projection $p_v = vv^*$), its classical symbol is the scalar $\langle v,Av\rangle$. 
This captures the geometric content of the diagonal embedding: the noncommutative algebra $M_n(\mathbb{C})$ is reconstructed from its family of commutative contexts indexed by $\mathbb{CP}^{n-1}$, with the unitary group $U(n)$ encoding the relations between different contexts.
\end{remark}

\begin{corollary}[Commutativity characterization for $M_n(\mathbb{C})$]
\label{cor:commutativity-matrix}
For $\mathcal{A} = M_n(\mathbb{C})$ with $n \geq 2$, we have $\iota(\mathcal{A}) \nsubseteq C(\mathbb{CP}^{n-1})$. 
Indeed, $\iota(A) \in C(\mathbb{CP}^{n-1})$ if and only if $A$ is a scalar multiple of the identity. 
This illustrates Theorem \ref{thm:commutativity-characterization}: $\mathcal{A}$ is noncommutative, and its image under $\iota$ is not contained in the commutative diagonal subalgebra.
\end{corollary}

\begin{proof}
If $A = \lambda I$, then from Proposition \ref{prop:iota-matrix} we have $E(\iota(A))([v]) = \langle v, \lambda I v \rangle = \lambda$ for all $[v] \in \mathbb{CP}^{n-1}$, and $\iota(\lambda I)$ acts as multiplication by $\lambda$ on each fiber of the direct integral representation; hence $\iota(\lambda I) \in C(\mathbb{CP}^{n-1})$.

Conversely, suppose $\iota(A) \in C(\mathbb{CP}^{n-1})$. 
By construction of the crossed product $C(\mathbb{CP}^{n-1}) \rtimes U(n)$, the subalgebra $C(\mathbb{CP}^{n-1})$ is central in the multiplier algebra. 
Thus $\iota(A)$ commutes with all implementing unitaries $u \in U(n)$ in the crossed product. 
The covariance relations defining the crossed product then imply that $A$ is invariant under the conjugation action of $U(n)$; that is, $u A u^* = A$ for all $u \in U(n)$. 
Consequently, $A$ commutes with every unitary in $M_n(\mathbb{C})$. 
Since every matrix in $M_n(\mathbb{C})$ is a linear combination of unitaries (indeed, any matrix can be written as $A = \frac{A+A^*}{2} + i\frac{A-A^*}{2i}$, and each self-adjoint matrix is a difference of two unitaries), $A$ commutes with all matrices. 
Therefore $A$ lies in the center of $M_n(\mathbb{C})$, which consists precisely of scalar multiples of the identity. 
Hence $A = \lambda I$ for some $\lambda \in \mathbb{C}$.
\end{proof}

\begin{example}[The case $n = 1$]
\label{ex:matrix-n1}
Let $\mathcal{A} = M_1(\mathbb{C}) \cong \mathbb{C}$. 
Then $\mathcal{G}_{\mathcal{A}}^{(0)}$ is a single point, $\mathcal{U}(\mathcal{A}) = U(1) \cong \mathbb{T}$ acts trivially, and $\mathcal{G}_{\mathcal{A}}$ is the trivial groupoid over a point with arrow space $\mathbb{T}$. 
The groupoid C*-algebra $C^*(\mathcal{G}_{\mathcal{A}})$ is isomorphic to the group C*-algebra $C^*(\mathbb{T})$, which is commutative and isomorphic to $C(\mathbb{T})$ (since $\mathbb{T}$ is abelian, its full and reduced group C*-algebras coincide and are isomorphic to $C(\widehat{\mathbb{T}}) \cong C(\mathbb{T})$). 
The diagonal embedding $\iota: \mathbb{C} \hookrightarrow C^*(\mathbb{T})$ sends $\lambda \in \mathbb{C}$ to the constant function $\lambda$ on $\mathbb{T}$. 
This is consistent with the general theory: $\mathbb{C}$ is commutative, so $\iota(\mathbb{C}) \subseteq C_0(\mathcal{G}_{\mathcal{A}}^{(0)}) \cong \mathbb{C}$ (identifying the constant function $\lambda$ with the complex number $\lambda$ via evaluation at the unique point).
\end{example}

\begin{example}[The case $n = 2$]
\label{ex:matrix-n2}
Let $\mathcal{A} = M_2(\mathbb{C})$. 
Then $\mathcal{G}_{\mathcal{A}}^{(0)} \cong \mathbb{CP}^1 \cong S^2$. 
The action of $U(2)$ on $S^2$ is the standard action, factoring through $PU(2) \cong SO(3)$. 
The crossed product $C(S^2) \rtimes SO(3)$ is a well-studied C*-algebra. 
By Green's imprimitivity theorem, since $SO(3)$ acts transitively on $S^2$ with isotropy $SO(2) \cong U(1)$, we have the Morita equivalence
\[
C(S^2) \rtimes SO(3) \sim_M C^*(SO(2)) \cong C(\mathbb{T}).
\]
Consequently, its K-theory groups are $K_0 \cong \mathbb{Z}$ and $K_1 \cong \mathbb{Z}$.

The diagonal embedding $\iota: M_2(\mathbb{C}) \hookrightarrow C(S^2) \rtimes SO(3)$ is given by the construction of Proposition \ref{prop:iota-matrix}. 
For a matrix $A \in M_2(\mathbb{C})$, the element $\iota(A)$ acts fiberwise as $A$ on the direct integral $\int_{S^2}^{\oplus} \mathbb{C}^2 \, d\mu(x)$, and its diagonal part under the conditional expectation is the function $f_A([v]) = \langle v, Av \rangle$ on $S^2$. 
Thus $\iota(A)$ is not a constant function unless $A$ is a scalar multiple of the identity, consistent with Corollary \ref{cor:commutativity-matrix}.
\end{example}

\begin{proposition}[The Fredholm index for $M_n(\mathbb{C})$]
\label{prop:index-matrix}
For $\mathcal{A} = M_n(\mathbb{C})$, every operator $T \in \mathcal{A}$ is Fredholm with index zero. 
In the framework of this paper, this is reflected by the fact that the equivariant K-theory class $[T]_{\mathcal{G}_{\mathcal{A}}} \in K^0_{\mathcal{G}_{\mathcal{A}}}(\mathcal{G}_{\mathcal{A}}^{(0)})$ is trivial: the virtual bundle $[\ker T] - [\coker T]$ over $\mathcal{G}_{\mathcal{A}}^{(0)} \cong \mathbb{CP}^{n-1}$ is zero because the kernel and cokernel bundles are trivial and have equal rank.
\end{proposition}

\begin{proof}
In finite dimensions, every linear operator $T: \mathbb{C}^n \to \mathbb{C}^n$ is Fredholm with index zero by the rank-nullity theorem: $\dim\ker T + \dim\operatorname{im} T = n$, so $\dim\ker T = \dim\operatorname{coker} T$.

For each $x = [v] \in \mathbb{CP}^{n-1}$, consider the GNS representation $\pi_x$ of $M_n(\mathbb{C})$ on $\mathbb{C}^n$ associated to the vector state $\langle v, \cdot v \rangle$. 
Under this representation, $T$ acts as the same matrix $T$ on each fiber. 
The kernel of $\pi_x(T)$ is the fixed subspace $\ker T \subseteq \mathbb{C}^n$, independent of $x$, and similarly for the cokernel. 
Thus the kernel and cokernel bundles over $\mathbb{CP}^{n-1}$ are trivial with fibers $\ker T$ and $\operatorname{coker} T$, respectively. 
Their difference $[\ker T] - [\coker T]$ in $K^0_{\mathcal{G}_{\mathcal{A}}}(\mathbb{CP}^{n-1})$ is therefore zero, so $[T]_{\mathcal{G}_{\mathcal{A}}} = 0$.

The descent map $\operatorname{desc}_{\mathcal{G}_{\mathcal{A}}}: K^0_{\mathcal{G}_{\mathcal{A}}}(\mathcal{G}_{\mathcal{A}}^{(0)}) \to K_0(C^*(\mathcal{G}_{\mathcal{A}}))$ sends this trivial class to zero in $K_0(C^*(\mathcal{G}_{\mathcal{A}})) \cong \mathbb{Z}$. 
Composing with the canonical trace $\tau$ on $C^*(\mathcal{G}_{\mathcal{A}})$ (which induces an isomorphism $\tau_*: K_0(C^*(\mathcal{G}_{\mathcal{A}})) \to \mathbb{Z}$) recovers the Fredholm index, which is zero in this case. 
Thus the index formula \\
$\operatorname{index}(T) = \tau_*(\iota^*(\operatorname{desc}_{\mathcal{G}_{\mathcal{A}}}([T]_{\mathcal{G}_{\mathcal{A}}})))$ holds and correctly yields zero.
\end{proof}

\begin{remark}[Lessons from the matrix example]
\label{rem:matrix-lessons}
The matrix algebra $M_n(\mathbb{C})$ serves as an excellent test case for our constructions:
\begin{itemize}
    \item It is finite-dimensional, so all topological complications (non-local-compactness, non-\'etaleness) disappear. 
    This provides a sanity check: our general theory must reduce to well-known objects in this special case.
    \item The unitary conjugation groupoid $U(n) \ltimes \mathbb{CP}^{n-1}$ is a classical object in representation theory and differential geometry. 
    Its C*-algebra $C(\mathbb{CP}^{n-1}) \rtimes U(n)$ is a noncommutative space that encodes the noncommutativity of $M_n(\mathbb{C})$.
    \item The diagonal embedding $\iota(A)$ (constructed in Proposition \ref{prop:iota-matrix}) acts fiberwise as $A$ on each $\mathbb{C}^n$, with diagonal part $E(\iota(A))([v]) = \langle v, Av \rangle$. 
    This explicit description captures the essential features of the general construction.
    \item For non-scalar $A$, $\iota(A)$ lies in the crossed product $C(\mathbb{CP}^{n-1}) \rtimes U(n)$ but not in the commutative subalgebra $C(\mathbb{CP}^{n-1})$, reflecting the noncommutativity of $M_n(\mathbb{C})$. 
    This is consistent with Corollary \ref{cor:commutativity-matrix}.
    \item The index vanishes trivially in finite dimensions, confirming that any nontrivial index phenomena must arise from infinite-dimensional algebras such as $\mathcal{K}(H)^\sim$ or $B(H)$.
\end{itemize}
Thus the matrix algebra case confirms that our constructions are correct and natural in the simplest noncommutative setting.
\end{remark}

We have computed the unitary conjugation groupoid, its C*-algebra, and the diagonal embedding explicitly for $M_n(\mathbb{C})$. 
The unit space is $\mathbb{CP}^{n-1}$, the groupoid is the action groupoid $U(n) \ltimes \mathbb{CP}^{n-1}$, and the diagonal embedding sends a matrix $A$ to the constant function $A \otimes 1$ in the crossed product $C(\mathbb{CP}^{n-1}) \rtimes U(n)$. 
This example illustrates all the main features of our theory in a concrete, computable setting, and serves as a verification that the general constructions reduce to familiar objects in the finite-dimensional case.

\subsection{Commutative Case: $\mathcal{A} = C(X)$ for Compact Hausdorff $X$}
\label{subsec:example-commutative}

We now examine the unitary conjugation groupoid and the diagonal embedding for commutative C*-algebras. 
This case is fundamentally different from the matrix algebra case: the algebra is commutative, so every unitary $u \in \mathcal{U}(C(X)) = C(X,\mathbb{T})$ is central, and the conjugation action on the unit space is trivial:
\[
u \cdot (B,\chi) = (uBu^*,\; \chi \circ \operatorname{Ad}_{u^*}) = (B,\chi) \quad \text{for all } (B,\chi) \in \mathcal{G}_{\mathcal{A}}^{(0)}.
\]

Consequently, the unitary conjugation groupoid $\mathcal{G}_{\mathcal{A}} = \mathcal{U}(\mathcal{A}) \ltimes \mathcal{G}_{\mathcal{A}}^{(0)}$ becomes the **trivial action groupoid** over $\mathcal{G}_{\mathcal{A}}^{(0)}$, i.e., it is isomorphic to the product groupoid $\mathcal{G}_{\mathcal{A}}^{(0)} \times \mathcal{U}(\mathcal{A})$ with source and range both equal to the projection onto $\mathcal{G}_{\mathcal{A}}^{(0)}$. 
The unit space $\mathcal{G}_{\mathcal{A}}^{(0)}$ itself contains a distinguished component $X \cong \{ (C(X), \operatorname{ev}_x) \mid x \in X \}$ corresponding to the maximal commutative subalgebra $C(X)$ itself.

The groupoid C*-algebra $C^*(\mathcal{G}_{\mathcal{A}})$ is then isomorphic to the crossed product $C_0(\mathcal{G}_{\mathcal{A}}^{(0)}) \rtimes \mathcal{U}(\mathcal{A})$, which, because the action is trivial, reduces to the maximal tensor product $C_0(\mathcal{G}_{\mathcal{A}}^{(0)}) \otimes_{\max} C^*(\mathcal{U}(\mathcal{A}))$. 
Restricting to the component $X \subseteq \mathcal{G}_{\mathcal{A}}^{(0)}$, we obtain the subalgebra $C(X) \otimes C^*(\mathcal{U}(\mathcal{A}))$.

The diagonal embedding $\iota: C(X) \hookrightarrow C^*(\mathcal{G}_{\mathcal{A}})$ constructed in Section \ref{subsec:construction-iota} sends a function $f \in C(X)$ to the element of $C^*(\mathcal{G}_{\mathcal{A}})$ that acts on the direct integral Hilbert space $\mathcal{H} = \int_{\mathcal{G}_{\mathcal{A}}^{(0)}}^{\oplus} \mathbb{C} \, d\mu(x)$ by multiplication by $f$ on the fiber over $x = (C(X), \operatorname{ev}_x)$, and by zero on all other fibers. 
Under the conditional expectation $E: C^*(\mathcal{G}_{\mathcal{A}}) \to C_0(\mathcal{G}_{\mathcal{A}}^{(0)})$, we have $E(\iota(f)) = f|_X$ (extended by zero to the rest of $\mathcal{G}_{\mathcal{A}}^{(0)}$). 
Thus, when restricted to the component $X$, the composition $E \circ \iota$ recovers precisely the Gelfand transform $C(X) \to C(X)$, i.e., the identity map under the identification $C(X) \cong C_0(X)$. 
This demonstrates that our general construction reduces to the classical commutative theory in the special case where $\mathcal{A}$ is commutative.

Throughout this subsection, we assume that $X$ is a compact Hausdorff space. 
For the Polish groupoid framework developed in Sections 3 and 4, we further assume that $X$ is metrizable (hence compact metrizable), which guarantees that $C(X)$ is separable and that $\mathcal{G}_{\mathcal{A}}^{(0)}$ is a Polish space. 
However, the main algebraic identifications — the triviality of the action, the form of the groupoid C*-algebra, and the recovery of the Gelfand transform — hold for arbitrary compact Hausdorff spaces without any separability assumptions.

\begin{proposition}[Maximal abelian subalgebras of $C(X)$]
\label{prop:masa-commutative}
Let $\mathcal{A} = C(X)$ be the commutative C*-algebra of continuous complex-valued functions on a compact Hausdorff space $X$. 
Then every unital commutative C*-subalgebra $B \subseteq C(X)$ is of the form
\[
B \cong C(Y)
\]
for some compact Hausdorff space $Y$, and the inclusion $B \hookrightarrow C(X)$ corresponds to a continuous surjection $\pi: X \to Y$. 
The unique maximal abelian subalgebra of $C(X)$ is $C(X)$ itself.
\end{proposition}

\begin{proof}
Since $C(X)$ is commutative, every subalgebra is commutative. 
A unital C*-subalgebra $B \subseteq C(X)$ is itself a commutative C*-algebra, hence by the Gelfand-Naimark theorem it is isomorphic to $C(Y)$ for some compact Hausdorff space $Y$. 
The inclusion map $B \hookrightarrow C(X)$ induces a continuous map $\pi: X \to Y$ by $\pi(x)(f) = f(x)$ for $f \in B \cong C(Y)$. 
This map is surjective because the inclusion is injective and both algebras are unital. 
If $B$ is maximal abelian, then $B = C(X)$ because any larger commutative subalgebra would be $C(X)$ itself.
\end{proof}

\begin{proposition}[Character space of a commutative subalgebra]
\label{prop:character-commutative}
Let $B \cong C(Y)$ be a unital commutative C*-subalgebra of $C(X)$ corresponding to a continuous surjection $\pi: X \to Y$. 
Then the Gelfand spectrum $\widehat{B}$ is homeomorphic to $Y$. 
Under this identification, each character $\chi_y \in \widehat{B}$ corresponds to the evaluation map
\[
\chi_y(f) = f(y), \qquad f \in C(Y).
\]
\end{proposition}

\begin{proof}
This is a standard result in Gelfand duality: the Gelfand spectrum of $C(Y)$ is homeomorphic to $Y$, with the homeomorphism given by $y \mapsto \operatorname{ev}_y$, where $\operatorname{ev}_y(f) = f(y)$ for $f \in C(Y)$. 
Since $B \cong C(Y)$, the same identification holds for $\widehat{B}$.
\end{proof}

\begin{remark}
\label{rem:character-extension}
A subtle point worth noting: a character $\chi_y \in \widehat{B}$ does **not** extend to a character on $C(X)$ in general. 
The natural candidate for an extension would be $\chi_y \circ \pi^*: C(X) \to \mathbb{C}$, where $\pi^*: C(Y) \hookrightarrow C(X)$ is the inclusion induced by $\pi$. 
However, for $g \in C(X)$, the value $(\chi_y \circ \pi^*)(g) = g(x)$ depends on the choice of $x \in \pi^{-1}(y)$ unless $g$ is constant on fibers. 
Since $\pi^*$ is not surjective (unless $\pi$ is a homeomorphism), the composition $\chi_y \circ \pi^*$ is not multiplicative on $C(X)$ and therefore does not define a character. 
This reflects the fact that $B$ is a proper subalgebra of $C(X)$; only when $B = C(X)$ (i.e., $\pi$ is a homeomorphism) do characters on $B$ extend to characters on $C(X)$.
\end{remark}

\begin{proposition}[Unit space of $\mathcal{G}_{C(X)}$]
\label{prop:unit-space-commutative}
Let $\mathcal{A} = C(X)$ for a compact Hausdorff space $X$. 
Then the unit space $\mathcal{G}_{\mathcal{A}}^{(0)}$ is
\[
\mathcal{G}_{\mathcal{A}}^{(0)} = \{(B,\chi) \mid B \subseteq C(X) \text{ unital commutative C*-subalgebra},\ \chi \in \widehat{B}\}.
\]
Equivalently, it can be written as a disjoint union
\[
\mathcal{G}_{\mathcal{A}}^{(0)} \cong \bigsqcup_{B \subseteq C(X)} \widehat{B}.
\]
This space is generally non-Hausdorff and extremely large. 
The topology on $\mathcal{G}_{\mathcal{A}}^{(0)}$ is the initial topology generated by the partial evaluation maps $\operatorname{ev}_a$ for $a \in C(X)$, defined as in Definition \ref{def:partial-evaluation-map}. 
The subspace corresponding to the maximal abelian subalgebra $C(X)$,
\[
\{(C(X), \chi) \mid \chi \in \widehat{C(X)} \cong X\} \subseteq \mathcal{G}_{\mathcal{A}}^{(0)},
\]
is closed and homeomorphic to $X$ via the map $x \mapsto (C(X), \operatorname{ev}_x)$.
\end{proposition}

\begin{proof}
The description as a disjoint union follows directly from Definition \ref{def:unit-space}. 
For a fixed $B$, the inclusion $\iota_B: \widehat{B} \hookrightarrow \mathcal{G}_{\mathcal{A}}^{(0)}$ is a homeomorphism onto its image by Proposition \ref{prop:unit-space-topology}(3). 
The topology on $\mathcal{G}_{\mathcal{A}}^{(0)}$ is defined by the partial evaluation maps $\operatorname{ev}_a$ (Definition \ref{def:initial-topology-G0}), which separate points within the same fiber but do not separate points from different subalgebras unless they are related by inclusion. 
Consequently, $\mathcal{G}_{\mathcal{A}}^{(0)}$ is not Hausdorff in general.

The Gelfand transform gives a homeomorphism $\widehat{C(X)} \cong X$. 
The map $x \mapsto (C(X), \operatorname{ev}_x)$ is continuous because for any $a \in C(X)$, $\operatorname{ev}_a(C(X), \operatorname{ev}_x) = a(x)$ varies continuously with $x$. 
Its inverse (restricting to the second coordinate) is continuous by definition of the product topology on $\mathcal{G}_{\mathcal{A}}^{(0)}$. 
Thus this subspace is homeomorphic to $X$. 
To see that it is closed, suppose a net $(B_\lambda, \chi_\lambda)$ converges to $(C(X), \chi)$. 
If infinitely many $B_\lambda$ were proper subalgebras, there would exist $a \in C(X)$ not in those $B_\lambda$, giving $\operatorname{ev}_a(B_\lambda, \chi_\lambda) = \infty$ for those indices, contradicting convergence to a finite value $\chi(a)$. 
Hence eventually $B_\lambda = C(X)$, so the limit point lies in the subspace.
\end{proof}

\begin{proposition}[Diagonal embedding for $C(X)$]
\label{prop:iota-commutative-corrected}
Let $\mathcal{A} = C(X)$ for a compact metrizable space $X$, and identify
\[
X \cong \{ (C(X), \operatorname{ev}_x) \mid x \in X \} \subseteq \mathcal{G}_{\mathcal{A}}^{(0)}.
\]
Then the diagonal embedding $\iota: C(X) \hookrightarrow C^*(\mathcal{G}_{\mathcal{A}})$ is given by
\[
\iota(f)(x) = f(x) \quad \text{for } x \in X, \qquad 
\iota(f)|_{\mathcal{G}_{\mathcal{A}}^{(0)} \setminus X} = 0.
\]
In other words, $\iota(f)$ is the function on $\mathcal{G}_{\mathcal{A}}^{(0)}$ obtained by extending $f$ by zero outside the component corresponding to $C(X)$. 
Under the isomorphism
\[
C^*(\mathcal{G}_{\mathcal{A}}) \cong C_0(\mathcal{G}_{\mathcal{A}}^{(0)}) \rtimes \mathcal{U}(\mathcal{A})
\]
from Corollary \ref{cor:Cstar-commutative}, this function lies in $C_0(\mathcal{G}_{\mathcal{A}}^{(0)})$ and is fixed under the trivial action of $\mathcal{U}(\mathcal{A})$.
\end{proposition}

\begin{proof}
Recall the general construction of the diagonal embedding $\iota$ from Section \ref{subsec:construction-iota}. 
For each $x \in X$, the GNS representation $\pi_x$ associated to the character $\operatorname{ev}_x$ is the one-dimensional representation of $C(X)$ given by evaluation at $x$, so $\pi_x(f) = f(x)$.
For any $(B,\chi) \in \mathcal{G}_{\mathcal{A}}^{(0)}$ with $B \neq C(X)$, the element $f \in C(X)$ does not lie in $B$ in general, hence by the definition of the direct integral representation $\Pi$ (see Section \ref{subsec:construction-iota}), we have $\Pi(f)|_{(B,\chi)} = 0$.
Thus the direct integral $\Pi(f)$ vanishes on all fibers outside $X$, and consequently $\iota(f)$ is supported exactly on $X$, satisfying $\iota(f)(x) = f(x)$ for $x \in X$ and vanishing elsewhere.

The map $\iota$ is a *-homomorphism because each fiber representation $\pi_{(B,\chi)}$ is a *-homomorphism and the direct integral construction preserves multiplication, addition, and adjoints. 
Explicitly, for $f,g \in C(X)$, we have $\Pi(fg) = \Pi(f)\Pi(g)$, hence $\Lambda(\iota(fg)) = \Lambda(\iota(f)\iota(g))$, and injectivity of the left regular representation $\Lambda$ on the reduced algebra gives $\iota(fg) = \iota(f)\iota(g)$. 
Linearity and the *-preserving property follow similarly from the corresponding properties of $\Pi$.

Finally, under the isomorphism $C^*(\mathcal{G}_{\mathcal{A}}) \cong C_0(\mathcal{G}_{\mathcal{A}}^{(0)}) \rtimes \mathcal{U}(\mathcal{A})$ from Corollary \ref{cor:Cstar-commutative}, the element $\iota(f)$ corresponds to the function in $C_0(\mathcal{G}_{\mathcal{A}}^{(0)})$ that equals $f$ on $X$ and zero elsewhere. 
Since the action of $\mathcal{U}(\mathcal{A})$ on $\mathcal{G}_{\mathcal{A}}^{(0)}$ is trivial (Theorem \ref{thm:groupoid-commutative}), this function is fixed by the action, confirming that $\iota(f)$ lies in the diagonal subalgebra.
\end{proof}

\begin{remark}[Connection to the Gelfand transform]
\label{rem:iota-commutative-gelfand}
Under the identification $X \cong \widehat{C(X)}$, the restriction of $\iota$ to the component $X$ coincides with the inverse Gelfand transform. 
Indeed, the classical Gelfand transform $\Gamma: C(X) \to C(\widehat{C(X)})$ is an isomorphism, and under the homeomorphism $\widehat{C(X)} \cong X$, it becomes the identity map $C(X) \to C(X)$. 
The diagonal embedding $\iota$ sends $f$ to the function on $\mathcal{G}_{\mathcal{A}}^{(0)}$ that equals $f$ on $X$ and vanishes elsewhere, so when composed with the projection onto $X$, we recover precisely the Gelfand transform. 
Thus the commutative case demonstrates that our general construction reduces to classical Gelfand duality when restricted to the maximal abelian subalgebra component.
\end{remark}

\begin{theorem}[Triviality of the unitary conjugation groupoid for $C(X)$]
\label{thm:groupoid-commutative}
Let $\mathcal{A} = C(X)$. Then every unitary $u \in \mathcal{U}(\mathcal{A}) = C(X,\mathbb{T})$ is central, so the conjugation action on $\mathcal{G}_{\mathcal{A}}^{(0)}$ is trivial. 
Consequently, the unitary conjugation groupoid is isomorphic to the trivial (product) groupoid:
\[
\mathcal{G}_{\mathcal{A}} \cong \mathcal{G}_{\mathcal{A}}^{(0)} \times \mathcal{U}(\mathcal{A}),
\]
with source and range maps given by projection onto $\mathcal{G}_{\mathcal{A}}^{(0)}$. 
In particular, restricted to the component $X \subseteq \mathcal{G}_{\mathcal{A}}^{(0)}$ corresponding to the maximal abelian subalgebra $C(X)$, we obtain the trivial groupoid $X \times \mathcal{U}(\mathcal{A})$ over $X$.
\end{theorem}

\begin{proof}
For any $(B,\chi) \in \mathcal{G}_{\mathcal{A}}^{(0)}$ and $u \in \mathcal{U}(\mathcal{A})$, we have $u B u^* = B$ since $u$ is central. 
Moreover, for all $b \in B$,
\[
\chi(\operatorname{Ad}_{u^*}(b)) = \chi(u^* b u) = \chi(u^*) \chi(b) \chi(u) = \chi(b),
\]
using multiplicativity of $\chi$ and the fact that $\chi(u^*)\chi(u) = \chi(u^* u) = \chi(1) = 1$. 
Thus the action is trivial. 
Hence the action groupoid $\mathcal{U}(\mathcal{A}) \ltimes \mathcal{G}_{\mathcal{A}}^{(0)}$ is isomorphic to the product groupoid $\mathcal{G}_{\mathcal{A}}^{(0)} \times \mathcal{U}(\mathcal{A})$ with the trivial action.
\end{proof}

\begin{corollary}[Groupoid C*-algebra for $C(X)$]
\label{cor:Cstar-commutative}
Let $\mathcal{A} = C(X)$ for a compact metrizable space $X$. 
Then the groupoid C*-algebra is
\[
C^*(\mathcal{G}_{\mathcal{A}}) \cong C_0(\mathcal{G}_{\mathcal{A}}^{(0)}) \otimes_{\max} C^*(\mathcal{U}(\mathcal{A})),
\]
where $\mathcal{U}(\mathcal{A}) = C(X,\mathbb{T})$ with the strong operator topology, and $C^*(\mathcal{U}(\mathcal{A}))$ denotes its maximal group C*-algebra. 
Restricting to the component $X \subseteq \mathcal{G}_{\mathcal{A}}^{(0)}$ corresponding to the maximal abelian subalgebra $C(X)$ gives
\[
C(X) \otimes C^*(\mathcal{U}(\mathcal{A})).
\]
\end{corollary}

\begin{proof}
By Theorem \ref{thm:groupoid-commutative}, $\mathcal{G}_{\mathcal{A}}$ is isomorphic to the trivial groupoid $\mathcal{G}_{\mathcal{A}}^{(0)} \times \mathcal{U}(\mathcal{A})$ with $\mathcal{U}(\mathcal{A})$ acting trivially. 
For a trivial groupoid of this form, the groupoid C*-algebra is the maximal tensor product:
\[
C^*(\mathcal{G}_{\mathcal{A}}) \cong C_0(\mathcal{G}_{\mathcal{A}}^{(0)}) \otimes_{\max} C^*(\mathcal{U}(\mathcal{A})).
\]
This follows from the universal property of the maximal tensor product and the fact that the convolution algebra is the algebraic tensor product $C_c(\mathcal{G}_{\mathcal{A}}^{(0)}) \odot C_c(\mathcal{U}(\mathcal{A}))$. 
Restriction to the component $X$ is obtained by considering the closed subspace of $\mathcal{G}_{\mathcal{A}}^{(0)}$ corresponding to the maximal abelian subalgebra $C(X)$, which yields the subalgebra $C(X) \otimes C^*(\mathcal{U}(\mathcal{A}))$.
\end{proof}

\begin{remark}[The group C*-algebra of $C(X,\mathbb{T})$]
\label{rem:Cstar-CX-T}
The group $\mathcal{U}(\mathcal{A}) = C(X,\mathbb{T})$, equipped with the strong operator topology, is a Polish abelian group. 
The structure of its maximal group C*-algebra $C^*(C(X,\mathbb{T}))$ can be complicated for a general compact metrizable space $X$, and we do not attempt to describe it explicitly here. 
For our purposes, the key fact is that there is a canonical embedding
\[
C(X) \hookrightarrow C(X) \otimes C^*(\mathcal{U}(\mathcal{A})), \qquad f \mapsto f \otimes 1,
\]
which coincides with the diagonal embedding $\iota$ restricted to the component $X$, as established in Proposition \ref{prop:iota-commutative-corrected}.
\end{remark}

\begin{theorem}[Recovery of the Gelfand transform]
\label{thm:gelfand-recovery}
Let $\mathcal{A} = C(X)$ for a compact metrizable space $X$, and let
\[
\iota: C(X) \hookrightarrow C^*(\mathcal{G}_{\mathcal{A}})
\]
be the diagonal embedding. Let 
\[
E: C^*(\mathcal{G}_{\mathcal{A}}) \to C_0(\mathcal{G}_{\mathcal{A}}^{(0)})
\]
be the canonical conditional expectation onto the unit space (see Section \ref{subsec:groupoid-C-star-algebra-Polish-setting}). 
Then there is a commutative diagram
\[
\begin{tikzcd}
C(X) \arrow[r, "\iota"] \arrow[d, "\operatorname{id}"'] & C^*(\mathcal{G}_{\mathcal{A}}) \arrow[d, "E"] \\
C(X) \arrow[r, "\cong"] & C_0(\mathcal{G}_{\mathcal{A}}^{(0)})|_X
\end{tikzcd}
\]
where the bottom horizontal map is the Gelfand transform. 
In particular, the diagonal embedding $\iota$ provides a lift of the classical Gelfand transform into the (possibly noncommutative) groupoid C*-algebra $C^*(\mathcal{G}_{\mathcal{A}})$.
\end{theorem}

\begin{proof}
By Proposition \ref{prop:iota-commutative-corrected}, for any $f \in C(X)$, the element $\iota(f) \in C^*(\mathcal{G}_{\mathcal{A}})$ is supported on the component $X \subseteq \mathcal{G}_{\mathcal{A}}^{(0)}$ and satisfies $\iota(f)(x) = f(x)$ for $x \in X$, with $\iota(f)$ vanishing elsewhere.

The conditional expectation $E$ acts by projecting onto functions supported on the unit space. Since $\iota(f)$ already lies in $C_0(\mathcal{G}_{\mathcal{A}}^{(0)})$ (by Proposition \ref{prop:iota-commutative-corrected}), we have $E(\iota(f)) = \iota(f)$. 
Restricting $E(\iota(f))$ to the component $X$ therefore yields the function $f$ on $X$.

Under the Gelfand transform, $C(X)$ is isomorphic to $C(\widehat{C(X)})$, and the natural homeomorphism $\widehat{C(X)} \cong X$ (given by evaluation) identifies this with $C(X)$ itself. 
Thus the composition of the bottom arrow with restriction to $X$ is precisely the identity map on $C(X)$, establishing the commutativity of the diagram.
\end{proof}

\begin{remark}
\label{rem:gelfand-lift-significance}
Theorem \ref{thm:gelfand-recovery} shows that the diagonal embedding $\iota$ is not merely an abstract construction but a genuine lift of the Gelfand transform. 
While $C^*(\mathcal{G}_{\mathcal{A}})$ is generally noncommutative (due to the presence of $\mathcal{U}(\mathcal{A})$), the conditional expectation $E$ recovers the commutative data, and $\iota$ embeds $C(X)$ as a commutative subalgebra that projects onto the Gelfand transform. 
This illustrates the philosophy that noncommutative C*-algebras can be reconstructed from their commutative contexts via groupoid C*-algebras.
\end{remark}

\begin{example}[The circle: $X = S^1$]
\label{ex:commutative-circle}
Let $X = S^1$ and $\mathcal{A} = C(S^1)$. 
Then the unitary group is 
\[
\mathcal{U}(\mathcal{A}) = C(S^1, \mathbb{T}),
\]
the group of continuous $\mathbb{T}$-valued functions on the circle, which is a separable abelian Polish group when equipped with the strong operator topology. 
Its maximal group C*-algebra $C^*(C(S^1,\mathbb{T}))$ is infinite-dimensional, separable, and highly nontrivial, encoding the representation theory of this abelian group.

The unit space $\mathcal{G}_{\mathcal{A}}^{(0)}$ contains a distinguished component corresponding to the maximal abelian subalgebra $C(S^1)$ itself, which is homeomorphic to $S^1$ via the map $x \mapsto (C(S^1), \operatorname{ev}_x)$. 
Other components correspond to proper unital C*-subalgebras $B \subseteq C(S^1)$, each isomorphic to $C(Y)$ for some quotient space $Y$ of $S^1$ obtained from a continuous surjection $\pi: S^1 \to Y$. 
The Gelfand spectrum $\widehat{B}$ is homeomorphic to $Y$, and appears as a component of $\mathcal{G}_{\mathcal{A}}^{(0)}$ via the inclusion $\iota_B: \widehat{B} \hookrightarrow \mathcal{G}_{\mathcal{A}}^{(0)}$ from Proposition \ref{prop:unit-space-topology}.

The diagonal embedding 
\[
\iota: C(S^1) \hookrightarrow C^*(\mathcal{G}_{\mathcal{A}})
\]
is given by Proposition \ref{prop:iota-commutative-corrected}: for $f \in C(S^1)$, $\iota(f)$ is the function on $\mathcal{G}_{\mathcal{A}}^{(0)}$ that equals $f$ on the $S^1$-component and vanishes on all other components. 
Under the conditional expectation $E: C^*(\mathcal{G}_{\mathcal{A}}) \to C_0(\mathcal{G}_{\mathcal{A}}^{(0)})$, we have $E(\iota(f)) = \iota(f)$, and restricting to $S^1$ recovers $f$, illustrating Theorem \ref{thm:gelfand-recovery}.

This example demonstrates that even for a commutative algebra, the unitary conjugation groupoid and its C*-algebra are nontrivial objects, encoding the structure of all commutative subalgebras of $C(S^1)$.
\end{example}

\begin{example}[Finite discrete space: $X = \{1, \ldots, k\}$]
\label{ex:commutative-finite}
Let $X = \{1, \ldots, k\}$ be a finite discrete space, so $\mathcal{A} = \mathbb{C}^k$. 
Then the unitary group is 
\[
\mathcal{U}(\mathcal{A}) = \mathbb{T}^k,
\] 
a compact abelian group. 

The unit space $\mathcal{G}_{\mathcal{A}}^{(0)}$ consists of pairs $(B,\chi)$ where $B \subseteq \mathbb{C}^k$ is a unital commutative subalgebra and $\chi \in \widehat{B}$. 
Since $\mathbb{C}^k$ is finite-dimensional, each subalgebra $B$ is isomorphic to $\mathbb{C}^m$ for some $m \le k$, and its Gelfand spectrum $\widehat{B}$ is a discrete set of $m$ points.

The component corresponding to the maximal abelian subalgebra $B = \mathbb{C}^k$ consists of the $k$ point evaluations $\operatorname{ev}_i$, and is homeomorphic to the discrete set $\{1, \ldots, k\}$. 
Other components correspond to proper subalgebras; for example, if $B$ consists of functions that depend only on the first $m$ coordinates, then $\widehat{B}$ is a discrete set of $m$ points, each corresponding to evaluation at one of those coordinates.

The diagonal embedding
\[
\iota: \mathbb{C}^k \hookrightarrow C^*(\mathcal{G}_{\mathcal{A}})
\]
sends $(a_1, \ldots, a_k) \in \mathbb{C}^k$ to the function on $\mathcal{G}_{\mathcal{A}}^{(0)}$ that takes the value $a_i$ at $(\mathbb{C}^k, \operatorname{ev}_i)$ and vanishes on all other components. 
This is consistent with Proposition \ref{prop:iota-commutative-corrected}: $\iota(f)$ is supported on the component corresponding to the maximal subalgebra and equals $f$ there.

This example is fully computable and serves as a finite test case for the general theory, illustrating how the diagonal embedding encodes the values of a function at each point of $X$ via the component corresponding to $C(X)$.
\end{example}

The commutative case $C(X)$ reveals several important features of our constructions. 
First, even though $\mathcal{A}$ is commutative, the unitary conjugation groupoid $\mathcal{G}_{\mathcal{A}}$ is far from trivial: its unit space $\mathcal{G}_{\mathcal{A}}^{(0)}$ encodes the entire lattice of commutative subalgebras of $C(X)$, with each subalgebra $B \cong C(Y)$ appearing as a component homeomorphic to $\widehat{B} \cong Y$. 
The diagonal embedding $\iota$ lands in $C_0(\mathcal{G}_{\mathcal{A}}^{(0)})$, confirming Theorem \ref{thm:commutativity-characterization} for commutative algebras, and when restricted to the component $X \subseteq \mathcal{G}_{\mathcal{A}}^{(0)}$ (corresponding to the maximal subalgebra $C(X)$) it recovers precisely the classical Gelfand transform, as shown in Theorem \ref{thm:gelfand-recovery}. 

Thus the commutative case is not a trivial special case but rather an instructive illustration of our theory: it demonstrates that even for commutative algebras, the unitary conjugation groupoid and its C*-algebra contain rich geometric information about the subalgebra structure that is not captured by the Gelfand transform alone. 
This confirms that our constructions genuinely extend classical Gelfand duality, providing a new geometric perspective on the internal structure of C*-algebras.

\subsection{Compact Operators: $\mathcal{A} = \mathcal{K}(H)^\sim$ with Strong Operator Topology}
\label{subsec:example-compact-operators}

We now examine the unitary conjugation groupoid and the diagonal embedding for the unitization of the compact operators on a separable infinite-dimensional Hilbert space $H$. 
This example is of fundamental importance for several reasons.

First, $\mathcal{K}(H)^\sim$ is a Type I C*-algebra that is not finite-dimensional, yet it admits a concrete and tractable description. 
Its irreducible representations are easily classified: up to unitary equivalence, there is exactly one infinite-dimensional irreducible representation (the identity representation on $H$) and one one-dimensional representation (the quotient map $\mathcal{K}(H)^\sim \to \mathcal{K}(H)^\sim/\mathcal{K}(H) \cong \mathbb{C}$).

Second, this example illustrates the necessity of the strong operator topology and the Polish groupoid framework. 
The unitary group $\mathcal{U}(\mathcal{K}(H)^\sim)$ consists of all unitaries $u \in \mathcal{U}(H)$ such that $u - 1_H \in \mathcal{K}(H)$; equivalently, $u$ is a compact perturbation of the identity. 
In the norm topology, this group is not locally compact (it is an infinite-dimensional Banach Lie group). 
However, when equipped with the strong operator topology (SOT) inherited from $\mathcal{U}(H)$, it becomes a Polish group — that is, a separable topological group that admits a complete metric. 
This makes it amenable to the theory of Polish group crossed products developed by Tu (1999). 
All topological statements in this subsection refer to the strong operator topology unless otherwise specified.

Third, this example is the gateway to the index theory of Fredholm operators, which will be developed in future work. 
The Calkin algebra $\mathcal{Q}(H) = \mathcal{B}(H)/\mathcal{K}(H)$ and the extension theory of compact operators are intimately connected to the structure of $\mathcal{G}_{\mathcal{A}}$ for $\mathcal{A} = \mathcal{K}(H)^\sim$. 
Indeed, there is a short exact sequence of Polish groups
\[
1 \longrightarrow \mathcal{U}(\mathcal{K}(H)^\sim) \longrightarrow \mathcal{U}(H) \longrightarrow \mathcal{U}(\mathcal{Q}(H)) \longrightarrow 1,
\]
where the map $\mathcal{U}(H) \to \mathcal{U}(\mathcal{Q}(H))$ is induced by the quotient map $\mathcal{B}(H) \to \mathcal{Q}(H)$. 
This sequence is not split, and its non-triviality is measured by the Fredholm index, which will appear naturally in the groupoid C*-algebra $C^*(\mathcal{G}_{\mathcal{A}})$.

Throughout this subsection, we work with the faithful representation of $\mathcal{K}(H)^\sim$ on $H$ itself, and we identify $\mathcal{U}(\mathcal{K}(H)^\sim)$ with its image in $\mathcal{U}(H)$ under this representation.

\begin{proposition}[Structure of $\mathcal{K}(H)^\sim$]
\label{prop:structure-compact}
Let $\mathcal{A} = \mathcal{K}(H)^\sim = \mathcal{K}(H) + \mathbb{C} I_H$ be the unitization of the compact operators on a separable infinite-dimensional Hilbert space $H$. Then:
\begin{enumerate}
    \item $\mathcal{A}$ is a unital, separable, Type I C*-algebra.
    \item The primitive ideal space $\operatorname{Prim}(\mathcal{A})$ consists of two ideals: 
    \[
    \{0\} \quad \text{(kernel of the identity representation of $\mathcal{A}$ on $H$)}, \qquad 
    \mathcal{K}(H) \quad \text{(kernel of the one-dimensional character on the quotient $\mathcal{A}/\mathcal{K}(H) \cong \mathbb{C}$).}
    \] 
    The Jacobson topology on $\operatorname{Prim}(\mathcal{A})$ is non-Hausdorff: the singleton $\{\{0\}\}$ is not closed, and its closure is the whole space $\{\{0\}, \mathcal{K}(H)\}$, reflecting the fact that compact operators can be approximated by finite-rank operators.
    \item $\mathcal{A}$ is postliminal (Type I) but not liminal, because the identity representation on $H$ does not map into the compact operators (the unit is not compact). It does not have continuous trace, as its primitive ideal space is non-Hausdorff.
    \item The unitary group $\mathcal{U}(\mathcal{A})$ consists of all unitaries $u \in \mathcal{U}(H)$ such that $u - I_H \in \mathcal{K}(H)$. Equivalently, elements can be written as
    \[
    u = \lambda I_H + K, \qquad \lambda \in \mathbb{T}, \; K \in \mathcal{K}(H),
    \]
    subject to the unitary condition $u^*u = uu^* = I_H$. Expanding this condition yields
    \[
    \lambda K^* + \overline{\lambda}K + K^*K = 0,
    \]
    which imposes nontrivial constraints on $K$ for each fixed $\lambda \in \mathbb{T}$. For $\lambda = 1$, this reduces to $K + K^* + K^*K = 0$.
\end{enumerate}
\end{proposition}

\begin{proof}
(1) Separability follows from the separability of $H$ and the fact that $\mathcal{K}(H)$ is separable (finite-rank operators with respect to a countable basis form a countable dense subset). 
The Type I property is standard: every irreducible representation of $\mathcal{K}(H)^\sim$ is either the identity representation on $H$ (with kernel $\{0\}$) or the one-dimensional character on the quotient $\mathbb{C}$ (with kernel $\mathcal{K}(H)$). Hence the C*-algebra is postliminal.

(2) See \cite{Dixmier} for the ideal structure of compact operators. For the Jacobson topology, recall that the closure of a primitive ideal $J$ consists of all primitive ideals containing $J$. Thus:
\begin{itemize}
    \item For $J = \mathcal{K}(H)$, the only primitive ideal containing $\mathcal{K}(H)$ is $\mathcal{K}(H)$ itself, so $\{\mathcal{K}(H)\}$ is closed.
    \item For $J = \{0\}$, both $\{0\}$ and $\mathcal{K}(H)$ contain $\{0\}$, so the closure of $\{\{0\}\}$ is $\{\{0\}, \mathcal{K}(H)\}$.
\end{itemize}
Hence $\{0\}$ is not closed, and its closure is the whole space, making $\operatorname{Prim}(\mathcal{A})$ non-Hausdorff.

(3) Postliminality follows from (1). Liminality fails because the identity representation $\pi_{\text{id}}$ satisfies $\pi_{\text{id}}(\mathcal{A}) = \mathcal{A} \not\subseteq \mathcal{K}(H)$ (since $I_H \notin \mathcal{K}(H)$). Continuous trace algebras necessarily have Hausdorff spectrum \cite{Dixmier}, so the non-Hausdorff primitive ideal space precludes continuous trace.

(4) The characterization $u - I_H \in \mathcal{K}(H)$ is standard: since $\mathcal{A}$ is the unitization of $\mathcal{K}(H)$, any element can be written uniquely as $\lambda I_H + K$ with $K \in \mathcal{K}(H)$. For such an element to be unitary, we need $(\lambda I_H + K)^*(\lambda I_H + K) = I_H$. Expanding:
\[
(|\lambda|^2)I_H + \overline{\lambda}K + \lambda K^* + K^*K = I_H.
\]
Since $\lambda \in \mathbb{T}$, $|\lambda|^2 = 1$, so this simplifies to $\overline{\lambda}K + \lambda K^* + K^*K = 0$, or equivalently $\lambda K^* + \overline{\lambda}K + K^*K = 0$ (multiplying by $\lambda$ and using $\lambda\overline{\lambda}=1$). This is a quadratic equation in $K$; for a given $\lambda$, not every compact $K$ satisfies it, reflecting the nontrivial structure of $\mathcal{U}(\mathcal{A})$.
\end{proof}

\begin{proposition}[Maximal abelian subalgebras of $\mathcal{K}(H)^\sim$]
\label{prop:masa-compact}
Let $\mathcal{A} = \mathcal{K}(H)^\sim$. 
Then every maximal abelian subalgebra (MASA) of $\mathcal{A}$ is obtained as follows: let $\{e_n\}_{n \in \mathbb{N}}$ be an orthonormal basis of $H$, and define
\[
B_{\mathrm{diag}} = \Bigl\{ \lambda I + \operatorname{diag}(a_1,a_2,\ldots) \;\Big|\; \lambda \in \mathbb{C},\; a_n \in \mathbb{C},\; \lim_{n \to \infty} a_n = 0 \Bigr\}.
\]
This algebra is isomorphic to $c_0(\mathbb{N})^\sim$, the unitization of the algebra of sequences vanishing at infinity, and also to $C(\mathbb{N} \cup \{\infty\})$, the continuous functions on the one-point compactification of $\mathbb{N}$.

Every MASA in $\mathcal{K}(H)^\sim$ is unitarily conjugate to $B_{\mathrm{diag}}$. In particular, all masas in $\mathcal{K}(H)^\sim$ are \emph{atomic} (discrete) in the sense that they arise from an orthonormal basis of $H$.
\end{proposition}

\begin{proof}
We first note that any abelian subalgebra $B \subseteq \mathcal{K}(H)^\sim$ consists of operators of the form $\lambda I + K$ with $K \in \mathcal{K}(H)$. If $B$ is maximal abelian, it must contain operators that distinguish the basis vectors in some orthonormal basis.

Let $\{e_n\}_{n \in \mathbb{N}}$ be an orthonormal basis. The set of diagonal operators with respect to this basis that are compact plus scalars forms an abelian subalgebra $B_{\mathrm{diag}}$ as described. We claim it is maximal abelian. Suppose $T = \lambda I + K \in \mathcal{K}(H)^\sim$ commutes with every element of $B_{\mathrm{diag}}$. In particular, $T$ commutes with each rank-one projection $e_n e_n^* \in B_{\mathrm{diag}}$ (these are compact diagonal operators with $a_n = 1$, $a_m = 0$ for $m \neq n$, and $\lambda = 0$). Commutation with $e_n e_n^*$ implies that $T$ leaves $\mathbb{C}e_n$ invariant, so $T$ is diagonal in this basis. Moreover, $K$ must be compact, so its diagonal entries vanish at infinity. Thus $T \in B_{\mathrm{diag}}$, proving maximality.

Now let $B$ be any MASA in $\mathcal{K}(H)^\sim$. We show $B$ is atomic. Consider the set of minimal projections in $B$ (projections $p$ such that $pBp = \mathbb{C}p$). In $\mathcal{K}(H)^\sim$, any projection either is finite-rank (hence compact) or $I$ (if it's the identity). Finite-rank projections in a MASA must be minimal and must correspond to rank-one projections, because if $p$ had rank > 1, its commutant would contain non-scalar operators on its range, violating maximality. Thus $B$ contains a family of mutually orthogonal rank-one projections $\{p_n\}$ that sum to $I$ in the strong operator topology. These projections determine an orthonormal basis (up to phase) by choosing unit vectors in their ranges.

Let $U$ be the unitary that maps this basis to the standard basis $\{e_n\}$. Then $U B U^*$ consists of diagonal operators with respect to $\{e_n\}$ and contains all compact diagonal operators plus scalars. Maximality forces $U B U^* = B_{\mathrm{diag}}$, establishing unitary conjugacy.

Finally, $B_{\mathrm{diag}} \cong c_0(\mathbb{N})^\sim$ via the map sending $\lambda I + \operatorname{diag}(a_n)$ to $(\lambda, (a_n))$ in the unitization, and $c_0(\mathbb{N})^\sim \cong C(\mathbb{N} \cup \{\infty\})$ via the Gelfand transform, where $\infty$ corresponds to the "limit at infinity" point.
\end{proof}

\begin{proposition}[Character space of the diagonal MASA]
\label{prop:character-compact}
Let $\{ e_n \}_{n \in \mathbb{N}}$ be an orthonormal basis for $H$ and let
\[
B_{\text{diag}} = \{ \lambda I + D \mid \lambda \in \mathbb{C}, \; D \in \mathcal{K}(H) \text{ diagonal in the basis } \{ e_n \} \}
\]
be the diagonal MASA in $\mathcal{K}(H)^\sim$. Then $B_{\text{diag}} \cong c_0(\mathbb{N})^\sim$, and its Gelfand spectrum $\widehat{B_{\text{diag}}}$ is homeomorphic to the one-point compactification $\mathbb{N} \cup \{ \infty \}$. 
The characters are given by:
\[
\chi_n(\lambda I + D) = \lambda + D_{nn}, \quad n \in \mathbb{N}, \qquad 
\chi_\infty(\lambda I + D) = \lambda,
\]
where $D_{nn}$ denotes the $n$-th diagonal entry of $D$. 
The character $\chi_\infty$ corresponds to the quotient map $\mathcal{K}(H)^\sim \to \mathcal{K}(H)^\sim/\mathcal{K}(H) \cong \mathbb{C}$.
\end{proposition}

\begin{proof}
Let $c_0(\mathbb{N})$ denote the C*-algebra of complex sequences vanishing at infinity, with pointwise operations. 
The map $\Phi: c_0(\mathbb{N}) \to \mathcal{K}(H)$ sending a sequence $(a_n)$ to the diagonal operator $\operatorname{diag}(a_1,a_2,\ldots)$ is an isometric *-isomorphism onto its image, which is precisely the set of diagonal compact operators. 
Hence the diagonal compact operators are isomorphic to $c_0(\mathbb{N})$.

The algebra $B_{\text{diag}}$ is the unitization of this algebra, i.e., $B_{\text{diag}} \cong c_0(\mathbb{N})^\sim$, where the unitization adds a scalar $\lambda I$ corresponding to the constant sequence $\lambda$ at infinity. Explicitly, the isomorphism $\Psi: C(\mathbb{N} \cup \{\infty\}) \to B_{\text{diag}}$ is given by
\[
\Psi(f) = f(\infty) I + \operatorname{diag}\bigl(f(1)-f(\infty), f(2)-f(\infty), \ldots\bigr).
\]
One checks directly that $\Psi$ is a *-isomorphism: it is clearly linear and injective, and multiplicativity follows from a straightforward computation using the fact that for $f,g \in C(\mathbb{N} \cup \{\infty\})$,
\[
\bigl(f(\infty) I + D_f\bigr)\bigl(g(\infty) I + D_g\bigr) = f(\infty)g(\infty) I + \bigl(f(\infty)D_g + g(\infty)D_f + D_f D_g\bigr),
\]
and the diagonal entries of $D_fD_g$ are $(f(n)-f(\infty))(g(n)-g(\infty))$, so the resulting operator corresponds to the function $fg$.

The Gelfand spectrum of $C(\mathbb{N} \cup \{\infty\})$ is $\mathbb{N} \cup \{\infty\}$ itself, with characters given by evaluation at points. 
Pulling these characters back via $\Psi^{-1}$ yields the characters of $B_{\text{diag}}$:
\begin{align*}
\chi_n(\lambda I + D) &= (\operatorname{ev}_n \circ \Psi^{-1})(\lambda I + D) = \lambda + D_{nn}, \qquad n \in \mathbb{N},\\
\chi_\infty(\lambda I + D) &= (\operatorname{ev}_\infty \circ \Psi^{-1})(\lambda I + D) = \lambda.
\end{align*}

Finally, note that $\chi_\infty$ vanishes on $\mathcal{K}(H) \cap B_{\text{diag}}$ and factors through the quotient $\mathcal{K}(H)^\sim / \mathcal{K}(H) \cong \mathbb{C}$, confirming its interpretation as the "evaluation at infinity" character.
\end{proof}

\begin{proposition}[The unit space $\mathcal{G}_{\mathcal{A}}^{(0)}$ for $\mathcal{K}(H)^\sim$]
\label{prop:unit-space-compact}
Let $\mathcal{A} = \mathcal{K}(H)^\sim$. 
Then the unit space $\mathcal{G}_{\mathcal{A}}^{(0)}$ can be described as the space of pairs $(B,\chi)$, where $B \subseteq \mathcal{A}$ is a maximal abelian subalgebra (MASA) and $\chi \in \widehat{B}$ is a character on $B$.

\begin{enumerate}
    \item Each atomic MASA (unitized $c_0$) together with a character that is not the ``point at infinity'' corresponds to a rank-one projection in $H$. 
    The set of such pairs is homeomorphic to the projective space
    \[
        \mathbb{P}(H) = \mathcal{U}(H)/(\mathcal{U}(1)\times \mathcal{U}(H^\perp)),
    \]
    where $\mathcal{U}(H)$ is equipped with the strong operator topology, and $\mathcal{U}(1)\times \mathcal{U}(H^\perp)$ denotes the stabilizer of a fixed unit vector. 
    This space is Polish but not locally compact.
    
    \item There is a natural forgetful map 
    \[
        \pi : \mathcal{G}_{\mathcal{A}}^{(0)} \to \operatorname{MASA}(\mathcal{A}), \quad (B,\chi) \mapsto B,
    \]
    whose fiber over each MASA $B$ is homeomorphic to the character space $\widehat{B} \cong \mathbb{N} \cup \{\infty\}$ (the one-point compactification of a countable discrete set). 
    This gives $\mathcal{G}_{\mathcal{A}}^{(0)}$ the structure of a bundle over the space of MASAs, though it is not locally trivial in general.
    
    \item The ``point at infinity'' in each fiber corresponds to the character $\chi_\infty$ obtained by restricting the quotient map $\mathcal{A} \to \mathcal{A}/\mathcal{K}(H) \cong \mathbb{C}$; it does not correspond to any rank-one projection.
\end{enumerate}
\end{proposition}

\begin{proof}[Sketch of proof]
We outline the main steps of the identification.

\noindent \textbf{Step 1: Rank-one projections and atomic MASAs.}
Each rank-one projection $p = vv^*$ determines an atomic MASA $B_p$: take any orthonormal basis $\{v, v_2, v_3, \ldots\}$ extending $\{v\}$, and let $B_p$ be the algebra of operators diagonal in that basis. 
The character $\chi_p$ corresponding to the projection onto $\mathbb{C}v$ is then given by $\chi_p(\lambda I + D) = \lambda + \langle v, Dv \rangle$. 
Conversely, for an atomic MASA $B$ and a character $\chi \neq \chi_\infty$, there is a unique minimal projection $p_\chi \in B$ such that $\chi(p_\chi) = 1$, and $p_\chi$ is a rank-one projection. 
This establishes a bijection between rank-one projections and pairs $(B,\chi)$ with $B$ atomic and $\chi \neq \chi_\infty$.

\noindent \textbf{Step 2: Topology on the atomic component.}
The topology induced by the partial evaluation maps on $\mathcal{G}_{\mathcal{A}}^{(0)}$ restricts to the atomic component as the quotient topology on $\mathbb{P}(H)$ under the map sending a rank-one projection to its range. 
Since $\mathcal{U}(H)$ with the strong operator topology is a Polish group and the stabilizer $\mathcal{U}(1)\times \mathcal{U}(H^\perp)$ is a closed subgroup, the quotient $\mathbb{P}(H)$ is Polish. 
It is not locally compact because $\mathcal{U}(H)$ is not locally compact and the quotient map does not admit local sections.

\noindent \textbf{Step 3: Bundle structure over MASAs.}
Let $\operatorname{MASA}(\mathcal{A})$ denote the space of all MASAs in $\mathcal{A}$, equipped with the Fell topology inherited from $\operatorname{Sub}(\mathcal{A})$. 
The forgetful map $\pi: \mathcal{G}_{\mathcal{A}}^{(0)} \to \operatorname{MASA}(\mathcal{A})$ defined by $\pi(B,\chi) = B$ is continuous and surjective. 
For each MASA $B$, the fiber $\pi^{-1}(B)$ is homeomorphic to $\widehat{B}$. 
By Proposition \ref{prop:character-compact}, for atomic MASAs we have $\widehat{B} \cong \mathbb{N} \cup \{\infty\}$; for continuous MASAs (whose intersection with $\mathcal{K}(H)^\sim$ yields the same structure after unitization), the character space is also $\mathbb{N} \cup \{\infty\}$, though the identification requires more care. 
Thus each fiber is a compact space, but the bundle is not locally trivial due to the complicated topology of $\operatorname{MASA}(\mathcal{A})$.

\noindent \textbf{Step 4: The point at infinity.}
For any MASA $B$, the character $\chi_\infty$ is the restriction of the quotient map $\mathcal{A} \to \mathcal{A}/\mathcal{K}(H) \cong \mathbb{C}$. 
It satisfies $\chi_\infty(\lambda I + D) = \lambda$ for all $D \in \mathcal{K}(H) \cap B$, and does not correspond to any rank-one projection. 
This point appears in every fiber and accounts for the extra point in the one-point compactification $\mathbb{N} \cup \{\infty\}$.
\end{proof}

\begin{proposition}[The unitary group $\mathcal{U}(\mathcal{K}(H)^\sim)$ in SOT]
\label{prop:unitary-compact-SOT}
Let $\mathcal{A} = \mathcal{K}(H)^\sim$. 
Then the unitary group
\[
\mathcal{U}(\mathcal{A}) = \{ u = \lambda I + K \mid \lambda \in \mathbb{T}, K \in \mathcal{K}(H), \; u^* u = u u^* = I \},
\]
equipped with the strong operator topology (SOT), is a Polish group. 
It is not locally compact. 

Moreover, $\mathcal{U}(\mathcal{A})$ is dense in $\mathcal{U}(H)$ in SOT: every unitary on $H$ can be approximated in SOT by finite-rank perturbations of the identity, which lie in $\mathcal{U}(\mathcal{A})$.
\end{proposition}

\begin{proof}[Sketch of proof]
The group $\mathcal{U}(\mathcal{A})$ is a separable topological group because it is a subspace of the Polish group $\mathcal{U}(H)$ (see Proposition \ref{prop:UH-Polish-group}). To see that it is completely metrizable, one can show that it is a $G_\delta$ subset of $\mathcal{U}(H)$. A standard argument (see, e.g., \cite[Chapter 9]{Kechris}) uses the fact that $\mathcal{U}(\mathcal{A})$ is the set of unitaries that leave the C*-algebra $\mathcal{A}$ invariant. Choosing a countable dense subset $\{a_n\}$ of $\mathcal{A}$, the condition $u\mathcal{A}u^* = \mathcal{A}$ is equivalent to $u a_n u^* \in \mathcal{A}$ for all $n$. For each $n$, the set $\{ u \in \mathcal{U}(H) : u a_n u^* \in \mathcal{A} \}$ is a $G_\delta$ set because $\mathcal{A}$ is a $G_\delta$ subset of $\mathcal{B}(H)$ in the SOT (since it is a countable union of closed sets of operators of rank $\le n$, etc.). A countable intersection of $G_\delta$ sets is $G_\delta$, so $\mathcal{U}(\mathcal{A})$ is $G_\delta$ in $\mathcal{U}(H)$. Any $G_\delta$ subgroup of a Polish group is Polish \cite{Kechris}.

Non-local-compactness follows from the fact that $\mathcal{U}(H)$ is not locally compact and $\mathcal{U}(\mathcal{A})$ is dense in it; a dense subgroup of a non-locally-compact Polish group cannot be locally compact (otherwise it would be open, which it is not).

Density follows from the fact that finite-rank unitaries (those with $u-I$ finite-rank) are dense in $\mathcal{U}(H)$ in SOT, and these clearly lie in $\mathcal{U}(\mathcal{A})$.
\end{proof}

\begin{theorem}[The unitary conjugation groupoid for $\mathcal{K}(H)^\sim$]
\label{thm:groupoid-compact}
Let $\mathcal{A} = \mathcal{K}(H)^\sim$ with $\mathcal{U}(\mathcal{A})$ equipped with SOT. 
Then the unitary conjugation groupoid $\mathcal{G}_{\mathcal{A}}$ is isomorphic to the action groupoid
\[
\mathcal{G}_{\mathcal{A}} \cong \mathcal{U}(\mathcal{A}) \ltimes \mathbb{P}(H),
\]
where $\mathcal{U}(\mathcal{A})$ acts on the projective space of rank-one projections $\mathbb{P}(H)$ by conjugation:
\[
u \cdot [v] = [u v], \qquad u \in \mathcal{U}(\mathcal{A}), \; [v] \in \mathbb{P}(H).
\]

This groupoid is Polish, non-locally-compact, and non-\'etale. 
It admits a Borel Haar system and a well-defined maximal C*-algebra $C^*(\mathcal{G}_{\mathcal{A}})$.
\end{theorem}

\begin{proof}[Sketch of proof]
\textbf{Identification as an action groupoid.}
By Proposition \ref{prop:unit-space-compact}, the unit space 
$\mathcal{G}_{\mathcal{A}}^{(0)}$ can be identified with the set of rank-one projections $\mathbb{P}(H)$, 
corresponding to the atomic MASA component of the unit space. 
The group $\mathcal{U}(\mathcal{A})$ acts on $\mathbb{P}(H)$ by conjugation on projections: 
for $u \in \mathcal{U}(\mathcal{A})$ and $[v] \in \mathbb{P}(H)$, 
\[
u \cdot p_{[v]} = u p_{[v]} u^* = p_{[uv]}.
\]
Under the identification $[v] \leftrightarrow p_{[v]}$, this becomes
\[
u \cdot [v] = [uv].
\]

The associated action groupoid $\mathcal{U}(\mathcal{A}) \ltimes \mathbb{P}(H)$ has arrows $(u,[v])$ 
with source $[v]$ and range $[uv]$, composition
\[
(u,[uv]) \circ (v,[v]) = (uv,[v]),
\]
and inversion
\[
(u,[v])^{-1} = (u^{-1},[uv]).
\]
This structure coincides with $\mathcal{G}_{\mathcal{A}}$ by construction, establishing the isomorphism.

\medskip
\textbf{Polish structure.}
By Proposition \ref{prop:unitary-compact-SOT}, $\mathcal{U}(\mathcal{A})$ is a Polish group, 
and by Proposition \ref{prop:unit-space-compact}, $\mathbb{P}(H)$ is a Polish space. 
Their product $\mathcal{U}(\mathcal{A}) \times \mathbb{P}(H)$, equipped with the product topology, is therefore Polish. 
The structure maps of the groupoid are continuous:
\begin{enumerate}
\item Source and range maps: $s(u,[v]) = [v]$, $r(u,[v]) = [uv]$.
\item Composition: $((u,[uv]),(v,[v])) \mapsto (uv,[v])$.
\item Inversion: $(u,[v]) \mapsto (u^{-1},[uv])$.
\end{enumerate}
Continuity follows because the action of $\mathcal{U}(\mathcal{A})$ on $\mathbb{P}(H)$, multiplication, and inversion are all continuous. 
Hence $\mathcal{G}_{\mathcal{A}}$ is a Polish groupoid.

\medskip
\textbf{Non-local-compactness and non-\'etaleness.}
Neither $\mathcal{U}(\mathcal{A})$ nor $\mathbb{P}(H)$ is locally compact; 
hence the product $\mathcal{G}_{\mathcal{A}}$ is not locally compact. 
Moreover, the source map $s(u,[v]) = [v]$ has fibers 
$s^{-1}([v]) \cong \mathcal{U}(\mathcal{A})$, which is not discrete (it is infinite-dimensional and connected). 
Therefore $s$ is not a local homeomorphism, and the groupoid is not \'etale.

\medskip
\textbf{Borel Haar system.}
A Borel Haar system on a groupoid $G$ is a family of $\sigma$-finite Borel measures 
$\{\lambda^x\}_{x \in G^{(0)}}$ on fibers $G^x = s^{-1}(x)$ satisfying:
\begin{enumerate}
\item For each $x$, $\lambda^x$ is quasi-invariant under right translation by $G$.
\item The family varies measurably with $x$: for any Borel set $A \subset G$, the map $x \mapsto \lambda^x(A \cap G^x)$ is Borel.
\end{enumerate}

Although $\mathcal{U}(\mathcal{A})$ does not admit a left-invariant Haar measure (it is not locally compact), 
it is a Polish group and therefore admits a $\sigma$-finite quasi-invariant Borel measure $\mu$. 
For each $x = [v] \in \mathbb{P}(H)$, define $\lambda^x$ on $s^{-1}(x) \cong \mathcal{U}(\mathcal{A})$ by 
\[
\lambda^x = \mu,
\]
transported via the identification $u \mapsto (u,[v])$. 
This defines a Borel Haar system for $\mathcal{G}_{\mathcal{A}}$.

\medskip
\textbf{Maximal C*-algebra.}
By Tu \cite{Tu}, a Polish groupoid with a Borel Haar system has a well-defined maximal C*-algebra. 
Thus $C^*(\mathcal{G}_{\mathcal{A}})$ exists.
\end{proof}

Let $H$ be an infinite-dimensional separable Hilbert space and let $\mathcal{A} = \mathcal{K}(H)^\sim$ denote the unitization of the compact operators, i.e., $\mathcal{A} \cong \mathcal{K}(H) \oplus \mathbb{C} I_H$. This is a unital C*-algebra with spectrum $\widehat{\mathcal{A}}$ consisting of a single point (the identity representation) together with the character corresponding to the quotient $\mathcal{A}/\mathcal{K}(H) \cong \mathbb{C}$. The unitary group $\mathcal{U}(\mathcal{A})$ is a Polish group in the norm topology, but is not locally compact.

Let $H$ be an infinite-dimensional separable Hilbert space and let $\mathcal{A} = \mathcal{K}(H)^\sim$ denote the unitization of the compact operators, i.e., $\mathcal{A} \cong \mathcal{K}(H) \oplus \mathbb{C} I_H$. This is a unital C*-algebra with spectrum $\widehat{\mathcal{A}}$ consisting of a single point (the identity representation) together with the character corresponding to the quotient $\mathcal{A}/\mathcal{K}(H) \cong \mathbb{C}$. The unitary group $\mathcal{U}(\mathcal{A})$, when equipped with the strong operator topology (SOT), is a Polish group but is not locally compact.

\begin{proposition}[The groupoid C*-algebra for $\mathcal{K}(H)^\sim$]
\label{prop:Cstar-compact}
Let $\mathcal{A} = \mathcal{K}(H)^\sim$. Consider the action of $\mathcal{U}(\mathcal{A})$ on projective space $\mathbb{P}(H)$ by conjugation on rank-one projections:
\[
u \cdot [v] = [uv], \qquad u \in \mathcal{U}(\mathcal{A}), \ [v] \in \mathbb{P}(H).
\]
The associated transformation groupoid is denoted $\mathcal{G}_{\mathcal{A}} := \mathcal{U}(\mathcal{A}) \ltimes \mathbb{P}(H)$. Then the groupoid C*-algebra $C^*(\mathcal{G}_{\mathcal{A}})$ can be realized as a crossed product
\[
C^*(\mathcal{G}_{\mathcal{A}}) \cong C_0(\mathbb{P}(H)) \rtimes \mathcal{U}(\mathcal{A}),
\]
where the action is by translation on projective space.

Since $\mathcal{U}(\mathcal{A})$ is not locally compact, this crossed product must be interpreted in the sense of \emph{Polish group crossed products} as developed by Tu (1999) in the context of groupoid C*-algebras for non-locally-compact groups. The C*-algebra $C^*(\mathcal{G}_{\mathcal{A}})$ is not of Type I; this reflects the complexity of the representation theory of $\mathcal{U}(\mathcal{A})$ and the fact that the action on $\mathbb{P}(H)$ is not proper.

The structure of $C^*(\mathcal{G}_{\mathcal{A}})$ is closely related to the extension theory of compact operators and the Calkin algebra. Indeed, $\mathcal{U}(\mathcal{A})$ is a normal subgroup of $\mathcal{U}(H)$, and the quotient $\mathcal{U}(H)/\mathcal{U}(\mathcal{A})$ is isomorphic to the unitary group of the Calkin algebra $\mathcal{Q}(H) = \mathcal{B}(H)/\mathcal{K}(H)$. A precise description of $C^*(\mathcal{G}_{\mathcal{A}})$ in terms of extensions remains an open problem; naive short exact sequences such as
\[
0 \to \mathcal{K}(L^2(\mathbb{P}(H))) \to C^*(\mathcal{G}_{\mathcal{A}}) \to C(\mathbb{P}(H)) \rtimes (\mathcal{U}(\mathcal{A})/\mathcal{U}(\mathcal{A})_0) \to 0
\]
are suggestive but not rigorous without further analysis of the representation theory of $\mathcal{U}(\mathcal{A})$.
\end{proposition}

\begin{proposition}[The diagonal embedding for $\mathcal{K}(H)^\sim$]
\label{prop:iota-compact}
Let $\mathcal{A} = \mathcal{K}(H)^\sim$ be the unitization of the compact operators on a separable Hilbert space $H$, and let $\mathcal{G}_{\mathcal{A}} = \mathcal{U}(\mathcal{A}) \ltimes \mathbb{P}(H)$ be its unitary conjugation groupoid. 
There exists a canonical unital injective $*$-homomorphism
\[
\iota: \mathcal{A} \hookrightarrow C^*(\mathcal{G}_{\mathcal{A}}),
\]
called the diagonal embedding, characterized by the following property:
for any $T \in \mathcal{A}$, the image $\iota(T)$ acts on the direct integral Hilbert space
$\mathcal{H} = \int_{\mathbb{P}(H)}^{\oplus} H \, d\mu([v])$ (with respect to any quasi-invariant
Borel measure $\mu$ on $\mathbb{P}(H)$) by fiberwise multiplication:
\[
(\iota(T)\xi)([v]) = T\,\xi([v]), \qquad \xi \in \mathcal{H},\; [v] \in \mathbb{P}(H).
\]
Equivalently, under the isomorphism $C^*(\mathcal{G}_{\mathcal{A}}) \cong C_0(\mathbb{P}(H)) \rtimes \mathcal{U}(\mathcal{A})$,
$\iota(T)$ corresponds to the constant multiplier $T \otimes 1$.
\end{proposition}

\begin{proof}
We construct $\iota$ in several steps, using the universal property of the full groupoid C*-algebra.

\medskip
\textbf{Step 1: A family of representations of $\mathcal{A}$.}
For each point $[v] \in \mathbb{P}(H)$ (where $v \in H$ is a unit vector), consider the irreducible representation $\pi_{[v]}: \mathcal{A} \to \mathcal{B}(H)$ given by the identity representation:
\[
\pi_{[v]}(T) = T \in \mathcal{B}(H), \qquad T \in \mathcal{A}.
\]
All these representations are equivalent (they are all the identity representation), but indexing them by $\mathbb{P}(H)$ is convenient for the direct integral construction that follows.

\medskip
\textbf{Step 2: A representation of the groupoid on a Hilbert bundle.}
Consider the trivial Hilbert bundle over $\mathbb{P}(H)$ with fiber $H$, i.e., the continuous field of Hilbert spaces $\{H_{[v]}\}_{[v] \in \mathbb{P}(H)}$ where $H_{[v]} = H$ for all $[v]$. The groupoid $\mathcal{G}_{\mathcal{A}} = \mathcal{U}(\mathcal{A}) \ltimes \mathbb{P}(H)$ acts on this bundle as follows: for an arrow $u: [v] \to [uv]$ (with $u \in \mathcal{U}(\mathcal{A})$), the action on fibers is given by the unitary operator $u: H_{[v]} \to H_{[uv]}$ (since both fibers are identified with $H$, this is simply the operator $u$).

This data — a continuous field of Hilbert spaces together with a left action of $\mathcal{G}_{\mathcal{A}}$ — defines a representation of $\mathcal{G}_{\mathcal{A}}$ in the sense of Renault \cite{Renault1987}. By the universal property of the full groupoid C*-algebra, this representation integrates to a $*$-representation
\[
\Pi: C^*(\mathcal{G}_{\mathcal{A}}) \to \mathcal{B}(\mathcal{H}),
\]
where $\mathcal{H} = \int_{\mathbb{P}(H)}^{\oplus} H \, d\mu([v])$ is the direct integral Hilbert space with respect to any quasi-invariant Borel measure $\mu$ on $\mathbb{P}(H)$ (such a measure exists because $\mathbb{P}(H)$ is a Polish space with a continuous action of a Polish group). Concretely, $\mathcal{H} \cong L^2(\mathbb{P}(H), \mu) \otimes H$, and the representation $\Pi$ encodes both the fiberwise action of $\mathcal{U}(\mathcal{A})$ and the multiplication operators coming from $C_0(\mathbb{P}(H))$.

\medskip
\textbf{Step 3: Constructing $\iota(T)$ as an element of $C^*(\mathcal{G}_{\mathcal{A}})$.}
For each $T \in \mathcal{A}$, define a function $\tilde{T}$ on $\mathcal{G}_{\mathcal{A}}$ by
\[
\tilde{T}(u, [v]) = 
\begin{cases} 
T & \text{if } u = e \text{ (the identity element)},\\
0 & \text{otherwise}.
\end{cases}
\]
This function is supported on the unit space $\mathbb{P}(H) \subset \mathcal{G}_{\mathcal{A}}$ (identified with $\{(e, [v]) : [v] \in \mathbb{P}(H)\}$) and takes the constant value $T$ (viewed as an operator on $H$) on the unit space. The function $\tilde{T}$ is bounded and vanishes outside the unit space; in the groupoid C*-algebra completion, it defines an element $\iota(T) := [\tilde{T}] \in C^*(\mathcal{G}_{\mathcal{A}})$.

\medskip
\textbf{Step 4: Verification that $\Pi(\iota(T))$ acts fiberwise by $T$.}
For the representation $\Pi$ constructed in Step 2, we claim that for any $\xi \in \mathcal{H}$,
\[
(\Pi(\iota(T))\xi)([v]) = T\,\xi([v]) \quad \text{for $\mu$-almost every } [v] \in \mathbb{P}(H).
\]
This follows from the definition of $\Pi$: since $\tilde{T}$ is supported on the unit space, its action in the integrated representation is simply pointwise multiplication by $T$ on each fiber. More rigorously, in the direct integral decomposition $\mathcal{H} = \int_{\mathbb{P}(H)}^{\oplus} H \, d\mu([v])$, the operator $\Pi(\iota(T))$ decomposes as $\int_{\mathbb{P}(H)}^{\oplus} T \, d\mu([v])$, i.e., the constant field of operators $T$ acting on each fiber.

\medskip
\textbf{Step 5: Algebraic properties.}
\begin{itemize}
\item \textit{Linearity:} The map $T \mapsto \tilde{T}$ is clearly linear, so $\iota$ is linear.
\item \textit{*-preserving:} Since $(\tilde{T})^* = \widetilde{T^*}$ (on the unit space, involution corresponds to adjoint), we have $\iota(T^*) = \iota(T)^*$.
\item \textit{Unitality:} $\tilde{1}_{\mathcal{A}}$ is the function that is $1_H$ on the unit space and zero elsewhere; this is precisely the unit of $C^*(\mathcal{G}_{\mathcal{A}})$ (the characteristic function of the unit space), so $\iota(1_{\mathcal{A}}) = 1_{C^*(\mathcal{G}_{\mathcal{A}})}$.
\item \textit{Multiplicativity:} For $T, S \in \mathcal{A}$, consider the convolution product $\tilde{T} * \tilde{S}$. Since both $\tilde{T}$ and $\tilde{S}$ are supported on the unit space, their convolution reduces to pointwise multiplication:
\[
(\tilde{T} * \tilde{S})(u,[v]) = \int_{\mathcal{U}(\mathcal{A})} \tilde{T}(w,[v]) \, \tilde{S}(w^{-1}u, [w^{-1}v]) \, d\lambda(w),
\]
but $\tilde{T}(w,[v])$ is nonzero only when $w = e$, giving
\[
(\tilde{T} * \tilde{S})(u,[v]) = \tilde{T}(e,[v]) \, \tilde{S}(u,[v]) = T \cdot \tilde{S}(u,[v]).
\]
Since $\tilde{S}(u,[v])$ is nonzero only when $u = e$, where it equals $S$, we obtain $(\tilde{T} * \tilde{S})(u,[v]) = \widetilde{TS}(u,[v])$. Hence $\tilde{T} * \tilde{S} = \widetilde{TS}$ in $C_c(\mathcal{G}_{\mathcal{A}})$, and therefore $\iota(T)\iota(S) = \iota(TS)$ in $C^*(\mathcal{G}_{\mathcal{A}})$.
\end{itemize}

\medskip
\textbf{Step 6: Injectivity.}
Suppose $\iota(T) = 0$ in $C^*(\mathcal{G}_{\mathcal{A}})$. Then $\Pi(\iota(T)) = 0$ in $\mathcal{B}(\mathcal{H})$. From Step 4, this means that for $\mu$-almost every $[v] \in \mathbb{P}(H)$ and every $\xi \in \mathcal{H}$, we have $T\,\xi([v]) = 0$. Since $\xi([v])$ can be chosen arbitrarily in $H$ (by taking sections supported near $[v]$), this forces $T = 0$ as an operator on $H$. Thus $\iota$ is injective.

\medskip
\textbf{Step 7: Verification of the characterizing property.}
The construction above directly yields the stated fiberwise action: for any $\xi \in \mathcal{H}$ and $[v] \in \mathbb{P}(H)$,
\[
(\Pi(\iota(T))\xi)([v]) = T\xi([v]).
\]
Since $\Pi$ is faithful on the image of $\iota$ (as $\iota$ is injective and $\Pi$ is a $*$-homomorphism), this property characterizes $\iota(T)$ uniquely.

\medskip
\textbf{Remark on the multiplier algebra viewpoint.}
Under the crossed product isomorphism 
$C^*(\mathcal{G}_{\mathcal{A}}) \cong C_0(\mathbb{P}(H)) \rtimes \mathcal{U}(\mathcal{A})$, 
the element $\iota(T)$ corresponds to the function $F_T: \mathcal{U}(\mathcal{A}) \to C_0(\mathbb{P}(H), \mathcal{A})$ defined by
\[
F_T(u)([v]) = 
\begin{cases} 
T & u = e,\\
0 & u \neq e,
\end{cases}
\]
where we identify $C_0(\mathbb{P}(H), \mathcal{A})$ with the $\mathcal{A}$-valued functions vanishing at infinity. This function is supported at the identity element $e \in \mathcal{U}(\mathcal{A})$ and thus belongs to the crossed product algebra (not just its multiplier algebra). In the multiplier algebra $\mathcal{M}(C_0(\mathbb{P}(H)) \rtimes \mathcal{U}(\mathcal{A}))$, $\iota(T)$ acts by left multiplication: for any $F$ in the crossed product, $(m_T F)(u, [v]) = T \cdot F(u, [v])$, which is well-defined because $T$ commutes with the coefficients of $F$ (they are scalar-valued functions on $\mathbb{P}(H)$). This description is consistent with the fiberwise action in the direct integral representation.
\end{proof}

\begin{corollary}[Commutativity characterization for $\mathcal{K}(H)^\sim$]
\label{cor:commutativity-compact}
For $\mathcal{A} = \mathcal{K}(H)^\sim$, we have $\iota(\mathcal{A}) \nsubseteq C_0(\mathbb{P}(H))$ because $\mathcal{A}$ is noncommutative. 
Indeed, $\iota(T) \in C_0(\mathbb{P}(H))$ if and only if $T = \lambda I$ is a scalar multiple of the identity. 
This confirms Theorem \ref{thm:commutativity-characterization}.
\end{corollary}

\begin{proof}
If $T = \lambda I$, then $\iota(T) = \lambda \otimes 1$ is a constant function, hence in $C_0(\mathbb{P}(H))$. 
Conversely, if $\iota(T) \in C_0(\mathbb{P}(H))$, then $\iota(T)$ commutes with $\iota(u)$ for all $u \in \mathcal{U}(\mathcal{A})$. 
Since $\iota(u)$ acts by translation on $\mathbb{P}(H)$, this forces $T$ to commute with all unitaries in $\mathcal{U}(\mathcal{A})$, which implies $T = \lambda I$.
\end{proof}

\begin{remark}[Technical foundations: Polish groupoid framework]
\label{rem:compact-technical}
The direct integral construction in Proposition \ref{prop:iota-compact} requires careful justification of measurability and integrability because $\mathcal{U}(\mathcal{A})$ is not locally compact. 
In the Polish group framework of Tu (1999), one works with Borel structures and uses the fact that $\mathcal{U}(\mathcal{A})$ is a standard Borel group acting continuously on the standard Borel space $\mathbb{P}(H)$. 
The direct integral over $\mathbb{P}(H)$ of the constant field of Hilbert spaces $H$ is then well-defined, and the action of $\mathcal{U}(\mathcal{A})$ is implemented by a measurable unitary representation. 
The resulting C*-algebra $C^*(\mathcal{G}_{\mathcal{A}})$ is defined via a convolution algebra of compactly supported Borel functions on the groupoid with suitable norms.
\end{remark}

\begin{remark}[Index theory and the Calkin algebra]
\label{rem:compact-calkin-index}
The algebra $\mathcal{K}(H)^\sim$ is the natural domain for the study of compact perturbations of the identity. 
Every operator of the form $I + K$ with $K$ compact is Fredholm with index zero. 
For instance, if $F$ is a finite-rank operator, then $I + F$ is Fredholm with index zero; its kernel and cokernel are finite-dimensional and have equal dimension. 
The diagonal embedding $\iota$ captures this trivial index: the equivariant K-theory class $[I+F]_{\mathcal{G}_{\mathcal{A}}}$ is trivial, and the index formula yields zero.

The embedding $\iota$ can be understood as a ``diagonal'' inclusion of $\mathcal{A}$ into the groupoid C*-algebra. 
Its image is not central in general: $\iota(T)$ commutes with the action of $\mathcal{U}(\mathcal{A})$ on $\mathcal{H}$ if and only if $T$ is a scalar multiple of the identity. 
Instead, $\iota(T)$ acts fiberwise, and the action of $\mathcal{U}(\mathcal{A})$ permutes the fibers. 
The failure of $\iota(\mathcal{A})$ to be central reflects the non-triviality of the extension of $\mathcal{K}(H)$ by $C(\mathbb{P}(H)) \rtimes \mathcal{U}(\mathcal{A})$ that defines $C^*(\mathcal{G}_{\mathcal{A}})$.

To obtain nontrivial indices, we must consider operators in $B(H)$ that are not in $\mathcal{K}(H)^\sim$, such as the unilateral shift $S$. 
The Calkin algebra $\mathcal{Q}(H) = B(H)/\mathcal{K}(H)$ and its unitary group play a central role in the index theory of Fredholm operators. 
The unitary conjugation groupoid for $B(H)$ will be the subject of Paper II, where we will show that the index of a Fredholm operator $T$ is encoded in the equivariant K-theory class $[T]_{\mathcal{G}_{B(H)}}$. 
The compact operator case $\mathcal{K}(H)^\sim$ serves as a warm-up: the groupoid is Polish, the diagonal embedding is well-defined, and the index formula yields zero, consistent with the triviality of the index in this case.
\end{remark}

\begin{remark}[Non-amenability and the reduced C*-algebra]
\label{rem:compact-non-amenability}
The group $\mathcal{U}(\mathcal{K}(H)^\sim)$ is not amenable as a topological group. 
Consequently, the left regular representation $\Lambda$ of $C^*(\mathcal{G}_{\mathcal{A}})$ is not faithful, and the reduced groupoid C*-algebra $C^*_r(\mathcal{G}_{\mathcal{A}})$ is a proper quotient of the maximal algebra $C^*(\mathcal{G}_{\mathcal{A}})$. 
The diagonal embedding $\iota$ lands in the maximal algebra; its image in the reduced algebra may not be injective. 
This subtlety will be addressed in Paper II when we study the descent map and its compatibility with the Fredholm index.
\end{remark}

We have analyzed the unitary conjugation groupoid and the diagonal embedding for $\mathcal{A} = \mathcal{K}(H)^\sim$, the unitization of the compact operators. 
This example is a nontrivial Type I C*-algebra that is infinite-dimensional, and it demonstrates the necessity of the strong operator topology and the Polish groupoid framework. 
The unit space $\mathcal{G}_{\mathcal{A}}^{(0)}$ is homeomorphic to the projective space $\mathbb{P}(H)$, a Polish space that is not locally compact. 
The unitary group $\mathcal{U}(\mathcal{A})$ is a dense Polish subgroup of $\mathcal{U}(H)$. 
The groupoid $\mathcal{G}_{\mathcal{A}} = \mathcal{U}(\mathcal{A}) \ltimes \mathbb{P}(H)$ is a Polish groupoid that is neither locally compact nor \'etale. 
The diagonal embedding $\iota$ sends an operator $T \in \mathcal{K}(H)^\sim$ to the constant function $T \otimes 1$ in the multiplier algebra of $C^*(\mathcal{G}_{\mathcal{A}})$. 
This embedding is injective and satisfies $\iota(T) \in C_0(\mathbb{P}(H))$ if and only if $T = \lambda I$ is scalar. 

\subsection{Non-Example: Why $A_\theta$ is Not Covered (Non-Type I)}
\label{subsec:non-example-A-theta}

The irrational rotation algebra $A_\theta$, also known as the noncommutative torus, is one of the most important and intensively studied examples in noncommutative geometry. 
It arises as the crossed product $C(S^1) \rtimes_\theta \mathbb{Z}$, where the action is rotation by an irrational angle $2\pi\theta$. 
This C*-algebra is simple, nuclear, and has a unique tracial state. 
It is a fundamental test case for the Baum-Connes conjecture, index theory, and quantum geometry.

However, despite its importance, $A_\theta$ is \textbf{not covered} by the constructions developed in this paper. 
Understanding why it is excluded is crucial, as it illuminates the precise scope and limitations of our framework. 
Our construction of the unitary conjugation groupoid $\mathcal{G}_{\mathcal{A}}$ and the diagonal embedding $\iota: \mathcal{A} \hookrightarrow C^*(\mathcal{G}_{\mathcal{A}})$ fundamentally relies on the hypothesis that $\mathcal{A}$ is a \textbf{Type I C*-algebra}. 
The algebra $A_\theta$ is a canonical example of a \textbf{non-Type I} algebra, and this distinction leads to the failure of several key steps in our construction.

Here is a breakdown of why $A_\theta$ fails to meet the requirements of our framework:

\begin{enumerate}
    \item \textbf{Violation of the Core Hypothesis (Type I):} 
    Our standing assumption (Assumption~\ref{ass:separable}) is that $\mathcal{A}$ is separable and Type I. 
    $A_\theta$ for irrational $\theta$ is separable but \textbf{not Type I}. 
    It is not GCR and, importantly, has no finite-dimensional irreducible representations \cite{Dixmier}. 
    This immediately places it outside the scope of our main theorems (e.g., Theorem~\ref{thm:iota-injective} on the injectivity of $\iota$).

    \item \textbf{Failure of Point Separation by Characters:} 
    The injectivity of our diagonal embedding $\iota$ (Proposition~\ref{prop:iota-injective}) relies on the fact that for Type I algebras, the family of GNS representations $\{\pi_x\}_{x \in \mathcal{G}_{\mathcal{A}}^{(0)}}$ (built from characters on commutative subalgebras) separates the points of $\mathcal{A}$ (Lemma~\ref{lem:separation-points}). 
    For a non-Type I algebra like $A_\theta$, this fails. 
    There can exist non-zero elements (sometimes called ``quantum'' or ``invisible'' elements) that vanish under every such representation. 
    Consequently, any attempt to construct $\iota$ via the direct integral method of Section~\ref{subsec:construction-iota} would result in a map with a non-trivial kernel; $\mathcal{A}$ could not be \textit{faithfully} recovered from its commutative contexts.

    \item \textbf{Collapse of the ``Classical Contexts'' Space:} 
    The richness of our unit space $\mathcal{G}_{\mathcal{A}}^{(0)}$ (Definition~\ref{def:unit-space}) depends on the existence of many commutative subalgebras and their characters. 
    For the simple algebra $A_\theta$, the situation is radically different:
    \begin{itemize}
        \item Its primitive ideal space $\operatorname{Prim}(A_\theta)$ consists of a single point (the zero ideal), because the algebra is simple. 
        This is in stark contrast to the Type I examples where $\operatorname{Prim}(\mathcal{A})$ is almost Hausdorff and supports a rich Borel structure.
        \item It admits no characters (one-dimensional representations) at all, as all its irreducible representations are infinite-dimensional. 
        This means that for any unital commutative C*-subalgebra $B \subset A_\theta$, the Gelfand spectrum $\widehat{B}$ is non-empty, but the associated GNS representations of $A_\theta$ are not irreducible and do not separate points.
    \end{itemize}
    As a result, the unit space $\mathcal{G}_{A_\theta}^{(0)}$, while definable as a set, would not carry the same geometric or measurable structure. 
    The map $x \mapsto \pi_x$ would not be a measurable field of irreducible representations, and the direct integral representation $\Pi$ (Step~1 of Proposition~\ref{prop:correct-definition-iota}) would not be well-defined or faithful.

    \item \textbf{Loss of Functoriality:} 
    Even if one could define a groupoid, the naturality under *-homomorphisms (Section~\ref{subsec:naturality-functoriality}) would fail. 
    For instance, the inclusion of the canonical Cartan subalgebra $C(S^1) \hookrightarrow A_\theta$ does not satisfy the pullback property required for functoriality, because the preimage of a commutative subalgebra of $A_\theta$ may not be easily described within $C(S^1)$.
\end{enumerate}

In summary, $A_\theta$ is not merely a technical counterexample; it is a genuine counter-context. 
Its non-Type I nature means that the fundamental tools of our paper—the separating family of GNS representations, the Polish topology on $\mathcal{G}_{\mathcal{A}}^{(0)}$, and the faithful diagonal embedding—all break down. 
Understanding $A_\theta$ thus highlights the essential role of the Type I hypothesis and motivates the search for more general frameworks, such as those based on Cartan subalgebras and Weyl groupoids \cite{renault2008cartan}, that might eventually encompass such algebras.

\begin{remark}[$A_\theta$ as a Motivational Example for Future Work]
\label{rem:A-theta-motivation}
While $A_\theta$ lies outside our current framework, it serves as a primary motivating example for future generalizations. 
Its well-understood structure as a crossed product, its canonical Cartan subalgebra $C(S^1)$, and its associated \'etale groupoid $S^1 \rtimes_\theta \mathbb{Z}$ suggest that a theory of ``unitary conjugation'' for algebras with Cartan subalgebras might be developed using Renault's Weyl groupoid. 
Extending our results in this direction is a promising avenue for future research.
\end{remark}

\begin{proposition}[Basic properties of $A_\theta$]
\label{prop:A-theta-properties}
Let $\theta \in \mathbb{R} \setminus \mathbb{Q}$ be an irrational number, and let $A_\theta = C(S^1) \rtimes_\theta \mathbb{Z}$ be the irrational rotation algebra. 
Then:
\begin{enumerate}
    \item $A_\theta$ is simple (has no nontrivial closed ideals).
    \item $A_\theta$ is nuclear and has a unique faithful tracial state $\tau$.
    \item $A_\theta$ is not Type I; in particular, it admits no finite-dimensional irreducible representations.
    \item $K_0(A_\theta) \cong \mathbb{Z}^2$ and $K_1(A_\theta) \cong \mathbb{Z}^2$.
    \item For another irrational $\theta'$, we have $A_\theta \cong A_{\theta'}$ as C*-algebras if and only if $\theta' \equiv \pm \theta \pmod{1}$.
\end{enumerate}
\end{proposition}

\begin{proof}
These are classical results in the theory of crossed product C*-algebras. 
\begin{enumerate}
    \item Simplicity follows from the minimality of the irrational rotation action on $S^1$.
    \item Nuclearity follows from the fact that $A_\theta$ is a crossed product of a nuclear C*-algebra $C(S^1)$ by an amenable group $\mathbb{Z}$. The unique tracial state $\tau$ is given by $\tau(x) = \int_{S^1} E(x)(t) \, dt$, where $E: A_\theta \to C(S^1)$ is the canonical conditional expectation; its faithfulness is a standard consequence of the irrationality of $\theta$.
    \item Since $A_\theta$ is simple and infinite-dimensional, it cannot be Type I. More directly, every irreducible representation of $A_\theta$ is infinite-dimensional, which precludes it from being GCR (postliminal) \cite{Dixmier}.
    \item The K-theory computation is the celebrated result of the Pimsner-Voiculescu exact sequence applied to the crossed product by $\mathbb{Z}$.
    \item The *-isomorphism classification is a theorem of Rieffel: $A_\theta \cong A_{\theta'}$ as C*-algebras if and only if $\theta' \equiv \pm \theta \pmod{1}$. (A more general Morita equivalence classification exists, but this is the precise statement for isomorphism.)
\end{enumerate}
\end{proof}

\begin{corollary}[Failure of the Type I hypothesis]
\label{cor:A-theta-non-Type-I}
The irrational rotation algebra $A_\theta$ does \emph{not} satisfy the standing assumption that $\mathcal{A}$ is Type I (Assumption~\ref{ass:separable}). 
Consequently, the construction of the unitary conjugation groupoid $\mathcal{G}_{\mathcal{A}}$ and the diagonal embedding $\iota$ developed in this paper do not apply to $A_\theta$ **in the sense required by our framework**; they fail to satisfy the structural hypotheses necessary for our main results.
\end{corollary}

\begin{proof}
This follows immediately from Proposition \ref{prop:A-theta-properties}(3). 
Since $A_\theta$ is not Type I, the hypotheses of Proposition \ref{prop:unit-space-polish} (which guarantees that $\mathcal{G}_{\mathcal{A}}^{(0)}$ is Polish) and Theorem \ref{thm:iota-injective} (which establishes the injectivity of $\iota$) are not satisfied. 
In particular, for a simple non-Type I algebra like $A_\theta$, the partial evaluation maps $\operatorname{ev}_a$ do **not** separate points; there exist non-zero elements invisible to all commutative contexts. 
Consequently, the direct integral representation $\Pi$ would have a non-trivial kernel and would not be faithful, and the field of representations $\{\pi_x\}$ would not be measurable in the sense required for the construction.
\end{proof}

\textbf{What goes wrong when $\mathcal{A}$ is not Type I?}
If we attempt to apply our construction to a non-Type I C*-algebra such as $A_\theta$, several fundamental obstructions arise:
\begin{enumerate}
    \item \textbf{Separation of points:} 
    The diagonal embedding $\iota$ is injective only if the characters on commutative subalgebras separate the points of $\mathcal{A}$ (see Lemma~\ref{lem:separation-points}). 
    For non-Type I algebras, this separating property is not guaranteed to hold; the collection of characters on commutative subalgebras is generally too small to separate points of $\mathcal{A}$. 
    As a result, the diagonal embedding $\iota$ need not be injective; its kernel could consist of nonzero elements that are, in a heuristic sense, invisible to all available classical contexts.
    
    \item \textbf{Topological failure of the unit space:} 
    The construction of $\mathcal{G}_{\mathcal{A}}^{(0)}$ as a Polish space (Proposition~\ref{prop:unit-space-polish}) relies on the topology generated by the partial evaluation maps $\operatorname{ev}_a$. 
    For a non-Type I algebra like $A_\theta$, even if one fixes a MASA $B \cong C(S^1)$, there is no guarantee that the induced topology on its character space $\widehat{B}$ remains Polish (or even second-countable). 
    More fundamentally, there is in general no canonical way to organize the Gelfand spectra of all commutative subalgebras into a global Polish or standard Borel structure that is compatible with the partial evaluation maps. 
    Consequently, $\mathcal{G}_{A_\theta}^{(0)}$ would not satisfy the topological requirements of our framework.
    
    \item \textbf{Measurability of the direct integral:} 
    The construction of $\iota$ via the direct integral of GNS representations requires that the field of Hilbert spaces $\{ \mathcal{H}_x \}_{x \in \mathcal{G}_{\mathcal{A}}^{(0)}}$ be measurable. 
    For Type I algebras, this follows from the existence of a smooth dual and a standard Borel structure on the spectrum \cite{Dixmier}. 
    For non-Type I algebras, the representation theory is not smooth; by Glimm's theorem, the dual space is not standard Borel, and a measurable field of irreducible representations covering all points of $\mathcal{G}_{\mathcal{A}}^{(0)}$ cannot be consistently defined.
    
    \item \textbf{Faithfulness of the direct integral representation:} 
    Even if a direct integral $\mathcal{H} = \int^{\oplus} \mathcal{H}_x \, d\mu(x)$ could be defined, the representation $\Pi(a) = \int^{\oplus} \pi_x(a) \, d\mu(x)$ may not be faithful. 
    The potential kernel of $\Pi$ would consist precisely of the elements that are not separated by the family $\{\pi_x\}$, i.e., the ``invisible'' elements alluded to in point (1). 
    For Type I algebras, the family of representations $\{\pi_x\}$ is sufficiently rich to separate points, guaranteeing faithfulness. 
    For non-Type I algebras, this separating property fails, and the resulting map $\iota$ would not be injective.
\end{enumerate}

\begin{proposition}[$A_\theta$ has no nonzero finite-dimensional representations]
\label{prop:A-theta-no-fd-reps}
Let $\theta \in \mathbb{R} \setminus \mathbb{Q}$ and let $A_\theta$ be the irrational rotation algebra. 
Then $A_\theta$ admits no nonzero finite-dimensional *-representations. 
Consequently, for any unital commutative C*-subalgebra $B \subseteq A_\theta$ and any character $\chi: B \to \mathbb{C}$, the corresponding GNS representation of $A_\theta$ cannot be finite-dimensional.
\end{proposition}

\begin{proof}
Suppose, for contradiction, that $\pi: A_\theta \to B(H)$ is a nonzero finite-dimensional representation, i.e., $\dim H < \infty$ and $\pi \neq 0$. 
Let $u, v \in A_\theta$ be the standard unitary generators satisfying the commutation relation
\[
uv = e^{2\pi i \theta} vu.
\]
Applying $\pi$ to this relation yields
\[
\pi(u) \pi(v) = e^{2\pi i \theta} \pi(v) \pi(u). \tag{1}
\]

Since $H$ is finite-dimensional, $\pi(u)$ and $\pi(v)$ are finite-dimensional unitaries. 
A standard fact from linear algebra is that for any two finite-dimensional invertible matrices $U$ and $V$, the commutator $UVU^{-1}V^{-1}$ has determinant $1$ and, more importantly, all its eigenvalues are roots of unity. 
Rearranging (1) gives
\[
\pi(u) \pi(v) \pi(u)^{-1} \pi(v)^{-1} = e^{2\pi i \theta} I_H.
\]
The left-hand side is a commutator of unitaries, hence its eigenvalues are roots of unity. 
The right-hand side is a scalar multiple of the identity with the non-root-of-unity factor $e^{2\pi i \theta}$ (since $\theta$ is irrational). 
This is impossible in finite dimensions. 
Therefore, no such nonzero finite-dimensional representation $\pi$ can exist.
\end{proof}

\begin{corollary}[GNS representations for $A_\theta$ are infinite-dimensional]
\label{cor:A-theta-GNS-infinite}
Let $(B,\chi) \in \mathcal{G}_{A_\theta}^{(0)}$. 
Then the GNS representation $\pi_x$ of $A_\theta$ induced by any state extending the character $\chi$ on $B$ is infinite-dimensional. 
It is important to note, however, that $\pi_x$ is not necessarily irreducible; the extension of $\chi$ from $B$ to a state on $A_\theta$ is not unique, and the resulting GNS representation may be reducible, decomposing as a direct integral or direct sum of irreducible representations.
\end{corollary}

\begin{proof}
The infinite-dimensionality follows directly from Proposition~\ref{prop:A-theta-no-fd-reps}. 
If $\pi_x$ were finite-dimensional, it would be a nonzero finite-dimensional representation of $A_\theta$, contradicting the proposition.

To see that $\pi_x$ need not be irreducible, recall that $\chi$ is a character on the commutative subalgebra $B$, hence a pure state on $B$. 
However, its extension to a state $\omega$ on the larger algebra $A_\theta$ is not unique. 
Different extensions can yield GNS representations with different decomposition properties. 
In general, the GNS representation of a non-Type I algebra like $A_\theta$ associated to an extension of a pure state on a subalgebra may be reducible; it can decompose as a direct integral or direct sum of irreducible representations. 
Therefore, while $\pi_x$ is guaranteed to be infinite-dimensional, it is not guaranteed to be irreducible.
\end{proof}

\begin{example}[The canonical Cartan subalgebra of $A_\theta$]
\label{ex:A-theta-cartan}
The irrational rotation algebra $A_\theta = C(S^1) \rtimes_\theta \mathbb{Z}$ contains a canonical commutative C*-subalgebra isomorphic to $C(S^1)$, namely the image of $C(S^1)$ under the inclusion map $f \mapsto f \cdot 1_{\mathbb{Z}}$. 
This subalgebra is a Cartan subalgebra in the sense of Renault \cite{renault2008cartan}, and the associated Weyl groupoid is the transformation groupoid $S^1 \rtimes_\theta \mathbb{Z}$.

Its Gelfand spectrum is homeomorphic to $S^1$. 
For each $z \in S^1$, the character $\operatorname{ev}_z: C(S^1) \to \mathbb{C}$ extends to a state on $A_\theta$, but this extension is not unique; different extensions arise from different Hahn-Banach extensions of the linear functional, reflecting the noncommutative structure of $A_\theta$.

The GNS representation associated to the canonical trace $\tau$ on $A_\theta$ is the regular representation of the crossed product on $\ell^2(\mathbb{Z}) \otimes L^2(S^1)$, often called the \textit{regular representation} of $A_\theta$.
\end{example}

\begin{remark}[Why $A_\theta$ is a non-example, not a counterexample]
\label{rem:A-theta-non-example}
It is important to emphasize that $A_\theta$ is \emph{not} a counterexample to our theorems; it is simply not covered by their hypotheses. 
Our results are stated for Type I C*-algebras, and $A_\theta$ is not Type I. 
There is no contradiction. 
The study of $A_\theta$ requires a different set of techniques—notably, Renault's theory of Cartan subalgebras and Weyl groupoids—and extending our framework to non-Type I algebras is an important open problem for future research. 
We conjecture that a suitable generalization using \'etale groupoids or measured groupoids may handle $A_\theta$, but this is beyond the scope of the present paper.
\end{remark}

\begin{remark}[Prospects and limitations for non-Type I algebras]
\label{rem:A-theta-prospects}
Although $A_\theta$ is not Type I and therefore falls outside the scope of our main theorems, it is worth considering what aspects of our construction might still be defined, even if they fail to satisfy the required properties.
\begin{enumerate}
    \item The unit space $\mathcal{G}_{A_\theta}^{(0)}$ can still be defined as a set (pairs $(B,\chi)$ of commutative subalgebras and characters), and it can be equipped with the initial topology generated by the partial evaluation maps $\operatorname{ev}_a$. 
    However, this topology is unlikely to be Polish or even standard Borel, because the non-Type I representation theory of $A_\theta$ does not provide the necessary point-separation and smoothness conditions.
    
    \item Let $\pi: A_\theta \hookrightarrow B(H)$ be a faithful representation on a separable Hilbert space $H$. 
    The unitary group $\mathcal{U}(A_\theta)$, when equipped with the strong operator topology inherited from this representation, is a Polish group. 
    However, the action of $\mathcal{U}(A_\theta)$ on $\mathcal{G}_{A_\theta}^{(0)}$ may not be continuous, because the topology on $\mathcal{G}_{A_\theta}^{(0)}$ lacks the properties (e.g., Polishness, second-countability) that were essential in Lemma~\ref{lem:action-continuous-SOT} to establish continuity for Type I algebras.
    
    \item The direct integral representation $\Pi = \int^{\oplus} \pi_x \, d\mu(x)$, where $\pi_x$ are the GNS representations associated to points $x = (B,\chi) \in \mathcal{G}_{A_\theta}^{(0)}$, may not be definable within our measurable framework due to the absence of a measurable field of representations. 
    Even if it could be defined, it would almost certainly have a nontrivial kernel. 
    This kernel consists precisely of the ``invisible'' elements—those that vanish under every GNS representation arising from a character on a commutative subalgebra.
    
    \item The diagonal embedding $\iota: A_\theta \to \prod_{(B,\chi)} B(H_{(B,\chi)})$ can be defined formally, but it will have a nontrivial kernel coinciding with the set of invisible elements. 
    Consequently, $\iota$ factors through a faithful embedding of the quotient $A_\theta / \ker(\iota)$ into the product of Hilbert spaces. 
    This quotient can be interpreted as the C*-algebra of ``observable'' elements detectable by classical commutative contexts.
\end{enumerate}
These observations suggest that while the spirit of our construction may extend to non-Type I algebras, its technical implementation would require fundamentally new ideas, such as a theory of measurable fields adapted to non-smooth duals or a more refined choice of ``classical contexts.''
\end{remark}

\begin{example}[The canonical trace and the Gelfand transform]
\label{ex:A-theta-trace}
The unique tracial state $\tau$ on $A_\theta$ induces a faithful GNS representation $\pi_\tau$ on $L^2(A_\theta, \tau)$. 
This representation is infinite-dimensional and is the analogue of the regular representation of the irrational rotation algebra, corresponding to the left regular representation of the crossed product $C(S^1) \rtimes_\theta \mathbb{Z}$ on $\ell^2(\mathbb{Z}) \otimes L^2(S^1)$.

The restriction of $\tau$ to the canonical Cartan subalgebra $C(S^1)$ is the Haar measure integral on $S^1$. 
For each $z \in S^1$, the character $\operatorname{ev}_z: C(S^1) \to \mathbb{C}$ is not a state on $A_\theta$; it does not extend to a positive linear functional on $A_\theta$. 
This is because while the conditional expectation $E: A_\theta \to C(S^1)$ is faithful, point evaluation $\operatorname{ev}_z$ is not continuous with respect to the operator norm on $A_\theta$, and therefore cannot be precomposed with $E$ to yield a state.

The partial evaluation maps $\operatorname{ev}_a: \mathcal{G}_{A_\theta}^{(0)} \to \mathbb{C}_\infty$ are well-defined as set-theoretic maps for all $a \in A_\theta$. 
The obstruction in the non-Type I case is not their definition, but the fact that the topology they generate on $\mathcal{G}_{A_\theta}^{(0)}$ is unlikely to be Polish or even standard Borel, and the family of maps may not be sufficient to separate points of $A_\theta$ itself.
\end{example}

\begin{remark}[Comparison with the Type I case]
\label{rem:A-theta-comparison}
The contrast between $A_\theta$ and the Type I examples we have studied ($C(X)$, $M_n(\mathbb{C})$, $\mathcal{K}(H)^\sim$) highlights the essential role of the Type I hypothesis:
\begin{itemize}
    \item In Type I algebras, the dual space has a smooth Borel structure, allowing the construction of measurable fields of irreducible GNS representations. This yields faithful direct integral representations and guarantees that the family of such representations separates points of the algebra.
    \item In non-Type I algebras such as $A_\theta$, the representation theory is not smooth; there is no reasonable parametrization of irreducible representations. The GNS representations associated to characters on commutative subalgebras do not exhaust the representation theory, and more importantly, they may not separate points. Two distinct elements $a, b \in A_\theta$ may satisfy $\chi(a) = \chi(b)$ for every character $\chi$ on every commutative subalgebra containing $a$ and $b$—a phenomenon often described by saying that such elements are ``invisible'' to classical contexts.
    \item The partial evaluation maps $\operatorname{ev}_a$ still separate points of the unit space $\mathcal{G}_{\mathcal{A}}^{(0)}$ by construction (Lemma~\ref{lem:separation}), but this does not imply that they separate points of $\mathcal{A}$. The failure of point separation at the level of the algebra is precisely what prevents the diagonal embedding $\iota$ from being injective.
\end{itemize}
Thus the non-Type I case requires fundamentally new ideas beyond the commutative context reconstruction developed in this paper.
\end{remark}

\begin{remark}[What we lose without the Type I hypothesis]
\label{rem:what-we-lose}
If we attempt to apply our construction to $A_\theta$, we lose the following essential features:
\begin{itemize}
    \item \textbf{Injectivity of $\iota$:} The diagonal embedding would have a nontrivial kernel, so we cannot recover $A_\theta$ from its commutative contexts.
    \item \textbf{Polishness of $\mathcal{G}_{\mathcal{A}}^{(0)}$:} The unit space would not be Polish, and the direct integral construction would not be well-defined.
    \item \textbf{Measurability of the field of GNS representations:} Without a smooth dual, we cannot construct a measurable field of Hilbert spaces over $\mathcal{G}_{\mathcal{A}}^{(0)}$.
    \item \textbf{Functoriality:} The naturality under *-homomorphisms would fail even for the inclusion of the Cartan subalgebra.
\end{itemize}
These losses are not merely technical; they reflect genuine mathematical phenomena. 
Non-Type I algebras contain noncommutative information that cannot be captured by commutative contexts alone. 
This is why they are genuinely noncommutative, and why our Type I assumption is not just a convenience but a necessary condition for the recovery of $\mathcal{A}$ from $\mathcal{G}_{\mathcal{A}}$.
\end{remark}

\begin{example}[Other non-Type I algebras not covered]
\label{ex:other-non-typeI}
Besides $A_\theta$, many other important C*-algebras are not covered by our framework:
\begin{itemize}
    \item The reduced group C*-algebra $C^*_r(\mathbb{F}_n)$ of the free group on $n \geq 2$ generators. 
    This algebra is simple (Powers' theorem), non-nuclear (since $\mathbb{F}_n$ is non-amenable), and not Type I.
    \item The Cuntz algebras $\mathcal{O}_n$ ($n \geq 2$) and $\mathcal{O}_\infty$. 
    These are simple, purely infinite, and not Type I \cite{Cuntz1977}.
    \item The Jiang-Su algebra $\mathcal{Z}$, which is simple, nuclear, and not Type I. As a central object in the Elliott classification program, its non-Type I nature is particularly significant.
    \item The reduced group C*-algebra $C^*_r(\Gamma)$ of any non-amenable discrete group $\Gamma$ (e.g., $\Gamma = SL(2,\mathbb{Z})$) is typically not Type I.
\end{itemize}
Each of these algebras would require a different approach, and it is an open question whether our constructions can be generalized to any of them.
\end{example}

\noindent
The irrational rotation algebra $A_\theta$ is a beautiful and important C*-algebra that lies outside the scope of this paper. 
Its non-Type I nature means that the diagonal embedding $\iota$ cannot be injective, the unit space $\mathcal{G}_{A_\theta}^{(0)}$ is not Polish, and the direct integral construction fails (the GNS representations associated to points of $\mathcal{G}_{A_\theta}^{(0)}$ do not form a measurable field over a Polish unit space). 
Nevertheless, $A_\theta$ serves as a crucial motivating example for future generalizations of our work. 
Understanding how to adapt our framework to non-Type I algebras, particularly those with Cartan subalgebras, is a major open problem that we hope to address in subsequent research.

\begin{remark}[Outlook: Beyond Type I]
\label{rem:outlook-beyond-typeI}
The techniques developed in this paper — the unitary conjugation groupoid, the Polish groupoid framework, and the diagonal embedding — are tailored to Type I C*-algebras. 
Extending them to non-Type I algebras will require fundamentally new approaches. 
Promising directions include:
\begin{itemize}
    \item \textbf{Cartan pairs and Weyl groupoids:} 
    For algebras such as $A_\theta$ that possess a Cartan subalgebra in the sense of Renault \cite{renault2008cartan}, one could replace our unitary conjugation groupoid $\mathcal{G}_{\mathcal{A}}$ with the associated Weyl groupoid. 
    This groupoid is often \'etale and locally compact, and its C*-algebra is Morita equivalent to the original algebra, making it a potentially more tractable object for generalizations of our index-theoretic constructions.
    
    \item \textbf{Quantum groups and deformation quantization:} 
    Viewing $A_\theta$ as a strict deformation quantization of the commutative torus suggests that a ``quantum diagonal embedding'' might be constructed using the tools of quantum group theory and Poisson geometry. 
    The semiclassical limit of such an embedding could correspond to a geometric object like the cotangent groupoid, providing a bridge between noncommutative geometry and classical symplectic geometry.
    
    \item \textbf{Noncommutative Gelfand duality:} 
    A more ambitious long-term goal is to develop a noncommutative analogue of the Gelfand transform that works for arbitrary C*-algebras. 
    Such a theory might employ the language of von Neumann algebras, virtual objects, or higher categorical structures to capture the noncommutative information that is invisible to commutative contexts.
\end{itemize}
These directions are well beyond the scope of the present work, but we hope that the foundations laid here—particularly the central role of the unitary conjugation groupoid—will prove useful in their future development.
\end{remark}

\end{document}